\documentclass[11pt]{article}
\usepackage{amssymb}
\usepackage{mathrsfs}
\usepackage{latexsym,amsmath,amssymb,amsfonts,epsfig,graphicx,cite,psfrag}
\usepackage{cases}
\usepackage{eepic,color,colordvi,amscd}
\usepackage{ebezier}
\usepackage{verbatim}
\usepackage{comment}

\newtheorem{theo}{Theorem}[section]

\newtheorem{lem}[theo]{Lemma}
\newtheorem{coro}[theo]{Corollary}

\newtheorem{conj}[theo]{Conjecture}
\newtheorem{defi}[theo]{Definition}
\newtheorem{nota}[theo]{Notation}

\newtheorem{rem}[theo]{Remark}

\def\qed{\hfill \rule{4pt}{7pt}}

\def\pf{\noindent {\it Proof. }}

\begin{document}

\title{Induced Forests in Bipartite Planar Graphs\footnote{Partially supported by NSF grants DMS-1265564 and CNS-1443894}}
\author{Yan Wang\footnote{Email: yanwang@gatech.edu}, Qiqin Xie\footnote{Email: qxie7@math.gatech.edu}, Xingxing Yu\footnote{Email: yu@math.gatech.edu}\\
School  of Mathematics\\
Georgia Institute of Technology\\
Atlanta, GA 30332-0160, USA}


\maketitle

\begin{abstract} 

Akiyama and Watanabe
conjectured that every simple planar bipartite graph on $n$ vertices contains 
an induced forest on at least $5n/8$ vertices. We apply the discharging method to show that every simple bipartite 
planar graph on $n$ vertices contains an induced forest on at least $\lceil (4n+3)/7 \rceil$ vertices.


\end{abstract}


\section{Introduction}

In this paper, we consider simple graphs only. Clearly, every bipartite graph contains an independent set of size at least half of its vertices. 
It is natural to ask under what conditions can we find  considerably larger sparse induced subgraphs, 
for example, induced forests? The study of the maximum size of induced forests was initiated 
by Erd\H{o}s, Saks, and S\'{o}s in 1986 \cite{Erdos}. Later, Matou\v{s}ek and \v{S}\'{a}mal \cite{Matousek}, and 
also Fox, Loh, and Sudakov \cite{Fox} studied large induced trees in triangle-free graphs and $K_r$-free graphs, respectively.

For a graph $G$, let $|G|=|V(G)|$ and let $a(G)$ denote the largest number of vertices of an induced forest in $G$. For later convenience, we use $A(G)$ to denote an induced forest in $G$ 
of size $a(G)$. Albertson and Berman \cite{Albertson1} (also see Albertson and Haas in \cite{Albertson2}) conjectured in 1979 that $a(G) \geq |G|/2$
for any planar graph $G$.
For bipartite planar graphs, Akiyama and Watanabe \cite{Akiyama} made the following in 1987

\begin{conj}
\label{conjecture2}
If $G$ is a bipartite planar graph, then $a(G) \geq 5|G|/8$.
\end{conj}

The bound in Conjecture \ref{conjecture2} is tight with $Q_3$ (the $3$-cube), and more examples can be constructed, for example, 
by adding a matching between two $4$-cycles in two $Q_3$'s.  

Planar graphs have average degree strictly less than $6$. Alon \cite{Alon1} considered bipartite graphs $G$ with average degree at most $d \geq 1$, 
and showed that $a(G) \geq (\frac{1}{2} + e^{-bd^2})|G|$, for some absolute constant $b > 0$. 
Conlon \textit{et al.} \cite{Conlon} improved Alon's bound to $(\frac{1}{2} + d^{-b'd})|G|$, for some constant $b'>0$.  
Since the average degree of any bipartite planar graph is less than 4, the above results  
give a nontrivial bound for Conjecture \ref{conjecture2}.

There has been some recent activities on  Conjecture \ref{conjecture2}. 
It is shown in \cite{Salavatipour} (also see \cite{Kowalik}) that if $G$ is a triangle-free planar graph then $a(G)\ge (17|G|+24)/32$, 
which is improved to $(6|G|+7)/11$ in \cite{Dross}. 
In this paper, we prove the following

\begin{theo}
\label{MainTheorem}
Let $G$ be a bipartite planar graph. Then $a(G) \geq \lceil (4|G|+3)/7 \rceil$.
\end{theo}

In our proof of Theorem \ref{MainTheorem}, we apply the discharging method. 
Suppose Theorem~\ref{MainTheorem} is false, and let $G$ be a counterexample with $|G|$ minimum. 
Using the discharging technique, we force some small configurations, which are reducible in the sense 
that after certain operations we can use an induced forest from a smaller graph to construct an induced forest in $G$. Often such operations involve the identification of vertices, which  may result in multiple edges; 
we remove all but one such edges after the identification.
Note that we always identify vertices in the same color class of the bipartite graph $G$. Hence, there will be no loop after the identification. 

We need some notations and terminologies. 
Let $v \in V(G)$ and $X,Y \subseteq V(G)$. $N(v)$ denotes the set of neighbors of $v$, and $G[X]$ 
denotes the induced subgraph of $G$ on $X$. We define $G - v := G[V(G) - \{v\}]$, $G - X := G[V(G) - X]$, $G[X + v] := G[X \cup \{v\}]$ and $G[X + Y] := G[X \cup Y]$. 
Let $n$ be a positive integer.
We denote $V_n$, $V_{\leq n}$, $V_{\geq n}$ the set of vertices of degree exactly $n$, at most $n$, and at least $n$, respectively.
We call a vertex $v$ in $G$ is a $n$-vertex ($n^+$-vertex, $n^-$-vertex, respectively) if $v \in V_n$ ($v \in V_{\geq n}$, $v \in V_{\leq n}$, respectively)
If $G$ is a planar graph and $v_1,v_2,...,v_k$ are vertices of $G$ incident with a common face $F$, then $G/v_1v_2...v_k$ 
denotes the simple plane graph obtained from $G$ by identifying $v_1,v_2,...,v_k$ in $F$ as a new vertex $w$. 
We define $G/\{v_1v_2,...,v_{k-1}v_k\} = (G/v_1v_2)/\{v_3v_4,...,v_{k-1}v_k\}$.
$G+v_1v_2$ denotes the simple plane graph obtained from $G$ by adding the edge $v_1v_2$ in $F$ if $v_1v_2 \not\in E(G)$.
$X \triangle Y$ denotes the symmetric difference between $X$ and $Y$. 
A {\it separation} in a graph $G$ consists of a pair of subgraphs $G_1, G_2$, denoted as $(G_1,G_2)$, such that
$E(G_1) \cup E(G_2)=E(G)$, $E(G_1\cap G_2)=\emptyset$, 
$G_1 \not\subseteq G_2$, and $G_2 \not\subseteq G_1$.
$e(X)$ denotes the number of edges in $G[X]$ and $e(X,Y)$ denotes the number of  edges of  $G$ between vertices in $X$ and vertices in $Y$.

The rest of the paper is organized as follows. 
In Section 2, we present some inequalities that we use, which can be established by considering remainders modular 7. 
We also set up some notation for a minimum counterexample $G$ of Theorem \ref{MainTheorem}, and prove some basic properties about $G$. 
In Section 3, we derive information about the structures around a vertex of degree 2 in $G$. In Section 4, we work on the neighbors of a degree 3 vertex. 
In Section 5 and 6, we deal with two forbidden configurations around a 3-vertex. 
In Section 7, we work with degree 5 and 6 vertices.   We prove Theorem \ref{MainTheorem} in Section 8 by giving discharging rules based on the structural information obtained in the previous sections.  


\section{Useful inequalities and the minimum counterexample}

We begin with some inequalities that will be used frequently throughout the paper. 

\begin{lem}\label{ineq1}
Let $a_1,a_2\ge 1$ be integers such that $a_1+a_2=n+3-k$, with $k\le 8$. Then 
$\max\{\lceil(4a_1+3)/7\rceil+ \lceil(4a_2+3)/7\rceil+2, \lceil(4a_1-1)/7\rceil+ \lceil(4a_2-1)/7\rceil+3\}\ge \lceil(4n+3)/7\rceil.$
\end{lem}
\pf Note the symmetry between $a_1$ and $a_2$. If $4a_1+3\equiv 0 \mod 7$ then 
$\lceil(4a_1-1)/7\rceil+ \lceil(4a_2-1)/7\rceil+3\ge (4a_1-1+4)/7 +(4a_2-1)/7+3=(4n+3-4k+32)/7 \geq (4n+3)/7. $

So we may assume  $4a_i+3\not\equiv 0 \mod 7$ for $i=1,2$. Let $4a_i+3\equiv r_i \mod 7$ with $1\le r_i\le 6$ for $i=1,2$. If $r_1\ne 6$ or $r_2\ne 6$ then 
$\lceil(4a_1+3)/7\rceil+ \lceil(4a_2+3)/7\rceil+2\ge (4a_1+3)/7+(4a_2+3)/7+2+3/7=(4n+3-4k+32)/7\ge (4n+3)/7.$

So assume $r_1=r_2=6$. Then 
$\lceil(4a_1-1)/7\rceil+ \lceil(4a_2-1)/7\rceil+3\ge (4a_1-1+5)/7 +(4a_2-1+5)/7+3=(4n+3-4k+38)/7>(4n+3)/7. $

Therefore, the conclusion holds since the left hand side of the inequality is an integer.
\qed

\medskip
With similar, but more involved arguments, we have the following inequalities. We leave out the details.

\begin{lem}\label{ineq2}
Let $a,a_1,a_2,...,a_k,c,n$ be positive integers where $k \geq 1$.
Let $L$ be a set of integers and $b_j$ be a positive integer for all $j \in L$.
\begin{itemize}
\item[(1)] If $(4a+3)/7 + \sum\limits_{i=1}^k (4a_i+3)/7 + \sum\limits_{j \in L} (4b_j+3)/7 + c - k  \geq (4n+3-3k)/7$, then 
$\max\limits_{A_i \in \{0,1\},\forall i \in [k]} \{ \lceil (4(a-\sum\limits_{i=1}^k A_i)+3)/7 \rceil +  \sum\limits_{i=1}^k \lceil (4(a_i-A_i)+3)/7 \rceil + \sum\limits_{j\in L} \lceil (4b_j+3)/7 \rceil + c - \sum\limits_{i=1}^k (1 - A_i) \} \geq \lceil (4n+3)/7 \rceil;$
\item[(2)] If $(4a+3)/7 + (4a_1+3)/7 + c - 1  \geq (4n-1)/7$ and $(4a+3,4a_1+3) \not\equiv (0,4),(4,0) \mod 7$, then 
$\max\limits_{A_1 \in \{0,1\}} \{ \lceil (4(a-A_1)+3)/7 \rceil +  \lceil (4(a_1-A_1)+3)/7 \rceil + c - (1 - A_1) \} \geq \lceil (4n+3)/7 \rceil;$

\item[(3)] If $(4a+3)/7 + (4a_1+3)/7 + c  \geq (4n-1)/7$, then $ \lceil (4a+3)/7 \rceil +  \lceil (4a_1+3)/7 \rceil + c  \geq \lceil (4n+3)/7 \rceil$ if $(4a+3,4a_1+3) \not\equiv (0,0),(0,6),(0,5),(0,4),(4,0),(6,5),(5,6),(5,0),(6,6),(6,0) \mod 7$;

\item[(4)] If $(4a+3)/7 + \sum\limits_{i=1}^2  (4a_i+3)/7 + c - 2  \geq (4n-4)/7$, then $\max\limits_{A_1,A_2 \in \{0,1\}} \{ \lceil (4(a-\sum\limits_{i=1}^2 A_i)+3)/7 \rceil +  \sum\limits_{i=1}^2 \lceil (4(a_i-A_i)+3)/7 \rceil + c - \sum\limits_{i=1}^2 (1-A_i) \} \geq \lceil (4n+3)/7 \rceil,$ unless $(4a+3,4a_1+3,4a_2+3) \equiv (1,0,0), (4,0,4), (4,4,0), (0,4,4) \mod 7$;

\item[(5)] If $(4a+3)/7 + \sum\limits_{i=1}^2  (4a_i+3)/7 + c - 2  \geq (4n-5)/7$, then $\max\limits_{A_1,A_2 \in \{0,1\}} \{ \lceil (4(a-\sum\limits_{i=1}^2 A_i)+3)/7 \rceil + \sum\limits_{i=1}^2 \lceil (4(a_i-A_i)+3)/7 \rceil + c - \sum\limits_{i=1}^2 (1-A_i) \} \geq \lceil (4n+3)/7 \rceil,$ unless $(4a+3,4a_1+3,4a_2+3) \equiv $ $(0, 0, 0)$, $(1, 0, 0)$, $(4, 0, 3), (4, 3, 0), (3, 0, 4), (4, 0, 4), (3, 4, 0), (4, 4, 0)$, $(1, 6, 0)$, $(1, 0, 6)$, $(0, 3, 4)$, $(0, 4, 3), (0, 4, 4), (6,4,4), (4,4,6), (4,6,4)  \mod 7$;

\item[(6)] If $\sum\limits_{i=1}^k (4a_i+3)/7 + c  \geq (4n+2)/7$, then $\sum\limits_{i=1}^k \lceil (4a_i+3)/7 \rceil + c \geq \lceil (4n+3)/7 \rceil$, unless $4a_i+3 \equiv 0 \mod 7$ for $i \in [k]$;

\item[(7)] If $\sum\limits_{i=1}^k (4a_i+3)/7 + c  \geq (4n+1)/7$, then $\sum\limits_{i=1}^k \lceil (4a_i+3)/7 \rceil + c \geq \lceil (4n+3)/7 \rceil$, unless there exists $j \in [n]$ such that $4a_j+3 \equiv 0,6 \mod 7$ and $4a_i+3 \equiv 0 \mod 7$ for $i \in [k] - \{j\}$;

\item[(8)] If $(4a+3)/7 + (4a_1+3)/7 + c  \geq 4n/7$, then $ \lceil (4a+3)/7 \rceil +  \lceil (4a_1+3)/7 \rceil + c  \geq \lceil (4n+3)/7 \rceil$ unless $(4a+3,4a_1+3) \equiv (0,0),(0,6),(0,5),(5,0),(6,6),(6,0) \mod 7$.
\end{itemize}
\end{lem}

Note that in applications $a_1,a_2,...,a_k,b_1,...,b_l$ are the numbers of vertices in some subgraphs of a given graph, and $A_i$ is the indicator function whether a vertex is included or not. Moreover, we have $k \leq 4$ and $l \leq 2$ in all applications.

We now set up some notation for the proof of  Theorem~\ref{MainTheorem}. Throughtout the remainder of this paper, 
let $G$ be a bipartite plane graph with $|G|=n$ such that 
\begin{itemize}
\item [(i)] $a(G)<\lceil (4n+3)/7\rceil$, 
\item[(ii)] subject to (i), $|G|$ is minimum, and 
\item[(iii)] subject to (ii), $|E(G)|$ is maximum. 
\end{itemize}

\begin{lem} \label{basic_lemma}
$G$ is a connected quadrangulation,  $\delta(G) \geq 2$, and for each  $v \in V_{\le 3}$ we may choose $A(G)$ so that $v\in A(G)$.
\end{lem}

\pf If $G$ is disconnected, let $G_1,...,G_k$ be the components of $G$ (hence $k \geq 2$). 
By the choice of $G$, $a(G_i) \geq \lceil (4|G_i|+3)/7 \rceil$ for $i\in [k]$. So $a(G) 
\geq \sum_{i=1}^k \lceil (4|G_i|+3)/7 \rceil \geq \lceil (4n+3)/7 \rceil$, a contradiction. So $G$ is connected.

If $G$ is not a quadrangulation, then $G$ has a facial walk $a_1a_2...a_ka_1$ with $k \geq 6$. By the choice of $G$, 
$a(G + a_1a_4) \geq \lceil (4n+3)/7 \rceil$. This implies that $a(G) \geq \lceil (4n+3)/7 \rceil$, a contradiction.
Thus $G$ is a quadrangulation, and hence, $\delta(G)\ge 2$. 


Now let $F=A(G)$ with $v\in V_{\le 3} - V(F)$. By the maximality of $A(G)$, $N(v)\cap V(F)\ne \emptyset$. 
If $|V(F) \cap N(v)| \leq 2$, then let $w \in V(F) \cap N(v)$; if $|V(F) \cap N(v)| =3$, then there exists  $w \in V(F)$ such that no two vertices in 
$V(F) \cap N(v)$ are contained in the same component of $F - w$. Now $G[F-w + v]$ is a maximum induced forest in $G$ containing $v$.
\qed

\medskip

The following notation will be convenient when performing graph operations. 

\begin{nota}
\label{DefR}
Let $v \in V(G)$ and $U \subseteq N(v)$. Define $R_{v,U} := R_{v,U}^1 \cup R_{v,U}^2$ where $R_{v,U}^1 = \{\{r\} \subseteq  N(v) -U: r \in V_{\leq 2}\}$ and 
$R_{v,U}^2 = \{\{r_1,r_2\} \subseteq N(v) -U: r_1,r_2 \in V_{3}$ and $r_1,r_2$
are cofacial $\}$. 
\end{nota}

\begin{lem}
\label{NoTwoRNeighbors} 
For any $v \in V(G)$ and $U\subseteq N(v)$, if $R_1,R_2 \in R_{v,U}$, then $R_1 \cap R_2 \neq \emptyset$.
\end{lem}

\pf First, assume that there exist distinct  $\{x\},\{y\} \in R_{v,U}^1$. Let $F'=A(G-\{v,x,y\})$. By the choice of $G$, 
$|F'| = a(G') \geq \lceil (4(n-3)+3)/7 \rceil$. Hence $G[F' + \{x,y\}]$ is an induced forest in $G$; so  
$a(G) \geq |F'|+2 \geq \lceil (4n+3)/7 \rceil$, a contradiction.

Now assume there exist $\{x\} \in R_{v,U}^1, \{y,z\} \in R_{v,U}^2$. Let $w\in V(G)$ such that $vywzv$ is a facial cycle. 
Let $F'= A(G-\{v,x,y,z,w\})$. Then $|F'| \geq \lceil (4(n-5)+3)/7 \rceil$ by the choice of $G$. 
Clearly, $G[F' + \{x,y,z\}]$ is an induced forest in $G$; so $a(G) \geq |F'|+3 \geq \lceil (4n+3)/7 \rceil$, a contradiction.

Finally, assume $\{x_1,x_2\}, \{y_1,y_2\} \in R_{v,U}^2$ with $\{x_1,x_2\} \cap \{y_1,y_2\} = \emptyset$. 
Let $x_3,y_3\in V(G)$ such that $vx_1x_3x_2v$ and $vy_1y_3y_2v$ are facial cycles. 
Let $F' = A(G-\{v,x_1,x_2,x_3,y_1,y_2,y_3\})$. By the choice of $G$, 
$|F'| \geq \lceil (4(n-7)+3)/7 \rceil$. Now $G[F' + \{x_1,x_2,y_1,y_2\}]$ is an induced forest in $G$, implying 
$a(G) \geq |F'|+4 \geq \lceil (4n+3)/7 \rceil$, a contradiction.
\qed

\begin{nota}
\label{RRemark} 
Let $v \in V(G)$ and $U \subseteq N(v)$, and let $R\in R_{v,U}$.  We define $G*R=G-\{v,r\}$ if $R=\{r\}$, and $G*R=(G-v)/r_1r_2$ if $R=\{r_1,r_2\}$. 
For $F\subseteq G*R$, define $F\cdot R=G[F + r]$ if $R=\{r\}$. If $R=\{r_1,r_2\}$ and $r \in F$ where $r$ denotes the identification of $r_1$ and $r_2$, then define $F\cdot R=G[F- r + \{r_1,r_2\}]$.
\end{nota}

\begin{rem}
Let $v \in V(G)$ and $U \subseteq N(v)$. If $R=\{r_1,r_2\} \in R_{v,U}^2$ and $r$ denotes the identification of $r_1$ and $r_2$, then by Lemma \ref{basic_lemma} there exists $F = A(G * R)$ such that $r \in F$.
\end{rem}



\section{Structure around $2$-vertices}

The objective of this section is to prove the following lemma about neighbors of a 2-vertex in $G$. This will be used later for 
discharging rules.
 
\begin{lem}
\label{2summary} 
For each $x\in V_2$, there exist $v_5,v_5'\in V_{\ge 5}\cap N(x)$
or there exist $v_4\in V_{\leq 4} \cap N(x)$ and $v_6\in V_{\ge 6} \cap N(x)$. 
\end{lem}

\begin{rem}
Apply Lemma \ref{NoTwoRNeighbors} with $v = v_5$  and $U=\emptyset$, we have $R_{v_5,\emptyset}=\{\{x\}\}$ because any two elements in $R_{v_5,\emptyset}$ intersect.
Similarly, $R_{v_5',\emptyset}=\{\{x\}\}$ and $R_{v_4,\emptyset}=R_{v_6,\emptyset}=\{\{x\}\}$.
\end{rem}

\pf First, $e(V_2)=0$. For, suppose there exists $xy\in E(G)$ with $x,y \in V_2$.   
Let $z \in N(y) - \{x\}$ and $F' = A(G-\{x,y,z\})$. Then $|F'| \geq \lceil (4(n-3)+3)/7 \rceil$. 
Clearly,  $G[F' + \{x,y\}]$ is an induced forest in $G$; so $a(G)\ge |F'|+2 \geq \lceil (4n+3)/7 \rceil$, a contradiction.

Next, we claim that  for each $y\in V_2$, 
it is impossible that $y$ has one neighbor of degree $3$ and the other neighbor of degree at most $5$.
For otherwise, there exists a path $xyz$ in $G$ with $x \in V_3, y \in V_2, z \in V_{\leq 5}$.  
Let $N(x)-\{y\}=\{x_1,x_2\}$. 
Note that $\{x_1,x_2\} \subseteq N(z)$ since $G$ is a quadrangulation.
Then, $d(z) = 5$; otherwise, with $F' = A(G-\{x,y,z,x_1,x_2\})$,  
$G[F' + \{x,y,z\}]$ is an induced forest in $G$ showing that 
$a(G)\ge |F'| +3\geq \lceil (4(n-5)+3)/7 \rceil+3\ge  \lceil (4n+3)/7 \rceil,$
a contradiction. So let $N(z)=\{x_1,y,x_2,z_2,z_1\}$ such that $x_i$ and $z_i$ are cofacial for $i=1,2$. If $|N(x_1) \cap N(z_1)| \leq 2$, then 
let $F' = A((G-\{x,y,z,x_2\})/x_1z_1)$ with $w$ as the identification of $x_1$ and $z_1$; now  
$G[F'+\{x,y,z\}]$ (if $w \not\in F'$) or $G[F'-w+\{x,y,x_1,z_1\}]$ (if $w\in F'$) is an induced forest in $G$ showing that 
$a(G)\ge |F'|+3  \ge \lceil (4(n-5)+3)/7 \rceil\geq \lceil (4n+3)/7 \rceil,$
a contradiction. Thus, let $|N(x_1) \cap N(z_1)|\ge 3$. Then there exist $u\in N(x_1)\cap N(z_1)-\{z\}$ and a separation 
$(G_1,G_2)$ in $G$ such that $V(G_1 \cap G_2) = \{x_1,z_1,u\}$, $\{x,y,z,x_2,z_2\} \subseteq V(G_1)$, 
and  $N(x_1)\cap N(z_1) - \{z\} \subseteq  V(G_2)$. Let $F_1^{(1)} = A(G_1-\{x_1,z_1,x,y,z,x_2\})$ and $F_2^{(1)} = A(G_2-\{x_1,z_1\})$.
Then $G[F_1^{(1)} \cup F_2^{(1)} +\{x,y,z\}- (\{u\}\cap (F_1^{(1)} \triangle F_2^{(1)} ))]$ is an induced forest  in $G$, which implies  that 
$$a(G) \geq |F_1^{(1)}| + |F_2^{(1)}| +2\ge  \lceil (4(|G_1|-6)+3)/7 \rceil+\lceil (4(|G_2|-2)+3)/7 \rceil+2.$$
Now let $F_1^{(2)} = A(G_1-\{x_1,z_1,x,y,z,x_2,u\})$ and $F_2^{(2)} = A(G_2-\{x_1,z_1,u\})$. Then  
$G[F_1^{(2)} \cup F_2^{(2)} +\{x,y,z\}]$ is an induced forest in $G$, showing that 
$$a(G)\geq |F_1^{(2)}| + |F_2^{(2)}| +3\ge \lceil (4(|G_1|-7)+3)/7 \rceil +  \lceil (4(|G_2|-3)+3)/7\rceil +3.$$
By Lemma~\ref{ineq1}, we have $a(G)\ge  \lceil (4n+3)/7\rceil$, a contradiction.

Thus, to complete the proof of Lemma~\ref{2summary}, it suffices to show that for each $y\in V_2$, 
it is impossible that $y$ has one neighbor of degree $4$ and the other neighbor of degree at most $5$. 
For otherwise, there exists a path $xyz$ such that $x \in V_4, y \in V_2$ and $z \in V_{\le 5}$. Thus, $z\in V_4\cup V_5$ by the above claims.  
Let $N(x) = \{x_1,x_2,x_3,y\}$ and $N(z) = \{z_1,z_2,x_2,x_3,y\}$ if $z\in V_5$ or $N(z) = \{z_1,x_2,x_3,y\}$ if $z\in V_4$. 

\medskip

Case 1.  $N(x_2) \cap N(x_3) = \{x,z\}$ and either $|N(z_1) \cap N(z_2)| \leq 2$ or $z\in V_4$. 

Let $F' = A((G-\{x,y,z\})/\{x_2x_3,z_1z_2\})$ (when $z\in V_5$) and 
$F'=A((G-\{x,y,z,z_1\})/x_2x_3)$ (when $z\in V_4$). Let $x'$ (respectively, $z'$ when $z\in V_5$) denote the identification of $x_2$ and $x_3$ (respectively, $z_1$ and $z_2$). 
Let $z'=z_1$ if $z\in V_4$. By the choice of $G$, $|F'| \geq \lceil (4(n-5)+3)/7 \rceil$. It is easy to see that one of the following is an induced forest in $G$: 
$G[F'+\{x,y,z\}]$ (if $x',z' \not\in F'$), or  $G[(F'-z')+\{x,y,z_1,z_2\}]$ (if $x' \not\in F'$ and $z' \in F'$), 
or $G[(F'-x')+ \{x_2,x_3,y,z\}]$ (if $x' \in F'$ and $z' \not\in F'$), or 
$G[(F' -\{x',z'\})+ \{x_2,x_3,y,z_1,z_2\}]$ (if $x',z' \in F'$). Therefore, 
$a(G)\ge |F'|+3 \geq \lceil (4n+3)/7 \rceil$, a contradiction.

\medskip

Case 2.  $|N(x_2) \cap N(x_3)|\ge 3$ and either $|N(z_1) \cap N(z_2)| \leq 2$ or $z \in V_4$.

Then there exist  $w\in N(x_2)\cap N(x_3)$ and a separation 
$(G_1,G_2)$ in $G$ such that $V(G_1 \cap G_2) = \{w,x_2,x_3,x\}$, $\{y,z,z_1,z_2\}\subseteq  V(G_1)$, and $N(x_2)\cap N(x_3) - \{z\} \subseteq V(G_2)$. Let  $F_1^{(1)} = A((G_1-\{w,x_2,x_3,x,y,z\})/z_1z_2)$ (when $z\in V_5$) or  $F_1^{(1)} = A(G_1-\{w,x_2,x_3,x,y,z,z_1\})$ (when $z\in V_4$),   and 
let $F_2^{(2)} = A(G_2-\{w,x_2,x_3,x\})$. Let $z'$ denote the identification of $z_1$ and $z_2$. Then $G[F_1^{(1)} \cup F_2^{(1)} + \{x,y,z\}]$ (if $z \in V_4$ or if  $z \in V_5$ and $z' \not\in F_1^{(1)}$), or $G[F_1^{(1)} \cup F_2^{(1)} - z' + \{x,y,z_1,z_2\}]$ (if $z \in V_5$ and $z' \in F_1^{(1)}$) is an induced forest in $G$, showing that 
$$a(G)\ge |F_1^{(1)}| + |F_2^{(1)}| +3\ge  \lceil (4(|G_1|-7)+3)/7 \rceil +\lceil (4(|G_2|-4)+3)/7 \rceil +3.$$
Let $F_1^{(2)} = A((G_1-\{x,x_2,x_3,x,y,z\})/z_1z_2)$ (when $z\in V_5$) with $z'$ as the identification of $z_1$ and $z_2$, or  $F_1^{(2)} = A(G_1-\{x_2,x_3,x,y,z,z_1\})$ (when $z\in V_4$),   and 
let $F_2^{(2)} = A(G_2-\{x_2,x_3,x\})$. Then $G[F_1^{(2)} \cup F_2^{(2)} + \{x,y,z\} - (\{w\}\cap(F_1^{(2)} \triangle F_2^{(2)}))]$ (if $z \in V_4$ or $z \in V_5$ and $z' \not\in F_1^{(2)}$), or $G[F_1^{(2)} \cup F_2^{(2)} - z' + \{x,y,z_1,z_2\} - (\{w\}\cap(F_1^{(2)} \triangle F_2^{(2)}))]$ (if $z \in V_5$ and $z' \in F_1^{(2)}$)  is an induced forest in $G$, giving
$$a(G)\ge |F_1^{(2)}| + |F_2^{(2)}| +2\ge  \lceil (4(|G_1|-6)+3)/7 \rceil +\lceil (4(|G_2|-3)+3)/7 \rceil +2.$$
Hence, by Lemma~\ref{ineq2}(1) (with $k=1$, $a = |G_1|-6, a_1 = |G_2|-3, c = 3, L = \emptyset$), $a(G)\ge \lceil (4n+3)/7\rceil$, a contradiction.

\medskip
Case 3.  $N(x_2) \cap N(x_3) = \{x,z\}$ and $|N(z_1) \cap N(z_2)|\ge 3$. 

Then there exist $u\in N(z_1)\cap N(z_2)$ and a separation $(G_1,G_2)$ in $G$ such that $V(G_1 \cap G_2) = \{z_1,z_2,u\}$, 
$\{x,y,z,x_2,x_3\}\subseteq V(G_1)$, and $N(z_1)\cap N(z_2)-\{z\}\subseteq V(G_2)$. 
Let $F_1^{(1)} = A((G_1-\{z_1,z_2,x,y,z\})/x_2x_3)$ with $x'$ as the identification of $x_2$ and $x_3$, and $F_2^{(1)} = A(G_2-\{z_1,z_2\})$. 
Then $G[F_1^{(1)} \cup F_2^{(1)} + \{x,y,z\} - (\{u\}\cap(F_1^{(2)} \triangle F_2^{(2)}))]$ (if $x' \not\in F_1^{(1)}$) or $G[(F_1^{(1)}-x') \cup F_2^{(1)} +\{x_2,x_3,y,z\} - (\{u\}\cap(F_1^{(2)} \triangle F_2^{(2)}))]$ (if $x' \in F_1^{(1)}$) 
is an induced forest in $G$, which, by the choice of $G$, implies 
$$a(G)\ge |F_1^{(1)}| + |F_2^{(1)}| + 2\ge \lceil (4(|G_1|-6)+3)/7 \rceil+\lceil (4(|G_2|-2)+3)/7\rceil +2.$$
Let $F_1^{(2)} = A((G_1-\{u,z_1,z_2,x,y,z\})/x_2x_3)$ with $x'$ as the identification of $x_2$ and $x_3$, and $F_2^{(2)} = A(G_2-\{u,z_1,z_2\})$. 
Then $G[F_1^{(2)} \cup F_2^{(2)} + \{x,y,z\}]$ (if $z' \not\in F_1^{(2)}$) or $G[(F_1^{(2)}-x') \cup F_2^{(2)} +\{x_2,x_3,y,z\}]$ (if $z' \in F_1^{(2)}$)
is an induced forest in $G$. So by the choice of $G$. 
$$a(G)\ge |F_1^{(2)}| + |F_2^{(2)}| + 3 \ge \lceil (4(|G_1|-7)+3)/7 \rceil+ \lceil (4(|G_2|-3)+3)/7 \rceil+3.$$
So by Lemma~\ref{ineq2}(1) (with $k=1$, $a = |G_1|-6, a_1 = |G_2|-2, c = 3, L = \emptyset$), $a(G)\ge  \lceil (4n+3)/7\rceil$, a contradiction. 

\medskip

Case 4.  $|N(x_2) \cap N(x_3)|\ge 3$ and  $|N(z_1) \cap N(z_2)|\ge 3$. 

Then there exist $w\in N(x_2)\cap N(x_3)-\{x,z\}$, $u\in N(z_1)\cap N(z_2)-\{z\}$, and  subgraphs $G_1,G_2,G_3$ of $G$  such that 
$G_2$ is the maximal subgraph of $G$ contained in the closed region of the plane bounded by the cycle $wx_2xx_3w$ containing $N(x_2)\cap N(x_3) - \{z\}$,  
$G_3$ is the maximal subgraph of $G$ contained in  the closed region of the plane bounded by the cycle $zz_1uz_2z$ containing $N(z_1)\cap N(z_2) - \{z\}$, 
and $G_1$ is obtained from $G$ by removing $G_2-\{w,x,x_2,x_3\}$ and $G_3-\{u,z,z_1,z_2\}$. 

Define $A_i = \{u\}$ for $i=1,3$, $A_i=\emptyset$ for $i=2,4$, and $\overline{A_i} = \{u\} - A_i$. 
Define $W_i = \{w\}$ for $i=3,4$, $W_i=\emptyset$ for $i=1,2$ and $\overline{W_i} = \{w\} - W_i$. 
For $i\in [4]$, let $F_1^{(i)} = A(G_1-\{x,y,z,x_2,x_3,z_1,z_2\}-A_i-W_i)$ and $F_2^{(i)} = A(G_2-\{x_2,x_3,x\} - W_i)$ and $F_3^{(i)} = A(G_3-\{z_1,z_2\}-A_i)$. 
Then $|F_1^{(i)}|\ge \lceil (4(|G_1|-7-|A_i|-|W_i|)+3)/7 \rceil$, $|F_2^{(i)}| \ge \lceil (4(|G_2|-3-|W_i|)+3)/7 \rceil$, and 
$|F_3^{(i)}|\ge \lceil (4(|G_3|-2-|A_i|)+3)/7 \rceil+3$. Since 
$G[F_1^{(i)} \cup F_2^{(i)} \cup F_3^{(i)} +\{x,y,z\} - \{u,w\}\cap (F_1^{(i)}\triangle (F_2^{(i)}\cup F_3^{(i)}))]$ is an induced forest in $G$, 
$a(G) \ge  |F_1^{(i)}| + |F_2^{(i)}| + |F_3^{(i)}| + 3 - (1-|A_i|) - (1-|W_i|).$ Let $(n_1,n_2,n_3) := (4(|G_1|-7)+3,4(|G_2|-3)+3,4(|G_3|-2)+3)$. 
So by Lemma~\ref{ineq2}(4) (with $a = |G_1|-7, a_1 = |G_2| - 3, a_2 = |G_3| -2, c=3$), 
$$(n_1,n_2,n_3) \equiv (1,0,0), (4,0,4), (4,4,0), (0,4,4) \mod 7.$$

{\it Subcase 4.1} . $(n_1,n_2,n_3) \equiv (1,0,0)$ (resp. $(4,4,0)) \mod 7$

Let $W_5=\overline{W_6}=\{w\}$ and $W_6=\overline{W_5}=\emptyset$. Let $i=5$ if $(n_1,n_2,n_3) \equiv (1,0,0) \mod 7$ and $i=6$ if $(n_1,n_2,n_3) \equiv (4,4,0) \mod 7$. 
Let $F_1^{(i)} = A((G_1-\{x,y,z,x_2,x_3\} - W_i)/z_1z_2)$ 
with $z'$ as the identification of $z_1$ and $z_2$,  $F_2^{(i)} = A(G_2-\{x_2,x_3,x\} - W_i)$, and $F_3^{(i)} = A(G_3)$. 
By the choice of $G$,  $|F_1^{(i)}| \geq \lceil (4(|G_1|-6-|W_i|)+3)/7 \rceil$, $|F_2^{(i)}| \geq \lceil (4(|G_2|-3-|W_i|)+3)/7 \rceil$,
and $|F_3^{(i)}| \geq \lceil (4|G_3|+3)/7 \rceil$. 
Then $G[F_1^{(i)} \cup F_2^{(i)} \cup F_3^{(i)} + \{x,y,z\} - \{z_1,z_2,u,w\}\cap (F_1^{(i)}\triangle (F_2^{(i)}\cup F_3^{(i)})]$ (if $z' \not\in F_1^{(i)}$) or 
$G[(F_1^{(i)}-z') \cup F_2^{(i)} \cup F_3^{(i)} +\{x,y,z_1,z_2\}- \{u,w,z_1,z_2\}\cap ((F_1^{(i)}\cup \{z_1,z_2\}) \triangle (F_2^{(i)}\cup F_3^{(i)})]$  (if $z' \in F_1^{(i)}$)
is an induced forest in $G$, showing that 
$a(G)\ge |F_1^{(i)}| + |F_2^{(i)}| + |F_3^{(i)}| + 3 - 3 - | \overline{W_i} | \ge \lceil (4n+3)/7 \rceil$, a contradiction.

\medskip

{\it Subcase 4.2} .  $(n_1,n_2,n_3) \equiv (4,0,4) \mod 7$.

Let $F_1^{(7)} = A(G_1-\{x,y,z,x_2,x_3,z_1,w\})$, $F_2^{(7)} = A(G_2-\{x_2,x_3,x,w\})$, and $F_3^{(7)} = A(G_3 - \{z_1\})$. 
Then  $|F_1^{(7)}| \geq \lceil (4(|G_1|-7)+3)/7 \rceil$,  $|F_2^{(7)}| \geq \lceil (4(|G_2|-4)+3)/7 \rceil$, 
and $|F_3^{(7)}| \geq \lceil (4(|G_3|-1)+3)/7 \rceil$. Clearly,  $G[F_1^{(7)} \cup F_2^{(7)} \cup F_3^{(7)}+ \{x,y,z\} - \{u,z_2\}\cap (F_1^{(7)}\triangle (F_2^{(7)}\cup F_3^{(7)})]$ 
is an induced forest in $G$, showing that
$a(G)\ge |F_1^{(7)}| + |F_2^{(7)}| + |F_3^{(7)}| + 1 \ge \lceil (4n+3)/7 \rceil$, a contradiction.

\medskip

{\it Subcase 4.3} . $(n_1,n_2,n_3) \equiv (0,4,4) \mod 7$. 

Let $F_1^{(8)} = A(G_1-\{y,z,x_2,x_3,z_1\} + xz_2)$,  $F_2^{(8)} = A(G_2-\{x_2,x_3\})$, and $F_3^{(8)} =  A(G_3 - \{z_1\})$. 
Then  $|F_1^{(8)}| \geq \lceil (4(|G_1|-5)+3)/7 \rceil$, $|F_2^{(8)}| \geq \lceil (4(|G_2|-2)+3)/7 \rceil$, and 
$|F_3^{(8)}|\ge  \lceil (4(|G_3|-1)+3)/7 \rceil$. 
Now $G[F_1^{(8)} \cup F_2^{(8)} \cup F_3^{(8)}+ \{y,z\} - (\{u,w,x,z_2\}\cap ( F_1^{(8)}\triangle (F_2^{(8)}\cup F_3^{(8)}))]$
is an induced forest in $G$, which implies that $a(G)\ge |F_1^{(8)}| + |F_2^{(8)}| + |F_3^{(8)}| -2
\ge \lceil (4n+3)/7 \rceil$, a contradiction. This completes the proof of Lemma~\ref{2summary}.


\section{Structure around $3$-vertices}

In this section, we derive useful information about strutures around a 3-vertex. 

\begin{lem}
\label{No434Edge}
Let $x_1\in V_3$ and $N(x_1) = \{x, y_1,z_1\}$, with $y_1,z_1 \in V_4$, $x_2 \in N(x) \cap N(y_1) -\{x_1\}$ and $xx_1y_1x_2x$ be a facial cycle in $G$.  Then $z_1x_2\notin E(G)$. 
\end{lem}

\pf
For, suppose $z_1x_2\in E(G)$. Then $G$ has a separation $(G_1,G_2)$ such that $V(G_1 \cap G_2) = \{x_1,x_2,z_1\}$, 
$y_1 \in V(G_1)$, and $x \in V(G_2)$. For $i=1,2$, let  $F_i^{(1)} = A(G_i - \{z_1,x_1,x_2\})$; so 
$|F_i^{(1)}| \geq \lceil (4(|G_i|-3)+3)/7 \rceil$. Now $G[F_1^{(1)} \cup F_2^{(1)}+x_1]$ 
is an induced forest in $G$, giving $a(G)\ge |F_1^{(1)}| + |F_2^{(1)}| + 1$. 

Let $F_1^{(2)} = A(G_1 - \{z_1,x_1,x_2,y_1\})$ and $F_2^{(2)} = A(G_2 - \{z_1,x_1,x_2,x\})$. Then 
 $|F_i^{(2)}| \geq \lceil (4(|G_i|-4)+3)/7 \rceil$ for $i=1,2$. 
If $N(z_1) \cap V(G_1) - \{x_1,x_2\} \neq \emptyset$ and $N(z_1) \cap V(G_2) - \{x_1,x_2\} \neq \emptyset$, then  
$G[F_1^{(2)} \cup F_2^{(2)} +\{x_1,z_1\}]$ is an induced forest in $G$, giving $a(G)\ge |F_1^{(2)}| + |F_2^{(2)}| + 2$. 
Thus, by Lemma~\ref{ineq1}, $a(G)\ge \lceil (4n+3)/7 \rceil$, a contradiction. 

If $N(z_1) \cap V(G_1) - \{x_1,x_2\} = \emptyset$, then since $G$ is a quadrangulation,  $y_1,x_1,z_1,x_2$ are incident to a common face. This is a contradiction since $|N(y_1)| = 4$.
So $N(z_1) \cap V(G_2) - \{x_1,x_2\} = \emptyset$.
Then since $G$ is a quadrangulation, $x,x_1,z_1,x_2$ are incident to a common face.
This implies that $|N(x)| = 2$. 
So $G[F_1^{(2)} \cup F_2^{(2)} +\{x_1,x\}]$ is an induced forest in $G$, giving $a(G)\ge |F_1^{(2)}| + |F_2^{(2)}| + 2$. 
Thus, by Lemma~\ref{ineq1}, $a(G)\ge \lceil (4n+3)/7 \rceil$, a contradiction. 
\qed

\begin{lem}
\label{No333} 
$\Delta(G[V_{\leq 3}])\le 1$. 
\end{lem}

\pf First, we claim $e(V_2)=0$. For, suppose there exists $xy\in E(G)$ with $x,y \in V_2$.   
Let $z \in N(y) - \{x\}$ and $F' = A(G-\{x,y,z\})$. Then $|F'| \geq \lceil (4(n-3)+3)/7 \rceil$. 
Clearly,  $G[F' + \{x,y\}]$ is an induced forest in $G$; so $a(G)\ge |F'|+2 \geq \lceil (4n+3)/7 \rceil$, a contradiction.

Suppose $G[V_{\leq 3}]$ contains a path, say  $xyz$. 
By the claim above and Lemma \ref{NoTwoRNeighbors}, we may assume that $|N(y)| = |N(z)| = 3$. 
Suppose $|N(x)| = 2$. 
Since every face of $G$ has length $4$, $x$ and $z$ have a common neighbor, say $s$.
Let $N(x) = \{s,y\}$, $N(y) = \{y_1,x,z\}$ and $N(z) = \{z_1,s,y\}$. 
Let $F' = A(G-\{x,y,z,s,y_1\})$. 
Then by the choice of $G$, $|F'|\ge  \lceil (4(n-5)+3)/7 \rceil$. Now $G[F' + \{x,y,z\}]$ is an induced forest in $G$ and, hence,
$a(G)\ge |F'|+3  \ge \lceil (4n+3)/7 \rceil,$ a contradiction. So $|N(x)| = 3$.

Since every face of $G$ has length $4$, $x$ and $z$ have a common neighbor, say $s$.
Let $N(x) = \{x_1,s,y\}$, $N(y) = \{y_1,x,z\}$ and $N(z) = \{z_1,s,y\}$. If $x_1 = z_1$, let $F' = A(G-\{x,y,z,s,x_1\})$. Then 
by the choice of $G$, $|F'|\ge  \lceil (4(n-5)+3)/7 \rceil$. Now $G[F' + \{x,y,z\}]$ is an induced forest in $G$ and, hence,
$a(G)\ge |F'|+3  \ge \lceil (4n+3)/7 \rceil,$ a contradiction. So $x_1 \neq z_1$.
 
If $N(x_1) \cap N(z_1) = \{y_1\}$, let $F' = A((G-\{x,y,z,s\})/x_1z_1)$ with $x'$ as the identification of $x_1$ and $z_1$. 
Then $|F'| \geq \lceil (4(n-5)+3)/7 \rceil$. Now $G[F'+ \{x,y,z\}]$ (if $x' \not\in F'$) or $G[(F'-x') +\{x,z,x_1,z_1\}] $ 
(if $x'\in F'$) is an induced forest in $G$. So $a(G)\ge |F'|+3 \geq \lceil (4n+3)/7 \rceil$, a contradiction.

So $|N(x_1) \cap N(z_1)|\ge 2$. 
Then there exist $w\in N(x_1)\cap N(z_1)-\{y_1\}$ and a separation $(G_1,G_2)$ in $G$ such that 
$V(G_1 \cap G_2) = \{w,x_1,y_1,z_1\}$, $\{x,y,z,s\}\subseteq  V(G_1)$, and $N(x_1)\cap N(z_1)\subseteq V(G_2)$. 
Let $W_1=\overline{W_2}=\{w\}$ and $\overline{W_1}=W_2=\emptyset$. 
For $i=1,2$, let $F_1^{(i)} = A(G_1-\{s,x,y,z,x_1,y_1,z_1\} - W_i)$ and $F_2^{(i)} = A(G_2-\{x_1,z_1\} - W_i)$. 
Then $|F_1^{(i)}| \geq \lceil (4(|G_1|-7-|W_i|)+3)/7 \rceil$ and $|F_2^{(i)}| \geq \lceil (4(|G_2|-2-|W_i|)+3)/7 \rceil$. 
Now $G[F_1^{(i)} \cup F_2^{(i)} + \{x,y,z\} - (\{w\}\cap (F_1^{(i)} \triangle F_2^{(i)}))]$ is an induced forest in $G$, giving 
$a(G)\ge |F_1^{(i)}| + |F_2^{(i)}| + 3 - |\overline{W_i}|$. By Lemma~\ref{ineq2}(1) (with $k=1$, $a = |G_1| - 7, a_1 = |G_2| - 2, L = \emptyset, c = 3$), $a(G)\ge \lceil (4n+3)/7 \rceil$, 
a contradiction.  \qed

\begin{lem}
\label{No3RR}
Let $x\in V_3$. If $y \in N(x)$ and $R_{y,\{x\}} \neq \emptyset$ then for any $z\in N(x)-\{y\}$, $R_{z,\{x\}}=\emptyset$. 
\end{lem}

\pf For otherwise, suppose $z\in N(x)-\{y\}$ and $R_{z,\{x\}} \neq \emptyset$. Let $R_1\in R_{y,\{x\}}$ and $R_2 \in R_{z,\{x\}}$. 

If $|R_1| = 1$ or $|R_2| = 1$, let $F' = A(((G - \{x,y,z\})*R_1)*R_2)$. Then $|F'| \geq \lceil (4(n-5)+3)/7 \rceil$. Now $G[((F' + x) \cdot R_1) \cdot R_2]$ is an induced forest in $G$, showing $a(G) \geq |F'| + 3 \geq \lceil (4n+3)/7 \rceil$, a contradiction.

So $|R_1| = |R_2| = 2$, let $R_1 = \{r_1,r_2\}$ and $yr_1y'r_2y$ bound a $4$-face. 
Suppose $y' = z$. 
Let $F' = A(G - \{x,y,z,r_1,r_2\})$.
Then $|F'| \geq \lceil (4(n-5)+3)/7 \rceil$.
Now $G[F' + \{x,r_1,r_2\} ]$ is an induced forest in $G$, showing $a(G) \geq |F'| + 3 \geq \lceil (4n+3)/7 \rceil$, a contradiction.

Now, we may assume $y' \neq z$. 
Suppose $R_1 \cap R_2 \neq \emptyset$. 
Without loss of generality, let $R_2 = \{r_2,r_3\}$. 
Since $G$ is a quadrangulation, $zr_2y'r_3z$ bounds a $4$-face. 
Let $F'' = A(G - \{x,y,z,r_1,r_2,r_3,y'\})$.
Then $|F''| \geq \lceil (4(n-7)+3)/7 \rceil$.
Now $G[F'' + \{x,r_1,r_2,r_3\} ]$ is an induced forest in $G$, showing $a(G) \geq |F''| + 4 \geq \lceil (4n+3)/7 \rceil$, a contradiction.

Finally, we may assume $R_1 \cap R_2 = \emptyset$.
Let $F''' = A( (G - \{x,y,z,r_1,r_2,y'\})*R_2 )$. Then $|F'''| \geq \lceil (4(n-7)+3)/7 \rceil$. Now $G[(F''' + \{x,r_1,r_2\}) \cdot R_2]$ is an induced forest in $G$, showing $a(G) \geq |F'''| + 4 \geq \lceil (4n+3)/7 \rceil$, a contradiction.

\qed

\begin{lem}
\label{No4-3R}
Let $x\in V_3$. If $y\in N(x)\cap V_{\le 4}$ then for any $z\in N(x)-\{y\}$, $R_{z,\{x\}}=\emptyset$. 
\end{lem}

\pf Let $N(x) = \{u,y,z\}$, $y\in V_{\le 4}$ and $R\in R_{z,\{x\}}$. Let $vyxzv$ be a facial cycle, 
$N(y) = \{y_1,x,v\}$ if $y\in V_3$ and $ N(y) = \{y_1,y_2,x,v\}$ if $y\in V_4$. 
In the proof below, we assume $y\in V_4$ as for $y\in V_3$. 
We simply delete $y_1$ instead of identifying $y_1$ and $y_2$. 
Define $W_i = \{v\}$ for $i=1,3,5,8$ and $W_i=\emptyset$ if $i=2,4,6,7$, and let $\overline{W_i} = \{v\} - W_i$ for $i\in [8]$.

Suppose $R = \{y_2\}$. This implies that $zy_2 \in E(G)$ and $|N(y_2)| = 2$.
Since $G$ is a plane graph, $uv \not \in E(G)$. 
Let $F = A(G - \{x,y,z,y_1,y_2\} + uv)$. 
By the choice of $G$, $|F| \geq \lceil (4(n-5)+3)/7 \rceil$.
Then $G[F + \{x,y,y_2\}]$ is an induced forest in $G$. 
So $a(G)\ge |F|+3 \geq \lceil (4n+3)/7 \rceil$, a contradiction.
So $R \neq \{y_2\}$.
Similarly,  $R \neq \{y_1\}$.

\medskip

Case 1. $|N(y_1) \cap N(y_2)| \leq 2$ and $uv \not\in E(G)$. 

Let $F' = A((G-\{x,y,z\}) * R)/y_1y_2 + uv)$ with $y'$ as the identification of $y_1$ and $y_2$. 
By the choice of $G$,  $|F'| \geq \lceil (4(n-5)+3)/7 \rceil$.
Then $G[(F' + \{x,y\})\cdot R]$ (if $y' \not\in F'$) or $G[(F'-\{y'\}+ \{x,y_1,y_2\}) \cdot R]$ (if $y' \in F'$)
is an induced forest in $G$. So $a(G)\ge |F'|+3 \geq \lceil (4n+3)/7 \rceil$, a contradiction.

\medskip

Case 2. $|N(y_1) \cap N(y_2)| \leq 2$ and $uv \in E(G)$. 

Then $G$ has a separation $(G_1,G_2)$ such that $V(G_1 \cap G_2) = \{u,v,x\}$, $\{y,y_1,y_2\}\subseteq V(G_1)$, $z \in V(G_2)$. 
For $i=1,2$, let $F_1^{(i)} = A((G_1-\{u,x,y\} - W_i)/y_1y_2))$ with $y'$ as the identification of $y_1$ and $y_2$, 
 and $F_2^{(i)} = A((G_2-\{u,x,z\}- W_i)* R )$. Then  $|F_j^{(i)}| \geq \lceil (4(|G_j|-4-|W_i|)+3)/7 \rceil$ for $j=1,2$. 
Now $G[(F_1^{(i)} \cup F_2^{(i)} + \{x,y\})\cdot  R - (\{v\}\cap (F_1^{(i)} \triangle F_2^{(i)}))]$ (if $y' \not\in F_1^{(i)})$ or  
$G[((F_1^{(i)}-y') \cup F_2^{(i)}+ \{x,y_1,y_2\})\cdot R - (\{v\}\cap (F_1^{(i)} \triangle F_2^{(i)}))]$ is an induced forest in $G$, showing 
that $a(G)\ge |F_1^{(i)}| + |F_2^{(i)}| +3 - |\overline{W_i}|$. By Lemma~\ref{ineq2}(1) (with $k=1, a = |G_1| - 4, a_1 = |G_2|-4, L = \emptyset, c=3$), $a(G)\ge  \lceil (4n+3)/7 \rceil$, a contradiction.

\medskip
Case 3. $|N(y_1) \cap N(y_2)|\ge 3$ and $uv \not\in E(G)$.  

There exist $w\in N(y_1) \cap N(y_2)$ and a separation $(G_1,G_2)$ in $G$ such that $V(G_1 \cap G_2) = \{y_1,y_2,w\}$, $\{x,y,z,u,v\} \subseteq V(G_1)$, 
and $N(y_1)\cap N(y_2)-\{y\}\subseteq  V(G_2)$. Define $A_i = \{w\}$ if $i=1,3,4$ and  $A_i=\emptyset$ if $i=2,5,6$,  and let $\overline{A_i} = \{w\} - A_i$. 
 For $i=1,2$, let $F_1^{(i)} = A((G_1-\{x,y,z,y_1,y_2\} - A_i) * R  + uv)$, and $F_2^{(i)} = A(G_2-\{y_1,y_2\} - A_i)$. 
Then $|F_1^{(i)}| \geq \lceil (4(|G_1|-6 - |A_i|)+3)/7 \rceil$, and $|F_2^{(i)}| \geq \lceil (4(|G_2|-2 - |A_i|)+3)/7 \rceil$. 
Now $G[(F_1^{(i)}\cup  F_2^{(i)} +\{x,y\})\cdot R - (\{w\}\cap (F_1^{(i)} \triangle F_2^{(i)}))]$ is an induced forest in $G$, implying 
$a(G)\ge |F_1^{(1)}| + |F_2^{(1)}| + 3 - |\overline{A_i}|$. So by Lemma~\ref{ineq2}(1) (with $k=1, a = |G_1| - 6, a_1 = |G_2|-2, L = \emptyset, c=3$),  $a(G)\ge  \lceil (4n+3)/7 \rceil$, a contradiction.

\medskip

Case 4. $|N(y_1) \cap N(y_2)|\ge 3$ and $uv \in E(G)$. 

There exist $w\in N(y_1) \cap N(y_2)$ and subgraphs $G_1,G_2,G_3$ of $G$ such that $G_2$ is the maximal subgraph of $G$ 
contained in the closed region of the plane bounded by $uxzvu$ and containing $R$, $G_3$ is obtained by deleting $y$ from the maximal subgraph of $G$ contained in the closed 
region bounded by $y_1yy_2wy_1$ and containing $N(y_1)\cap N(y_2)$, and $G_1$ is obtained from $G$ by removing $G_2-\{u,v,x\}$ and $G_3-\{w,y_1,y_2\}$.  
For $i=3,4,5,6$, let $F_1^{(i)} = A(G_1-\{x,u,y,y_1,y_2\}-A_i-W_i)$, $F_2^{(i)} = A((G_2-\{u,x,z\}-W_i) * R)$, and $F_3^{(i)} = A(G_3-\{y_1,y_2\}-A_1)$. 
Then $|F_1^{(i)}| \geq \lceil (4(|G_1|-5-|W_i|-|A_i|)+3)/7 \rceil$, $|F_2^{(i)}| \geq \lceil (4(|G_2|-4-|W_i|)+3)/7 \rceil$, and 
$|F_3^{(i)}| \geq \lceil (4(|G_3|-2-|A_i|)+3)/7 \rceil$. Now $G[(F_1^{(i)} \cup F_2^{(i)} \cup F_3^{(i)} + \{x,y\})\cdot R - \{v,w\}\cap (F_1^{(i)} \triangle (F_2^{(i)} \cup F_3^{(i)}))]$
is an induced forest in $G$, showing $a(G)\ge |F_1^{(i)}| + |F_2^{(i)}| + |F_3^{(i)}| + 3 - |\overline{W_i}| - |\overline{A_i}|$. 
Let $(n_1,n_2,n_3):=(4(|G_1|-5)+3,4(|G_2|-4)+3,4(|G_3|-2)+3)$. 
By Lemma~\ref{ineq2}(4) (with $a = |G_1|-5, a_1 = |G_2|-4, a_2 = |G_3|-2$),  $ (n_1,n_2,n_3) \equiv (1,0,0), (4,0,4), (4,4,0), (0,4,4) \mod 7$. 

\medskip

{\it Subcase 4.1} . $(n_1,n_2,n_3) \equiv (1,0,0)$ (resp. $(4,4,0) \mod 7$). 

For $i=7$ (resp. $i=8$), let $F_1^{(i)} = A(G_1-\{u,x,y\} - W_i)$, $F_2^{(i)} = A(G_2-\{u,x\}-W_i)$ and $F_3^{(i)} = A(G_3)$. 
Then $|F_1^{(i)}| \geq \lceil (4(|G_1|-3)+3)/7 \rceil$, $|F_2^{(i)}| \geq \lceil (4(|G_2|-2)+3)/7 \rceil$ and $|F_3^{(i)}| \geq \lceil (4|G_3|+3)/7 \rceil$. 
Now $G[F_1^{(i)} \cup F_2^{(i)} \cup F_3^{(i)}+\{x\} -\overline{W_i} - \{y_1,y_2,w\} \cap (F_1^{(i)} \triangle F_3^{(i)} )]$ is an induced forest in $G$, 
showing $a(G)\ge |F_1^{(i)}| + |F_2^{(i)}| + |F_3^{(i)}| + 1 - 3 - |\overline{W_i}| \geq \lceil (4n+3)/7 \rceil$, a contradiction.
 
\medskip

{\it Subcase 4.2} . $(n_1,n_2,n_3) \equiv (4,0,4),(0,4,4) \mod 7$. 

Let $F_1^{(9)} = A(G_1-\{u,x,y,y_1,v\})$, $F_2^{(9)} = A(G_2-\{u,x,v\})$ and $F_3^{(9)} = A(G_3 - \{y_1\})$. 
Then $|F_1^{(9)}| \geq \lceil (4(|G_1|-5)+3)/7 \rceil$, $|F_2^{(9)}| \geq \lceil (4(|G_2|-3)+3)/7 \rceil$, and 
$|F_3^{(9)}| \geq \lceil (4(|G_3|-1)+3)/7 \rceil$. Now $G[F_1^{(9)} \cup F_2^{(9)} \cup F_3^{(9)}+ \{x,y\} - \{y_2,w\}\cap (F_1^{(9)} \triangle F_3^{(9)}) ]$ is 
an induced forest in $G$, showing $a(G)\ge |F_1^{(9)}| + |F_2^{(9)}| + |F_3^{(9)}| + 2 - 2 \geq \lceil (4n+3)/7 \rceil$, a contradiction. 
\qed

\begin{lem}
\label{No34-4-4-}
For each $x\in V_3$, $N(x)\not\subseteq V_{\le 4}$. 
\end{lem}

\pf Let $x\in V_3$ with $N(x) = \{w,y,z\}\subseteq V_{\leq 4}$. 
By Lemma \ref{No333}, $|N(x)\cap V_{\leq 3}|\le 1$; so let $N(z)=\{x,z_1,z_2,w_1\}$ and $N(w)=\{x,w_1,w_2,y_1\}$. 
Suppose $y \in V_2$. Let $N(y) = \{x,y_1\}$. 
Since $G$ is a quadrangulation, we may assume $z_1 = y_1$. 
Let $F = A(G - \{x,y,z,w,y_1,w_1,z_2\})$.
Then $|F| \geq \lceil (4(n-7)+3)/7 \rceil$.
Therefore, $G[F + \{x,y,z,w\}]$ is an induced forest in $G$, showing that $a(G)\ge  
|F|+4 \geq \lceil (4n+3)/7 \rceil$, a contradiction.

Now let $N(y)=\{x,y_1,z_1\}$ if $y\in V_3$  and $N(y)=\{x,y_1,y_2,z_1\}$ if $y\in V_4$. 
In the argument to follow, we treat the case  $y\in V_4$, as the proof for $y\in V_3$ is the same by replacing identification of $y_1$ and $y_2$ with 
the deletion of $y_1$.

\medskip

Case 1. $|N(y_1) \cap N(y_2)| \leq 2$, $|N(z_1) \cap N(z_2)| \leq 2$ and $|N(w_1) \cap N(w_2)| \leq 2$.

Let $F' = A(G-\{x,y,z,w\}/\{y_1y_2,z_1z_2,w_1w_2\})$ with $y',z',w'$ as the identifications of  $y_1$ and $y_2$, $z_1$ and $z_2$, and $w_1$ and $w_2$, respectively. 
Then $|F'| \geq \lceil (4(n-7)+3)/7 \rceil$. Let $F=F' + \{x,y,z,w\}$  if $w',y',z' \not\in F'$, and otherwise, let 
$F$ be obtained from $F' + \{x,y,z,w\}$ by deleting  $w,w'$ (respectively, $y,y'$, $z,z'$) adding  $\{w_1,w_2\}$ (respectively, $\{y_1,y_2\}$, $\{z_1,z_2\}$) 
if $w'\in F'$ (respectively, $y'\in F'$, $z'\in F'$). 
Then $G[F]$ is an induced forest in $G$, showing that $a(G)\ge  
|F'|+4 \geq \lceil (4n+3)/7 \rceil$, a contradiction.

\medskip

Case 2. Exactly one of $|N(y_1) \cap N(y_2)|$, $|N(z_1) \cap N(z_2)|$, $|N(w_1) \cap N(w_2)|$ is greater than $2$.

By symmetry, assume $|N(z_1) \cap N(z_2)|\ge 3$. Then there exist $z'\in N(z_1)\cap N(z_2)$ and a separation $(G_1,G_2)$ in $G$ 
such that $V(G_1 \cap G_2) = \{z_1,z_2,z'\}$, $\{x,y,z,w,y_1,y_2,w_1,w_2\}\subseteq V(G_1)$, $N(z_1)\cap N(z_2)-\{z\}\subseteq V(G_2)$. 
Define $A_i = \{z'\}$ for $i=1,5$ or $A_i=\emptyset$ for $i=2,6$,  and let $\overline{A_i} = \{z'\} - A_i$. 
For $i=1,2$, let $F_1^{(i)} = A((G_1-\{x,y,z,w,z_1,z_2\}-A_i)/\{y_1y_2,w_1w_2\})$ with $y',w'$ as the identifications of $y_1$ and $y_2$, $w_1$ and $w_2$, 
respectively,  and let $F_2^{(i)} = A(G_2-\{z_1,z_2\}-A_i)$. Then 
$|F_1^{(i)}| \geq \lceil (4(|G_1|-8 - |A_i|)+3)/7 \rceil$ and $|F_2^{(i)}| \geq \lceil (4(|G_2|-2-|A_i|)+3)/7 \rceil$. 
Let $F^{(i)}=F_1^{(i)} \cup F_2^{(i)} +\{x,y,z,w\} - (\{z'\} \cap (F_1^{(i)} \triangle F_2^{(i)}))$ if $w',y' \not\in F_1^{(i)}$, and otherwise, let  
$F^{(i)}$ be obtained from $F_1^{(i)} \cup F_2^{(i)} +\{x,y,z,w\} -  (\{z'\} \cap (F_1^{(i)} \triangle F_2^{(i)}))$ by deleting $\{y,y'\}$ (respectively, $\{w,w'\}$) and adding $\{y_1,y_2\}$ (respectively, 
$\{w_1,w_2\}$) when $y'\in F_1^{(i)}$ (respectively, $w'\in  F_1^{(i)}$).  Then $G[F^{(i)}]$ is an induced forest in $G$, giving $a(G)\ge 
|F_1^{(i)}| + |F_2^{(i)}| + 4 - |\overline{A_i}|$. By Lemma~\ref{ineq2}(2)
(with $a=|G_1|-8, a_1=|G_2|-2, c=4$)
,  $(4(|G_1|-8)+3,4(|G_2|-2)+3) \equiv (4,0), (0,4) \mod 7$. 

\medskip

{\it Subcase 2.1}. $(4(|G_1|-8)+3,4(|G_2|-2)+3) \equiv (4,0) \mod 7$. 

Let $F_1^{(3)} = A((G_1- \{x,y, z,w\}) /\{y_1y_2,w_1w_2,z_1z_2\})$ with $y',w',z''$ as the identification of $y_1$ and $y_2$, $w_1$ and $w_2$, and $z_1$ and $z_2$, 
respectively,  and let $F_2^{(3)} = A(G_2)$. Then  $|F_1^{(3)}| \geq \lceil (4(|G_1|-7)+3)/7 \rceil$ and $|F_2^{(3)}| \geq \lceil (4|G_2|+3)/7 \rceil$. 
Let $F^{(3)}= \overline{F_1}^{(3)} \cup F_2^{(3)}  - (\{z',z_1,z_2\} \cap (\overline{F_1}^{(3)} \triangle F_2^{(3)}))$ where $\overline{F_1}^{(3)} = {F_1}^{(3)}+\{x,y,z,w\}$ if  $w',y',z'' \not\in F_1^{(3)}$; otherwise, let $\overline{F_1}^{(3)}$ be obtained from 
${F_1}^{(3)}+\{x,y,z,w\}$ by deleting $y,y'$ (respectviely, $w,w'$, $z,z''$) and adding $\{y_1,y_2\}$ (respectively, 
$\{w_1,w_2\}$, $\{z_1,z_2\}$) when $y'\in F_1^{(3)}$ (respectively, $w'\in  F_1^{(3)}$, $z''\in F_1^{(3)}$). 
Therefore, $G[F^{(3)}]$ is an induced forest in $G$, showing that $a(G)\ge |F_1^{(3)}| + |F_2^{(3)}| + 4 - 3 \geq \lceil (4n+3)/7 \rceil$, a contradiction.

\medskip
{\it Subcase} 2.2. $(4(|G_1|-8)+3,4(|G_2|-2)+3) \equiv (0,4) \mod 7$. 

If $wz_2 \not\in E(G)$, then let $F_1^{(4)} = A((G_1-\{x,y,z,z_1,w_1\})/y_1y_2 + wz_2)$ with $y'$ as the identification of  $y_1$ and $y_2$, 
and $F_2^{(4)} = A(G_2 - \{z_1\})$. Then $|F_1^{(4)}| \geq \lceil (4(|G_1|-6)+3)/7 \rceil$ and $|F_2^{(4)}| \geq \lceil (4(|G_2|-1)+3)/7 \rceil$. 
Now $G[F_1^{(4)} \cup F_2^{(4)} + \{x,y,z\} - (\{z',z_2\} \cap (F_1^{(4)} \triangle F_2^{(4)}))]$ (if $y' \not\in F_1^{(4)}$) or 
$G[(F_1^{(4)} - y') \cup F_2^{(4)} + \{x,y_1,y_2,z\} - (\{z',z_2\} \cap (F_1^{(4)} \triangle F_2^{(4)}))]$ (if $y' \in F_1^{(4)}$) is an induced forest in $G$, 
giving $a(G)\ge |F_1^{(4)}| + |F_2^{(4)}| + 3 - 2 \geq \lceil (4n+3)/7 \rceil$, a contradiction. 

So $wz_2 \in E(G)$. Then there exist subgraphs 
$G_1',G_2',G_3'$ of $G$ such that 
$G_2'=G_2$,
$G_3'$ is the maximal subgraph of $G$ contained in the closed region of the plane bounded by the cycle $wxzz_2w$ and containing $N(w)\cap N(z) - \{x\}$,
and $G_1'$ is obtained from $G$ by removing $G_2'-\{z_1,z_2,z'\}$ and $G_3'-\{w,z,z_2\}$.
For $i=5,6$, let $F_1^{(i)} = A(G_1'-\{w,x,z,z_1,z_2\}-A_i)$,
 $F_2^{(i)} = A(G_2' - \{z_1,z_2\}-A_i)$,
 and  $F_3^{(i)} = A(G_3' - \{w,z,z_2\})$. Then  
$|F_1^{(i)}| \geq \lceil (4(|G_1'|-5 - |A_i|)+3)/7 \rceil$, $|F_2^{(i)}| \geq \lceil (4(|G_2'|-2 - |A_i|)+3)/7 \rceil$ and $|F_3^{(i)}| \geq \lceil (4(|G_3'|-3)+3)/7 \rceil$. 
So $G[F_1^{(i)} \cup F_2^{(i)} \cup F_3^{(i)} + \{x,z\} - (\{z'\}\cap (F_1^{(i)} \triangle F_2^{(i)}))]$ is an induced forest in $G$, giving 
$a(G)\ge |F_1^{(i)}| + |F_2^{(i)}| + |F_3^{(i)}| + 2 - |\overline{A_1}|$. By Lemma~\ref{ineq2}(1) (with $k=1, L = \{1\}, a = |G_1'|-5, a_1 = |G_2'| - 2, b_1 = |G_3'| - 3, c=2$), $a(G)\ge  \lceil (4n+3)/7 \rceil$, a contradiction.


\medskip
Thus, by symmetry, we have  
Case 3. At least two of $|N(y_1) \cap N(y_2)|$, $|N(z_1) \cap N(z_2)|$ and $|N(w_1) \cap N(w_2)|$ are greater than $2$, and at least
 two of $|N(z_1) \cap N(y_2)|$, $|N(w_1) \cap N(z_2)|$ and $|N(w_2) \cap N(y_1)|$ are greater than $2$. 

First, suppose $|N(y_1) \cap N(y_2)|>2$, $|N(y_1) \cap N(w_2)|> 2$, $|N(w_1) \cap N(w_2)|>2$ and $|N(w_1) \cap N(z_2)|>2$. 
Then there exist $y'\in N(y_1) \cap N(y_2)-\{y\}$, $w'\in N(y_1) \cap N(w_2) -\{w\}$, $w''\in N(w_1) \cap N(w_2)-\{w\}$, 
$z'\in N(w_1) \cap N(z_2) -\{z\}$, and subgraphs $G_1,G_2,G_3,G_4,G_5$ of $G$ such that 
$G_2$ is the maximal subgraph of $G$ contained in the closed region of the plane bounded by the cycle $yy_1y'y_2y$ and containing $N(y_1)\cap N(y_2) - \{y\}$, $G_3$ is the maximal subgraph of $G$ contained in the closed region of the plane bounded by the cycle $wy_1w'w_2w$ and containing $N(y_1)\cap N(w_2) - \{w\}$, $G_4$ is the maximal subgraph of $G$ contained in the closed region of the plane bounded by the cycle $ww_1w''w_2w$ and containing $N(w_1)\cap N(w_2) - \{w\}$, $G_5$ is the maximal subgraph of $G$ contained in the closed region of the plane bounded by the cycle $zw_1w''z_2z$ and containing $N(z_2)\cap N(w_1) - \{z\}$, and $G_1$ is obtained from $G$ by removing $G_2-\{y_1,y_2,y'\}$, $G_3-\{y_1,w_2,w'\}$, $G_4-\{w_1,w_2,w''\}$ and $G_5-\{w_1,z_2,z'\}$. 
Let $A_1 \subseteq \{y'\}$, $B_1 \subseteq \{w'\}$, $C_1 \subseteq \{w''\}$, $D_1 \subseteq \{z'\}$. 
Let $\overline{A_1} = \{y'\} - A_1$, $\overline{B_1} = \{w'\} - B_1$, $\overline{C_1} = \{w''\} - C_1$, $\overline{D_1} = \{z'\} - D_1$. 
For all choices of $A_1,B_1,C_1,D_1$, let $F_1^{(i)} = A(G_1-\{w,y,x,z,z_1,z_2,y_1,y_2,w_1,w_2\} - A_1 - B_1 - C_1 - D_1)$, $F_2^{(i)} = A(G_2 - \{y_1,y_2\} - A_1)$, 
$F_3^{(i)} = A(G_3 - \{y_1,w_2\} - B_1)$, $F_4^{(i)} = A(G_4 - \{w_1,w_2\} - C_1)$, and $F_5^{(i)} = A(G_5 - \{w_1,z_2\} - D_1)$. 
Then $|F_1^{(i)}| \geq \lceil (4(|G_1|-10-|A_1|-|B_1|-|C_1|-|D_1|)+3)/7 \rceil$, $|F_2^{(i)}| \geq \lceil (4(|G_2|-2-|A_1|)+3)/7 \rceil$, 
$|F_3^{(i)}| \geq \lceil (4(|G_3|-2-|B_1|)+3)/7 \rceil$, $|F_4^{(i)}| \geq \lceil (4(|G_4|-2-|C_1|)+3)/7 \rceil$, and $|F_5^{(i)}| \geq \lceil (4(|G_5|-2-|D_1|)+3)/7 \rceil$. 
Now $G[F_1^{(i)} \cup F_2^{(i)} \cup F_3^{(i)} \cup F_4^{(i)} \cup F_5^{(i)} + \{w,x,y,z\}  - (\{y'\} \cap (F_1^{(i)} \triangle F_2^{(i)})) - (\{w'\} \cap (F_1^{(i)} \triangle F_3^{(i)})) - (\{w''\} \cap (F_1^{(i)} \triangle F_4^{(i)})) - (\{z'\} \cap (F_1^{(i)} \triangle F_5^{(i)})) ]$ 
is an induced forest in $G$. Hence, $a(G)\ge |F_1^{(i)}| + |F_2^{(i)}| + |F_3^{(i)}| + |F_4^{(i)}| + |F_5^{(i)}| + 4 - |\overline{A_1}| - |\overline{B_1}| - |\overline{C_1}| - |\overline{D_1}|$. 
By Lemma~\ref{ineq2}(1) 
(with $k=4, a = |G_1|-10, a_j = |G_{j+1}|-2$ for $j=1,2,3,4$, $L = \emptyset, c=4$), 
$a(G)\ge  \lceil (4n+3)/7 \rceil$, a contradiction. 

Thus, by symmetry, we may assume that $|N(y_1) \cap N(y_2)|>2$, $|N(y_1) \cap N(w_2)|> 2$, $|N(z_1) \cap N(z_2)|>2$ and $|N(w_1) \cap N(z_2)|>2$. 
Then there exist $y'\in N(y_1) \cap N(y_2)-\{y\}$, $w'\in N(y_1) \cap N(w_2) -\{w\}$, $z'\in N(z_1) \cap N(z_2)-\{w\}$, 
$z''\in N(w_1) \cap N(z_2) -\{z\}$, and subgraphs $G_1,G_2,G_3,G_4,G_5$ of $G$ such that 
$G_2$ is the maximal subgraph of $G$ contained in the closed region of the plane bounded by the cycle $yy_1y'y_2y$ and containing $N(y_1)\cap N(y_2) - \{y\}$, $G_3$ is the maximal subgraph of $G$ contained in the closed region of the plane bounded by the cycle $wy_1w'w_2w$ and containing $N(y_1)\cap N(w_2) - \{w\}$, $G_4$ is the maximal subgraph of $G$ contained in the closed region of the plane bounded by the cycle $zz_1z'z_2z$ and containing $N(z_1)\cap N(z_2) - \{z\}$, $G_5$ is the maximal subgraph of $G$ contained in the closed region of the plane bounded by the cycle $zw_1z''z_2z$ and containing $N(z_2)\cap N(w_1) - \{z\}$, and $G_1$ is obtained from $G$ by removing $G_2-\{y_1,y_2,y'\}$, $G_3-\{y_1,w_2,w'\}$, $G_4-\{z_1,z_2,z'\}$ and $G_5-\{w_1,z_2,z''\}$. 
Let $A_1 \subseteq  \{y'\}$, $B_1 \subseteq  \{w'\}$, $C_1 \subseteq  \{z'\}$, $D_1 \subseteq  \{z''\}$. 
Let $\overline{A_1} = \{y'\} - A_1$, $\overline{B_1} = \{w'\} - B_1$, $\overline{C_1} = \{z'\} - C_1$, $\overline{D_1} = \{z''\} - D_1$. 
For all choices of $A_1,B_1,C_1,D_1$, let $F_1^{(i)} = A(G_1-\{w,y,x,z,z_1,z_2,y_1,y_2,w_1,w_2\} - A_1 - B_1 - C_1 - D_1)$, $F_2^{(i)} = A(G_2 - \{y_1,y_2\} - A_1)$, 
$F_3^{(i)} = A(G_3 - \{y_1,w_2\} - B_1)$,  $F_4^{(i)} = A(G_4 - \{z_1,z_2\} - C_1)$, and $F_5^{(i)} = A(G_5 - \{w_1,z_2\} - D_1)$. 
Then $|F_1^{(i)}| \geq \lceil (4(|G_1|-10-|A_1|-|B_1|-|C_1|-|D_1|)+3)/7 \rceil$, 
$|F_2^{(i)}| \geq \lceil (4(|G_2|-2-|A_1|)+3)/7 \rceil$, $|F_3^{(i)}| \geq \lceil (4(|G_3|-2-|B_1|)+3)/7 \rceil$, $|F_4^{(i)}| \geq \lceil (4(|G_4|-2-|C_1|)+3)/7 \rceil$, 
and $|F_5^{(i)}| \geq \lceil (4(|G_5|-2-|D_1|)+3)/7 \rceil$. Now $G[F_1^{(i)} \cup F_2^{(i)} \cup F_3^{(i)} \cup F_4^{(i)} \cup F_5^{(i)} + \{w,x,y,z\} - (\{y'\} \cap (F_1^{(i)} \triangle F_2^{(i)})) - (\{w'\} \cap (F_1^{(i)} \triangle F_3^{(i)})) - (\{z'\} \cap (F_1^{(i)} \triangle F_4^{(i)})) - (\{z''\} \cap (F_1^{(i)} \triangle F_5^{(i)})) ]$
is an induced forest in $G$, showing $a(G)\ge |F_1^{(i)}| + |F_2^{(i)}| + |F_3^{(i)}| + |F_4^{(i)}| + |F_5^{(i)}| + 4 - |\overline{A_1}| - |\overline{B_1}| - |\overline{C_1}| - |\overline{D_1}|$. 
By Lemma~\ref{ineq2}(1) 
(with $k=4, a = |G_1|-10, a_j = |G_{j+1}|-2$ for $j=1,2,3,4$, $L = \emptyset, c=4$),
 $a(G)\ge  \lceil (4n+3)/7 \rceil$, a contradiction. 
\qed

By Lemmas \ref{No3RR}, \ref{No4-3R}, \ref{No34-4-4-}, we have the following: 

\begin{coro}
\label{3summary} 
Let $x\in V_3$. Then there exists $v\in N(x)\cap V_{\ge 5}$ such that $R_{v,\{x\}}=\emptyset$. 
\end{coro}

\section{A forbidden configuration around a 3-vertex}

We prove the following, which eliminates two configurations around a 3-vertex. 


\begin{lem}
\label{5-2-B-combined}
Let $x \in V_3$, $\{y,z\} \subseteq V_{\leq 4}$, $N(x) = \{w,y,z\}$. Suppose $xzvwx$ is a facial cycle and $w\in V_5$. Then $R_{v,\{w,z\}}=\emptyset$ and 
$v \not\in V_{\leq 4}$.
\end{lem}




\pf
We may assume $\{y,z\}\subseteq V_4$ because the case when $y \in V_3$ or $z \in V_3$ is identical by replacing identifying neighbors of $4$-vertex with deleting a neighbor of $3$-vertex. 

In the first part, we prove $R_{v,\{w,z\}}=\emptyset$. For, suppose $R\in R_{v,\{w,z\}}$. 
Let $N(y) = \{x,y_1,y_2,z_1\}$, $N(z) = \{x,v,z_1,z_2\}$ and $y_2w \in E(G)$. 

First, we claim that $wz_1 \not\in E(G)$. For, suppose $wz_1 \in E(G)$. 
There exists a separation $(G_1,G_2)$ such that $V(G_1 \cap G_2) = \{w,z,z_1\}$, $\{x,y,y_1,y_2\} \subseteq V(G_1)$, and $v \in V(G_2)$. 
Let $F_1^{(1)} = A(G_1 - \{z_1,z,w,x\})$, and $F_2^{(1)} = A((G_2 - \{z_1,z,w,v\}) * R)$. 
Then $|F_1^{(1)}| \geq \lceil (4(|G_1|-4)+3)/7 \rceil$, and $|F_2^{(1)}| \geq \lceil (4(|G_2|-5)+3)/7 \rceil$. 
Now $G[(F_1^{(1)} \cup F_2^{(1)} + \{x,z\}) \cdot R]$ is an induced forest in $G$, showing $a(G) \geq |F_1^{(1)}| + |F_2^{(1)}| + 3 \geq \lceil (4n+3)/7 \rceil$, a contradiction. 

Secondly, we claim that $wz_2 \not\in E(G)$. For otherwise,
there exists a separation $(G_1,G_2)$ such that $V(G_1 \cap G_2) = \{w,v,z_2\}$, $\{x,y,z,z_1,y_1,y_2\} \subseteq V(G_1)$, and $R \subseteq V(G_2)$. 
Let $F_1^{(2)} = A(G_1 - \{x,z,z_1,z_2,w,v\})$, and $F_2^{(2)} = A(G_2 - \{z_2,w\})$.
Then $|F_1^{(2)}| \geq \lceil (4(|G_1|-6)+3)/7 \rceil$, and $|F_2^{(2)}| \geq \lceil (4(|G_2|-2)+3)/7 \rceil$. 
Now $G[F_1^{(2)} \cup F_2^{(2)} + \{x,z\}]$ is an induced forest in $G$, showing $a(G) \geq |F_1^{(2)}| + |F_2^{(2)}| + 2$. 
This implies $4(|G_2|-2)+3 \equiv 0,5,6 \mod 7$. 
If $|N(y_1) \cap N(y_2)| \leq 2$, 
 let $F_1^{(3)} = A((G_1 - \{x,y,z,z_1,w,v\})/y_1y_2)$ with $y'$ as the identification of $\{y_1,y_2\}$, and $F_2^{(3)} = A((G_2 - \{w,v\}) * R)$. 
Then $|F_1^{(3)}| \geq \lceil (4(|G_1|-7)+3)/7 \rceil$, and $|F_2^{(3)}| \geq \lceil (4(|G_2|-3)+3)/7 \rceil$. 
Now $F^{(3)} := G[(F_1^{(3)} \cup F_2^{(3)} + \{x,y,z\}) \cdot R - (\{z_2\} \cap (F_1^{(3)} \triangle F_2^{(3)}))]$ (if $y' \not\in F_1^{(3)}$) or $G[( (F_1^{(3)} - y') \cup F_2^{(3)} + \{x,y_1,y_2,z\} ) \cdot R - (\{z_2\} \cap (F_1^{(3)} \triangle F_2^{(3)}))]$ (if $y' \in F_1^{(3)}$)  is an induced forest in $G$, showing $a(G) \geq |F_1^{(3)}| + |F_2^{(3)}| + 4 - 1$. By Lemma \ref{ineq2}(1) 
(with $k=1, a = |G_1| - 6, a_1 = |G_2| - 2, L = \emptyset, c = 3$), 
$a(G) \geq \lceil (4n + 3)/7 \rceil$, a contradiction.
So  $|N(y_1) \cap N(y_2)| > 2$. 
Then there exist $a_1 \in N(y_1) \cap N(y_2)$ and subgraphs $G_1',G_2',G_3'$ of $G$ such that 
$G_2' = G_2$, 
$G_3'$ is the maximal subgraph of $G$ contained in the closed region of the plane bounded by the cycle $yy_1a_1y_2y$ and containing $N(y_1)\cap N(y_2) - \{y\}$,
and $G_1$ is obtained from $G$ by removing $G_2'-\{w,v,z_2\}$ and $G_3'-\{a_1,y_2,y_1\}$. 
Let $A_4 = \emptyset$ and $A_5 = \{a_1\}$.
For $i=4,5$, 
let $F_1^{(i)} = A(G_1' - \{x,y,z,z_1,w,v,y_1,y_2,z_2\} - A_i)$, 
$F_2^{(i)} = A((G_2' - \{w,v,z_2\}) * R)$, 
and $F_3^{(i)} = A(G_3' - \{y_1,y_2\} - A_i)$. 
Then $|F_1^{(i)}| \geq \lceil (4(|G_1'|-9-|A_i|)+3)/7 \rceil$, 
$|F_2^{(i)}| \geq \lceil (4(|G_2'|-4)+3)/7 \rceil$, 
and $|F_3^{(i)}| \geq \lceil (4(|G_3'|-2-|A_i|)+3)/7 \rceil$. 
Now $G[(F_1^{(i)} \cup F_2^{(i)} \cup F_3^{(i)} + \{x,y,z,w\}) \cdot R - (\{a_1\} \cap (F_1^{(i)} \triangle F_3^{(i)})) ]$ is an induced forest in $G$, showing $a(G) \geq |F_1^{(i)}| + |F_2^{(i)}| + |F_3^{(i)}| + 5 - (1-|A_i|)$. 
By Lemma \ref{ineq2}(1) 
(with $k=1,a=|G_1'|-8,a_1 = |G_3'|-2,L=\{1\},b_1 = |G_2'|-4,c=5$),
  $a(G) \geq \lceil (4n+3)/7 \rceil$, a contradiction.

\medskip
Case 1: $|N(y_1) \cap N(y_2)| \leq 2$ and $|N(z_1) \cap N(z_2)| \leq 2$.

Let $F' = A(((G - \{x,y,z,v\}) * R)/\{y_1y_2,z_1z_2\} + wz')$ with $y'$ (respectively, $z'$) as the identifications of $\{y_1,y_2\}$ (respectively, $\{z_1,z_2\}$). 
Then $|F'| \geq \lceil (4(n-7)+3)/7 \rceil$.
Let $F := (F' + \{x,y,z\}) \cdot R$ if $y',z' \not\in F'$, and otherwise $F'$ obtained from $(F' + \{x,y,z\}) \cdot R$ by deleting $\{y,y'\}$ (respectively, $\{z',z\}$) and adding $\{y_2,y_1\}$ (respectively, $\{z_1,z_2\}$) when $y' \in F'$ (respectively, $z' \in F'$).
Therefore, $G[F']$ is an induced forest in $G$, showing $a(G) \geq |F'| + 4 \geq \lceil (4n+3)/7 \rceil$, a contradiction. 

\medskip
Case 2: $|N(y_1) \cap N(y_2)|  > 2$. 

There exist $a_1 \in N(y_1) \cap N(y_2)$ and a separation $(G_1,G_2)$ such that $V(G_1 \cap G_2) = \{y_1,y_2,a_1\}$, and $\{x,y,z,w,v\} \subseteq V(G_1)$, $N(y_1) \cap N(y_2) - \{y\} \subseteq V(G_2)$. 
Define $A_1 = \overline{A_2} = \{a_1\}$ and $A_2 = \overline{A_1} = \emptyset$. 
For $i=1,2$, let $F_1^{(i)} = A((G_1 - \{x,y,z,z_1,y_1,y_2,v\})* R - A_i + wz_2)$, and $F_2^{(i)} = A(G_2 - \{y_1,y_2\} - A_i)$. 
Then $|F_1^{(i)}| \geq \lceil (4(|G_1|-8-|A_i|)+3)/7 \rceil$, and $|F_2^{(i)}| \geq \lceil (4(|G_2|-2-|A_i|)+3)/7 \rceil$. 
Now $G[(F_1^{(i)} \cup F_2^{(i)} + \{x,y,z\}) \cdot R - (\{a_1\} \cap (F_1^{(i)} \triangle F_2^{(i)}))]$ is an induced forest in $G$, showing $a(G) \geq |F_1^{(i)}| + |F_2^{(i)}| + 4 - (1-|A_i|)$. 
By Lemma \ref{ineq2}(2), 
$(4(|G_1|-8)+3,4(|G_2|-2)+3) \equiv (4,0),(0,4) \mod 7$.

\medskip
{\it Subcase 2.1:} $(4(|G_1|-8)+3,4(|G_2|-2)+3) \equiv (4,0) \mod 7$.

Let $F_1^{(5)} = A(((G_1 - \{x,y,z,z_1,v\}) * R)/y_1y_2 + wz_2)$ with $y'$ as the identification of $\{y_1,y_2\}$, and $F_2^{(5)} = A(G_2)$. 
Then $|F_1^{(5)}| \geq \lceil (4(|G_1|-7)+3)/7 \rceil$, and $|F_2^{(i)}| \geq \lceil (4|G_2|+3)/7 \rceil$. 
Now 
$G[(F_1^{(5)} \cup F_2^{(5)} + \{x,y,z\}) \cdot R - (\{y_1,y_2,a_1\} \cap (F_1^{(5)} \triangle F_2^{(5)}))]$ (if $y' \not\in F_1^{(5)}$) or $G[(F_1^{(5)} \cup F_2^{(5)} + \{x,y_1,y_2,z\} - \{y'\}) \cdot R - (\{a_1\} \cap (F_1^{(5)} \triangle F_2^{(5)})) - (\{y_1,y_2\} - F_2^{(5)})]$ (if $y' \in F_1^{(5)}$) is an induced forest in $G$, showing $a(G) \geq |F_1^{(5)}| + |F_2^{(5)}| + 4 - 3 \geq \lceil (4n+3)/7 \rceil$, a contradiction.

\medskip
{\it Subcase 2.2:} $(4(|G_1|-8)+3,4(|G_2|-2)+3) \equiv (0,4) \mod 7$.
 
If $wy_1 \not\in E(G)$ and $|N(v) \cap N(z_2)| \leq 2$, then 
let $F_1^{(6)} = A(G_1 - \{x,y,z,z_1,y_2\})/vz_2 + wy_1)$ with $z'$ as the identification of $\{v,z_2\}$, and $F_2^{(6)} = A(G_2 - y_2)$. 
Then $|F_1^{(6)}| \geq \lceil (4(|G_1|-6)+3)/7 \rceil$, and $|F_2^{(6)}| \geq \lceil (4(|G_2|-1)+3)/7 \rceil$. 
Let $F = F_1^{(6)} \cup F_2^{(6)} + \{x,y,z\} - (\{y_1,a_1\} \cap (F_1^{(6)} \triangle F_2^{(6)}))$.
Now $G[F]$ (if $z' \not\in F_1^{(6)}$) or $G[F - \{z,z'\} + \{v,z_2\}]$ (if $z' \in F_1^{(6)}$)
is an induced forest in $G$, showing $a(G) \geq |F_1^{(6)}| + |F_2^{(6)}| + 3 - 2 \geq \lceil (4n+3)/7 \rceil$, a contradiction. So we have $wy_1 \in E(G)$ or $|N(v) \cap N(z_2)| \geq 3$. 


If $wy_1 \in E(G)$, then there exists a separation $(G_1',G_2')$ such that $V(G_1' \cap G_2') = \{y_1,y_2,w\}$, $\{x,y,z,v\} \subseteq V(G_1')$, and $N(y_1) \cap N(y_2) - \{y\} \subseteq V(G_2')$. 
Let $F_1^{(7)} = A((G_1 - \{w,x,y,z,y_1,y_2,z_1,v\}) * R)$ and $F_2^{(7)} = A(G_2 - \{y_1,w\})$. 
Then $|F_1^{(7)}| \geq \lceil (4(|G_1'|-9)+3)/7 \rceil$ and $|F_2^{(7)}| \geq \lceil (4(|G_2'|-2)+3)/7 \rceil$. 
Now $G[(F_1^{(7)} \cup F_2^{(7)} + \{x,y,z\}) \cdot R]$ is an induced forest in $G$, showing $a(G) \geq |F_1^{(7)}| + |F_2^{(7)}| + 4$. 
Let $F_1^{(8)} = A((G_1 - \{w,x,y,z,y_1,y_2,z_1,z_2,v\}) * R)$, and $F_2^{(8)} = A(G_2 - \{y_1,w,y_2\})$. 
Then $|F_1^{(8)}| \geq \lceil (4(|G_1'|-10)+3)/7 \rceil$, and $|F_2^{(8)}| \geq \lceil (4(|G_2'|-3)+3)/7 \rceil$. 
Now $F^{(8)} := G[(F_1^{(8)} \cup F_2^{(8)} + \{x,y,z,w\}) \cdot R]$  is an induced forest in $G$, showing $a(G)\geq |F_1^{(8)}| + |F_2^{(8)}| + 5$. By Lemma \ref{ineq2}(1) 
(with $k=1, a=|G_1'|-9, a_1 = |G_2'|-2, L = \emptyset, c=4$), 
$a(G) \geq \lceil (4n+3)/7 \rceil$, a contradiction. 


If $|N(v) \cap N(z_2)| >2$, then there exist $c_1 \in N(v) \cap N(z_2)$ and subgraphs $G_1'',G_2'',G_3''$ of $G$ such that
$G_2'' = G_2$,
$G_3''$ is the maximal subgraph of $G$ contained in the closed region of the plane bounded by the cycle $zvc_1z_2$ and containing $N(v)\cap N(z_2) - \{z\}$,
and $G_1''$ is obtained from $G$ by removing $G_2''-\{y_1,y_2,a_1\}$ and $G_3''-\{v,z_2,c_1\}$. 
By symmetry, assume $R \subseteq G_1''$.
Define $C_8 = \overline{C_9} = \{c_1\}$ and $C_9 = \overline{C_8} =\emptyset$.  
For $i=8,9$, let $F_1^{(i)} = A((G_1'' - \{x,y,z,y_2,z_1,z_2,v\} - C_i) * R  + wy_1)$, $F_2^{(i)} = A(G_2'' - y_2)$, and $F_3^{(i)} = A(G_3'' - \{c,z_2\} - C_i)$.
Then $|F_1^{(i)}| \geq \lceil (4(|G_1''|-8-|C_i|)+3)/7 \rceil$, $|F_2^{(i)}| \geq \lceil (4(|G_2''|-1)+3)/7 \rceil$, and $|F_3^{(i)}| \geq \lceil (4(|G_3''|-2-|C_i|)+3)/7 \rceil$. 
Now $G[(F_1^{(i)} \cup F_2^{(i)} \cup F_3^{(i)} + \{x,y,z\}) \cdot R - (\{y_1,a_1\} \cap (F_1^{(i)} \triangle F_2^{(i)})) - (\{c_1\} \cap (F_1^{(i)} \triangle F_3^{(i)}))]$ is an induced forest in $G$, showing $a(G) \geq |F_1^{(i)}| + |F_2^{(i)}| + |F_3^{(i)}| + 4 - 2 - |\overline{C_1}|$.
By Lemma \ref{ineq2}(1) 
(with $k=1, a=|G_1''|-8, a_1 = |G_3''|-2, L = \{1\}, b_1 = |G_2''|-1, c=2$),
 $a(G) \geq \lceil (4n+3)/7 \rceil$, a contradiction.

\medskip
Case 3: $|N(z_1) \cap N(z_2)| > 2$. 

There exist $b_1 \in N(z_1) \cap N(z_2)$ and a separation $(G_1,G_2)$ such that $V(G_1 \cap G_2) = \{z_1,z_2,b_1\}$, $\{x,y,z,w,v\} \subseteq V(G_1)$, and $N(z_1) \cap N(z_2) - \{z\} \subseteq V(G_2)$. 
Let $B_1 = \{b_1\}$ and $B_2 = \emptyset$. 
For $i=1,2$, let $F_1^{(i)} = A(((G_1 - \{x,y,z,z_1,z_2,v\} - B_i) * R )/y_1y_2)$ with $y'$ as the identification of $\{y_1,y_2\}$, and $F_2^{(i)} = A(G_2 - \{z_1,z_2\} - B_i)$. 
Then $|F_1^{(i)}| \geq \lceil (4(|G_1|-8-|B_i|)+3)/7 \rceil$, and $|F_2^{(i)}| \geq \lceil (4(|G_2|-2-|B_i|)+3)/7 \rceil$. 
Let $F = (F_1^{(i)} \cup F_2^{(i)} + \{x,y,z\}) \cdot R - (\{b_1\} \cap (F_1^{(i)} \triangle F_2^{(i)}))$.
Now $G[F]$ (if $y' \not\in F_1^{(i)}$) or $G[F - \{y,y'\} + \{y_1,y_2\}]$ (if $y' \in F_1^{(i)}$)
is an induced forest in $G$, showing $a(G) \geq |F_1^{(i)}| + |F_2^{(i)}| + 4 - ( 1- |B_i|)$. 
By Lemma \ref{ineq2}(2), 
$(4(|G_1|-8)+3,4(|G_2|-2)+3) \equiv (0,4),(4,0) \mod 7$.

\medskip
{\it Subcase 3.1:} $(4(|G_1|-8)+3,4(|G_2|-2)+3) \equiv (0,4) \mod 7$.

Let $F_1^{(3)} = A(((G_1 - \{x,y,z,v,z_1\}) * R)/y_1y_2 + wz_2)$ with $y'$ as the identification of $\{y_1,y_2\}$, and $F_2^{(3)} = A(G_2 - z_1)$. 
Then $|F_1^{(3)}| \geq \lceil (4(|G_1|-7)+3)/7 \rceil$, and $|F_2^{(3)}| \geq \lceil (4(|G_2|-1)+3)/7 \rceil$. 
Let $F = (F_1^{(3)} \cup F_2^{(3)} + \{x,y,z\}) \cdot R - (\{z_2,b_1\} \cap (F_1^{(3)} \triangle F_2^{(3)}))$. 
Now $G[F]$ (if $y' \not\in F_1^{(3)}$) or $G[F - \{y,y'\}+\{y_1,y_2\}]$ (if $y' \in F_1^{(3)}$)
is an induced forest in $G$, showing  $a(G) \geq |F_1^{(3)}| + |F_2^{(3)}| + 4 - 2 \geq \lceil (4n+3)/7 \rceil$, a contradiction.

\medskip
{\it Subcase 3.2:} $(4(|G_1|-8)+3,4(|G_2|-2)+3) \equiv (4,0) \mod 7$.

Let $F_1^{(4)} = A(((G_1 - \{x,y,z,v\}) * R)/\{y_1y_2,z_1z_2\} + wz')$ with $y'$ (respectively, $z'$) as the identification of $\{y_1,y_2\}$ (respectively, $\{z_1,z_2\}$), and $F_2^{(9)} = A(G_2)$.
Then $|F_1^{(4)}| \geq \lceil (4(|G_1|-7)+3)/7 \rceil$, and $|F_2^{(4)}| \geq \lceil (4|G_2|+3)/7 \rceil$. 
Now $F^{(4)} := (\overline{F_1}^{(4)} \cup F_2^{(4)}) \cdot R - (\{z_1,z_2,b_1\} \cap (\overline{F_1}^{(4)} \triangle F_2^{(4)})) $
where $\overline{F_1}^{(4)} = {F_1}^{(4)} + \{x,y,z\}$ if $y',z' \not\in F_1^{(4)}$;
or obtained from ${F_1}^{(4)} + \{x,y,z\}$ by deleting $\{z,z'\}$ ($\{y,y'\}$ respectively) and adding $\{z_1,z_2\}$ ($\{y_1,y_2\}$ respectively) when $z' \in {F_1}^{(4)}$ ($y' \in {F_1}^{(4)}$ respectively).
Therefore, $G[F^{(4)}]$ is an induced forest of size $|F_1^{(4)}| + |F_2^{(4)}| + 4 - 3 \geq \lceil (4n+3)/7 \rceil$, a contradiction.

\bigskip


We now prove $v \not\in V_{\leq 4}$. 
By Lemma \ref{2summary}, $v \not V_2$.
For otherwise, $v \in V_4$. 
The case $v \in V_3$ is identical by replacing identification of neighbors of $v$ with deletion of a neighbor of $v$.
Let $N(y) = \{x,y_1,y_2,z_1\}$ and $N(z) = \{x,v,z_1,z_2\}$ and $y_2w \in E(G)$. Let $N(v) = \{z,w,v_1,v_2\}$ and $v_1vzz_2v_1$ be a facial cycle. 

\medskip
\textit{Claim 1: $wz_1 \not\in E(G)$. }


For, suppose $wz_1 \in E(G)$. There exists a separation $(G_1,G_2)$ such that $V(G_1 \cap G_2) = \{w,x,z_1\}$, $\{y,y_1,y_2\} \subseteq V(G_1)$, and $\{z,v\} \subseteq V(G_2)$. 
For $i=1,2$, let $F_i^{(1)} = A(G_i - \{z_1,w,x\})$. 
Then $|F_i^{(1)}| \geq \lceil (4(|G_i|-3)+3)/7 \rceil$. 
Now $G[F_1^{(1)} \cup F_2^{(1)} + x]$ is an induced forest in $G$, showing  $a(G) \geq |F_1^{(1)}| + |F_2^{(1)}| + 1$. 
By Lemma \ref{ineq2}(7) (with $k=1, a_i = |G_i| - 3$ for $i=1,2$), 
$(4(|G_1|-3)+3,4(|G_2|-3)+3) \equiv (0,0), (0,6), (6,0) \mod 7$. 

If $wz_2 \not\in E(G)$ and $wy_1 \not\in E(G)$, let $F_1^{(2)} = A(G_1 - \{z_1,x,y,y_2\} + wy_1)$, and $F_2^{(2)} = A(G_2 - \{z_1,z,x,v\} + wz_2)$. 
For $i=1,2$ $|F_i^{(2)}| \geq \lceil (4(|G_i|-4)+3)/7 \rceil$.
Now $G[F_1^{(2)} \cup F_2^{(2)} + \{x,y,z\} - (\{w\} \cap (F_1^{(2)} \triangle F_2^{(2)}))]$  is an induced forest in $G$, showing $a(G) \geq |F_1^{(2)}| + |F_2^{(2)}| + 3 - 1 \geq \lceil (4n+3)/7 \rceil$, a contradiction.

If $wz_2 \in E(G)$,  let $F_1^{(3)} = A(G_1 - \{z_1,w,x,y,y_2\})$, and $F_2^{(3)} = A(G_2 - \{z_1,w,x,z,v,z_2\})$. 
Then $|F_1^{(3)}| \geq \lceil (4(|G_1|-5)+3)/7 \rceil$, and $|F_2^{(3)}| \geq \lceil (4(|G_2|-6)+3)/7 \rceil$.
Now $G[F_1^{(3)} \cup F_2^{(3)} + \{x,y,z,w\}]$ is an induced forest in $G$, showing $a(G) \geq |F_1^{(3)}| + |F_2^{(3)}| + 4 \geq \lceil (4n+3)/7 \rceil$, a contradiction.

So $wy_1 \in E(G)$. 
Let $F_1^{(4)} = A(G_1 - \{z_1,w,x,y,y_2,y_1\})$, and $F_2^{(4)} = A(G_2 - \{z_1,w,x,z,v\})$. 
Then $|F_1^{(4)}| \geq \lceil (4(|G_1|-6)+3)/7 \rceil$, and $|F_2^{(4)}| \geq \lceil (4(|G_2|-5)+3)/7 \rceil$.
Now $G[F_1^{(4)} \cup F_2^{(4)} + \{x,y,z,w\}]$ is an induced forest in $G$, showing $a(G) \geq |F_1^{(4)}| + |F_2^{(4)}| + 4 \geq \lceil (4n+3)/7 \rceil$, a contradiction. 



\medskip
\textit{Claim 2: $wz_2 \not\in E(G)$. (By symmetry,  $wy_1 \not\in E(G)$) }



For, suppose $wz_2 \in E(G)$. 
There exists a separation $(G_1,G_2)$ such that $V(G_1 \cap G_2) = \{w,z,z_2\}$, $\{x,y,y_1,y_2\} \subseteq V(G_1)$, and $v \in V(G_2)$. 
Let $F_1^{(1)} = A(G_1 - \{w,z,z_2,x,z_1\})$, and $F_2^{(1)} = A(G_2 - \{w,z,z_2\})$. Then $|F_1^{(1)}| \geq \lceil (4(|G_1|-5)+3)/7 \rceil$, and $|F_2^{(1)}| \geq \lceil (4(|G_2|-3)+3)/7 \rceil$. 
Now $G[F_1^{(1)} \cup F_2^{(1)} + \{x,z\}]$ is an induced forest in $G$, showing $a(G) \geq |F_1^{(1)}| + |F_2^{(1)}| + 2$. 
By Lemma \ref{ineq2}(8)
 (with $a = |G_1|-5, a_1 = |G_2|-3$), 
$(4(|G_1|-5)+3,4(|G_2|-3)+3) \equiv (0,0),(0,6), (0,5), (5,0), (6,6), (6,0) \mod 7$.

If $wy_1 \not\in E(G)$, let $w' \in N(w) - \{y_2,v,z_2,x,y_1\}$. Let $e = w'y$ if $w' \in G_1$ and otherwise $e = \emptyset$.
Let $F_1^{(2)} = A(G_1 - \{w,z,z_2,x,y_2,z_1\} + e)$,  and $F_2^{(2)} = A(G_2 - \{w,z,z_2,v\})$. 
Then $|F_1^{(2)}| \geq \lceil (4(|G_1|-6)+3)/7 \rceil$, and $|F_2^{(2)}| \geq \lceil (4(|G_2|-4)+3)/7 \rceil$. 
Now $G[F_1^{(2)} \cup F_2^{(2)} + \{w,x,z\}]$  is an induced forest of size $|F_1^{(2)}| + |F_2^{(2)}| + 3 \geq \lceil (4n+3)/7 \rceil$, a contradiction. 

So $wy_1 \in E(G)$. 
By Lemma \ref{2summary}, $|N(y_2)| > 2$.
There exist subgraphs $G_1',G_2',G_3'$ of $G$ such that 
$G_2' = G_2$,
$G_3'$ is the maximal subgraph of $G$ contained in the closed region of the plane bounded by the cycle $wy_2yy_1w$ and containing $N(y_2)$,
and $G_1'$ is obtained from $G$ by removing $G_2'-\{w,z,z_2\}$ and $G_3'-\{w,y_2,y,y_1\}$. 
Note $G_3' - \{y,y_1,y_2,w\} \neq \emptyset$ since $|N(y_2)| > 2$. 
Let $F_1^{(4)} = A(G_1' - \{z_1,w,x,y,z,y_1,y_2,z_2\})$, $F_2^{(4)} = A(G_2' - \{z_2,w,z,v\})$, and $F_3^{(3)} = A(G_3' - \{y_1,y_2,w,y\})$. 
Then $|F_1^{(4)}| \geq \lceil (4(|G_1'|-8)+3)/7 \rceil$, $|F_i^{(4)}| \geq \lceil (4(|G_i'|-4)+3)/7 \rceil$ for $i=2,3$.
Note $|F_2^{(4)}| \geq \lceil (4(|G_2|-4)+3)/7 \rceil \geq (4(|G_2|-4)+3)/7 + 4/7$.
Now $G[F_1^{(4)} \cup F_2^{(4)} \cup F_3^{(4)} + \{x,y,z,w\}]$  is an induced forest in $G$, showing $a(G) \geq |F_1^{(4)}| + |F_2^{(4)}| + |F_3^{(4)}| + 4 \geq \lceil (4n+3)/7 \rceil$, a contradiction.

Note that we did not use the information on $v$ in the above proof. So by symmetry,  $wy_1 \not\in E(G)$.



\medskip
\textit{Claim 3: $v_2z_2 \not\in E(G)$. }


For, suppose $v_2z_2 \in E(G)$. There exists a separation $(G_1,G_2)$ of $G$ such that $V(G_1 \cap G_2) = \{v,v_2,z_2\}$, $\{x,y,y_1,y_2\} \subseteq V(G_1)$, and $v_1 \in V(G_2)$. 
Let $F_1^{(1)} = A(G_1 - \{v,v_2,z_2,z,x\} + wz_1)$, and $F_2^{(1)} = A(G_2 - \{z_2,v,v_2\})$. 
Then $|F_1^{(1)}| \geq \lceil (4(|G_1|-5)+3)/7 \rceil$, and $|F_2^{(1)}| \geq \lceil (4(|G_2|-3)+3)/7 \rceil$. 
Now $G[F_1^{(1)} \cup F_2^{(1)} + \{v,z\}]$ is an induced forest in $G$, showing $a(G) \geq |F_1^{(1)}| + |F_2^{(1)}| + 2$. 
By Lemma \ref{ineq2}(8)
(with $a=|G_1|-5, a_1 = |G_2|-3$), 
$(4(|G_1|-5)+3,4(|G_2|-3)+3) \equiv (0,0),(0,6),(0,5),(5,0),(6,6),(6,0) \mod 7$.

If $|N(z_1) \cap N(z_2)| \leq 2$, then let $F_1^{(2)} = A((G_1 - \{x,z,v,w,v_2\})/z_1z_2)$ with $z'$ as the identification of $\{z_1,z_2\}$, and $F_2^{(2)} = A(G_2 - \{v,v_2\})$. 
Then $|F_1^{(2)}| \geq \lceil (4(|G_1|-6)+3)/7 \rceil$, and $|F_2^{(2)}| \geq \lceil (4(|G_2|-2)+3)/7 \rceil$. 
Now $G[F_1^{(2)} \cup F_2^{(2)} + \{x,v,z\} - \{z_2\}]$ (if $z' \not\in F_1^{(2)}$) or  $G[(F_1^{(2)} - z') \cup F_2^{(2)}  + \{x,v,z_2,z_1\}  - (\{z_2\} - F_2^{(2)})] $ (if $z' \in F_1^{(2)}$) is an induced forest of size $|F_1^{(2)}| + |F_2^{(2)}| + 3 - 1$, which implies $a(G) \geq \lceil (4n+3)/7 \rceil$, a contradiction.

So $|N(z_1) \cap N(z_2)| > 2$. Then there exist $a_1 \in N(z_1) \cap N(z_2)$ and subgraphs  $G_1',G_2',G_3'$ such that 
$G_2' = G_2$,
$G_3'$ is the maximal subgraph of $G$ contained in the closed region of the plane bounded by the cycle $zz_1a_1z_2z$ and containing $N(z_1) \cap N(z_2) - \{z\}$,
and $G_1'$ is obtained from $G$ by removing $G_2'-\{v_2,v,z_2\}$ and $G_3'-\{z_2,a_1,z_1\}$. 
Let $A_3 = \{a_1\}$ and $A_4 = \emptyset$. 
For $i=3,4$, let $F_1^{(i)} = A(G_1' - \{x,w,z,v,z_1,z_2,v_2\} - A_i)$, $F_2^{(i)} = A(G_2' - \{z_2,v_2,v\})$, and $F_3^{(i)} = A(G_3' - \{z_1,z_2\} - A_i)$. 
Then $|F_1^{(i)}| \geq \lceil (4(|G_1'|-7 - |A_i|)+3)/7 \rceil$, $|F_2^{(i)}| \geq \lceil (4(|G_2'|-3)+3)/7 \rceil$, and $|F_3^{(i)}| \geq \lceil (4(|G_3'|-2 - |A_i|)+3)/7 \rceil$. 
Now $F^{(i)} := G[F_1^{(i)} \cup F_2^{(i)} \cup F_3^{(i)} + \{x,v,z\} - (\{a_1\} \cap (F_1^{(i)} \triangle F_2^{(i)})) ]$ is an induced forest in $G$ showing $a(G) \geq |F_1^{(3)}| + |F_2^{(3)}| + |F_3^{(3)}| + 3 - (1 - |A_i|)$. 
Let $(n_1,n_2,n_3):=(4(|G_1'|-7)+3,4(|G_2'|-3)+3,4(|G_3'|-2)+3)$. 
By Lemma \ref{ineq2}(2), 
$(n_1,n_2,n_3) \equiv (0,0,4), (4,0,0) \mod 7$.

If $(n_1,n_2,n_3) \equiv (0,0,4) \mod 7$, let $F_1^{(4)} = A(G_1' - \{x,y,z,v,z_1\} + wz_2)$, $F_2^{(4)} = A(G_2' - v)$, and $F_3^{(4)} = A(G_3' - z_1)$. 
Then $|F_1^{(4)}| \geq \lceil (4(|G_1'|-5)+3)/7 \rceil$, $|F_2^{(4)}| \geq \lceil (4(|G_2'|-1)+3)/7 \rceil$, and $|F_3^{(4)}| \geq \lceil (4(|G_3'|-1)+3)/7 \rceil$. Now $G[F_1^{(4)} \cup F_2^{(4)} \cup F_3^{(4)} + \{x,z\} - (\{a_1,z_2\} \cap (F_1^{(4)} \triangle F_3^{(4)})) -  (\{v_2,z_2\} \cap (F_1^{(4)} \triangle F_2^{(4)})) ]$ is an induced forest in $G$, showing $a(G) \geq |F_1^{(4)}| + |F_2^{(4)}| + |F_3^{(4)}| + 2 - 4 \geq \lceil (4n+3)/7 \rceil$, a contradiction. 

If $(n_1,n_2,n_3) \equiv (4,0,0) \mod 7$, then by Lemma \ref{No434Edge}, $yv \not\in E(G)$. 
Let $F_1^{(5)} = A(G_1' - \{w,x,z,z_1,z_2,a_1\} + yv)$, $F_2^{(5)} = A(G_2' - z_2)$, and $F_3^{(5)} = A(G_3' - \{z_1,z_2,a_1\})$. 
Then $|F_1^{(5)}| \geq \lceil (4(|G_1'|-6)+3)/7 \rceil$, $|F_2^{(5)}| \geq \lceil (4(|G_2'|-1)+3)/7 \rceil$, and $|F_3^{(5)}| \geq \lceil (4(|G_3'|-3)+3)/7 \rceil$. Now $G[F_1^{(5)} \cup F_2^{(5)} \cup F_3^{(5)} + \{x,z\} - (\{v,v_2\} \cap (F_1^{(5)} \triangle F_2^{(5)}))]$ is an induced forest in $G$, showing $a(G) \geq |F_1^{(5)}| + |F_2^{(5)}| + |F_3^{(5)}| + 2 - 2 \geq \lceil (4n+3)/7 \rceil$, a contradiction. 

\medskip
\textit{Claim 4: $v_2z_1 \not\in E(G)$. }


For, suppose $v_2z_1 \in E(G)$. 
By Lemma \ref{No434Edge}, $y \not\in \{v_1,v_2\}$. 
There exists a separation $(G_1,G_2)$ such that $V(G_1 \cap G_2) = \{v,v_2,z,z_1\}$, $\{x,y,y_1,y_2\} \subseteq V(G_1)$, and $\{z_2,v_1\} \subseteq V(G_2)$. 
Let $F_1^{(1)} = A(G_1 - \{v,v_2,z,z_1,w,x\})$ and $F_2^{(1)} = A(G_2 - \{v,v_2,z,z_1,v_1\})$. 
Then $|F_1^{(1)}| \geq \lceil (4(|G_1|-6)+3)/7 \rceil$, and $|F_2^{(1)}| \geq \lceil (4(|G_2|-5)+3)/7 \rceil$. 
Now $G[F_1^{(1)} \cup F_2^{(1)} + \{v,z,x\}]$ is an induced forest of size $|F_1^{(1)}| + |F_2^{(1)}| + 3$. 
Let $(n_1,n_2):=(4(|G_1|-6)+3,4(|G_2|-5)+3)$. 
By Lemma \ref{ineq2}(3), 
$(n_1,n_2) \equiv (0,0), (0,6), (0,5), (0,4), (4,0), (6,5), (5,6), (5,0), (6,6), (6,0) \mod 7$. 

If $(n_1,n_2) \equiv (0,0), (0,6), (0,5), (0,4) \mod 7$, then for $i=1,2$, let $F_i^{(2)} = A(G_i - \{v,v_2,z,z_1\})$. 
Then $|F_i^{(2)}| \geq \lceil (4(|G_i|-4)+3)/7 \rceil$. 
Now $G[F_1^{(2)} \cup F_2^{(2)} + \{v\}]$ is an induced forest of size $|F_1^{(2)}| + |F_2^{(2)}| + 1 \geq \lceil (4n+3)/7 \rceil$, a contradiction.

If $(n_1,n_2) \equiv (5,0), (6,0), (6,5), (5,6), (6,6) \mod 7$, 
then let $F_1^{(3)} = A(G_1 - \{x,y,z,v,z_1,$ $v_2,y_2\} + wy_1)$ and $F_2^{(3)} = A(G_2 - \{v,v_2,z,z_1\})$
Then $|F_1^{(3)}| \geq \lceil (4(|G_1|-7)+3)/7 \rceil$ and $|F_2^{(3)}| \geq \lceil (4(|G_2|-4)+3)/7 \rceil$. 
Now $G[F_1^{(3)} \cup F_2^{(3)} + \{x,y,z\}]$ is an induced forest of size $|F_1^{(3)}| + |F_2^{(3)}| + 3 \geq \lceil (4n+3)/7 \rceil$, a contradiction.

So $(n_1,n_2) \equiv (4,0) \mod 7$, 
then let $F_1^{(4)} = A(G_1 - \{z,z_1,v\} + xv_2)$ and $F_2^{(4)} = A(G_2 - \{z,z_1,v\})$
Then $|F_1^{(4)}| \geq \lceil (4(|G_1|-3)+3)/7 \rceil$ and $|F_2^{(4)}| \geq \lceil (4(|G_2|-3)+3)/7 \rceil$. 
Now $G[F_1^{(4)} \cup F_2^{(4)} + \{z\} - \{v_2\} \cap (F_1^{(4)} \triangle F_2^{(4)}) ]$ is an induced forest of size $|F_1^{(4)}| + |F_2^{(4)}| + 1 - 1 \geq \lceil (4n+3)/7 \rceil$, a contradiction.

\medskip
\textit{Claim 5: $yz_2 \not\in E(G)$, $y_1z \not\in E(G)$. } 


By symmetry, suppose that $y_1 = z_2$. 
By Lemma \ref{2summary}, $|N(z_1)| \geq 3$.
Then there exists a separation $(G_1,G_2)$ such that $V(G_1 \cap G_2) = \{z_1,y_1\}$, $\{x,y,z\} \subseteq V(G_1)$, and $N(z_1) - \{y,z\} \subseteq V(G_2)$. 
Let $F_1^{(1)} = A(G_1 - \{x,y,z,y_1,z_1\})$, and $F_2^{(1)} = A(G_2 - \{y_1,z_1\})$. 
Then $|F_1^{(1)}| \geq \lceil (4(|G_1|-5)+3)/7 \rceil$, and $|F_2^{(1)}| \geq \lceil (4(|G_2|-2)+3)/7 \rceil$. 
Now $G[F_1^{(1)} \cup F_2^{(1)} + \{y,z\}]$  is an induced forest in $G$, showing $a(G) \geq |F_1^{(1)}| + | F_2^{(1)}| + 2$. 
By Lemma \ref{ineq2}(8) (with $a = |G_1|-5, a_1 = |G_2|-2, c=2$), 
$(4(|G_1|-5)+3,4(|G_2|-2)+3) \equiv (0,0),(0,6),(0,5),(5,0),(6,6),(6,0) \mod 7$.
Let $F_1^{(2)} = A(G_1 - \{x,y,z,y_1,z_1,v\})$ and $F_2^{(2)} = A(G_2 - y_1)$. Then $|F_1^{(2)}| \geq \lceil (4(|G_1|-6)+3)/7 \rceil$ and $|F_2^{(2)}| \geq \lceil (4(|G_2|-1)+3)/7 \rceil$. 
Now $G[F_1^{(2)} \cup F_2^{(2)} + \{y,z\}]$  is an induced forest in $G$, showing $a(G) \geq |F_1^{(2)}| + | F_2^{(2)}| + 2 \geq \lceil (4n+3)/7 \rceil$, a contradiction.

\medskip

Next, we distinguish several cases.

\medskip
Case 1: $|N(z_1) \cap N(z_2)| \leq 2$, $|N(y_1) \cap N(y_2)| \leq 2$ and $|N(w) \cap N(v_2)| \leq 2$.

Let $F' = A((G - \{x,y,z,v\})/\{z_1z_2,y_1y_2,wv_2\} + v'z')$ with $z'$ (respectively, $y',v'$) as the identification of $\{z_1,z_2\}$ (respectively, $\{y_1,y_2\},\{w,v_2\}$). 
Then $|F'| \geq \lceil (4(n-7)+3)/7 \rceil$. 
Let $F: = F' + \{x,y,v,z\}$ if $z',y',v' \not\in F'$ and otherwise, let $F$  be obtained by $F' + \{x,y,v,z\}$ by deleting $\{z,z'\}$ (respectively, $\{y,y'\},\{v,v'\}$) and adding $\{z_1,z_2\}$ (respectively, $\{y_1,y_2\},\{v_2,w\}$) when $z' \in F'$ (respectively, when $y'\in F' ,v' \in F'$). 
Therefore, $G[F]$ is an induced forest in $G$, showing $|F'| + 4 \geq \lceil (4n+3)/7 \rceil$, a contradiction.

\medskip

Case 2: $|N(z_1) \cap N(z_2)| > 2$, $|N(y_1) \cap N(y_2)| \leq 2$ and $|N(w) \cap N(v_2)| \leq 2$. 

There exist $a_1 \in N(z_1) \cap N(z_2)$ and a separation $(G_1,G_2)$ of $G$ such that $V(G_1 \cap G_2) = \{z_1,z_2,a_1\}$, $\{x,y,y_1,y_2\} \subseteq V(G_1)$, and $N(z_1) \cap N(z_2) - \{z\} \subseteq V(G_2)$. 
Let $A_1 = \overline{A_2} = \{a_1\}$ and  $A_2 = \overline{A_1} =\emptyset$. 
For $i=1,2$, let $F_1^{(i)} = A((G_1 - \{x,y,z,v,z_1,z_2\} - A_i)/\{y_1y_2,wv_2\})$ with $y'$ (respectively, $v'$) as the identification of $\{y_1,y_2\}$ (respectively, $\{w,v_2\}$), and $F_2^{(i)} = A(G_2 - \{z_1,z_2\} - A_i)$. 
Then $|F_1^{(i)}| \geq \lceil (4(|G_1|-8- |A_i|)+3)/7 \rceil$, and $|F_2^{(i)}| \geq \lceil (4(|G_2|-2-|A_i|)+3)/7 \rceil$. 
Let $F^{(i)} = F_1^{(i)} \cup F_2^{(i)} + \{x,y,v,z\} - (\{a_1\} \cap (F_1^{(i)} \triangle F_2^{(i)}))$ if $y',v' \not\in F'$ and otherwise, let $F^{(i)}$  be obtained by $F_1^{(i)} \cup F_2^{(i)} + \{x,y,v,z\} - (\{a_1\} \cap (F_1^{(i)} \triangle F_2^{(i)}))$ by deleting $\{y,y'\}$ (respectively, $\{v,v'\}$) and adding $\{y_1,y_2\}$ (respectively, $\{v_2,w\}$) when $y' \in F'$ (respectively, $v' \in F'$).
Therefore, $G[F^{(i)}]$ is an induced forest in $G$, showing $a(G) \geq |F_1^{(i)}| + |F_2^{(i)}| + 4 - |\overline{A_i}|$. 
By Lemma \ref{ineq2}(2) (with $a = |G_1|-8, a_1 = |G_2|-2$), $(4(|G_1|-8)+3,4(|G_2|-2)+3) \equiv (0,4), (4,0) \mod 7$.

\medskip
{\it Subcase 2.1}: $(4(|G_1|-8)+3,4(|G_2|-2)+3) \equiv (4,0) \mod 7$. 

Let $F_1^{(1)} = A((G_1 - \{x,y,z,v\})/\{z_1z_2,y_1y_2,wv_2\} + v'z')$ with $z'$ (respectively, $y',v'$) as the identification of $\{z_1,z_2\}$ (respectively, $\{y_1,y_2\},\{w,v_2\}$), and $F_2^{(1)} = A(G_2)$. 
Then $|F_1^{(1)}| \geq \lceil (4(|G_1|-7)+3)/7 \rceil$, and $|F_2^{(1)}| \geq \lceil (4|G_2|+3)/7 \rceil$. 
Let $F^{(1)} = \overline{F_1}^{(1)} \cup F_2^{(1)} - (\{z_1,z_2,a_1\} \cap (\overline{F_1}^{(1)} \triangle F_2^{(1)}))$,
where $\overline{F_1}^{(1)} = {F_1}^{(1)}  + \{x,y,v,z\} $ if $z',y',v' \not\in F'$ and otherwise, let $\overline{F_1}^{(1)}$  be obtained by ${F_1}^{(1)}  + \{x,y,v,z\}$ by deleting $\{z,z'\}$ (respectively, $\{y,y'\},\{v,v'\}$) and adding $\{z_1,z_2\}$ (respectively, $\{y_1,y_2\},\{v_2,w\}$) when $z' \in F'$ (respectively, when $y' \in F',v' \in F'$). 
Therefore, $G[F^{(1)}]$ is an induced forest in $G$, showing $a(G) \geq |F_1^{(1)}| + | F_2^{(1)}| + 4 - 3 \geq \lceil (4n+3)/7 \rceil$, a contradiction.

\medskip
{\it Subcase 2.2}: $(4(|G_1|-8)+3,4(|G_2|-2)+3) \equiv (0,4) \mod 7$. 

Let $F_1^{(2)} = A((G_1 - \{x,y,z,v,z_1\})/y_1y_2 + wz_2)$ with $y'$ as the identification of $\{y_1,y_2\}$, and $F_2^{(2)} = F(G_2 - z_1)$. 
Then $|F_1^{(2)}| \geq \lceil (4(|G_1|-6)+3)/7 \rceil$, and $|F_2^{(2)}| \geq \lceil (4(|G_2|-1)+3)/7 \rceil$.
Let $F = F_1^{(2)} \cup F_2^{(2)} + \{x,y,z\} - (\{z_2,a_1\} \cap (F_1^{(2)} \triangle F_2^{(2)}))$.
Now $G[F]$ (if $y' \not\in F_1^{(2)}$) or $G[F - \{y,y'\} + \{y_1,y_2\}]$ (if $y' \in F_1^{(2)}$)
is an induced forest in $G$, showing $a(G) \geq |F_1^{(2)}| + | F_2^{(2)}| + 3 - 2 \geq \lceil (4n+3)/7 \rceil$, a contradiction. 

\medskip

Case 3: $|N(z_1) \cap N(z_2)| \leq 2$, $|N(y_1) \cap N(y_2)| \leq 2$ and $|N(w) \cap N(v_2)| > 2$. 

There exist $c_1 \in N(z_1) \cap N(z_2)$ and a separation $(G_1,G_2)$ of $G$ such that $V(G_1 \cap G_2) = \{w,v_2,c_1\}$, $\{x,y,y_1,y_2\} \subseteq V(G_1)$, and $N(w) \cap N(v_2) - \{v\} \subseteq V(G_2)$. 
Let $C_1 = \overline{C_2} =  \{c_1\}$ and $\overline{C_1} = C_2 = \emptyset$. 
For $i=1,2$, let $F_1^{(i)} = A((G_1 - \{x,y,z,v,w,v_2\} - C_i)/\{y_1y_2,z_1z_2\})$ with $y'$ (respectively, $z'$) as the identification of $\{y_1,y_2\}$ (respectively, $\{z_1,z_2\}$) and $F_2^{(i)} = A(G_2 - \{w,v_2\} - C_i)$. 
Then $|F_1^{(i)}| \geq \lceil (4(|G_1|-8- |C_i|)+3)/7 \rceil$ and $|F_2^{(i)}| \geq \lceil (4(|G_2|-2-|C_i|)+3)/7 \rceil$. 
Let $F^{(i)} = F_1^{(i)} \cup F_2^{(i)} + \{x,y,v,z\} - (\{c_1\} \cap (F_1^{(i)} \triangle F_2^{(i)}))$ if $y',z' \not\in F_1^{(i)}$ and otherwise, let $F^{(i)}$  be obtained by $F_1^{(i)} \cup F_2^{(i)} + \{x,y,v,z\} - (\{c_1\} \cap (F_1^{(i)} \triangle F_2^{(i)}))$ by deleting $\{y,y'\}$ (respectively, $\{z,z'\}$) and adding $\{y_1,y_2\}$ (respectively, $\{z_1,z_2\}$) when $y' \in F'$ (respectively, $z' \in F'$).
Therefore, $G[F^{(i)}]$ is an induced forest in $G$, showing $a(G)\geq |F_1^{(i)}| + |F_2^{(i)}| + 4 - |\overline{C_i}|$. 
By Lemma \ref{ineq2}(2) (with $a =|G_1|-8, a_1 = |G_2|-2$), $(4(|G_1|-8)+3,4(|G_2|-2)+3) \equiv (0,4), (4,0) \mod 7$.

\medskip
{\it Subcase 3.1}: $(4(|G_1|-8)+3,4(|G_2|-2)+3) \equiv (4,0) \mod 7$. 

Let $F_1^{(1)} = A((G_1 - \{x,y,z,v\})/\{z_1z_2,y_1y_2,wv_2\} + v'z')$ with $z'$ (respectively, $y',v'$) as the identification of $\{z_1,z_2\}$ (respectively, $\{y_1,y_2\},\{w,v_2\}$), and $F_2^{(1)} = A(G_2)$. 
Then $|F_1^{(1)}| \geq \lceil (4(|G_1|-7)+3)/7 \rceil$ and $|F_2^{(1)}| \geq \lceil (4|G_2|+3)/7 \rceil$. 
Let $F^{(1)} = \overline{F_1}^{(1)} \cup F_2^{(1)} - (\{w,v_2,c_1\} \cap (\overline{F_1}^{(1)} \triangle F_2^{(1)}))$,
where $\overline{F_1}^{(1)} = {F_1}^{(1)}  + \{x,y,v,z\} $ if $z',y',v' \not\in F'$ and otherwise, let $\overline{F_1}^{(1)}$  be obtained by ${F_1}^{(1)}  + \{x,y,v,z\}$ by deleting $\{z,z'\}$ (respectively, $\{y,y'\},\{v,v'\}$) and adding $\{z_1,z_2\}$ (respectively, $\{y_1,y_2\},\{v_2,w\}$) when $z' \in F'$ (respectively, $y' \in F',v' \in F'$). 
Therefore, $G[F^{(1)}]$ is an induced forest in $G$, showing $a(G) \geq |F_1^{(1)}| + | F_2^{(1)}| + 4 - 3 \geq \lceil (4n+3)/7 \rceil$, a contradiction.

\medskip
{\it Subcase 3.2}: $(4(|G_1|-8)+3,4(|G_2|-2)+3) \equiv (0,4) \mod 7$. 

Let $F_1^{(2)} = A(G_1 - \{w,x,y,z\}/\{y_1y_2,z_1z_2\})$ with $y'$ (respectively, $z'$) as the identification of $\{y_1,y_2\}$ (respectively, $\{z_1,z_2\}$), and $F_2^{(2)} = A(G_2 - w)$. 
Then $|F_1^{(2)}| \geq \lceil (4(|G_1|-6)+3)/7 \rceil$, and $|F_2^{(2)}| \geq \lceil (4(|G_2|-1)+3)/7 \rceil$. 
Let $F^{(2)} := \overline{F_1}^{(2)} \cup F_2^{(2)}  - (\{v_2,c_1\} \cap (\overline{F_1}^{(2)} \triangle F_2^{(2)}))$, where $\overline{F_1}^{(2)} = {F_1}^{(2)} + \{x,y,z\}$ if $y',z' \not\in F_1^{(2)}$, and otherwise, $\overline{F_1}^{(2)}$ obtained from $\overline{F_1}^{(2)} = {F_1}^{(2)} + \{x,y,z\}$ by deleting $y,y'$ (respectively, $z,z'$) and adding $\{y_1,y_2\}$ (respectively, $\{z_1,z_2\}$) when $y' \in F_1^{(2)}$ (respectively, $z' \in F_1^{(2)}$). 
Therefore, $G[F^{(2)}]$ is an induced forest in $G$, showing $a(G) \geq |F_1^{(2)}| + | F_2^{(2)}| + 3 - 2 \geq \lceil (4n+3)/7 \rceil$, a contradiction.

\medskip

Case 4: $|N(z_1) \cap N(z_2)| \leq 2$, $|N(y_1) \cap N(y_2)| > 2$ and $|N(w) \cap N(v_2)| \leq 2$. 

There exist $b_1 \in N(y_1) \cap N(y_2)$ and a separation $(G_1,G_2)$ of $G$ such that $V(G_1 \cap G_2) = \{y_1,y_2,b_1\}$, $\{x,z,y_1,y_2\} \subseteq V(G_1)$, and $N(y_1) \cap N(y_2) - \{y\} \subseteq V(G_2)$. 
Let $B_1 = \overline{B_2} =  \{b_1\}$, and let $B_2 = \overline{B_1} = \emptyset$.  
For $i=1,2$, let $F_1^{(i)} = A((G_1 - \{x,y,z,v,y_1,y_2\} - B_i)/\{wv_2,z_1z_2\} + v'z')$ with $v'$ (respectively, $z'$) as the identification of $\{w,v_2\},\{z_1,z_2\}$ and $F_2^{(i)} = A(G_2 - \{y_1,y_2\} - B_i)$. 
Then $|F_1^{(i)}| \geq \lceil (4(|G_1|-8- |B_i|)+3)/7 \rceil$ and $|F_2^{(i)}| \geq \lceil (4(|G_2|-2-|B_i|)+3)/7 \rceil$. 
Let $F^{(i)} := G[\overline{F_1}^{(i)} \cup F_2^{(i)} - (\{b_1\} \cap (\overline{F_1}^{(i)} \triangle F_2^{(i)}))]$,
where $\overline{F_1}^{(i)} := {F_1}^{(i)} + \{x,y,v,z\}$ if $v',z' \not\in F_1^{(i)}$, and otherwise let $\overline{F_1}^{(i)}$ be obtained from ${F_1}^{(i)} + \{x,y,v,z\}$ by deleting $\{z,z'\}$ (respectively, $\{v,v'\}$) and adding $\{z_1,z_2\}$ (respectively, $\{v_2,w\}$) when $z' \in {F_1}^{(i)}$ (respectively,  $y' \in {F_1}^{(i)}$). 
Therefore, $F^{(i)}$ is an induced forest in $G$, showing $a(G) \geq |F_1^{(i)}| + |F_2^{(i)}| + 4 - |\overline{B_i}|$. 
By Lemma \ref{ineq2}(2) (with $a = |G_1|-8, a_1=|G_2|-2$), 
$(4(|G_1|-8)+3,4(|G_2|-2)+3) \equiv (4,0),(0,4) \mod 7$.

\medskip
{\it Subcase 4.1:} $(4(|G_1|-8)+3,4(|G_2|-2)+3) \equiv (4,0) \mod 7$.

Let $F_1^{(1)} = A((G_1 - \{x,y,z,v\})/\{z_1z_2,y_1y_2,wv_2\} + v'z')$ with $z'$ (respectively, $y',v'$) as the identification of $\{z_1,z_2\}$ (respectively, $\{y_1,y_2\},\{w,v_2\}$), and $F_2^{(1)} = A(G_2)$. 
Then $|F_1^{(1)}| \geq \lceil (4(|G_1|-7)+3)/7 \rceil$, and $|F_2^{(1)}| \geq \lceil (4|G_2|+3)/7 \rceil$. 
Let $F^{(1)} = \overline{F_1}^{(1)} \cup F_2^{(1)} - (\{y_1,y_2,b_1\} \cap (\overline{F_1}^{(1)} \triangle F_2^{(1)}))$,
where $\overline{F_1}^{(1)} = {F_1}^{(1)}  + \{x,y,v,z\} $ if $z',y',v' \not\in F'$, and otherwise, $\overline{F_1}^{(1)}$  obtained from ${F_1}^{(1)}  + \{x,y,v,z\}$ by deleting $\{z,z'\}$ (respectively, $\{y,y'\},\{v,v'\}$) and adding $\{z_1,z_2\}$ (respectively, $\{y_1,y_2\},\{v_2,w\}$) when $z' \in F'$ (respectively, $y' \in F', v' \in F'$). 
Therefore, $G[F^{(1)}]$ is an induced forest in $G$, showing $a(G) \geq |F_1^{(1)}| + | F_2^{(1)}| + 4 - 3 \geq \lceil (4n+3)/7 \rceil$, a contradiction.

\medskip
{\it Subcase 4.2:} $(4(|G_1|-8)+3,4(|G_2|-2)+3) \equiv (0,4) \mod 7$.

By Claim 3, $v_2z_2 \not\in E(G)$; so $|N(v) \cap N(z_2)| \leq 2$. 
Let $F_1^{(2)} = A((G_1 - \{z_1,x,y,z,y_2\})/vz_2 $ $ + wy_1)$ with $z'$ as the identification of $\{v,z_2\}$, and $F_2^{(2)} = A(G_2 - y_2)$. 
Then $|F_1^{(2)}| \geq \lceil (4(|G_1|-6)+3)/7 \rceil$, and $|F_2^{(2)}| \geq \lceil (4(|G_2|-1)+3)/7 \rceil$. 
Let $F = F_1^{(2)} \cup F_2^{(2)} + \{x,y,z\} - (\{y_1,b_1\} \cap (F_1^{(2)} \triangle F_2^{(2)}))$. 
Now $G[F]$ (if $z' \not\in F_1^{(2)} $) or $G[F-\{z,z'\}+\{v,z_2\}]$ (if $z' \in F_1^{(2)} $)
is an induced forest of size $|F_1^{(2)}| + | F_2^{(2)}| + 3 - 2 \geq \lceil (4n+3)/7 \rceil$, a contradiction.

\medskip
Case 5: $|N(z_1) \cap N(z_2)| > 2$, $|N(y_1) \cap N(y_2)| > 2$.

There exist $a_1 \in N(z_1) \cap N(z_2), b_1 \in N(y_1) \cap N(y_2)$ and subgraphs $G_1,G_2,G_3$ of $G$ such that 
$G_2$ is the maximal subgraph of $G$ contained in the closed region of the plane bounded by the cycle $zz_1a_1z_2z$ and containing $N(z_1) \cap N(z_2) - \{z\}$,
$G_3$ is the maximal subgraph of $G$ contained in the closed region of the plane bounded by the cycle $yy_1b_1y_2y$ and containing $N(y_1) \cap N(y_2) - \{y\}$,
and $G_1$ is obtained from $G$ by removing $G_2-\{z_1,z_2,a_1\}$ and $G_3-\{y_1,b_1,y_2\}$.
Let $A_i = \{a_1\}$ if $i=1,2$ and $A_i = \emptyset$ if $i=3,4$.
Let $B_i = \{b_1\}$ if $i=1,3$ and $B_i = \emptyset$ if $i=2,4$.
For $i \in [4]$, let $F_1^{(i)} = A((G_1 - \{x,y,z,v,y_1,y_2,z_1,z_2\} - A_i - B_i)/wv_2)$ with $v'$ as the identification of $\{w,v_2\}$,
$F_2^{(i)} = A(G_2 - \{z_1,z_2\} - A_i)$,
and $F_3^{(i)} = A(G_3 - \{y_1,y_2\} - B_i)$. 
Then $|F_1^{(i)}| \geq \lceil (4(|G_1|-9- |A_i| - |B_i|)+3)/7 \rceil$, $|F_2^{(i)}| \geq \lceil (4(|G_2|-2-|A_i|)+3)/7 \rceil$, and $|F_3^{(i)}| \geq \lceil (4(|G_3|-2-|B_i|)+3)/7 \rceil$.
Let $F = F_1^{(i)} \cup F_2^{(i)} \cup F_3^{(i)} + \{x,y,v,z\} - (\{a_1\} \cap (F_1^{(i)} \triangle F_2^{(i)})) - (\{b_1\} \cap (F_1^{(i)} \triangle F_3^{(i)}))$.  
Now $G[F]$ (if $v' \not\in F_1^{(i)}$) or $G[F - \{v',v\} + \{w,v_2\}]$ (if $v' \in F_1^{(i)}$)
 is an induced forest of size $|F_1^{(i)}| + |F_2^{(i)}| + |F_3^{(i)}| + 4 - (1-|A_i|) - (1-|B_i|)$. 
 Let $(n_1,n_2,n_3):=(4(|G_1|-9)+3,4(|G_2|-2)+3,4(|G_3|-2)+3)$. 
By Lemma \ref{ineq2}(5) (with $a = |G_1|-9, a_1 = |G_2|-2, a_2 = |G_3|-2$), 
$(n_1,n_2,n_3) \equiv (0,0,0), (1,0,0), (4,0,3), (4,3,0), (3,0,4), (4,0,4), (3,4,0), (4,4,0), (1,6,0), (1,0,6), (0,3,4)$, 
$(0,4,3), (0,4,4), (6,4,4), (4,4,6), (4,6,4) \mod 7$. 

\medskip
{\it Subcase 5.1: } $(n_1,n_2,n_3) \equiv (0,0,0), (1,0,0) \mod 7$. 

If $|N(w) \cap N(v_2)| \leq 2$, 
let $F_1^{(1)} = A((G_1 - \{x,y,z,v\})/\{z_1z_2,y_1y_2,wv_2\} + v'z')$ with $z'$ (respectively, $y',v'$) as the identification of $\{z_1,z_2\}$ (respectively, $\{y_1,y_2\},\{w,v_2\}$) and $F_2^{(1)} = A(G_2)$, and $F_3^{(1)} = A(G_3)$. 
Then $|F_1^{(1)}| \geq \lceil (4(|G_1|-7)+3)/7 \rceil$, $|F_2^{(1)}| \geq \lceil (4|G_2|+3)/7 \rceil$, and $|F_3^{(1)}| \geq \lceil (4|G_3|+3)/7 \rceil$. 
Let $F^{(1)} := G[\overline{F_1}^{(1)} \cup F_2^{(1)} \cup F_3^{(1)} - (\{z_1,z_2,a_1\} \cap (\overline{F_1}^{(1)} \triangle F_2^{(1)})) - (\{y_1,y_2,b_1\} \cap (\overline{F_1}^{(1)} \triangle F_3^{(1)}))]$
where $\overline{F_1}^{(1)} := {F_1}^{(1)} + \{x,y,v,z\}$ if $v',y',z' \not\in F_1^{(1)}$ and otherwise $\overline{F_1}^{(1)}$ obtained from ${F_1}^{(1)} + \{x,y,v,z\}$ by deleting $\{z,z'\}$ (respectively, $\{v,v'\},\{y,y'\}$) and adding $\{z_1,z_2\}$ (respectively, $\{v_2,w\}, \{y_1,y_2\}$) when $z' \in {F_1}^{(1)}$ (respectively, $v' \in {F_1}^{(1)}$, $y' \in {F_1}^{(1)}$). 
Therefore, $F^{(1)}$ is an induced forest in $G$, showing $a(G) \geq |F_1^{(1)}| + | F_2^{(1)}| + | F_3^{(1)}| + 4 - 6 \geq \lceil (4n+3)/7 \rceil$, a contradiction.

If $|N(w) \cap N(v_2)| > 2$, there exist $c_1 \in N(w) \cap N(v_2)$ and subgraphs $G_1',G_2',G_3',G_4'$ of $G$ such that
$G_2' = G_2, G_3' = G_3,$
$G_4'$ is the maximal subgraph of $G$ contained in the closed region of the plane bounded by the cycle $vwc_1v_2v$ and containing $N(w) \cap N(v_2) - \{v\}$,
and $G_1'$ is obtained from $G$ by removing $G_2-\{z_1,z_2,a_1\}$, $G_3'-\{y_1,b_1,y_2\}$ and $G_4' - \{w,v_2,v\}$.
Let $C_9 = \emptyset$ and $C_{10} = \{c_1\}$. 
For $i=9,10$, let $F_1^{(i)} = A((G_1' - \{x,y,w,z,v,v_2\} - C_i)/\{y_1y_2,z_1z_2\})$ with $y'$ (respectively, $z'$) as identification of $\{y_1,y_2\}$ (respectively, $\{z_1,z_2\}$), $F_2^{(i)} = A(G_2')$, $F_3^{(i)} = A(G_3')$, and $F_4^{(i)} = A(G_4-\{w,v_2\} - C_i)$. 
Then $|F_1^{(i)}| \geq \lceil (4(|G_1'|-8-|C_i|)+3)/7 \rceil$, $|F_2^{(i)}| \geq \lceil (4|G_2'|+3)/7 \rceil$, $|F_3^{(i)}| \geq \lceil (4|G_3'|+3)/7 \rceil$, and $|F_4^{(i)}| \geq \lceil (4(|G_4'|-2-|C_i|)+3)/7 \rceil$. 
Let $F^{(i)} := G[\overline{F_1}^{(i)} \cup F_2^{(i)} \cup F_3^{(i)} \cup F_4^{(i)}  - (\{z_1,z_2,a_1\} \cap (\overline{F_1}^{(i)} \triangle F_2^{(i)}) ) - (\{y_1,y_2,b_1\} \cap (\overline{F_1}^{(i)} \triangle F_2^{(i)}) ) - (\{c_1\} \cap (\overline{F_1}^{(i)} \triangle F_4^{(i)}) )]$
where $\overline{F_1}^{(i)} = {F_1}^{(i)} + \{x,y,z,v\}$ if $y',z' \not\in F_1^{(i)}$ and $\overline{F_1}^{(i)}$ obtained from ${F_1}^{(i)} + \{x,y,z,v\}$ by deleting $y,y'$ (respectively, $\{z,z'\}$) and adding $\{y_1,y_2\}$ (respectively, $\{z_1,z_2\}$) when $y' \in {F_1}^{(i)}$ (respectively, $z' \in {F_1}^{(i)}$). 
Therefore, $F^{(i)}$ is an induced forest in $G$, showing $a(G) \geq |F_1^{(i)}| + | F_2^{(i)}| + | F_3^{(i)}| + | F_4^{(i)}| + 4 - 6 - (1-|C_i|)$.
Let $(n_1',n_2',n_3',n_4'):= (4(|G_1'|-8)+3,4(|G_2'|-2)+3,4(|G_3'|-2)+3,4(|G_4'|-2)+3)$.
By Lemma \ref{ineq2}(2), $(n_1',n_2',n_3',n_4') \equiv (4,0,0,0),(0,0,0,4) \mod 7$. 

If $(n_1',n_2',n_3',n_4') \equiv (0,0,0,4) \mod 7$,
 let $F_1^{(11)} = A((G_1' - \{x,y,z,w\})/\{z_1z_2,y_1y_2\})$ with $z'$ (respectively, $y'$) as the identification of $\{z_1,z_2\}$ (respectively, $\{y_1,y_2\}$), $F_2^{(11)} = A(G_2')$, $F_3^{(11)} = A(G_3')$,  and $F_4^{(11)} = A(G_4' - w)$. 
 Then $|F_1^{(11)}| \geq \lceil (4(|G_1'|-6)+3)/7 \rceil$, $|F_2^{(11)}| \geq \lceil (4|G_2'|+3)/7 \rceil$, $|F_3^{(11)}| \geq \lceil (4|G_3'|+3)/7 \rceil$ and $|F_4^{(11)}| \geq \lceil (4(|G_4'|-1)+3)/7 \rceil$. 
 Let $F^{(11)} := G[ \overline{F_1}^{(11)} \cup F_2^{(11)} \cup F_3^{(11)} \cup F_4^{(11)}  - (\{z_1,z_2,a_1\} \cap (\overline{F_1}^{(11)} \triangle F_2^{(11)})) - (\{y_1,y_2,b_1\} \cap (\overline{F_1}^{(11)} \triangle F_3^{(11)}))   - (\{v_2,c_1\} \cap (\overline{F_1}^{(11)} \triangle F_4^{(11)})) ]$
 where $\overline{F_1}^{(11)} = {F_1}^{(11)}+ \{x,y,z\} $ when $z',y' \not\in F_1^{(11)}$, and otherwise, let $\overline{F_1}^{(11)}$ be obtained from ${F_1}^{(11)}$ by deleting $\{z,z'\}$ (respectively, $\{y,y'\}$) and adding $\{z_1,z_2\}$ (respectively, $\{y_1,y_2\}$) when $z' \in {F_1}^{(11)}$ (respectively, $y' \in {F_1}^{(11)}$). 
 Therefore, $F^{(11)}$ is an induced forest in $G$, showing $a(G) \geq |F_1^{(11)}| + | F_2^{(11)}| + | F_3^{(11)}| + | F_4^{(11)}| + 3 - 8 \geq \lceil (4n+3)/7 \rceil$, a contradiction. 

If $(n_1',n_2',n_3',n_4') \equiv (4,0,0,0) \mod 7$,
let $F_1^{(12)} = A((G_1' - \{x,y,z,v\})/\{z_1z_2,y_1y_2,wv_2\} + v'z')$ with $z'$ (respectively, $y',v'$) as the identification of $\{z_1,z_2\}$ (respectively, $\{y_1,y_2\},\{w,$ $v_2\}$), $F_2^{(12)} = A(G_2')$, $F_3^{(12)} = A(G_3')$, and $F_4^{(12)} =A(G_4')$. 
Then $|F_1^{(12)}| \geq \lceil (4(|G_1'|-7)+3)/7 \rceil$, $|F_2^{(12)}| \geq \lceil (4|G_2'|+3)/7 \rceil$, $|F_3^{(12)}| \geq \lceil (4|G_3'|+3)/7 \rceil$ and $|F_4^{(12)}| \geq \lceil (4|G_4'|+3)/7 \rceil$. 
Let $F^{(12)} := G[ \overline{F_1}^{(12)} \cup F_2^{(12)} \cup F_3^{(12)} \cup F_4^{(12)}  - (\{z_1,z_2,a_1\} \cap (\overline{F_1}^{(12)} \triangle F_2^{(12)})) - (\{y_1,y_2,b_1\} \cap (\overline{F_1}^{(12)} \triangle F_3^{(12)}))   - (\{w,v_2,c_1\} \cap (\overline{F_1}^{(12)} \triangle F_4^{(12)})) ]$
 where $\overline{F_1}^{(12)} = {F_1}^{(12)}+ \{x,y,z,v\} $ when $z',y',v' \not\in F_1^{(12)} $, and otherwise, let $\overline{F_1}^{(12)}$ be  obtained from ${F_1}^{(12)}$ by deleting $\{z,z'\}$ (respectively, $\{y,y'\},\{v,v'\}$) and adding $\{z_1,z_2\}$ (respectively, $\{y_1,y_2\},\{w,v_2\}$) when $z' \in {F_1}^{(12)}$ (respectively, $y',v' \in {F_1}^{(12)}$).
Therefore, $F^{(12)}$ is an induced forest in $G$, showing $a(G) \geq |F_1^{(12)}| + | F_2^{(12)}| + | F_3^{(12)}| + | F_4^{(12)}| + 4 - 9 \geq \lceil (4n+3)/7 \rceil$, a contradiction.

\medskip
{\it Subcase 5.2: } $(n_1,n_2,n_3) \equiv (3,0,4),(4,0,4),(3,4,0),(4,4,0),(0,4,4),(6,4,4),(4,4,6)$, $(4,6,4) \mod 7$. 

Let $F_1^{(2)} = A((G_1 - \{x,y,z,y_2,z_1\})/vz_2 + wy_1)$ with $v'$ as the identification of $\{v,z_2\}$, $F_2^{(2)} = A(G_2 - z_1)$, and $F_3^{(2)} = A(G_3 - y_2)$. 
Then $|F_1^{(2)}| \geq \lceil (4(|G_1|-6)+3)/7 \rceil$, $|F_2^{(2)}| \geq \lceil (4(|G_2|-1)+3)/7 \rceil$, and $|F_3^{(2)}| \geq \lceil (4(|G_3|-1)+3)/7 \rceil$.
Now $G[F_1^{(2)} \cup F_2^{(2)} \cup F_3^{(2)} + \{x,y,z\} - (\{z_2,a_1\} \cap (F_1^{(2)}  \triangle F_2^{(2)})) - (\{y_1,b_1\} \cap (F_1^{(2)} \triangle F_3^{(2)}))]$ (if $v' \not\in F_1^{(2)}$) or $G[(F_1^{(2)} - v') \cup F_2^{(2)} \cup F_3^{(2)} + \{x,y,v,z_2\} - (\{z_2,a_1\} \cap ((F_1^{(2)} \cup \{z_2\}) \triangle F_2^{(2)})) - (\{y_1,b_1\} \cap (F_1^{(2)} \triangle F_3^{(2)}))]$ (if $v' \in F_1^{(2)}$)
 is an induced forest in $G$, showing $a(G) \geq |F_1^{(2)}| + | F_2^{(2)}| + | F_3^{(2)}| + 3 - 4 \geq \lceil (4n+3)/7 \rceil$, a contradiction.

\medskip
{\it Subcase 5.3: } $(n_1,n_2,n_3) \equiv (0,4,3) \mod 7$. 

Let $F_1^{(3)} = A(G_1 - \{x,y,z,v,y_1,y_2,z_1\} + wz_2)$, $F_2^{(3)} = A(G_2-z_1)$, and $F_3^{(3)} = A(G_3 - \{y_1,y_2\})$. 
Then $|F_1^{(3)}| \geq \lceil (4(|G_1|-7)+3)/7 \rceil$, $|F_2^{(3)}| \geq \lceil (4(|G_2|-1)+3)/7 \rceil$ and $|F_3^{(3)}| \geq \lceil (4(|G_3|-2)+3)/7 \rceil$. 
Now $G[F_1^{(3)} \cup F_2^{(3)} \cup F_3^{(3)} + \{x,y,z\} - ( \{z_2,a_1\} \cap (F_1^{(3)} \triangle F_2^{(3)})) -  ( \{b_1\} \cap (F_1^{(3)} \triangle F_3^{(3)}))]$ is an induced forest in $G$, showing $a(G) \geq |F_1^{(3)}| + | F_2^{(3)}| + | F_3^{(3)}| + 3 - 3 \geq \lceil (4n+3)/7 \rceil$, a contradiction.

\medskip
{\it Subcase 5.4: } $(n_1,n_2,n_3) \equiv (0,3,4) \mod 7$. 

Let $F_1^{(4)} = A(G_1 - \{x,y,z,v,y_2,z_1,z_2\} + wy_1)$, $F_2^{(4)} = A(G_2 - \{z_1,z_2\})$, and $F_3^{(4)} = A(G_3 - y_2)$. 
Then $|F_1^{(4)}| \geq \lceil (4(|G_1|-7)+3)/7 \rceil$, $|F_2^{(4)}| \geq \lceil (4(|G_2|-2)+3)/7 \rceil$ and $|F_3^{(4)}| \geq \lceil (4(|G_3|-1)+3)/7 \rceil$. 
Now $G[F_1^{(4)} \cup F_2^{(4)} \cup F_3^{(4)} + \{x,y,z\} - (\{a_1\} \cap (F_1^{(4)} \triangle F_2^{(4)})) - (\{b_1,y_1\} \cap (F_1^{(4)} \triangle F_3^{(4)})) ]$  is an induced forest in $G$, showing $a(G) \geq |F_1^{(4)}| + | F_2^{(4)}| + | F_3^{(4)}| + 3 - 3 \geq \lceil (4n+3)/7 \rceil$, a contradiction.


\medskip
{\it Subcase 5.5: } $(n_1,n_2,n_3) \equiv (4,3,0) \mod 7$ (respectively, $(1,6,0) \mod 7$).

If $|N(w) \cap N(v_2)| \leq 2$, 
let $A_5 = \emptyset$ and $A_6 = \{a_1\}$. 
For $i=5$ (respectively, $i=6$), let $F_1^{(i)} = A((G_1 - \{x,y,z,v,z_1,z_2\} - A_i)/\{y_1y_2,wv_2\})$ with $y'$ (respectively, $v'$) as the identification of  $\{y_1,y_2\}$ (respectively, $\{w,v_2\}$),
 $F_2^{(i)} = A(G_2 - \{z_1,z_2\} - A_i)$ 
 and $F_3^{(i)} = A(G_3)$.
Then $|F_1^{(i)}| \geq \lceil (4(|G_1|-8- |A_i|)+3)/7 \rceil$, $|F_2^{(i)}| \geq \lceil (4(|G_2|-2-|A_i|)+3)/7 \rceil$ and $|F_3^{(i)}| \geq \lceil (4|G_3|+3)/7 \rceil$. 
Let $F^{(i)} := G[\overline{F_1}^{(i)} \cup F_2^{(i)} \cup F_3^{(i)} - (\{a_1\} \cap (\overline{F_1}^{(i)} \triangle F_2^{(i)})) - (\{y_1,y_2,b_1\} \cap (\overline{F_1}^{(i)} \triangle F_3^{(i)}))]$,
where $\overline{F_1}^{(i)} := {F_1}^{(i)} + \{x,y,v,z\}$ if $v',y' \not\in F_1^{(1)}$, and otherwise, $\overline{F_1}^{(1)}$ obtained from ${F_1}^{(i)} + \{x,y,v,z\}$ by deleting $\{y,y'\}$ (respectively, $\{v,v'\}$) and adding $\{y_1,y_2\}$ (respectively, $\{v_2,w\}$) when $y' \in {F_1}^{(i)}$ (respectively, $v' \in {F_1}^{(i)}$). 
Therefore, $F^{(i)}$ is an induced forest in $G$, showing $a(G) \geq |F_1^{(i)}| + | F_2^{(i)}| + | F_3^{(i)}| + 4 - 3 - (1-|A_i|) \geq \lceil (4n+3)/7 \rceil$, a contradiction.

If $|N(w) \cap N(v_2)| > 2$,
there exist $c_1 \in N(w) \cap N(v_2)$ and subgraphs $G_1',G_2',G_3',G_4'$ of $G$ such that
$G_2' = G_2, G_3' = G_3,$
$G_4'$ is the maximal subgraph of $G$ contained in the closed region of the plane bounded by the cycle $vwc_1v_2v$ and containing $N(w) \cap N(v_2) - \{v\}$,
and $G_1'$ is obtained from $G$ by removing $G_2-\{z_1,z_2,a_1\}$, $G_3'-\{y_1,b_1,y_2\}$ and $G_4' - \{w,v_2,v\}$.
Let $A_i = \{a_1\}$ if $i=13,14,17,19$ and $A_i =\emptyset$ if $i=15,16,18,20$. 
Let $C_i = \{c_1\}$ if $i=13,15$ and $C_i= \emptyset$ if $i=14,16$. 
For $i=13,14,15,16$, let $F_1^{(i)} = A((G_1' - \{x,y,w,z,v,v_2,z_1,z_2\} - A_i - C_i)/y_1y_2)$ with $y'$ as the identification of $\{y_1,y_2\}$, $F_2^{(i)} = A(G_2' - \{z_1,z_2\} - A_i)$, $F_3^{(i)} = A(G_3')$, and $F_4^{(i)} = A(G_4'-\{w,v_2\} - C_i)$. 
Note $|F_1^{(i)}| \geq \lceil (4(|G_1'|-9-|A_i|-|C_i|)+3)/7 \rceil$, $|F_2^{(i)}| \geq \lceil (4(|G_2'|-2-|A_i|)+3)/7 \rceil$, $|F_3^{(i)}| \geq \lceil (4|G_3'|+3)/7 \rceil$, and $|F_4^{(i)}| \geq \lceil (4(|G_4'|-2-|C_i|)+3)/7 \rceil$. 
Now $G[F_1^{(i)} \cup F_2^{(i)} \cup F_3^{(i)} \cup F_4^{(i)} + \{x,y,z,v\} - (\{y_1,y_2,b_1\} \cap ({F_1}^{(i)} \triangle F_3^{(i)}))  - (\{a_1\} \cap ({F_1}^{(i)} \triangle F_2^{(i)}))  - (\{c_1\} \cap ({F_1}^{(i)} \triangle F_4^{(i)}))]$ (if $y' \not\in F_1^{(i)}$) 
or $G[(F_1^{(i)} - y') \cup F_2^{(i)} \cup F_3^{(i)} \cup F_4^{(i)} + \{x,y_1,y_2,z,v\} - (\{y_1,y_2,b_1\} \cap (({F_1}^{(i)}+\{y_1,y_2\}) \triangle F_3^{(i)}))  - (\{a_1\} \cap ({F_1}^{(i)} \triangle F_2^{(i)}))  - (\{c_1\} \cap ({F_1}^{(i)} \triangle F_4^{(i)}))]$ (if $y' \in F_1^{(i)}$)
is an induced forest in $G$, showing $a(G) \geq |F_1^{(i)}| + | F_2^{(i)}| + | F_3^{(i)}| + | F_4^{(i)}| + 4 - 3  - (1-|A_i|) - (1-|C_i|)$. 
Let $(n_1',n_2',n_3',n_4'):= (4(|G_1'|-9)+3,4(|G_2'|-2)+3,4(|G_3'|-2)+3,4(|G_4'|-2)+3)$.
By Lemma \ref{ineq2}(2), $(n_1',n_2',n_3',n_4') \equiv (4,3,0,0),(4,6,0,0),(0,3,0,4),(0,6,0,4) \mod 7$.

If $(n_1',n_2',n_3',n_4') \equiv (4,6,0,0) \mod 7$ (respectively, $(4,3,0,0)\mod 7$),
for $i=17$ (respectively, $i=18$), let $F_1^{(i)} = A((G_1' - \{x,y,z,v,z_1,z_2\} - A_i)/\{y_1y_2,wv_2\})$ with $y'$ (respectively, $v'$) as the identification of $\{y_1,y_2\}$ (respectively $\{w,v_2\}$), $F_2^{(i)} = A(G_2'-\{z_1,z_2\} - A_i)$, $F_3^{(i)} = A(G_3')$, and $F_4^{(i)} = A(G_4')$. 
Then $|F_1^{(i)}| \geq \lceil (4(|G_1'|-8-|A_i|)+3)/7 \rceil$, $|F_2^{(i)}| \geq \lceil (4(|G_2'|-2-|A_i|)+3)/7 \rceil$, $|F_3^{(i)}| \geq \lceil (4|G_3'|+3)/7 \rceil$, and $|F_4^{(i)}| \geq \lceil (4|G_4'|+3)/7 \rceil$. 
Let $F^{(i)} := G[\overline{F_1}^{(i)} \cup F_2^{(i)} \cup F_3^{(i)} \cup F_4^{(i)}  - (\{y_1,y_2,b_1\} \cap (\overline{F_1}^{(i)} \triangle F_3^{(i)})) - (\{w,v_2,c_1\} \cap (\overline{F_1}^{(i)} \triangle F_3^{(i)})) - (\{a_1\} \cap (\overline{F_1}^{(i)} \triangle F_4^{(i)}))]$,
where $\overline{F_1}^{(i)} = {F_1}^{(i)} + \{x,y,z,v\}$ if $y',v' \not\in F_1^{(i)}$, and otherwise, $\overline{F_1}^{(i)}$ obtained from ${F_1}^{(i)} + \{x,y,z,v\}$ by deleting $\{y,y'\}$ (respectively, $\{v,v'\}$) and adding $\{y_1,y_2\}$ (respectively,  $\{v_2,w\}$) when $y' \in {F_1}^{(i)}$ (respectively, $v' \in {F_1}^{(i)}$).
Therefore, $F^{(i)}$ is an induced forest in $G$, showing $a(G) \geq |F_1^{(i)}| + | F_2^{(i)}| + | F_3^{(i)}| + | F_4^{(i)}| + 4 - 6  - (1-|A_i|) \geq \lceil (4n+3)/7 \rceil$, a contradiction. 

If $(n_1',n_2',n_3',n_4') \equiv (0,6,0,4) \mod 7$ (respectively, $(0,3,0,4)\mod 7$),
for $i=19$ (respectively, $i=20$), let $F_1^{(i)} = A((G_1' - \{x,y,z,w,z_1,z_2\} - A_i)/y_1y_2)$ with $y'$ as the identification of $\{y_1,y_2\}$, $F_2^{(i)} = A(G_2'-\{z_1,z_2\} - A_i)$, $F_3^{(i)} = A(G_3')$, and $F_4^{(i)} = A(G_4' - w)$. 
Then $|F_1^{(i)}| \geq \lceil (4(|G_1'|-7-|A_i|)+3)/7 \rceil$, $|F_2^{(i)}| \geq \lceil (4(|G_2'|-2-|A_i|)+3)/7 \rceil$, $|F_3^{(i)}| \geq \lceil (4|G_3'|+3)/7 \rceil$, and $|F_4^{(i)}| \geq \lceil (4(|G_4'|-1)+3)/7 \rceil$. 
Now $G[F_1^{(i)} \cup F_2^{(i)} \cup F_3^{(i)} \cup F_4^{(i)} + \{x,y,z\} - (\{a_1\} \cap (F_1^{(i)} \triangle F_2^{(i)})) - (\{y_1,y_2,b_1\} \cap (F_1^{(i)} \triangle F_3^{(i)})) - (\{v_2,c_1\} \cap (F_1^{(i)} \triangle F_4^{(i)})) ]$ (if $y' \not\in F_1^{(i)}$)
or $G[(F_1^{(i)} - y') \cup F_2^{(i)} \cup F_3^{(i)} \cup F_4^{(i)} + \{x,y_1,y_2,z\} - (\{a_1\} \cap (F_1^{(i)} \triangle F_2^{(i)})) - (\{y_1,y_2,b_1\} \cap ((F_1^{(i)} + \{y_1,y_2\}) \triangle F_3^{(i)})) - (\{v_2,c_1\} \cap (F_1^{(i)} \triangle F_4^{(i)})) ]$ (if $y' \in F_1^{(i)}$)
is an induced forest in $G$, showing $a(G) \geq |F_1^{(i)}| + | F_2^{(i)}| + | F_3^{(i)}| + | F_4^{(i)}| + 3 - 5  - (1-|A_i|) \geq \lceil (4n+3)/7 \rceil$, a contradiction.

\medskip
{\it Subcase 5.6: } $(n_1,n_2,n_3) \equiv (4,0,3) \mod 7$ (respectively, $(1,0,6) \mod 7$).

If $|N(v_2) \cap N(w)| \leq 2$, 
let $B_7 = \emptyset$ and $B_8 = \{b_1\}$. 
For $i=7$ (respectively, $i=8$), let $F_1^{(i)} = A((G_1 - \{x,y,z,v,y_1,y_2\} - B_i)/\{z_1z_2,wv_2\} + z'v')$ with $z'$ (respectively, $v'$) as the identification of $\{z_1,z_2\}$ (respectively, $\{w,v_2\}$), $F_2^{(i)} = A(G_2)$, and $F_3^{(i)} = A(G_3 - \{y_1,y_2\} - B_1)$. 
Then $|F_1^{(i)}| \geq \lceil (4(|G_1|-8- |B_i|)+3)/7 \rceil$, $|F_2^{(i)}| \geq \lceil (4|G_2|+3)/7 \rceil$, and $|F_3^{(i)}| \geq \lceil (4(|G_3|-2-|B_1|)+3)/7 \rceil$. 
Let $F^{(i)} := G[\overline{F_1}^{(i)} \cup F_2^{(i)} \cup F_3^{(i)} - (\{b_1\} \cap (\overline{F_1}^{(i)} \triangle F_3^{(i)})) - (\{z_1,z_2,a_1\} \cap (\overline{F_1}^{(i)} \triangle F_2^{(i)}))]$,
where $\overline{F_1}^{(i)} := {F_1}^{(i)} + \{x,y,v,z\}$ if $z',v' \not\in F_1^{(1)}$, and otherwise, $\overline{F_1}^{(1)}$ obtained from ${F_1}^{(i)} + \{x,y,v,z\}$ by deleting $\{z,z'\}$ (respectively, $\{v,v'\}$) and adding $\{z_1,z_2\}$ (respectively, $\{v_2,w\}$) when $z' \in {F_1}^{(i)}$ (respectively, $v' \in {F_1}^{(i)}$). 
Therefore, $F^{(i)}$ is an induced forest in $G$, showing $a(G) \geq |F_1^{(i)}| + | F_2^{(i)}| + | F_3^{(i)}| + 4 - 3 - (1-|B_i|) \geq \lceil (4n+3)/7 \rceil$, a contradiction.

If $|N(v_2) \cap N(w)| > 2$,
there exist $c_1 \in N(w) \cap N(v_2)$ and subgraphs $G_1',G_2',G_3',G_4'$ of $G$ such that
$G_2' = G_2, G_3' = G_3,$
$G_4'$ is the maximal subgraph of $G$ contained in the closed region of the plane bounded by the cycle $vwc_1v_2v$ and containing $N(w) \cap N(v_2) - \{v\}$,
and $G_1'$ is obtained from $G$ by removing $G_2-\{z_1,z_2,a_1\}$, $G_3'-\{y_1,b_1,y_2\}$ and $G_4' - \{w,v_2,v\}$.
Let $B_1 = \{b_1\}$ if $i=21,22,25,27$ and $\emptyset$ if $i=23,24,26,28$.
Let $C_1 = \{c_1\}$ if $i=21,23$ and $\emptyset$ if $i=22,24$.
 For $i=21,22,23,24$, let $F_1^{(i)} = A((G_1' - \{x,y,w,z,v,v_2,y_1,y_2\} - B_i - C_i)/z_1z_2)$ with $z'$ as the identification of $\{z_1,z_2\}$, $F_2^{(i)} = A(G_2')$, $F_3^{(i)} = A(G_3' - \{y_1,y_2\} - B_i)$, and $F_4^{(i)} = A(G_4'-\{w,v_2\} - C_i)$. 
 Then $|F_1^{(i)}| \geq \lceil (4(|G_1'|-9-|A_i|-|C_i|)+3)/7 \rceil$, $|F_2^{(i)}| \geq \lceil (4|G_2'|+3)/7 \rceil$, $|F_3^{(i)}| \geq \lceil (4(|G_3'|-2-|B_i|)+3)/7 \rceil$ and $|F_4^{(i)}| \geq \lceil (4(|G_4'|-2-|C_i|)+3)/7 \rceil$. 
 Now $G[F_1^{(i)} \cup F_2^{(i)} \cup F_3^{(i)} \cup F_4^{(i)} + \{x,y,z,v\} - (\{z_1,z_2,a_1\} \cap (F_1^{(i)} \triangle F_2^{(i)})) - (\{b_1\} \cap (F_1^{(i)} \triangle F_3^{(i)})) - (\{c_1\} \cap (F_1^{(i)} \triangle F_4^{(i)}))$ (if $z' \not\in F_1^{(i)}$ )
 or  $G[(F_1^{(i)}-z') \cup F_2^{(i)} \cup F_3^{(i)} \cup F_4^{(i)} + \{x,y,z_1,z_2,v\} - (\{z_1,z_2,a_1\} \cap ( (F_1^{(i)} + \{z_1,z_2\}) \triangle F_2^{(i)})) - (\{b_1\} \cap (F_1^{(i)} \triangle F_3^{(i)})) - (\{c_1\} \cap (F_1^{(i)} \triangle F_4^{(i)}))$ (if $z' \in F_1^{(i)}$ )
  is an induced forest of size $|F_1^{(i)}| + | F_2^{(i)}| + | F_3^{(i)}| + | F_4^{(i)}| + 4 - 3  - (1-|B_i|) - (1-|C_i|)$. 
 Let $(n_1',n_2',n_3',n_4') := (4(|G_1'|-2)+3,4(|G_2'|-2)+3,4(|G_3'|-2)+3,4(|G_4'|-2)+3)$.
By Lemma \ref{ineq2}(2), $(n_1',n_2',n_3',n_4') \equiv (4,0,3,0),(4,0,6,0),(0,0,3,4),(0,0,6,4) \mod 7$.

If $(n_1',n_2',n_3',n_4') \equiv (4,0,6,0) \mod 7$ (respectively, $(4,0,3,0)\mod 7$),
 for $i=25$ (respectively, $i=26$), let $F_1^{(i)} = A((G_1' - \{x,y,z,v,y_1,y_2\} - B_i)/\{z_1z_2,wv_2\} + z'v')$ with $z'$ (respectively, $v'$) as the identification of $\{z_1z_2\}$ (respectively $\{w,v_2\}$), $F_2^{(i)} = A(G_2')$, $F_3^{(i)} = A(G_3'-\{z_1,z_2\} - B_i)$, and $F_4^{(i)} = A(G_4')$. 
 Then $|F_1^{(i)}| \geq \lceil (4(|G_1'|-8-|B_i|)+3)/7 \rceil$, $|F_2^{(i)}| \geq \lceil (4|G_2'|+3)/7 \rceil$, $|F_3^{(i)}| \geq \lceil (4(|G_3'|-2-|B_i|)+3)/7 \rceil$ and $|F_4^{(i)}| \geq \lceil (4|G_4'|+3)/7 \rceil$. 
 Now $F^{(i)} := G[\overline{F_1}^{(i)} \cup F_2^{(i)} \cup F_3^{(i)} \cup F_4^{(i)}   - (\{z_1,z_2,a_1\} \cap (\overline{F_1}^{(i)} \triangle F_3^{(i)}))  - (\{b_1\} \cap (\overline{F_1}^{(i)} \triangle F_3^{(i)})) - (\{w,v_2,c_1\} \cap (\overline{F_1}^{(i)} \triangle F_4^{(i)}))]$,
 where $\overline{F_1}^{(i)} = F_1^{(i)} + \{x,y,z,v\}$ if $z',v' \not\in F_1^{(i)}$, and otherwise, let $\overline{F_1}^{(i)}$ be obtained from $F_1^{(i)} + \{x,y,z,v\}$ by deleting $\{z,z'\}$ (respectively, $\{v,v'\}$) and adding $\{z_1,z_2\}$ (respectively, $\{v_2,w\}$) when $z' \in F_1^{(i)}$ (respectively, $v' \in F_1^{(i)}$). 
 Therefore, $F^{(i)}$ is an induced forest in $G$, showing $a(G) \geq |F_1^{(i)}| + | F_2^{(i)}| + | F_3^{(i)}| + | F_4^{(i)}| + 4 - 6  - (1-|B_i|) \geq \lceil (4n+3)/7 \rceil$, a contradiction.
 
 If $(n_1',n_2',n_3',n_4') \equiv (0,0,6,4) \mod 7$ (respectively, $(0,0,3,4)\mod 7$),
for $i=27,28$, let $F_1^{(i)} = A(G_1' - \{x,y,z,w,y_1,y_2\} - B_i)/z_1z_2)$ with $z'$ as the identification of $\{z_1,z_2\}$, $F_2^{(i)} = A(G_2')$, $F_3^{(i)} = A(G_3'-\{z_1,z_2\} - B_i)$, and $F_4^{(i)} = A(G_4' - w)$. 
Then $|F_1^{(i)}| \geq \lceil (4(|G_1'|-7-|B_i|)+3)/7 \rceil$, $|F_2^{(i)}| \geq \lceil (4|G_2'|+3)/7 \rceil$, $|F_3^{(i)}| \geq \lceil (4(|G_3'|-2-|B_i|)+3)/7 \rceil$ and $|F_4^{(i)}| \geq \lceil (4(|G_4'|-1)+3)/7 \rceil$. 
Now $G[F_1^{(i)} \cup F_2^{(i)} \cup F_3^{(i)} \cup F_4^{(i)} + \{x,y,z\} - (\{z_1,z_2,a_1\} \cap (F_1^{(i)} \triangle F_2^{(i)})) - (\{b_1\} \cap (F_1^{(i)} \triangle F_3^{(i)})) - (\{v_2,c_1\} \cap (F_1^{(i)} \triangle F_4^{(i)}))]$  (if $z' \not\in F_1^{(i)}$ )
or $G[(F_1^{(i)} - z') \cup F_2^{(i)} \cup F_3^{(i)} \cup F_4^{(i)} + \{x,y,z_1,z_2\} - (\{z_1,z_2,a_1\} \cap ((F_1^{(i)} + \{z_1,z_2\}) \triangle F_2^{(i)})) - (\{b_1\} \cap (F_1^{(i)} \triangle F_3^{(i)})) - (\{v_2,c_1\} \cap (F_1^{(i)} \triangle F_4^{(i)}))]$  (if $z' \in F_1^{(i)}$ )
is an induced forest in $G$, showing $a(G) \geq |F_1^{(i)}| + | F_2^{(i)}| + | F_3^{(i)}| + | F_4^{(i)}| + 3 - 5  - (1-|B_i|) \geq \lceil (4n+3)/7 \rceil$, a contradiction.

In the second case, for $i=19,20$, let $F_1^{(i)} = A(G_1 - \{x,y,z,w,y_1,y_2\} - B_1)/z_1z_2)$ with $z'$ as the identification of $\{z_1,z_2\}$, $F_2^{(i)} = F_2^{(3)}$, $F_3^{(i)} = F_3^{(i-4)}$, and $F_4^{(i)} = F_4^{(2)}$. Note $|F_1^{(i)}| \geq \lceil (4(|G_1|-7-|B_1|)+3)/7 \rceil$. Define $F^{(i)} := G[F_1^{(i)} \cup_m F_2^{(i)} \cup_m F_3^{(i)} \cup_m F_4^{(i)} \cup_m \{x,y,z\} -_m \{z_1,z_2,a_1,v_2,c_1\} -_m \overline{B_1}]$ if $z' \not\in F_1^{(i)}$ and else $F^{(i)} := G[F_1^{(i)} \cup_m F_2^{(i)} \cup_m F_3^{(i)} \cup_m F_4^{(i)} \cup_m \{x,y,z_1,z_2\} -_m \{z'\} -_m \{z_1,z_2,a_1,v_2,c_1\} -_m \overline{B_1}]$. So $F^{(i)}$ is an induced forest of size $|F_1^{(i)}| + | F_2^{(i)}| + | F_3^{(i)}| + | F_4^{(i)}| + 3 - 5  - |\overline{B_1}| \geq \lceil (4n+3)/7 \rceil$, a contradiction.

\medskip

Case 6: $|N(z_1) \cap N(z_2)| > 2$, $|N(y_1) \cap N(y_2)| \leq 2$ and $|N(w) \cap N(v_2)| > 2$. 

There exist $a_1 \in N(z_1) \cap N(z_2)$, $c_1 \in N(w) \cap N(v_2)$ and subgraphs $G_1,G_2,G_3$ of $G$ such that 
$G_2$ is the maximal subgraph of $G$ contained in the closed region of the plane bounded by the cycle $zz_1a_1z_2z$ and containing $N(z_1) \cap N(z_2) - \{z\}$,
$G_3$ is the maximal subgraph of $G$ contained in the closed region of the plane bounded by the cycle $vwc_1v_2v$ and containing $N(w) \cap N(v_2) - \{v\}$,
and $G_1$ is obtained from $G$ by removing $G_2-\{z_1,z_2,a_1\}$ and $G_3-\{w,c_1,v_2\}$.
Let $A_i = \{a_1\}$ if $i=1,2$ and $A_i = \emptyset$ if $i=3,4$.
Let $C_i = \{c_1\}$ if $i=1,3$ and $C_i = \emptyset$ if $i=2,4$.
For $i \in [4]$, let $F_1^{(i)} = A((G_1 - \{x,y,z,v,w,v_2,z_1,z_2\} - A_i - C_i)/y_1y_2)$ with $y'$ as the identification of $\{y_1,y_2\}$, $F_2^{(i)} = A(G_2 - \{z_1,z_2\} - A_i)$, and $F_3^{(i)} = A(G_3 - \{w,v_2\} - C_i)$. 
Then $|F_1^{(i)}| \geq \lceil (4(|G_1|-9- |A_i| - |C_i|)+3)/7 \rceil$, $|F_2^{(i)}| \geq \lceil (4(|G_2|-2-|A_i|)+3)/7 \rceil$, and $|F_3^{(i)}| \geq \lceil (4(|G_3|-2-|C_i|)+3)/7 \rceil$. 
Let $F = F_1^{(i)} \cup F_2^{(i)} \cup F_3^{(i)} + \{x,y,v,z\} - (\{a_1\} \cap (F_1^{(i)} \triangle F_2^{(i)})) - (\{c_1\} \cap (F_1^{(i)} \triangle F_3^{(i)}))$. 
Now $G[F]$ (if $y' \not\in F_1^{(i)}$) or $G[F - \{y,y'\}+\{y_1,y_2\}]$ (if $y' \in F_1^{(i)}$)
is an induced forest of size $a(G) \geq |F_1^{(i)}| + |F_2^{(i)}| + |F_3^{(i)}| + 4 - (1-|A_1|) - (1-|C_1|)$.
Let $(n_1,n_2,n_3):= (4(|G_1|-9)+3,4(|G_2|-2)+3,4(|G_3|-2)+3)$.
By Lemma \ref{ineq2}(5) (with $a = |G_1|-9, a_1 = |G_2|-2, a_2 = |G_3|-2, c = 4$), 
$(n_1,n_2,n_3) \equiv (0,0,0), (1,0,0), (4,0,3), (4,3,0), (3,0,4), (4,0,4), (3,4,0), (4,4,0), (1,6,0), (1,0,6), (0,3,4)$, 
$(0,4,3), (0,4,4), (6,4,4), (4,4,6), (4,6,4) \mod 7$. 

\medskip




{\it Subcase 6.1: } $|N(v_1) \cap N(v_2)| > 2$. 

There exist $d_1 \in N(v_1) \cap N(v_2)$ and subgraphs $G_1',G_3',G_4'$ such that
$G_3' = G_3$,
$G_4'$ is the maximal subgraph of $G$ contained in the closed region of the plane bounded by the cycle $vv_1d_1v_2v$ and containing $N(v_1) \cap N(v_2) - \{v\}$,
and $G_1'$ is obtained from $G$ by removing $G_3'-\{w,c_1,v_2\}$ and $G_4'-\{v_1,v_2,v\}$.
Let $C_i = \{c_1\}$ if $i=1,2$ and $C_i = \emptyset$ if $i=3,4$.
Let $D_i = \{d_1\}$ if $i=1,3$ and $D_i = \emptyset$ if $i=2,4$.
For $i \in [4]$, let $F_1^{(i)} = A(G_1' - \{x,z,z_1,v,w,v_1,v_2\} - C_i - D_i + yz_2)$, $F_2^{(i)} = A(G_3' - \{w,v_2\} - C_i)$, and $F_3^{(i)} = A(G_4' - \{v_1,v_2\} - D_i)$. 
Then $|F_1^{(i)}| \geq \lceil (4(|G_1'|-7 - |C_i| - |D_i|)+3)/7 \rceil$, $|F_2^{(i)}| \geq \lceil (4(|G_3'|-2-|C_i|)+3)/7 \rceil$, and $|F_3^{(i)}| \geq \lceil (4(|G_4'|-2-|D_i|)+3)/7 \rceil$. 
Now $G[F_1^{(i)} \cup F_2^{(i)} \cup F_3^{(i)} + \{x,z,v\} - (\{c_1\} \cap (F_1^{(i)} \triangle F_2^{(i)})) - (\{d_1\} \cap (F_1^{(i)} \triangle F_3^{(i)}))]$ is an induced forest in $G$, showing $a(G)\geq |F_1^{(i)}| + | F_2^{(i)}| + | F_3^{(i)}| + 3 - (1-|C_i|) - (1-|D_i|)$. 
Let $(n_1',n_2',n_3'):=(4(|G_1'|-7)+3,4(|G_3'|-2)+3,4(|G_4'|-2)+3)$.
By Lemma \ref{ineq2}(4) (with $a = |G_1'|-7, a_1 = |G_3'|-2, a_2 = |G_4'|-2, c=3$),
 $(n_1',n_2',n_3') \equiv (1,0,0),(0,4,4),(4,4,0),(4,0,4) \mod 7$. 

If $(n_1',n_2',n_3') \equiv (0,4,4) \mod 7$, 
let $F_1^{(5)} = A(G_1' - \{x,z,w,v,v_1\})$, $F_2^{(5)} = A(G_3' - w)$, and $F_3^{(5)} = A(G_4' - v_1)$. 
Then $|F_1^{(5)}| \geq \lceil (4(|G_1'|-5)+3)/7 \rceil$, $|F_2^{(5)}| \geq \lceil (4(|G_3'|-1)+3)/7 \rceil$, and $|F_3^{(5)}| \geq \lceil (4(|G_4'|-1)+3)/7 \rceil$. Now $G[F_1^{(5)} \cup F_2^{(5)} \cup F_3^{(5)} + \{x,v\} - (\{c_1,v_2\} \cap (F_1^{(5)} \triangle F_2^{(5)})) - (\{d_1,v_2\} \cap (F_1^{(5)} \triangle F_3^{(5)}))]$  is an induced forest in $G$, showing $a(G) \geq |F_1^{(5)}| + | F_2^{(5)}| + | F_3^{(5)}| + 2 - 4 \geq \lceil (4n+3)/7 \rceil$, a contradiction. 

If $(n_1',n_2',n_3') \equiv (4,4,0) \mod 7$, 
let $F_1^{(6)} = A((G_1' - \{x,z,w,v,z_1\})/v_1v_2 + yz_2)$ with $v'$ as the identification of $\{v_1,v_2\}$, $F_2^{(6)} = A(G_3' - w)$, and $F_3^{(6)} = A(G_4')$.
Then $|F_1^{(6)}| \geq \lceil (4(|G_1'|-6)+3)/7 \rceil$, $|F_2^{(6)}| \geq \lceil (4(|G_3'|-1)+3)/7 \rceil$, and $|F_3^{(6)}| \geq \lceil (4|G_4'|+3)/7 \rceil$. 
Now $G[F_1^{(6)} \cup F_2^{(6)} \cup F_3^{(6)} + \{x,v,z\} - (\{v_2,c_1\} \cap (F_1^{(6)} \triangle F_2^{(6)})) - (\{v_1,v_2,d_1\} \cap (F_1^{(6)} \triangle F_3^{(6)}))]$ (if $v' \not\in F_1^{(6)}$) or $G[(F_1^{(6)} - v') \cup F_2^{(6)} \cup F_3^{(6)} + \{x,v_1,v_2,z\} - (\{v_2,c_1\} \cap ( (F_1^{(6)} + v_2) \triangle F_2^{(6)})) - (\{v_1,v_2,d_1\} \cap ( (F_1^{(6)} + \{v_1,v_2\}) \triangle F_3^{(6)}))]$ (if $v' \not\in F_1^{(6)}$)
is an induced forest in $G$, showing $a(G) \geq |F_1^{(6)}| + | F_2^{(6)}| + | F_3^{(6)}| + 3 - 5 \geq \lceil (4n+3)/7 \rceil$, a contradiction. 

If $(n_1',n_2',n_3') \equiv (4,0,4) \mod 7$, 
let $F_1^{(7)} = A(G_1' - \{x,z,w,v,v_2,c_1\})$, $F_2^{(7)} = A(G_3' - \{w,v_2,c_1\})$, and $F_3^{(7)} = A(G_4' - \{v_2\})$. 
Then $|F_1^{(7)}| \geq \lceil (4(|G_1'|-6)+3)/7 \rceil$, $|F_2^{(7)}| \geq \lceil (4(|G_3'|-3)+3)/7 \rceil$, and $|F_3^{(7)}| \geq \lceil (4(|G_4'|-1)+3)/7 \rceil$.
Now $G[F_1^{(7)} \cup F_2^{(7)} \cup F_3^{(7)} + \{x,v\} - (\{v_1,d_1\} \cap (F_1^{(7)} \triangle F_3^{(7)}))]$ is an induced forest in $G$, showing $|F_1^{(7)}| + | F_2^{(7)}| + | F_3^{(7)}| + 2 - 2 \geq \lceil (4n+3)/7 \rceil$, a contradiction. 




If $(n_1',n_2',n_3') \equiv (1,0,0) \mod 7$, 
then there exist $a_1 \in N(z_1) \cap N(z_2)$ and subgraphs $G_1'',G_2'',G_3'',G_4''$ of $G$ such that 
$G_3'' = G_3', G_4'' = G_4'$,
$G_2''$ is the maximal subgraph of $G$ contained in the closed region of the plane bounded by the cycle $zz_1a_1z_2z$ and containing $N(v_1) \cap N(v_2) - \{v\}$,
and $G_1''$ is obtained from $G$ by removing $G_2'' - \{z_1,a_1,z_2\}$, $G_3''-\{w,c_1,v_2\}$ and $G_4''-\{v_1,v_2,v\}$.
Let $A_8 = \{a_1\}$ and $A_9 = \emptyset$.
For $i=8,9$, let $F_1^{(i)} = A((G_1'' - \{x,y,z,v,z_1,z_2\} - A_i)/\{y_1y_2,wv_2\})$ with $y'$ (respectively, $v'$) as the identification of $\{y_1,y_2\}$ (respectively, $\{w,v_2\}$), $F_2^{(i)} = A(G_2'' - \{z_1,z_2\} - A_i)$, $F_3^{(i)} = A(G_3'')$, and $F_4^{(i)} = A(G_4'')$. 
Then $|F_1^{(i)}| \geq \lceil (4(|G_1''|-8-|A_i|)+3)/7 \rceil$, $|F_2^{(i)}| \geq \lceil (4(|G_2''|-2-|A_i|)+3)/7 \rceil$, $|F_3^{(i)}| \geq \lceil (4|G_3''|+3)/7 \rceil$, and $|F_4^{(i)}| \geq \lceil (4|G_4''|+3)/7 \rceil$. 
Let $F^{(i)} := G[\overline{F_1}^{(i)} \cup F_2^{(i)} \cup F_3^{(i)} \cup F_4^{(i)}     - (\{a_1\} \cap (\overline{F_1}^{(i)} \triangle F_2^{(i)})) - (\{w,v_2,c_1\} \cap (\overline{F_1}^{(i)} \triangle F_3^{(i)})) - (\{v_1,v_2,d_1\} \cap (\overline{F_1}^{(i)} \triangle F_4^{(i)}))]$,
where $\overline{F_1}^{(i)} = {F_1}^{(i)} + \{x,y,v,z\}$ if $y',v' \not\in F_1^{(i)}$, and otherwise, $\overline{F_1}^{(i)}$ obtained from ${F_1}^{(i)} + \{x,y,v,z\}$ by deleting $\{y,y'\}$ (respectively, $\{v,v'\}$) and adding $\{y_1,y_2\}$ (respectively, $\{w,v_2\}$) when $y' \in {F_1}^{(i)}$ (respectively, $v' \in {F_1}^{(i)}$). 
Therefore, $F^{(i)}$ is an induced forest in $G$, showing $a(G) \geq |F_1^{(i)}| + | F_2^{(i)}| + | F_3^{(i)}| + 4 - 6 - (1-|A_i|)$. 
By Lemma \ref{ineq2}(2), $(4(|G_1''|-8)+3,4(|G_2''|-2)+3,4(|G_3''|-2)+3,4(|G_4''|-2)+3) \equiv (4,0,0,0),(0,4,0,0) \mod 7$. 

If  $(4(|G_1''|-8)+3,4(|G_2''|-2)+3,4(|G_3''|-2)+3,4(|G_4''|-2)+3) \equiv (4,0,0,0) \mod 7$, 
let $F_1^{(10)} = A((G_1'' - \{x,y,z,v\})/\{z_1z_2,y_1y_2,wv_2\} + v'z')$ with $z'$ (respectively, $y',v'$) as the identification of $\{z_1,z_2\}$ (respectively, $\{y_1,y_2\},\{w,v_2\}$) and $F_2^{(10)} = A(G_2'')$, $F_3^{(10)} = A(G_3'')$, and $F_4^{(10)} = A(G_4'')$. 
Then $|F_1^{(10)}| \geq \lceil (4(|G_1''|-7)+3)/7 \rceil$, $|F_2^{(10)}| \geq \lceil (4|G_2''|+3)/7 \rceil$, $|F_3^{(10)}| \geq \lceil (4|G_3''|+3)/7 \rceil$, and $|F_4^{(10)}| \geq \lceil (4|G_4''|+3)/7 \rceil$. 
Let $F^{(10)} := G[\overline{F_1}^{(10)} \cup F_2^{(10)} \cup F_3^{(10)} - (\{z_1,z_2,a_1\} \cap (\overline{F_1}^{(10)} \triangle F_2^{(10)})) - (\{w,v_2,c_1\} \cap (\overline{F_1}^{(10)} \triangle F_3^{(10)})) - (\{v_1,v_2,d_1\} \cap (\overline{F_1}^{(10)} \triangle F_4^{(10)}))]$,
where $\overline{F_1}^{(1)} := {F_1}^{(10)} + \{x,y,v,z\}$ if $v',y',z' \not\in F_1^{(10)}$, and otherwise, $\overline{F_1}^{(10)}$ obtained from ${F_1}^{(10)} + \{x,y,v,z\}$ by deleting $\{z,z'\}$ (respectively, $\{v,v'\},\{y,y'\}$) and adding $\{z_1,z_2\}$ (respectively, $\{v_2,w\}, \{y_1,y_2\}$) when $z' \in {F_1}^{(1)}$ (respectively, $v' \in {F_1}^{(10)}$, $y' \in {F_1}^{(10)}$). 
Therefore, $F^{(10)}$ is an induced forest in $G$, showing $a(G) \geq |F_1^{(10)}| + | F_2^{(10)}| + | F_3^{(10)}| + | F_4^{(10)}| + 4 - 9 \geq \lceil (4n+3)/7 \rceil$, a contradiction.

If  $(4(|G_1''|-8)+3,4(|G_2''|-2)+3,4(|G_3''|-2)+3,4(|G_4''|-2)+3) \equiv (0,4,0,0) \mod 7$, 
let $F_1^{(11)} = A(G_1'' - \{x,z,w,v,v_1,v_2,z_1,c_1,d_1\} + yz_2)$, $F_2^{(11)} = A(G_2'' - \{z_1\})$, $F_3^{(11)} = A(G_3'' - \{w,v_2,c_1\})$, and $F_4^{(11)} = A(G_4'' - \{v_1,v_2,d_1\})$. 
Then $|F_1^{(11)}| \geq \lceil (4(|G_1''|-9)+3)/7 \rceil$, $|F_2^{(11)}| \geq \lceil (4(|G_2''|-1)+3)/7 \rceil$, $|F_3^{(11)}| \geq \lceil (4(|G_3''|-3)+3)/7 \rceil$, and $|F_4^{(11)}| \geq \lceil (4(|G_4''|-3)+3)/7 \rceil$. 
Now $G[F_1^{(11)} \cup F_2^{(11)} \cup F_3^{(11)} \cup F_4^{(11)} + \{x,z,v\}] - (\{z_2,a_1\} \cap (F_1^{(11)} \triangle F_2^{(11)}))$. Therefore, $F^{(16)}$ is an induced forest in $G$, showing $a(G) \geq |F_1^{(11)}| + | F_2^{(11)}| + | F_3^{(11)}| + | F_4^{(11)}| + 3 - 2 \geq \lceil (4n+3)/7 \rceil$, a contradiction. 

\medskip
{\it Subcase 6.2: } $|N(v_1) \cap N(v_2)| \leq 2$.

\medskip
{\it Subcase 6.2.1: } $(n_1,n_2,n_3) \equiv (0,0,0), (1,0,0) \mod 7$. 

Let $F_1^{(1)} = A((G_1 - \{x,y,z,v\})/\{z_1z_2,y_1y_2,wv_2\} + v'z')$ with $z'$ (respectively, $y',v'$) as the identification of $\{z_1,z_2\}$ (respectively, $\{y_1,y_2\},\{w,v_2\}$), $F_2^{(1)} = A(G_2)$, and $F_3^{(1)} = A(G_3)$. 
Then $|F_1^{(1)}| \geq \lceil (4(|G_1|-7)+3)/7 \rceil$, $|F_2^{(1)}| \geq \lceil (4|G_2|+3)/7 \rceil$, and $|F_3^{(1)}| \geq \lceil (4|G_3|+3)/7 \rceil$. 
Let $F^{(1)} := G[\overline{F_1}^{(1)} \cup F_2^{(1)} \cup F_3^{(1)} - (\{z_1,z_2,a_1\} \cap (\overline{F_1}^{(1)} \triangle F_2^{(1)})) - (\{w,v_2,c_1\} \cap (\overline{F_1}^{(1)} \triangle F_3^{(1)}))]$,
where $\overline{F_1}^{(1)} := {F_1}^{(1)} + \{x,y,v,z\}$ if $v',y',z' \not\in F_1^{(1)}$, and otherwise, $\overline{F_1}^{(1)}$  obtained from ${F_1}^{(1)} + \{x,y,v,z\}$ by deleting $\{z,z'\}$ (respectively, $\{v,v'\},\{y,y'\}$) and adding $\{z_1,z_2\}$ (respectively, $\{v_2,w\}, \{y_1,y_2\}$) when $z' \in {F_1}^{(1)}$ (respectively, $v' \in {F_1}^{(1)}$, $y' \in {F_1}^{(1)}$). 
Therefore, $F^{(1)}$ is an induced forest in $G$, showing $a(G) \geq |F_1^{(1)}| + | F_2^{(1)}| + | F_3^{(1)}| + 4 - 6 \geq \lceil (4n+3)/7 \rceil$, a contradiction.

\medskip
{\it Subcase 6.2.2: } $(n_1,n_2,n_3) \equiv (3,0,4),(4,0,4),(3,4,0),(4,4,0),(0,4,4),(6,4,4), (4,4,6),$ $(4,6,4) \mod 7$. 

Let $F_1^{(2)} = A(G_1 - \{x,z,z_1,w,v\}/v_1v_2 + yz_2)$ with $v'$ as the identification of $\{v_1,v_2\}$, $F_2^{(2)} = A(G_2 - \{z_1\})$, and $F_3^{(2)} = A(G_3 - \{w\})$. 
Then $|F_1^{(2)}| \geq \lceil (4(|G_1|-6)+3)/7 \rceil$, $|F_2^{(2)}| \geq \lceil (4(|G_2|-1)+3)/7 \rceil$, and $|F_3^{(2)}| \geq \lceil (4(|G_3|-1)+3)/7 \rceil$.
Now $G[F_1^{(2)} \cup F_2^{(2)} \cup F_3^{(2)} + \{x,v,z\} - (\{z_2,a_1\} \cap (F_1^{(2)} \triangle F_2^{(2)}))  - (\{c_1,v_2\} \cap (F_1^{(2)} \triangle F_2^{(2)}))]$ (if $v' \not\in F_1^{(2)}$) or $G[(F_1^{(2)} - v') \cup F_2^{(2)} \cup F_3^{(2)} + \{x,v,z\} - (\{z_2,a_1\} \cap (F_1^{(2)} \triangle F_2^{(2)}))  - (\{c_1,v_2\} \cap ((F_1^{(2)} + v_2) \triangle F_2^{(2)}))]$ (if $v' \in F_1^{(2)}$)  
is an induced forest in $G$, showing $a(G) \geq |F_1^{(2)}| + | F_2^{(2)}| + | F_3^{(2)}| + 3 - 4 \geq \lceil (4n+3)/7 \rceil$, a contradiction.

\medskip
{\it Subcase 6.2.3: } $(n_1,n_2,n_3) \equiv (0,4,3) \mod 7$.

Let $F_1^{(3)} = A(G_1 - \{x,z,v,w,z_1,v_1,v_2\} + yz_2)$, $F_2^{(3)} = A(G_2 - z_1)$, and $F_3^{(3)} = A(G_3 - \{w,v_2\})$. 
Then $|F_1^{(3)}| \geq \lceil (4(|G_1|-7)+3)/7 \rceil$, $|F_2^{(3)}| \geq \lceil (4(|G_2|-1)+3)/7 \rceil$ and $|F_3^{(3)}| \geq \lceil (4(|G_3|-2)+3)/7 \rceil$. 
Then $G[F_1^{(3)} \cup F_2^{(3)} \cup F_3^{(3)} + \{x,v,z\} - (\{z_2,a_1\} \cap (F_1^{(3)} \triangle F_2^{(3)}) )  - (\{c_1\} \cap (F_1^{(3)} \triangle F_3^{(3)}) )]$ is an induced forest in $G$, showing $a(G) \geq |F_1^{(3)}| + | F_2^{(3)}| + | F_3^{(3)}| + 3 - 3 \geq \lceil (4n+3)/7 \rceil$, a contradiction.

\medskip
{\it Subcase 6.2.4: } $(n_1,n_2,n_3) \equiv (0,3,4) \mod 7$.

Let $F_1^{(4)} = A((G_1 - \{x,y,z,w,z_1,z_2\})/y_1y_2)$ with $y'$ as the identification of $\{y_1,y_2\}$, $F_2^{(4)} = A(G_2 - \{z_1,z_2\})$, and $F_3^{(4)} = A(G_3 - w)$. 
Then $|F_1^{(4)}| \geq \lceil (4(|G_1|-7)+3)/7 \rceil$, $|F_2^{(4)}| \geq \lceil (4(|G_2|-2)+3)/7 \rceil$ and $|F_3^{(4)}| \geq \lceil (4(|G_3|-1)+3)/7 \rceil$. 
Let $F = F_1^{(4)} \cup F_2^{(4)} \cup F_3^{(4)} + \{x,y,z\} - (\{a_1\} \cap (F_1^{(4)} \triangle F_2^{(4)})) - (\{c_1,v_2\} \cap (F_1^{(4)} \triangle F_3^{(4)}))$.
Now $G[F]$ (if $y' \not\in F_1^{(4)}$) or $G[F - \{y',y\} + \{y_1,y_2\}]$ (if $y' \in F_1^{(4)}$)
is an induced forest in $G$, showing $a(G) \geq |F_1^{(4)}| + | F_2^{(4)}| + | F_3^{(4)}| + 3 - 3 \geq \lceil (4n+3)/7 \rceil$, a contradiction.

\medskip
{\it Subcase 6.2.5: } $(n_1,n_2,n_3) \equiv (4,3,0) \mod 7$ (respectively $(1,6,0) \mod 7$).

Let $A_5 = \emptyset$ and $A_6 = \{a_1\}$. 
For $i=5,6$, let $F_1^{(i)} = A((G_1 - \{x,y,z,v,z_1,z_2\} - A_i)/\{y_1y_2,wv_2\})$ with $y'$ (respectively, $v'$) as the identification of $\{y_1,y_2\}$ (respectively, $\{w,v_2\}$), $F_2^{(i)} = A(G_2 - \{z_1,z_2\} - A_i)$, and $F_3^{(i)} = A(G_3)$.
Then $|F_1^{(i)}| \geq \lceil (4(|G_1|-8- |A_i|)+3)/7 \rceil$, $|F_2^{(i)}| \geq \lceil (4(|G_2|-2-|A_i|)+3)/7 \rceil$ and $|F_3^{(i)}| \geq \lceil (4|G_3|+3)/7 \rceil$. 
Let $F^{(i)} := G[\overline{F_1}^{(i)} \cup F_2^{(i)} \cup F_3^{(i)}  - (\{a_1\} \cap (\overline{F_1}^{(i)} \triangle F_2^{(i)})) - (\{w,v_2,c_1\} \cap (\overline{F_1}^{(i)} \triangle F_3^{(i)}))]$,
where $\overline{F_1}^{(i)} = {F_1}^{(i)} + \{x,y,v,z\}$ if $y',v' \not\in F_1^{(i)}$, and otherwise, $\overline{F_1}^{(i)}$ obtained from ${F_1}^{(i)}$ by deleting $\{y,y'\}$ (respectively, $\{v,v'\}$) and adding $\{y_1,y_2\}$ (respectively, $\{w,v_2\}$) when $y' \in {F_1}^{(i)}$ (respectively, $v' \in {F_1}^{(i)}$).
Therefore, $F^{(i)}$ is an induced forest in $G$, showing $a(G) \geq |F_1^{(i)}| + | F_2^{(i)}| + | F_3^{(i)}| + 4 - 3 - (1-|A_i|) \geq \lceil (4n+3)/7 \rceil$, a contradiction.

\medskip
{\it Subcase 6.2.6: } $(n_1,n_2,n_3) \equiv (4,0,3) \mod 7$ (respectively, $(1,0,6) \mod 7$).

Let $C_7 = \emptyset$ and $C_8 = \{c_1\}$. 
For $i=7,8$, let $F_1^{(i)} = A((G_1 - \{x,y,z,v,w,v_2\} - C_i)/\{z_1z_2,y_1y_2\})$ with $z'$ (respectively, $y'$) as the identification of  $\{z_1,z_2\}$ (respectively, $\{y_1,y_2\}$), $F_2^{(i)} = A(G_2)$, and $F_3^{(i)} = A(G_3 - \{w,v_2\} - C_i)$.
Then $|F_1^{(i)}| \geq \lceil (4(|G_1|-8- |C_i|)+3)/7 \rceil$, $|F_2^{(i)}| \geq \lceil (4|G_2|+3)/7 \rceil$ and $|F_3^{(i)}| \geq \lceil (4(|G_3|-2-|C_i|)+3)/7 \rceil$. 
Let $F^{(i)} := G[\overline{F_1}^{(i)} \cup F_2^{(i)} \cup F_3^{(i)}  - (\{c_1\} \cap (\overline{F_1}^{(i)} \triangle F_3^{(i)})) - (\{z_1,z_2,a_1\} \cap (\overline{F_1}^{(i)} \triangle F_2^{(i)}))]$,
where $\overline{F_1}^{(i)} = {F_1}^{(i)} + \{x,y,v,z\}$ if $y',z' \not\in F_1^{(i)}$, and otherwise $\overline{F_1}^{(i)}$ obtained from ${F_1}^{(i)}$ by deleting $\{y,y'\}$ (respectively, $\{z,z'\}$) and adding $\{y_1,y_2\}$ (respectively, $\{z_1,z_2\}$) when $y' \in {F_1}^{(i)}$ (respectively, $z' \in {F_1}^{(i)}$).
Therefore, $F^{(i)}$ is an induced forest in $G$, showing $a(G) \geq |F_1^{(i)}| + | F_2^{(i)}| + | F_3^{(i)}| + 4 - 3 - (1-|C_i|) \geq \lceil (4n+3)/7 \rceil$, a contradiction.

\medskip

Case 7: $|N(z_1) \cap N(z_2)| \leq 2$, $|N(y_1) \cap N(y_2)| > 2$ and $|N(w) \cap N(v_2)| > 2$. 

There exist $b_1 \in N(y_1) \cap N(y_2)$, $c_1 \in N(w) \cap N(v_2)$ and subgraphs $G_1,G_2,G_3$ of $G$ such that 
$G_2$ is the maximal subgraph of $G$ contained in the closed region of the plane bounded by the cycle $yy_1b_1y_2y$ and containing $N(y_1) \cap N(y_2) - \{y\}$,
$G_3$ is the maximal subgraph of $G$ contained in the closed region of the plane bounded by the cycle $vwc_1v_2v$ and containing $N(w) \cap N(v_2) - \{v\}$,
and $G_1$ is obtained from $G$ by removing $G_2-\{y_1,y_2,b_1\}$ and $G_3-\{w,c_1,v_2\}$.
Let $B_i = \{b_1\}$ if $i=1,2$ and $B_i = \emptyset$ if $i=3,4$.
Let $C_i = \{c_1\}$ if $i=1,3$ and $C_i = \emptyset$ if $i=2,4$.
For $i \in [4]$, let $F_1^{(i)} = A((G_1 - \{x,y,z,v,y_1,y_2,w,v_2\} - B_i - C_i)/z_1z_2)$ with $z'$ as the identification of $\{z_1,z_2\}$, 
$F_2^{(i)} = A(G_2 - \{y_1,y_2\} - B_i)$,
and $F_3^{(i)} = A(G_3 - \{w,v_2\} - C_i)$. 
Then $|F_1^{(i)}| \geq \lceil (4(|G_1|-9- |B_i| - |C_i|)+3)/7 \rceil$, $|F_2^{(i)}| \geq \lceil (4(|G_2|-2-|B_i|)+3)/7 \rceil$, and $|F_3^{(i)}| \geq \lceil (4(|G_3|-2-|C_i|)+3)/7 \rceil$. 
Let $F = F_1^{(i)} \cup F_2^{(i)} \cup F_3^{(i)} + \{x,y,v,z\} - (\{b_1\} \cap (F_1^{(i)} \triangle F_2^{(i)})) - (\{c_1\} \cap (F_1^{(i)} \triangle F_3^{(i)}))$.
Now $G[F]$ (if $z' \not\in F_1^{(i)}$) or $G[F - \{z,z'\}+\{z_1,z_2\}]$ (if $z' \in F_1^{(i)}$)
is an induced forest of size $|F_1^{(i)}| + |F_2^{(i)}| + |F_3^{(i)}| + 4 - (1-|B_i|) - (1-|C_i|)$. 
Let $(n_1,n_2,n_3) := (4(|G_1|-9)+3,4(|G_2|-2)+3,4(|G_3|-2)+3)$.
By Lemma \ref{ineq2}(5) (with $a = |G_1| - 9, a_1 = |G_2| - 2, a_2 = |G_3| - 2, c=4$),
 $(n_1,n_2,n_3) \equiv (0,0,0), (1,0,0), (4,0,3), (4,3,0), (3,0,4), (4,0,4), (3,4,0), (4,4,0), (1,6,0), (1,0,6), (0,3,4),$ 
$(0,4,3), (0,4,4), (6,4,4), (4,4,6), (4,6,4) \mod 7$.
Let $B_i = \emptyset$ if $i=1,4,8$ and $B_i = \{b_1\}$ if $i=3,5,9$.

\medskip
{\it Subcase 7.1: } $(n_1,n_2,n_3) \equiv (0,0,0),(1,0,0) \mod 7$ (respectively $(4,4,0) \mod 7$).

For $i=1$ (respectively, $i=5$), let $F_1^{(i)} = A((G_1 - \{x,y,z,v\} - B_i)/\{z_1z_2,y_1y_2,wv_2\} + v'z')$ with $z'$ (respectively, $y',v'$) as the identification of $\{z_1,z_2\}$ (respectively, $\{y_1,y_2\},\{w,v_2\}$) and $F_2^{(i)} = A(G_2 - B_i)$, and $F_3^{(i)} = A(G_3)$. 
Then $|F_1^{(i)}| \geq \lceil (4(|G_1|-7-|B_i|)+3)/7 \rceil$, $|F_2^{(i)}| \geq \lceil (4(|G_2|-|B_i|)+3)/7 \rceil$, and $|F_3^{(i)}| \geq \lceil (4|G_3|+3)/7 \rceil$. 
Let $F^{(i)} := G[\overline{F_1}^{(i)} \cup F_2^{(i)} \cup F_3^{(i)} - (\{y_1,y_2,b_1\} \cap (\overline{F_1}^{(i)} \triangle F_2^{(i)})) - (\{w,v_2,c_1\} \cap (\overline{F_1}^{(i)} \triangle F_3^{(i)}))]$,
where $\overline{F_1}^{(i)} := {F_1}^{(i)} + \{x,y,v,z\}$ if $v',y',z' \not\in F_1^{(i)}$, and otherwise, $\overline{F_1}^{(i)}$ obtained from ${F_1}^{(i)} + \{x,y,v,z\}$ by deleting $\{z,z'\}$ (respectively, $\{v,v'\},\{y,y'\}$) and adding $\{z_1,z_2\}$ (respectively, $\{v_2,w\}, \{y_1,y_2\}$) when $z' \in {F_1}^{(i)}$ (respectively, $v' \in {F_1}^{(i)}$, $y' \in {F_1}^{(i)}$). 
Therefore, $F^{(i)}$ is an induced forest in $G$, showing $a(G) \geq |F_1^{(i)}| + | F_2^{(i)}| + | F_3^{(i)}| + 4 - 5 - (1-|B_i|) \geq \lceil (4n+3)/7 \rceil$, a contradiction.

\medskip
{\it Subcase 7.2: } $(n_1,n_2,n_3) \equiv (0,4,4),(6,4,4) \mod 7$. 

Let $F_1^{(2)} = A(G_1 - \{x,y,w,z_1,y_2\} + y_1z)$, $F_2^{(2)} = A(G_2 - y_2)$, and $F_3^{(2)} = A(G_3 - w)$. 
Then $|F_1^{(2)}| \geq \lceil (4(|G_1|-5)+3)/7 \rceil$, $|F_2^{(2)}| \geq \lceil (4(|G_2|-1)+3)/7 \rceil$, and $|F_3^{(2)}| \geq \lceil (4(|G_3|-1)+3)/7 \rceil$.
Now $G[F_1^{(2)} \cup F_2^{(2)} \cup F_3^{(2)} + \{x,y\} - (\{y_1,b_1\} \cap (F_1^{(2)} \triangle F_2^{(2)})) - (\{v_2,c_1\} \cap (F_1^{(2)} \triangle F_3^{(2)}))]$  is an induced forest in $G$, showing $a(G) \geq |F_1^{(2)}| + | F_2^{(2)}| + | F_3^{(2)}| + 2 - 4 \geq \lceil (4n+3)/7 \rceil$, a contradiction.

\medskip
{\it Subcase 7.3: } $(n_1,n_2,n_3) \equiv (4,0,4),(4,6,4) \mod 7$ (respectively, $(0,3,4) \mod 7$).

For $i=3$ (respectively, $i=4$), let $F_1^{(i)} = A((G_1 - \{x,y,z,w,y_1,y_2\} - B_i)/z_1z_2)$ with $z'$ as the identification of $\{z_1,z_2\}$, $F_2^{(i)} = A(G_2 - \{y_1,y_2\} - B_i)$, and $F_3^{(i)} = A(G_3 - w)$. Then $|F_1^{(i)}| \geq \lceil (4(|G_1|-7-|B_i|)+3)/7 \rceil$, $|F_2^{(i)}| \geq \lceil (4(|G_2|-2-|B_i|)+3)/7 \rceil$ and $|F_3^{(i)}| \geq \lceil (4(|G_3|-1)+3)/7 \rceil$. 
Let $F = F_1^{(i)} \cup F_2^{(i)} \cup F_3^{(i)} + \{x,y,z\} - (\{v_2,c_1\} \cap (F_1^{(i)} \triangle F_3^{(i)})) - (\{b_1\} \cap (F_1^{(i)} \triangle F_2^{(i)}))$. 
Now $G[F]$ (if $z' \not\in F_1^{(i)}$) or $G[F - \{z,z'\} + \{z_1,z_2\}]$ (if $z' \in F_1^{(i)}$)
is an induced forest in $G$, showing $a(G) \geq |F_1^{(i)}| + | F_2^{(i)}| + | F_3^{(i)}| + 3 - 2 - (1-|B_i|) \geq \lceil (4n+3)/7 \rceil$, a contradiction.

\medskip
{\it Subcase 7.4: } $(n_1,n_2,n_3) \equiv (4,0,3) \mod 7$ (respectively, $(1,0,6) \mod 7$).

Let $C_6 = \emptyset$ and $C_7 = \{c_1\}$. 
For $i=6,7$, let $F_1^{(i)} = A((G_1 - \{x,y,z,w,v,v_2\} - C_i)/\{y_1y_2,z_1z_2\})$ with $y'$ (respectively, $z'$) as the identification of  $\{y_1,y_2\}$ (respectively, $\{z_1,z_2\}$), $F_2^{(i)} = A(G_2)$, and $F_3^{(i)} = A(G_3 - \{w,v_2\} - C_i)$. Then $|F_1^{(i)}| \geq \lceil (4(|G_1|-8-|C_i|)+3)/7 \rceil$, $|F_2^{(i)}| \geq \lceil (4|G_2|+3)/7 \rceil$, and $|F_3^{(i)}| \geq \lceil (4(|G_3|-2-|C_i|)+3)/7 \rceil$. 
Let $F^{(i)} := G[\overline{F_1}^{(i)} \cup F_2^{(i)} \cup F_3^{(i)} - (\{y_1,y_2,b_1\} \cap (\overline{F_1}^{(i)} \triangle F_2^{(i)})) - (\{c_1\} \cap (\overline{F_1}^{(i)} \triangle F_3^{(i)}))]$,
where $\overline{F_1}^{(i)} := {F_1}^{(i)} + \{x,y,v,z\}$ if $y',z' \not\in F_1^{(i)}$, and otherwise, $\overline{F_1}^{(i)}$ obtained from ${F_1}^{(i)} + \{x,y,v,z\}$ by deleting $\{z,z'\}$ (respectively, $\{y,y'\}$) and adding $\{z_1,z_2\}$ (respectively, $\{y_1,y_2\}$) when $z' \in {F_1}^{(i)}$ (respectively, $y' \in {F_1}^{(i)}$). 
Therefore, $F^{(i)}$ is an induced forest in $G$, showing $a(G) \geq |F_1^{(i)}| + | F_2^{(i)}| + | F_3^{(i)}| + 4 - 3 - (1-|C_i|) \geq \lceil (4n+3)/7 \rceil$, a contradiction.

\medskip
{\it Subcase 7.5: } $(n_1,n_2,n_3) \equiv (4,3,0) \mod 7$ (respectively, $(1,6,0) \mod 7$).

For $i=8$ (respectively, $i=9$), let $F_1^{(i)} = A((G_1 - \{x,y,z,v,y_1,y_2\} - B_i)/\{wv_2,z_1z_2\} + v'z')$ with $v'$ (respectively, $z'$) as the identification of $\{w,v_2\}$ (respectively, $\{z_1,z_2\}$), $F_2^{(i)} = A(G_2 - \{y_1,y_2\} - B_i)$, and $F_3^{(i)}=A(G_3)$. 
Then $|F_1^{(i)}| \geq \lceil (4(|G_1|-8-|B_i|)+3)/7 \rceil$,  $|F_2^{(i)}| \geq \lceil (4(|G_2| - 2 - |B_i|)+3)/7 \rceil$ and  $|F_3^{(i)}| \geq \lceil (4|G_3|+3)/7 \rceil$. 
Let $F^{(i)} := G[\overline{F_1}^{(i)} \cup F_2^{(i)} \cup F_3^{(i)} - (\{b_1\} \cap (\overline{F_1}^{(i)} \triangle F_2^{(i)})) - (\{w,v_2,c_1\} \cap (\overline{F_1}^{(i)} \triangle F_3^{(i)}))]$,
where $\overline{F_1}^{(i)} := {F_1}^{(i)} + \{x,y,v,z\}$ if $v',z' \not\in F_1^{(i)}$, and otherwise, $\overline{F_1}^{(i)}$ obtained from ${F_1}^{(i)} + \{x,y,v,z\}$ by deleting $\{z,z'\}$ (respectively, $\{v,v'\}$) and adding $\{z_1,z_2\}$ (respectively, $\{w,v_2\}$) when $z' \in {F_1}^{(i)}$ (respectively, $v' \in {F_1}^{(i)}$). 
Therefore, $F^{(i)}$ is an induced forest in $G$, showing $a(G) \geq |F_1^{(i)}| + | F_2^{(i)}| + | F_3^{(i)}| + 4 - 3 - (1-|B_i|) \geq \lceil (4n+3)/7 \rceil$, a contradiction.

\medskip
{\it Subcase 7.6: } $(n_1,n_2,n_3) \equiv (3,0,4) \mod 7$. 

Let $F_1^{(10)} = A((G_1 - \{x,y,z,w\})/\{y_1y_2,z_1z_2\})$ with $y'$ (respectively, $z'$) as the identification of  $\{y_1,y_2\}$ (respectively, $\{z_1,z_2\}$), $F_2^{(10)} = A(G_2)$, and $F_3^{(10)} = A(G_3 - w)$. 
Then $|F_1^{(10)}| \geq \lceil (4(|G_1|-6)+3)/7 \rceil$, $|F_2^{(10)}| \geq \lceil (4|G_2|+3)/7 \rceil$ and $|F_3^{(10)}| \geq \lceil (4(|G_3|-1)+3)/7 \rceil$. 
Let $F^{(10)} := G[\overline{F_1}^{(10)} \cup F_2^{(10)} \cup F_3^{(10)} - (\{y_1,y_2,b_1\} \cap (\overline{F_1}^{(10)} \triangle F_2^{(10)})) - (\{v_2,c_1\} \cap (\overline{F_1}^{(10)} \triangle F_3^{(10)}))]$,
where $\overline{F_1}^{(10)} := {F_1}^{(10)} + \{x,y,z\}$ if $y',z' \not\in F_1^{(10)}$, and otherwise, $\overline{F_1}^{(10)}$ obtained from ${F_1}^{(10)} + \{x,y,z\}$ by deleting $\{z,z'\}$ (respectively, $\{y,y'\}$) and adding $\{z_1,z_2\}$ (respectively, $\{y_1,y_2\}$) when $z' \in {F_1}^{(10)}$ (respectively, $y' \in {F_1}^{(10)}$). 
Therefore, $F^{(10)}$ is an induced forest in $G$, showing $a(G) \geq |F_1^{(10)}| + | F_2^{(10)}| + | F_3^{(10)}| + 3 - 5 \geq \lceil (4n+3)/7 \rceil$, a contradiction.

\medskip
{\it Subcase 7.7: } $(n_1,n_2,n_3) \equiv (4,4,6) \mod 7$. 

By Claim 5, $wy_1 \not\in E(G)$. Let $w' \in N(w) - \{v,c_1,x,y_2\}$. 
Let $F_1^{(11)} = A(G_1 - \{x,z,w,v,y_2,v_1,v_2,c_1\} + w'y)$,
 $F_2^{(11)} = A(G_2 - y_2)$, 
 and $F_3^{(11)} = A(G_3 - \{w,v_2,c_1\})$. 
 Then $|F_1^{(11)}| \geq \lceil (4(|G_1|-8)+3)/7 \rceil$, $|F_2^{(11)}| \geq \lceil (4(|G_2|-1)+3)/7 \rceil$, and $|F_3^{(11)}| \geq \lceil (4(|G_3|-3)+3)/7 \rceil$. 
 Now $G[F_1^{(11)} \cup F_2^{(11)} \cup F_3^{(11)} + \{x,w,v\} - (\{y_1,b_1\} \cap (F_1^{(11)} \triangle F_2^{(11)}))]$
 is an induced forest in $G$, showing $a(G) \geq |F_1^{(11)}| + | F_2^{(11)}| + | F_3^{(11)}| + 3 - 2 \geq \lceil (4n+3)/7 \rceil$, a contradiction.

\medskip
{\it Subcase 7.8: } $(n_1,n_2,n_3) \equiv (3,4,0) \mod 7$. 


We claim that $|N(y_1) \cap N(z_1)| \leq 2$. Otherwise, 
there exist $d_1 \in N(y_1) \cap N(z_1)$ and subgraphs $G_1',G_2',G_3',G_4'$ of $G$ such that 
$G_2' = G_2, G_3' = G_3,$
$G_4'$ is the maximal subgraph of $G$ contained in the closed region of the plane bounded by the cycle $zy_1d_1z_1z$ and containing $N(y_1) \cap N(z_1) - \{z\}$,
and $G_1'$ is obtained from $G$ by removing $G_2'-\{y_1,y_2,b_1\}$, $G_3'-\{w,c_1,v_2\}$ and $G_4'-\{y_1,d_1,z_1\}$.
Let $D_{14} = \{d_1\}$ and $D_{15} = \emptyset$. 
For $i=14,15$, let $F_1^{(i)} = A(G_1' - \{x,y,z,w,v,y_1,y_2,v_2,c_1,z_1\} - D_i)$, $F_2^{(i)} = A(G_2' - \{y_1,y_2\})$, $F_3^{(i)} = A(G_3'-\{w,v_2,c_1\})$, and $F_4^{(i)} = A(G_4'-\{y_1,z_1\} - D_i)$. 
Then $|F_1^{(i)}| \geq \lceil (4(|G_1'|-10-|D_i|)+3)/7 \rceil$, 
$|F_2^{(i)}| \geq \lceil (4(|G_2'| - 2)+3)/7 \rceil = \lceil (4(|G_2| - 2)+3)/7 \rceil = (4(|G_2| - 2)+3)/7 + 3/7$, 
$|F_3^{(i)}| \geq \lceil (4(|G_3'|-3)+3)/7 \rceil = \lceil (4(|G_3|-3)+3)/7 \rceil = (4(|G_3'|-3)+3)/7 + 4/7$, 
and $|F_4^{(i)}| \geq \lceil (4(|G_4'|-2-|D_i|)+3)/7 \rceil$. 
Note $N(w) - \{y_1,x,v\} \subseteq V(G_3')$.
Now $G[F_1^{(i)} \cup F_2^{(i)} \cup F_3^{(i)} \cup F_4^{(i)} + \{x,y,z,w\} - (\{b_1\} \cap (F_1^{(i)} \triangle F_2^{(i)})) - (\{d_1\} \cap (F_1^{(i)} \triangle F_4^{(i)}))]$ is an induced forest of size $|F_1^{(i)}| + | F_2^{(i)}| + | F_3^{(i)}| + | F_4^{(i)}| + 4 - 1 - (1-|D_i|)$. 
By Lemma \ref{ineq2}(1) (with $k=1, a = |G_1'|-10, a_1 = |G_4'|-2, L = \{1,2\}, b_1 = |G_2'| -2, b_2 = |G_3'|-3, c = 3$),
 $a(G) \geq \lceil (4n+3)/7 \rceil$, a contradiction.

Let $F_1^{(12)} = A((G_1 - \{x,y,z,y_2,w,v,v_2,c_1\})/z_1y_1)$ with $y'$ as the identification of $\{z_1,y_1\}$, $F_2^{(12)} = A(G_2 - y_2)$ and $F_3^{(12)} = A(G_3 - \{w,v_2,c_1\})$. 
Then $|F_1^{(12)}| \geq \lceil (4(|G_1|-9)+3)/7 \rceil$, $|F_2^{(12)}| \geq \lceil (4(|G_2|-1)+3)/7 \rceil$ and $|F_3^{(12)}| \geq \lceil (4(|G_3|-3)+3)/7 \rceil$.
Note $N(w) - \{y_2,x,v\} \subseteq V(G_3)$.
Now $G[F_1^{(12)} \cup F_2^{(12)} \cup F_3^{(12)} + \{x,y,w,v\} - (\{y_1,b_1\} \cap (F_1^{(12)} \triangle F_2^{(12)})) ]$ (if $y' \not\in F_1^{(12)}$) or $G[(F_1^{(12)} - y') \cup F_2^{(12)} \cup F_3^{(12)} + \{x,y_1,z_1,w,v\} - (\{y_1,b_1\} \cap ( (F_1^{(12)}+y_1) \triangle F_2^{(12)})) ]$ (if $y' \in F_1^{(12)}$) is an induced forest in $G$, showing $a(G) \geq |F_1^{(12)}| + | F_2^{(12)}| + | F_3^{(12)}| + 4 - 2 \geq \lceil (4n+3)/7 \rceil$, a contradiction.

\medskip
{\it Subcase 7.9: } $(n_1,n_2,n_3) \equiv (0,4,3) \mod 7$.

We claim that $|N(v_1) \cap N(v_2)| \leq 2$. Otherwise, 
there exist $e_1 \in N(v_1) \cap N(v_2)$ and subgraphs $G_1',G_3',G_4'$ such that
$G_3' = G_3$,
$G_4'$ is the maximal subgraph of $G$ contained in the closed region of the plane bounded by the cycle $vv_1e_1v_2v$ and containing $N(v_1) \cap N(v_2) - \{v\}$,
and $G_1'$ is obtained from $G$ by removing $G_3'-\{w,c_1,v_2\}$ and $G_4'-\{v_1,e_1,v_2\}$.
 Let $E_{16} = \{e_1\}$ and $E_{17} = \emptyset$.
 For $i=16,17$, let $F_1^{(i)} = A((G_1' - \{x,z,w,v,v_1,v_2\} - E_i)/z_1z_2)$ with $z'$ as the identification of $\{z_1,z_2\}$, $F_3^{(i)} = A(G_3'-\{w,v_2\})$, and $F_4^{(i)} = A(G_4'-\{v_1,v_2\} - E_i)$. 
 Then $|F_1^{(i)}| \geq \lceil (4(|G_1'|-7-|E_i|)+3)/7 \rceil$, 
 $|F_3^{(i)}| \geq \lceil (4(|G_3'|-2)+3)/7 \rceil = \lceil (4(|G_3|-2)+3)/7 \rceil= (4(|G_3|-2)+3)/7 + 4/7$, 
 and $|F_4^{(i)}| \geq \lceil (4(|G_4'|-2-|E_i|)+3)/7 \rceil$. 
 Let $F = F_1^{(i)} \cup F_3^{(i)} \cup F_4^{(i)} + \{x,z,v\} - (\{c_1\} \cap (F_1^{(i)} \triangle F_3^{(i)})) - (\{e_1\} \cap (F_1^{(i)} \triangle F_4^{(i)}))$.
Now $G[F]$ (if $z' \not\in F_1^{(i)}$) or $G[F -\{z,z'\}+\{z_1,z_2\}]$ (if $z' \in F_1^{(i)}$)
   is an induced forest in $G$, showing $a(G) \geq |F_1^{(i)}| + | F_3^{(i)}| + | F_4^{(i)}| + 3 - 1 - (1-|E_i|)$. 
   By Lemma \ref{ineq2}(1) (with $k=1, a = |G_1'| - 7, a_1 = |G_4'|-2, L = \{1\}, b_1 = |G_3'|-2, c=2$),
    $a(G) \geq \lceil (4n+3)/7 \rceil$, a contradiction.

Let $F_1^{(13)} = A((G_1 - \{x,w,v,y_2,c_1\})/\{yz,v_1v_2\})$ with $x'$ (respectively, $v'$) as the identification of $\{y,z\}$ (respectively, $\{v_1,v_2\}$), $F_2^{(13)} = A(G_2 - y_2)$, and $F_3^{(13)} = A(G_3 - \{w,c_1\})$. 
Then $|F_1^{(13)}| \geq \lceil (4(|G_1|-7)+3)/7 \rceil$, $|F_2^{(13)}| \geq \lceil (4(|G_2|-1)+3)/7 \rceil$ and $|F_3^{(13)}| \geq \lceil (4(|G_3|-2)+3)/7 \rceil$. 
Let $F^{(13)} := G[\overline{F_1}^{(13)} \cup F_2^{(13)} \cup F_3^{(13)}  - (\{y_1,b_1\} \cap (\overline{F_1}^{(13)} \triangle F_2^{(13)}) ) - (\{v_2\} \cap (\overline{F_1}^{(13)} \triangle F_3^{(13)}) )] $,
where $\overline{F_1}^{(13)} := {F_1}^{(13)} + \{x,w,v\}$ if $x',v' \not\in F_1^{(13)}$, and $\overline{F_1}^{(13)}$ obtained from ${F_1}^{(13)} + \{x,w,v\}$ by deleting $\{x,x'\}$ (respectively, $\{v,v'\}$) and adding $\{y,z\}$ (respectively, $\{v_1,v_2\}$) when $x' \in {F_1}^{(13)}$ (respectively, $v' \in {F_1}^{(13)}$).
Note $N(w) - \{y_2,x,v\} \subseteq V(G_3)$.
Therefore, $F^{(13)}$ is an induced forest in $G$, showing $a(G) \geq |F_1^{(13)}| + | F_2^{(13)}| + | F_3^{(13)}| + 3 -3 \geq \lceil (4n+3)/7 \rceil$, a contradiction.
\qed

\section{Another forbidden configuration at a 3-vertex}

In this section we prove that for any $x\in V_3$, $N(x)\cap V_i=\emptyset$ for some $i\in \{3,4,5\}$.  

\begin{lem}
\label{No3345}
Let $x\in V_3$. Then $N(x)\cap V_3=\emptyset$, or $N(x)\cap V_4=\emptyset$, or $N(x)\cap V_5=\emptyset$. 
\end{lem}


\pf
We begin the proof by assuming that $N(x)=\{w,y,z\}$ with $y\in V_3$, $z \in V_4$ and $w \in V_5$. 
Let  $N(y) = \{y_1,x,z_1\}$, $N(z) = \{x,z_1,z_2,z_3\}$, and  $N(w) = \{x,y_1,w_1,w_2,z_3\}$ where $w_1$ is co-facial with $y_1$.

\medskip
\textit{Claim 1: $N(z_1) \cap N(z_3) = \{z\}$. }

For, suppose $|N(z_1) \cap N(z_3)| \geq 2$. 
First, we claim that $N(z_1) \cap N(z_3) \cap N(w_2) = \emptyset$.
Otherwise, there exist $a_1 \in N(z_1) \cap N(z_3)$ and subgraphs $G_1,G_2,G_4$ of $G$ such that 
$G_2$ is the maximal subgraph of $G$ contained in the closed region of the plane bounded by the cycle $zz_1a_1z_3z$ and containing $N(z_1)\cap N(z_3) - \{z\}$, $G_4$ is the maximal subgraph of $G$ contained in the closed region of the plane bounded by the cycle $ww_2a_1z_3w$, and $G_1$ is obtained from $G$ by removing $G_2-\{z,z_1,a_1,z_3\}$ and $G_4-\{w,w_2,a_1,z_3\}$. 
Let $F_1^{(1)} = A(G_1 - \{w,x,y,z,y_1,z_1,z_3,a_1,w_2\})$, $F_2^{(1)} = A(G_2 - \{z_1,z,z_3,a_1\})$, and $F_4^{(1)} = A(G_4 - \{w,z_3,a_1,w_2\})$. Then $|F_1^{(1)}| \geq \lceil (4(|G_1|-9)+3)/7 \rceil$, $|F_2^{(1)}| \geq \lceil (4(|G_2|-4)+3)/7 \rceil$, and $|F_4^{(1)}| \geq \lceil (4(|G_4|-4)+3)/7 \rceil$. Now $G[F_1^{(1)} \cup F_2^{(1)} \cup F_4^{(1)} + \{w,x,y,z\}]$ is an induced forest in $G$, showing that $a(G) \geq |F_1^{(1)}| + |F_2^{(1)}| + |F_4^{(1)}| + 4$. 
By Lemma~\ref{ineq2}(7) 
(with $k=3, a_1 = |G_1| - 9, a_2 = |G_2| - 4, a_3 = |G_4| - 4, c=4$), 
$a(G) \geq \lceil (4n+3)/7 \rceil$ unless $(4(|G_1|-9)+3,4(|G_2|-4)+3),4(|G_4|-4)+3) \equiv (0,0,0),(0,6,0),(0,0,6),(6,0,0) \mod 7$. %
In first three cases, let $F_1^{(2)} = A(G_1 - \{w,x,y,z,z_1,z_3,a_1\})$, $F_2^{(2)} = A(G_2 - \{z_1,z,z_3,a_1\})$, and $F_4^{(2)} = A(G_4 - \{w,z_3,a_1\})$. Then $|F_1^{(2)}| \geq \lceil (4(|G_1|-7)+3)/7 \rceil$, $|F_2^{(2)}| \geq \lceil (4(|G_2|-4)+3)/7 \rceil$, and $|F_4^{(2)}| \geq \lceil (4(|G_4|-3)+3)/7 \rceil$. Now $G[F_1^{(2)} \cup F_2^{(2)} \cup F_4^{(2)} + \{x,y,z\} - \{w_2\}\cap(F_1^{(2)} \triangle F_4^{(2)}) ]$ is an induced forest in $G$, showing $a(G) \geq |F_1^{(2)}| + |F_2^{(2)}| + |F_4^{(2)}| + 3-1 \geq \lceil (4n+3)/7 \rceil$, a contradiction. Now, assume $(4(|G_1|-9)+3,4(|G_2|-4)+3),4(|G_4|-4)+3) \equiv (6,0,0) \mod 7$. 
If $y_1a_1 \not\in E(G)$, let $F_1^{(3)} = A(G_1 - \{w,x,y,z_1,z_3\} + y_1a_1)$, $F_2^{(3)} = A(G_2 - \{z_1,z_3\})$, and $F_4^{(3)} = A(G_4 - \{w,z_3\})$. Then $|F_1^{(3)}| \geq \lceil (4(|G_1|-5)+3)/7 \rceil$, $|F_2^{(3)}| \geq \lceil (4(|G_2|-2)+3)/7 \rceil$ and $|F_4^{(3)}| \geq \lceil (4(|G_4|-2)+3)/7 \rceil$. Now $G[F_1^{(3)} \cup F_2^{(3)} \cup F_4^{(3)} + \{x,y\} - (\{w_2,a_1\}\cap (F_1^{(3)} \triangle F_4^{(3)})) - (\{z,a_1\}\cap (F_1^{(3)} \triangle F_2^{(3)}))]$ an induced forest in $G$, showing $a(G) \geq |F_1^{(3)}| + |F_2^{(3)}| + |F_4^{(3)}| + 2 - 4 \geq \lceil (4n+3)/7 \rceil$, a contradiction. 
So $y_1a_1 \in E(G)$. Then there exist subgraphs $G_1',G_2',G_4',G_5'$ of $G$ such that $G_2'=G_2, G_4'=G_4$, 
$G_5'$ is the maximal subgraph of $G$ contained in the closed region of the plane bounded by the cycle $yy_1a_1z_1y$, and $G_1$ is obtained from $G$ by removing $G_2'-\{z,z_1,a_1,z_3\}$, $G_4'-\{w,w_2,a_1,z_3\}$, and $G_5'-\{y,y_1,a_1,z_1\}$. 
Let $F_1^{(4)} = A(G_1' - \{w,x,y,z,y_1,z_1,z_3,w_2,a_1\})$, 
$F_2^{(4)} = A(G_2'-\{z_1,z,z_3,a_1\})$, 
$F_4^{(4)} = A(G_4' - \{w,z_3,a_1,w_2\})$ 
and $F_5^{(4)} = A(G_5' - \{y_1,y,z_1,a_1\})$. Then $|F_1^{(4)}| \geq \lceil (4(|G_1'|-9)+3)/7 \rceil$, 
$|F_2^{(4)}| \geq \lceil (4(|G_2'|-2)+3)/7 \rceil$,
$|F_4^{(4)}| \geq \lceil (4(|G_4'|-2)+3)/7 \rceil$,
 and $|F_5^{(4)}| \geq \lceil (4(|G_5'|-4)+3)/7 \rceil$. Now $G[F_1^{(4)} \cup F_2^{(4)} \cup F_4^{(4)} \cup F_5^{(4)} + \{w,x,y,z\}]$ is an induced forest in $G$, showing $a(G) \geq |F_1^{(4)}| + |F_2^{(4)}| + |F_4^{(4)}| + |F_5^{(4)}|  + 4 \geq \lceil (4n+3)/7 \rceil$, a contradiction.

Secondly, we claim that $N(z_1) \cap N(z_3) \cap N(y_1) = \emptyset$. 
For otherwise, there exist $a_1 \in N(z_1) \cap N(z_3)$ and subgraphs $G_1,G_2,G_5$ of $G$ such that 
$G_2$ is the maximal subgraph of $G$ contained in the closed region of the plane bounded by the cycle $zz_1a_1z_3z$ and containing $N(z_1)\cap N(z_3) - \{z\}$, $G_5$ is the maximal subgraph of $G$ contained in the closed region of the plane bounded by the cycle $yy_1a_1z_1y$, and $G_1$ is obtained from $G$ by removing $G_2-\{z,z_1,a_1,z_3\}$ and $G_5-\{y,y_1,a_1,z_1\}$. 
Let $F_1^{(1)} = A(G_1 - \{x,y,z,y_1,z_1,z_3,a_1\})$, $F_2^{(1)} = A(G_2 - \{z_1,z,z_3,a_1\})$ and $F_5^{(1)} = A(G_5 - \{y_1,y,z_1,a_1\})$. Then $|F_1^{(1)}| \geq \lceil (4(|G_1|-7)+3)/7 \rceil$, $|F_2^{(1)}| \geq \lceil (4(|G_2|-4)+3)/7 \rceil$ and $|F_5^{(1)}| \geq \lceil (4(|G_5|-4)+3)/7 \rceil$. Now $G[F_1^{(1)} \cup F_2^{(1)} \cup F_5^{(1)} + \{x,y,z\}]$ is an induced forest in $G$, showing that $a(G) \geq |F_1^{(1)}| + |F_2^{(1)}| + |F_5^{(1)}| + 3$. 
By Lemma~\ref{ineq2}(6) (with $k=3, a_1 = |G_1| - 7, a_2 = |G_2| - 4, a_3 = |G_5| - 4, c=3$), 
$(4(|G_1|-7)+3,4(|G_2|-4)+3),4(|G_5|-4)+3) \equiv (0,0,0) \mod 7$. 
Let $F_1^{(2)} = A(G_1 - \{w,x,y,z,z_1\})$, $F_2^{(2)} = A(G_2 - \{z_1,z\})$, and $F_5^{(2)} = A(G_5 - \{y,z_1\})$. Then $|F_1^{(2)}| \geq \lceil (4(|G_1|-5)+3)/7 \rceil$, $|F_2^{(2)}| \geq \lceil (4(|G_2|-2)+3)/7 \rceil$, and $|F_5^{(2)}| \geq \lceil (4(|G_5|-2)+3)/7 \rceil$. Now $G[F_1^{(2)} \cup F_2^{(2)} \cup F_5^{(2)} + \{x,y\} - (\{z_3,a_1\} \cap (F_1^{(2)} \triangle F_2^{(2)})) - (\{y_1,a_1\} \cap (F_1^{(2)} \triangle F_5^{(2)})) ]$ is an induced forest in $G$, showing that $a(G) \geq |F_1^{(2)}| + |F_2^{(2)}| + |F_5^{(2)}| + 2 - 4 \geq \lceil (4n+3)/7 \rceil$, a contradiction.


Thirdly, we claim that $|N(y_1) \cap N(w_1)| \leq 2$. For otherwise, there exist $b_1 \in N(y_1) \cap N(w_1)$, 
$a_1 \in N(z_1) \cap N(z_3)$
 and subgraphs $G_1,G_2,G_3$ of $G$ such that 
$G_2$ is the maximal subgraph of $G$ contained in the closed region of the plane bounded by the cycle $zz_1a_1z_3z$ and containing $N(z_1)\cap N(z_3) - \{z\}$, $G_3$ is the maximal subgraph of $G$ contained in the closed region of the plane bounded by the cycle $wy_1b_1w_1w$ and containing $N(y_1)\cap N(w_1) - \{w\}$, and $G_1$ is obtained from $G$ by removing $G_2-\{z,z_1,a_1,z_3\}$ and $G_3-\{y_1,b_1,w_1\}$. 
Let $B_1 = B_3 = \overline{B_2} = \overline{B_4} = \{b_1\}$ and $B_2 = \overline{B_1} = B_4 = \overline{B_3} = \emptyset$. For $i=1,2$, let $F_1^{(i)} = A(G_1 - \{w,x,y,z,y_1,z_1,z_3,a_1,w_1\} - B_i)$, $F_2^{(i)} = A(G_2 - \{z_1,z,z_3,a_1\})$, and $F_3^{(i)} = A(G_3 - \{y_1,w_1\} - B_i)$. Then $|F_1^{(i)}| \geq \lceil (4(|G_1|-9 - |B_i|)+3)/7 \rceil$, $|F_2^{(i)}| \geq \lceil (4(|G_2|-4)+3)/7 \rceil$, and $|F_3^{(i)}| \geq \lceil (4(|G_3|-2-|B_i|)+3)/7 \rceil$. Now $G[F_1^{(i)} \cup F_2^{(i)} \cup F_3^{(i)} \cup \{w,x,y,z\} - \{b_1\} \cap (F_1^{(i)} \triangle F_3^{(i)})]$ is an induced forest in $G$, showing that $a(G) \geq |F_1^{(i)}| + |F_2^{(i)}| + |F_3^{(i)}| + 4 - |\overline{B_i}|$. 
For $j=3,4$, let $F_1^{(j)} = A(G_1 - \{w,x,y,z,y_1,z_1,z_3,w_1\} - B_j + w_2a_1)$, $F_2^{(j)} = A(G_2 - \{z_1,z,z_3\})$, and $F_3^{(j)} = A(G_3 - \{y_1,w_1\} - B_1)$. Then $|F_1^{(j)}| \geq \lceil (4(|G_1| -8 - |B_1|)+3)/7 \rceil$, $|F_2^{(j)}| \geq \lceil (4(|G_2|-3)+3)/7 \rceil$, and $|F_3^{(j)}| \geq \lceil (4(|G_3|-2-|B_j|)+3)/7 \rceil$. Now $G[F_1^{(j)} \cup F_2^{(j)} \cup F_3^{(j)} + \{w,x,y,z\} - (\{b_1\} \cap (F_1^{(j)} \triangle F_3^{(j)})) - (\{a_1\} \cap (F_1^{(j)} \triangle F_2^{(j)}))]$ is an induced forest in $G$, showing $a(G) \geq |F_1^{(j)}| + |F_2^{(j)}| + |F_3^{(j)}| + 4 - 1 - |\overline{B_1}|$. By Lemma~\ref{ineq2}(1) (with $k=1$), $a(G) \geq \lceil (4n+3)/7 \rceil$, a contradiction.

Since $|N(z_1) \cap N(z_3)| \geq 2$, there exist $a_1 \in N(z_1) \cap N(z_3)$ and subgraphs $G_1,G_2$ of $G$ such that $G_2$ is the maximal subgraph of $G$ contained in the closed region of the plane bounded by the cycle $zz_1a_1z_3z$ containing $N(z_1)\cap N(z_3) - \{z\}$, and $G_1$ is obtained from $G$ by removing $G_2-\{z,z_1,a_1,z_3\}$.
Let $F_1^{(1)} = A(G_1 - \{w,x,y,z,z_1,z_3,a_1\}/w_1y_1)$ with $w'$ as the identification of $w_1$ and $y_1$, and $F_2^{(1)} = A(G_2 - \{z_1,z,z_3,a_1\})$. Then $|F_1^{(1)}| \geq \lceil (4(|G_1|-8)+3)/7 \rceil$, and $|F_2^{(1)}| \geq \lceil (4(|G_2|-4)+3)/7 \rceil$. Now $G[F_1^{(1)} \cup F_2^{(1)} + \{w,x,y,z\}]$ (if $w' \not\in F_1^{(1)}$) or $G[(F_1^{(1)} -w') \cup F_2^{(1)} + \{w_1,y_1,x,y,z\} ]$ (if $w' \in F_1^{(1)}$) is an induced forest in $G$, showing $a(G) \geq |F_1^{(1)}| + |F_2^{(1)}|+ 4$. 
By Lemma~\ref{ineq2}(6) (with $k=2, a_1 = |G_1|-8, a_2 =|G_2|-4, c=4 $), 
$(4(|G_1|-8)+3,4(|G_2|-4)+3)) \equiv (0,0) \mod 7$. 
Let $F_1^{(2)} = A(G_1 - \{w,x,y,z,z_1,z_3\} + y_1a_1)$, and $F_2^{(2)} = A(G_2 - \{z_1,z,z_3\})$. Then $|F_1^{(2)}| \geq \lceil (4(|G_1|-6)+3)/7 \rceil$, and $|F_2^{(2)}| \geq \lceil (4(|G_2|-3)+3)/7 \rceil$. Now $G[F_1^{(2)} \cup F_2^{(2)} + \{x,y,z\} - (\{a_1\} \cap (F_1^{(2)} \triangle F_2^{(2)})) ]$ is an induced forest in $G$, showing $a(G) \geq |F_1^{(2)}| + |F_2^{(2)}|+ 3 - 1 \geq \lceil (4n+3)/7 \rceil$. This completes the proof of Claim 1. 

\medskip
\textit{Claim 2: $wz_2 \not\in E(G)$. }

Otherwise, $wz_2 \in E(G)$, there exists a separation $(G_1,G_2)$ such that $V(G_1 \cap G_2) = \{w,x,z,z_2\}$, $y \in V(G_1)$, and $z_3 \in V(G_2)$. Let $F_1^{(1)} = A(G_1 - \{w,x,z,z_2,y,z_1\})$, and $F_2^{(1)} = A(G_2 - \{w,x,z,z_2\})$. Then $|F_1^{(1)}| \geq \lceil (4(|G_1|-6)+3)/7 \rceil$, and $|F_2^{(1)}| \geq \lceil (4(|G_2|-4)+3)/7 \rceil$. Now $G[F_1^{(1)} \cup F_2^{(1)} \cup \{x,y,z\}]$ is an induced forest in $G$, showing $a(G) \geq |F_1^{(1)}| + |F_2^{(1)}|+ 3 \geq \lceil (4n+3)/7 \rceil$, a contradiction. This completes the proof of Claim 2. 

\medskip

\textit{Claim 3: $wz_1 \not\in E(G)$. }

Otherwise, $wz_1 \in E(G)$, there exists a separation $(G_1,G_2)$ in $G$ such that $V(G_1 \cap G_2) = \{w,x,y,z_1\}$, $y_1 \in V(G_1)$, and $z \in V(G_2)$. 
For $i = 1,2$, let $F_i^{(1)} = A(G_i - \{w,x,y,z_1\})$; 
so $|F_i^{(1)}| \geq \lceil (4(|G_i|-4)+3)/7 \rceil$.
Now $G[F_1^{(1)} \cup F_2^{(1)} + \{x,y\}]$ is an induced forest in $G$, showing $a(G) $ $\geq |F_1^{(1)}| + |F_2^{(1)}| + 2 \geq \lceil (4n+3)/7 \rceil$, a contradiction. This completes the proof of Claim 3.

\medskip

We now distinguish several cases. 

\medskip
Case 1: $|N(y_1) \cap N(w_1)| \leq 2$, $|N(z_1) \cap N(z_2)| \leq 2$ and $|N(z_2) \cap N(z_3)| \leq 2$. 

Let $F' = A(G - \{w,x,y,z\}/\{y_1w_1,z_1z_2z_3\})$ with $w'$ (respectively, $z'$) as identifications of $\{y_1,w_1\}$ (respectively, $\{z_1,z_2,z_3\}$). Then $|F'| \geq \lceil (4(n-7)+3)/7 \rceil$. 
Let $F = F' + \{w,x,y,z\}$ if $w',z' \not\in F'$; $F = F' + \{w_1,y_1,x,y,z\}$ if $z' \not\in F', w' \in F'$; $F = F' + \{x,y,z_1,z_2,z_3\} - \{z'\}$ if $w' \not\in F', z' \in F'$; and $F = F' + \{w_1,y_1,x,z_1,z_2,z_3\} - \{w',z'\}$ if $w', z' \in F'$.
Therefore, $G[F]$ is an induced forest in $G$, giving $a(G) \geq |F'| + 4 \geq \lceil (4n+3)/7 \rceil$, a contradiction.





\medskip
Case 2: $|N(y_1) \cap N(w_1)| \geq 3$, $|N(z_1) \cap N(z_2)| \leq 2$ and $|N(z_2) \cap N(z_3)| \leq 2$. 

There exist $b_1 \in N(y_1) \cap N(w_1)$ and a separation $(G_1,G_2)$ such that $V(G_1 \cap G_2) = \{w_1,y_1,b_1\}$, $x \in V(G_1)$ and $N(y_1) \cap N(w_1) - \{w\} \subseteq V(G_2)$. 
Let $B_1 = \overline{B_2} = \{b_1\}$ and $B_2 = \overline{B_1} =\emptyset$. 
For $i=1,2$, let $F_1^{(i)} = A((G_1 - \{w,x,y,z,y_1,w_1\} - B_i)/\{z_1z_2z_3\})$ with $z'$ as the identification of $\{z_1,z_2,z_3\}$, and $F_2^{(i)} = A(G_2 - \{y_1,w_1\} - B_i)$. Then $|F_1^{(i)}| \geq \lceil (4(|G_1| -8- |B_i|)+3)/7 \rceil$, and $|F_2^{(i)}| \geq \lceil (4(|G_2|-2-|B_i|)+3)/7 \rceil$. Now $G[F_1^{(i)} \cup F_2^{(i)} + \{w,x,y,z\} - (\{b_1\}\cap(F_1^{(i)} \triangle F_2^{(i)}))]$ (if $z' \not\in F_1^{(i)}$) or $G[ (F_1^{(i)} - z') \cup F_2^{(i)} + \{x,y,z_1,z_2,z_3\}  - (\{b_1\}\cap(F_1^{(i)} \triangle F_2^{(i)}))]$ (if $z' \in F_1^{(i)}$) is an induced forest in $G$, showing $a(G) \geq |F_1^{(i)}| + |F_2^{(i)}| + 4 - |\overline{B_i}|$. By Lemma~\ref{ineq2}(2) (with $a = |G_1|-8, a_1 = |G_2|-2, c=4$), 
$(4(|G_1|-8)+3,4(|G_2|-2)+3) \equiv (4,0),(0,4) \mod 7$. 

\medskip
{\it Subcase 2.1}: $(4(|G_1|-8)+3,4(|G_2|-2)+3) \equiv (4,0) \mod 7$. 

Let $F_1^{(1)} = A((G_1 - \{w,x,y,z\})/\{y_1w_1,z_1z_2z_3\})$ with $w'$ (respectively $z'$) as the identification of $\{y_1,w_1\}$ (respectively $\{z_1,z_2,z_3\}$), and $F_2^{(1)} = A(G_2)$. 
Then $|F_1^{(1)}| \geq \lceil (4(|G_1|-7)+3)/7 \rceil$, and $|F_2^{(1)}| \geq \lceil (4|G_2|+3)/7 \rceil$.  
Let $F^{(1)} := \overline{F_1}^{(1)} \cup F_2^{(1)} - \{y_1,w_1,b_1\} \cap (\overline{F_1}^{(1)} \triangle F_2^{(1)})$ 
where $\overline{F_1}^{(1)} = F_1^{(1)} + \{w,x,y,z\}$ if $w',z' \not\in F_1^{(1)}$; $\overline{F_1}^{(1)} = F_1^{(1)} + \{w_1,y_1,x,y,z\}$ if $z' \not\in F_1^{(1)}, w' \in F_1^{(1)}$; $\overline{F_1}^{(1)} = F_1^{(1)} + \{x,y,z_1,z_2,z_3\} - \{z'\}$ if $w' \not\in F_1^{(1)}, z' \in F_1^{(1)}$; and $\overline{F_1}^{(1)} = F_1^{(1)} + \{w_1,y_1,x,z_1,z_2,z_3\} - \{w',z'\}$ if $w', z' \in F_1^{(1)}$.
Therefore, $G[F^{(1)}]$ is an induced forest in $G$, showing $a(G) \geq |F_1^{(1)}| + |F_2^{(1)}| + 4 - 3 \geq \lceil (4n+3)/7 \rceil$, a contradiction.

\medskip
{\it Subcase 2.2}: $(4(|G_1|-8)+3,4(|G_2|-2)+3) \equiv (0,4) \mod 7$. 
Let $F_1^{(2)} = A(G_1 - \{y_1,x,y,z\}/\{z_1z_2z_3\})$ with $z'$ as the identification of $\{z_1,z_2,z_3\}$, and $F_2^{(2)} = A(G_2 -y_1)$. Then $|F_1^{(1)}| \geq \lceil (4(|G_1|-6)+3)/7 \rceil$ and $|F_1^{(2)}| \geq \lceil (4(|G_2|-1)+3)/7 \rceil$. Now $G[F_1^{(2)} \cup F_2^{(2)} + \{x,y,z\} - (\{w_1,b_1\} \cap (F_1^{(2)} \triangle F_2^{(2)}))]$ (if $z' \not\in F_1^{(2)}$) or $G[(F_1^{(2)} - z') \cup F_2^{(2)} + \{x,z_1,z_2,z_3\} - (\{w_1,b_1\} \cap (F_1^{(2)} \triangle F_2^{(2)}))]$ (if $z' \in F_1^{(2)}$) is an induced forest in $G$, showing $a(G) \geq |F_1^{(2)}| + |F_2^{(2)}| + 3 - 2 \geq \lceil (4n+3)/7 \rceil$, a contradiction.

\medskip
Case 3: $|N(y_1) \cap N(w_1)| \leq 2$, $|N(z_1) \cap N(z_2)| > 2$, $|N(z_2) \cap N(z_3)| \leq 2$. 

There exist $a_1 \in N(z_1) \cap N(z_2)$ and a separation $(G_1,G_2)$ such that $V(G_1 \cap G_2) = \{z_1,z_2,a_1\}$, $x \in V(G_1)$, and $N(z_1) \cap N(z_2) - \{z\} \subseteq V(G_2)$. 
Let $A_1 = \overline{A_2} = \{a_1\}$ and $A_2 = \overline{A_1} =\emptyset$. 
For $i=1,2$, let $F_1^{(i)} = A((G_1 - \{w,x,y,z,z_1,z_2,z_3\} - A_i)/y_1w_1)$ with $w'$ as the identification of $\{y_1,w_1\}$, and $F_2^{(i)} = A(G_2 - \{z_1,z_2\} - A_i)$. 
Then $|F_1^{(i)}| \geq \lceil (4(|G_1|-8- |A_i|)+3)/7 \rceil$ and $|F_2^{(i)}| \geq \lceil (4(|G_2|-2-|A_i|)+3)/7 \rceil$. 
Now $G[F_1^{(i)} \cup F_2^{(i)} + \{w,x,y,z\}  - (\{a_1\} \cap (F_1^{(i)} \triangle F_2^{(i)}))] $ (if $w' \not\in F_1^{(i)}$) or $G[(F_1^{(i)} - w') \cup F_2^{(i)} + \{w_1,y_1,x,y,z\}  - (\{a_1\} \cap (F_1^{(i)} \triangle F_2^{(i)}))]$ (if $w' \in F_1^{(i)}$) is an induced forest in $G$, showing $a(G) \geq |F_1^{(i)}| + |F_2^{(i)}| + 4 - |\overline{A_i}|$. By Lemma~\ref{ineq2}(2) (with $a = |G_1|-8, a_1 = |G_2|-2, c=4$), 
$(4(|G_1|-8)+3,4(|G_2|-2)+3) \equiv (4,0),(0,4) \mod 7$. 

\medskip
{\it Subcase 3.1}: $(4(|G_1|-8)+3,4(|G_2|-2)+3) \equiv (4,0) \mod 7$. 

Let $F_1^{(3)} = A((G_1 - \{w,x,y,z\})/\{y_1w_1,z_1z_2z_3\})$ with $w'$ (respectively, $z'$) as the identification of $\{y_1,w_1\}$ (respectively, $\{z_1,z_2,z_3\}$), and $F_2^{(3)} = A(G_2)$. 
Then $|F_1^{(3)}| \geq \lceil (4(|G_1|-7)+3)/7 \rceil$, and $|F_2^{(3)}| \geq \lceil (4|G_2|+3)/7 \rceil$.  
Let $F^{(3)} := \overline{F_1}^{(3)} \cup F_2^{(3)} - (\{z_1,z_2,a_1\} \cap (\overline{F_1}^{(3)} \triangle F_2^{(3)}))$ 
where $\overline{F_1}^{(3)} = F_1^{(3)} + \{w,x,y,z\}$ if $w',z' \not\in F_1^{(3)}$; $\overline{F_1}^{(3)} = F_1^{(3)} + \{w_1,y_1,x,y,z\}$ if $z' \not\in F_1^{(3)}, w' \in F_1^{(3)}$; $\overline{F_1}^{(3)} = F_1^{(3)} + \{x,y,z_1,z_2,z_3\} - \{z'\}$ if $w' \not\in F_1^{(3)}, z' \in F_1^{(3)}$; and $\overline{F_1}^{(3)} = F_1^{(3)} + \{w_1,y_1,x,z_1,z_2,z_3\} - \{w',z'\}$ if $w', z' \in F_1^{(3)}$.
Therefore, $G[F^{(3)}]$ is an induced forest in $G$, showing $a(G) \geq |F_1^{(3)}| + |F_2^{(3)}| + 4 - 3 \geq \lceil (4n+3)/7 \rceil$, a contradiction.


\medskip
{\it Subcase 3.2}: $(4(|G_1|-8)+3,4(|G_2|-2)+3) \equiv (0,4) \mod 7$. 

Let $F_1^{(4)} = A((G_1 - \{z_1,x,y,z,w\})/z_2z_3)$ with $z'$ as the identification of $\{z_2,z_3\}$, and $F_2^{(4)} = A(G_2 -z_1)$. Then $|F_1^{(4)}| \geq \lceil (4(|G_1|-6)+3)/7 \rceil$, and $|F_2^{(4)}| \geq \lceil (4(|G_2|-1)+3)/7 \rceil$. Now $G[F_1^{(4)} \cup F_2^{(4)} + \{x,y,z\} - (\{z_2,a_1\}\cap (F_1^{(4)} \triangle F_2^{(4)}))]$ (if $z' \not\in F_1^{(2)}$) or $G[(F_1^{(4)} - z') \cup F_2^{(4)} + \{x,y,z_2,z_3\}  - (\{z_2,a_1\}\cap ((F_1^{(4)} \cup \{z_2\}) \triangle F_2^{(4)}))]$ (if $z' \in F_1^{(2)}$) is an induced forest in $G$, showing $a(G) \geq |F_1^{(2)}| + |F_2^{(2)}| + 3 - 2 \geq \lceil (4n+3)/7 \rceil$, a contradiction.

\medskip
Case 4: $|N(y_1) \cap N(w_1)| \leq 2$, $|N(z_1) \cap N(z_2)| \leq 2$, $|N(z_2) \cap N(z_3)| > 2$. 

There exist $c_1 \in N(z_2) \cap N(z_3)$ and a separation $(G_1,G_2)$ such that $V(G_1 \cap G_2) = \{z_2,z_3,c_1\}$, $x \in V(G_1)$, and $N(z_2) \cap N(z_3) - \{z\} \subseteq V(G_2)$. 
Let $C_1 = \overline{C_2} = \{c_1\}$ and $C_2 = \overline{C_1} =\emptyset$.
For $i=1,2$, let $F_1^{(i)} = A((G_1 - \{w,x,y,z,z_1,z_2,z_3\} - C_i)/y_1w_1)$ with $w'$ as the identification of $\{y_1,w_1\}$, and $F_2^{(i)} = A(G_2 - \{z_2,z_3\} - C_i)$. Then $|F_1^{(i)}| \geq \lceil (4(|G_1|-8- |C_i|)+3)/7 \rceil$, and $|F_2^{(i)}| \geq \lceil (4(|G_2|-2-|C_i|)+3)/7 \rceil$. 
Now $G[F_1^{(i)} \cup F_2^{(i)} + \{w,x,y,z\} - (\{c_1\} \cap (F_1^{(i)} \triangle F_2^{(i)}))]$ (if $w' \not\in F_1^{(i)}$) or $G[(F_1^{(i)} - w') \cup F_2^{(i)} \cup \{w_1,y_1,x,y,z\} - (\{c_1\} \cap (F_1^{(i)} \triangle F_2^{(i)}))]$ (if $w' \in F_1^{(i)}$) is an induced forest in $G$, showing $a(G) \geq |F_1^{(i)}| + |F_2^{(i)}| + 4 - \overline{C_i}$. 
By Lemma~\ref{ineq2}(2), 
$(4(|G_1|-8)+3,4(|G_2|-2)+3) \equiv (4,0),(0,4) \mod 7$. 

\medskip
{\it Subcase 4.1}: $(4(|G_1|-8)+3,4(|G_2|-2)+3) \equiv (4,0) \mod 7$. 

Let $F_1^{(3)} = A((G_1 - \{w,x,y,z\})/\{y_1w_1,z_1z_2z_3\})$ with $w'$ (respectively $z'$) as the identification of $\{y_1,w_1\}$ (respectively $\{z_1,z_2,z_3\}$) and $F_2^{(3)} = A(G_2)$. 
Then $|F_1^{(3)}| \geq \lceil (4(|G_1|-7)+3)/7 \rceil$ and $|F_2^{(3)}| \geq \lceil (4|G_2|+3)/7 \rceil$.  
Let $F^{(3)} := \overline{F_1}^{(3)} \cup F_2^{(3)}  - \{z_2,z_3,c_1\} \cap (\overline{F_1}^{(3)} \triangle F_2^{(3)})$
where $\overline{F_1}^{(3)} = F_1^{(3)} + \{w,x,y,z\}$ if $w',z' \not\in F_1^{(3)}$; $\overline{F_1}^{(3)} = F_1^{(3)} + \{w_1,y_1,x,y,z\}$ if $z' \not\in F_1^{(3)}, w' \in F_1^{(3)}$; $\overline{F_1}^{(3)} = F_1^{(3)} + \{x,y,z_1,z_2,z_3\} - \{z'\}$ if $w' \not\in F_1^{(3)}, z' \in F_1^{(3)}$; and $\overline{F_1}^{(3)} = F_1^{(3)} + \{w_1,y_1,x,z_1,z_2,z_3\} - \{w',z'\}$ if $w', z' \in F_1^{(3)}$.
Therefore, $G[F^{(3)}]$ is an induced forest in $G$, showing $a(G) \geq |F_1^{(3)}| + |F_2^{(3)}| + 4 - 3 \geq \lceil (4n+3)/7 \rceil$, a contradiction.


\medskip
{\it Subcase 4.2}: $(4(|G_1|-8)+3,4(|G_2|-2)+3) \equiv (0,4) \mod 7$. 

Let $F_1^{(4)} = A((G_1 - \{y_1,x,y,z,z_3\})/z_1z_2 + wz')$ with $z'$ as the identification of $\{z_1,z_2\}$, and $F_2^{(4)} = A(G_2 -z_3)$. 
Then $|F_1^{(4)}| \geq \lceil (4(|G_1|-6)+3)/7 \rceil$, and $|F_1^{(4)}| \geq \lceil (4(|G_2|-1)+3)/7 \rceil$. 
Now $G[F_1^{(4)} \cup F_2^{(4)} + \{x,y,z\} - (\{z_2,c_1\} \cap (F_1^{(4)} \triangle F_2^{(4)}))]$ (if $z' \not\in F_1^{(4)}$) or $G[(F_1^{(4)}-z') \cup F_2^{(4)} + \{x,y,z_1,z_2\}  - (\{z_2,c_1\} \cap ( (F_1^{(4)} \cup \{z_2\}) \triangle F_2^{(4)}))]$ (if $z' \in F_1^{(4)}$) is an induced forest in $G$, showing $a(G) \geq |F_1^{(4)}| + |F_2^{(4)}| + 3 - 2 \geq \lceil (4n+3)/7 \rceil$, a contradiction.

\medskip
Case 5: $|N(z_1) \cap N(z_2)|  > 2$, $|N(z_2) \cap N(z_3)| > 2$. 

\medskip
{\it Subcase 5.1}: $|N(y_1) \cap N(w_1)| \leq 2$.

There exist $a_1 \in N(z_1) \cap N(z_2)$, $c_1 \in N(z_2) \cap N(z_3)$ and subgraphs $G_1,G_2,G_3$ such that 
$G_2$ is the maximal subgraph of $G$ contained in the closed region of the plane bounded by the cycle $zz_1a_1z_2z$ and containing $N(z_1)\cap N(z_2) - \{z\}$, 
$G_3$ is the maximal subgraph of $G$ contained in the closed region of the plane bounded by the cycle $zz_3c_1z_2z$ and containing $N(z_3)\cap N(z_2) - \{z\}$, 
and $G_1$ is obtained from $G$ by removing $G_2-\{z_1,a_1,z_2\}$ and $G_3-\{z_3,c_1,z_2\}$. 
Let $A_i = \{a_1\}$ if $i=1,2$ and $\emptyset$ if $i=3,4$ and $\overline{A_i} = \{a_1\} - A_i$. Let $C_i = \{c_1\}$ if $i=1,3$ and $\emptyset$ if $i=2,4$ and $\overline{C_i} = \{c_1\} - C_i$.
For $i \in [4]$, let $F_1^{(i)} = A((G_1 - \{w,x,y,z,z_1,z_2,z_3\} - A_i - C_i)/y_1w_1)$ with $w'$ as the identification of $\{y_1,w_1\}$, $F_2^{(i)} = A(G_2 - \{z_1,z_2\} - A_i)$, and $F_3^{(i)} = A(G_3 - \{z_2,z_3\} - C_i)$. Then $|F_1^{(i)}| \geq \lceil (4(|G_1|-8 - |A_i|- |C_i|)+3)/7 \rceil$, $|F_2^{(i)}| \geq \lceil (4(|G_2|-2-|A_i|)+3)/7 \rceil$, and $|F_3^{(i)}| \geq \lceil (4(|G_3|-2-|C_i|)+3)/7 \rceil$. 
Now $G[F_1^{(i)} \cup F_2^{(i)} \cup F_3^{(i)} + \{w,x,y,z\} - (\{a_1\} \cap (F_1^{(i)} \triangle F_2^{(i)})) - (\{c_1\} \cap (F_1^{(i)} \triangle F_3^{(i)}))]$ (if $w' \not\in F_1^{(i)}$) or $G[(F_1^{(i)} - w') \cup F_2^{(i)} \cup F_3^{(i)} + \{w_1,y_1,x,y,z\} - (\{a_1\} \cap (F_1^{(i)} \triangle F_2^{(i)})) - (\{c_1\} \cap (F_1^{(i)} \triangle F_3^{(i)}))]$ (if $w' \in F_1^{(i)}$) is an induced forest in $G$, showing $a(G) \geq |F_1^{(i)}| + |F_2^{(i)}| + |F_3^{(i)}|  + 4 - |\overline{A_i}| - |\overline{C_i}|$. 
By Lemma~\ref{ineq2}(1) (with $k=2, a = |G_1| - 8, a_1 = |G_2| -2, a_2 = |G_3|-2, L = \emptyset, c=4$), 
$a(G) \geq \lceil (4n+3)/7 \rceil$, a contradiction.

\medskip
{\it Subcase 5.2}:  $|N(y_1) \cap N(w_1)| \geq 3$. 

There exist $a_1 \in N(z_1) \cap N(z_2)$, $b_1 \in N(y_1) \cap N(w_1)$, $c_1 \in N(z_2) \cap N(z_3)$ and subgraphs $G_1,G_2,G_3,G_4$ of $G$ such that 
$G_2$ is the maximal subgraph of $G$ contained in the closed region of the plane bounded by the cycle $zz_1a_1z_2z$ and containing $N(z_1)\cap N(z_2) - \{z\}$, 
$G_3$ is the maximal subgraph of $G$ contained in the closed region of the plane bounded by the cycle $zz_3c_1z_2z$ and containing $N(z_3)\cap N(z_2) - \{z\}$, 
$G_4$ is the maximal subgraph of $G$ contained in the closed region of the plane bounded by the cycle $ww_1b_1y_1w$ and containing $N(y_1)\cap N(w_1) - \{w\}$,
and $G_1$ is obtained from $G$ by removing $G_2-\{z_1,a_1,z_2\}$, $G_3-\{z_3,c_1,z_2\}$ and $G_4-\{w_1,b_1,y_1\}$. 
Let $A_i \subseteq \{a_1\}$ and $\overline{A_i} = \{a_1\} - A_i$. Let $B_i \subseteq \{b_1\}$ and $\overline{B_i} = \{b_1\} - B_i$. Let $C_i \subseteq \{c_1\}$  and $\overline{C_i} = \{c_1\} - C_i$. 
For each choice of $A_i,B_i,C_i$, let $F_1^{(i)} = A(G_1 - \{w,x,y,z,z_1,z_2,z_3,y_1,w_1\} - A_i - B_i - C_i)$, $F_2^{(i)} = A(G_2 - \{z_1,z_2\} - A_i)$, and $F_3^{(i)} = A(G_3 - \{z_2,z_3\} - C_i)$ and $F_4^{(i)} = A(G_4 - \{y_1,w_1\} - B_i)$. Then $|F_1^{(i)}| \geq \lceil (4(|G_1|-9 - |A_i|-|B_i|-|C_i|)+3)/7 \rceil$, $|F_2^{(i)}| \geq \lceil (4(|G_2|-2-|A_i|)+3)/7 \rceil$, $|F_3^{(i)}| \geq \lceil (4(|G_3|-2-|C_i|)+3)/7 \rceil$, and $|F_4^{(i)}| \geq \lceil (4(|G_4|-2-|B_i|)+3)/7 \rceil$. 
Now $G[F_1^{(i)} \cup F_2^{(i)} \cup F_3^{(i)} \cup F_4^{(i)} + \{w,x,y,z\} - (\{a_1\} \cap (F_1^{(i)} \triangle F_2^{(i)})) - (\{c_1\} \cap (F_1^{(i)} \triangle F_3^{(i)})) - (\{b_1\} \cap (F_1^{(i)} \triangle F_4^{(i)}))]$  is an induced forest in $G$, showing $|F_1^{(i)}| + |F_2^{(i)}| + |F_3^{(i)}| + |F_4^{(i)}| + 4 - |\overline{A_i}| -|\overline{B_i}| -|\overline{C_i}|$. 
By Lemma~\ref{ineq2}(1) 
(with $k=3, a = |G_1| - 9, a_1 = |G_2| - 2, a_2 = |G_3| - 2, a_3 = |G_4| - 2, L = \emptyset, c=4$), 
$a(G) \geq \lceil (4n+3)/7 \rceil$, a contradiction.


\medskip
Case 6: $|N(z_1) \cap N(z_2)| > 2$, $|N(z_2) \cap N(z_3)| \leq 2$ and $|N(y_1) \cap N(w_1)| > 2$.
 
There exist $a_1 \in N(z_1) \cap N(z_2), b_1 \in N(y_1) \cap N(w_1)$ and subgraphs $G_1,G_2,G_3$ of $G$ such that 
$G_2$ is the maximal subgraph of $G$ contained in the closed region of the plane bounded by the cycle $zz_1a_1z_2z$ and containing $N(z_1)\cap N(z_2) - \{z\}$,  
$G_3$ is the maximal subgraph of $G$ contained in the closed region of the plane bounded by the cycle $ww_1b_1y_1w$ and containing $N(y_1)\cap N(w_1) - \{w\}$,
and $G_1$ is obtained from $G$ by removing $G_2-\{z_1,a_1,z_2\}$ and $G_3-\{w_1,b_1,y_1\}$. 
Let $A_i = \{a_1\}$ if $i=1,2,5,6,9,10$ and $A_i = \emptyset$ if $i=3,4,7,8,11,12$. Let $\overline{A_i} = \{a_1\} - A_i$. Let $B_i = \{b_1\}$ if $i=1,3,5,7$ and $B_i = \emptyset$ if $i=2,4,6,8$, and $\overline{B_i} = \{b_1\} - B_i$. 
For $i = 1,2,3,4$, let $F_1^{(i)} = A(G_1 - \{w,x,y,z,z_1,z_2,z_3,y_1,w_1\} - A_i - B_i)$, $F_2^{(i)} = A(G_2 - \{z_1,z_2\} - A_i)$, and $F_3^{(i)} = A(G_3 - \{y_1,w_1\} - B_i)$. 
Note $|F_1^{(i)}| \geq \lceil (4(|G_1|-9 - |A_i|- |B_i|)+3)/7 \rceil$, $|F_2^{(i)}| \geq \lceil (4(|G_2|-2-|A_i|)+3)/7 \rceil$, and $|F_3^{(i)}| \geq \lceil (4(|G_3|-2-|B_i|)+3)/7 \rceil$. 
Now $G[F_1^{(i)} \cup F_2^{(i)} \cup F_3^{(i)} + \{w,x,y,z\} - (\{a_1\} \cap (F_1^{(i)} \triangle F_2^{(i)})) - (\{b_1\} \cap (F_1^{(i)} \triangle F_3^{(i)}))]$ is an induced forest in $G$, showing $a(G) \geq |F_1^{(i)}| + |F_2^{(i)}| + |F_3^{(i)}|  + 4 - |\overline{A_i}| - |\overline{B_i}|$. 
By Lemma \ref{ineq2}(5) 
(with $a = |G_1|-9, a_1 = |G_2|-2, a_2 = |G_3|-2, c=4$),
$(n_1,n_2,n_3):=(4(|G_1|-9)+3,4(|G_2|-2)+3,4(|G_3|-2)+3) \equiv (0,0,0), (1,0,0), (4,0,3), (4,3,0), (3,0,4), (4,0,4), (3,4,0), (4,4,0), (1,6,0), (1,0,6), (0,3,4),$ $ (0,4,3), (0,4,4), (6,4,4), (4,4,6), (4,6,4) \mod 7$.


We claim that $4(|G_3|-2)+3 \not\equiv 4 \mod 7$. For, suppose that $4(|G_3|-2)+3 \equiv 4 \mod 7$. 
If $|N(w_2) \cap N(z_3)| \leq 2$, then for $i=5,7$, let $F_1^{(i)} = A((G_1 - \{w,x,y,z,z_1,z_2,y_1\} - A_i)/w_2z_3)$ with $w'$ as the identification of $\{w_2,z_3\}$, $F_2^{(i)} = A(G_2 - \{z_1,z_2\} - A_i)$, and $F_3^{(i)} = A(G_3 - \{y_1\})$. 
Then $|F_1^{(i)}| \geq \lceil (4(|G_1|-8-|A_i|)+3)/7 \rceil$, $|F_2^{(i)}| \geq \lceil (4(|G_2|-2-|A_i|)+3)/7 \rceil$ and $|F_3^{(i)}| \geq \lceil (4(|G_3|-1)+3)/7 \rceil = (4(|G_3|-1)+3)/7 + 6/7$. 
Now $G[F_1^{(i)} \cup F_2^{(i)} \cup F_3^{(i)} + \{w,x,y,z\} - (\{w_1,b_1\} \cap (F_1^{(i)} \triangle F_3^{(i)})) - (\{a_1\} \cap (F_1^{(i)} \triangle F_2^{(i)}))]$ (if $w' \not\in F_1^{(i)}$) or $G[(F_1^{(i)} - w') \cup F_2^{(i)} \cup F_3^{(i)} + \{w_2,z_3,x,y,z\}  - (\{w_1,b_1\} \cap (F_1^{(i)} \triangle F_3^{(i)})) - (\{a_1\} \cap (F_1^{(i)} \triangle F_2^{(i)}))]$  is an induced forest in $G$, showing $a(G) \geq |F_1^{(i)}| + |F_2^{(i)}| + |F_3^{(i)}|  + 4 - 2 - (1-|A_i|)$.
By Lemma \ref{ineq2}(1)
(with $k = 1, a = |G_1| - 8, a_1 = |G_2| - 2, L = \{1\}, b_1 = |G_3| - 1, c= 2$),
  $a(G) \geq \lceil (4n+3)/7 \rceil$, a contradiction.
So $|N(w_2) \cap N(z_3)| > 2$. 
Then there exist $a_1 \in N(z_1) \cap N(z_2), b_1 \in N(y_1) \cap N(w_1), d_1 \in N(w_2) \cap N(z_3)$ and subgraphs $G_1',G_2',G_3',G_4'$ of $G$ such that $G_2' = G_2, G_3' = G_3$,   
$G_4'$ is the maximal subgraph of $G$ contained in the closed region of the plane bounded by the cycle $ww_2d_1z_3w$ containing $N(w_2)\cap N(z_3) - \{w\}$,
and $G_1$ is obtained from $G$ by removing $G_2-\{z_1,a_1,z_2\}$, $G_3-\{w_1,b_1,y_1\}$ and $G_4-\{w_2,d_1,z_3\}$. 
Let $D_i = \{d_1\}$ if $i=9,11$ and $D_i = \emptyset$ if $i=10,12$,  and let $\overline{D_i} = \{d_1\} - D_i$. 
For $i=9,10,11,12$, let $F_1^{(i)} = A(G_1' - \{w,x,y,z,z_1,z_2,y_1,w_2,z_3\} - A_i - D_i)$, $F_2^{(i)} = A(G_2' - \{z_1,z_2\} - A_i)$, $F_3^{(i)} = A(G_3' - \{y_1\})$, and $F_4^{(i)} = A(G_4' - \{w_2,z_3\} - D_i)$. 
Then $|F_1^{(i)}| \geq \lceil (4(|G_1'|-9-|A_i|-|D_i|)+3)/7 \rceil$, $|F_2^{(i)}| \geq \lceil (4(|G_2'|-2-|A_i|)+3)/7 \rceil$, 
$|F_3^{(i)}| \geq \lceil (4(|G_3'|-1)+3)/7 \rceil = \lceil (4(|G_3|-1)+3)/7 \rceil = (4(|G_3|-1)+3)/7 + 6/7$, 
and $|F_4^{(i)}| \geq \lceil (4(|G_4'|-2-|D_i|)+3)/7 \rceil$. 
Now $G[F_1^{(i)} \cup F_2^{(i)} \cup F_3^{(i)} \cup F_4^{(i)} + \{w,x,y,z\} - (\{w_1,b_1\} \cap (F_1^{(i)} \triangle F_3^{(i)}))  - (\{a_1\} \cap (F_1^{(i)} \triangle F_2^{(i)})) - (\{d_1\} \cap (F_1^{(i)} \triangle F_4^{(i)}))]$ is an induced forest in $G$, showing $a(G) \geq |F_1^{(i)}| + |F_2^{(i)}| + |F_3^{(i)}| + |F_4^{(i)}| + 4 - 2 - |\overline{A_i}| - |\overline{D_i}|$ 
By Lemma \ref{ineq2}(1) (with $k=2,a=|G_1'|-9,a_1=|G_2'|-2,a_2=|G_4'|-2,L=\{1\},b_1=|G_3'|-1,c=2$)
$a(G) \geq \lceil (4n+3)/7 \rceil$, a contradiction. 

Hence, $4(|G_3|-2)+3 \not\equiv 4 \mod 7$. 
Therefore, $(n_1,n_2,n_3) \equiv (0,0,0), (1,0,0), (4,0,3),$  $(4,3,0), (3,4,0), (4,4,0),$ 
$(1,6,0), (1,0,6), (0,4,3), (4,4,6) \mod 7$.

\medskip
{\it Subcase 6.1}: $(n_1,n_2,n_3) \equiv (0,0,0), (1,0,0) \mod 7$. 

Let $F_1^{(1)} = A((G_1 - \{w,x,y,z\})/\{y_1w_1,z_1z_2z_3\})$ with $w'$ (respectively $z'$) as the identifications of $\{y_1,w_1\}$ (respectively, $\{z_1,z_2,z_3\}$), $F_2^{(1)} = A(G_2)$, and $F_3^{(1)} = A(G_3)$. 
Then $|F_1^{(1)}| \geq \lceil (4(|G_1|-7)+3)/7 \rceil$, $|F_2^{(1)}| \geq \lceil (4|G_2|+3)/7 \rceil$ and $|F_3^{(1)}| \geq \lceil (4|G_3|+3)/7 \rceil$. 
Let $F^{(1)} := \overline{F_1}^{(1)} \cup F_2^{(1)} \cup F_3^{(1)}  - (\{z_1,z_2,a_1\} \cap (\overline{F_1}^{(1)} \triangle F_2^{(1)})) - (\{y_1,w_1,b_1\} \cap (\overline{F_1}^{(1)} \triangle F_3^{(1)}))$ 
where $\overline{F_1}^{(1)} = F_1^{(1)} + \{w,x,y,z\}$ if $w',z' \not\in F_1^{(1)}$; $\overline{F_1}^{(1)} = F_1^{(1)} + \{w_1,y_1,x,y,z\}$ if $z' \not\in F_1^{(1)}, w' \in F_1^{(1)}$; $\overline{F_1}^{(1)} = F_1^{(1)} + \{x,y,z_1,z_2,z_3\} - \{z'\}$ if $w' \not\in F_1^{(1)}, z' \in F_1^{(1)}$; and $\overline{F_1}^{(1)} = F_1^{(1)} + \{w_1,y_1,x,z_1,z_2,z_3\} - \{w',z'\}$ if $w', z' \in F_1^{(1)}$.
Therefore, $G[F^{(1)}]$ is an induced forest in $G$, showing $a(G) \geq |F_1^{(1)}| + |F_2^{(1)}| + 4 - 6 \geq \lceil (4n+3)/7 \rceil$, a contradiction.

\medskip
{\it Subcase 6.2}: $(n_1,n_2,n_3) \equiv (4,3,0), (4,4,0)\mod 7$ (respectively, $(1,6,0) \mod 7$). 

Let $A_2 = \overline{A_3} = \emptyset$ and $A_3 = \overline{A_2} = \{a_1\}$.
For $i=2$ (respectively, $i=3$), let $F_1^{(i)} = A((G_1 - \{w,x,y,z,z_1,z_2,z_3\} - A_i)/y_1w_1)$ with $w'$ as the identification of $\{y_1,w_1\}$, $F_2^{(i)} = A(G_2 - \{z_1,z_2\} - A_i)$, and $F_3^{(i)} = A(G_3)$. 
Then $|F_1^{(i)}| \geq \lceil (4(|G_1|-8-|A_i|)+3)/7 \rceil$, $|F_2^{(i)}| \geq \lceil (4(|G_2|-2-|A_i|)+3)/7 \rceil$ and $|F_3^{(i)}| \geq \lceil (4|G_3|+3)/7 \rceil$. 
Now $G[F_1^{(i)} \cup F_2^{(i)} \cup F_3^{(i)} + \{w,x,y,z\} - (\{y_1,w_1,b_1\} \cap (F_1^{(i)} \triangle F_3^{(i)})) - (\{a_1\} \cap (F_1^{(i)} \triangle F_2^{(i)}))]$ (if $w' \not\in F_1^{(i)}$) or $G[(F_1^{(i)} - w') \cup F_2^{(i)} \cup F_3^{(i)} + \{w_1,y_1,x,y,z\}  - (\{y_1,w_1,b_1\} \cap ((F_1^{(i)} \cup  \{y_1,w_1\}) \triangle F_3^{(i)})) - (\{a_1\} \cap (F_1^{(i)} \triangle F_2^{(i)}))]$ (if $w' \in F_1^{(i)}$) is an induced forest in $G$, showing $a(G) \geq |F_1^{(i)}| + |F_2^{(i)}| + |F_3^{(i)}| + 4 - 3 - |\overline{A_i}|\geq \lceil (4n+3)/7 \rceil$, a contradiction.

\medskip
{\it Subcase 6.3}: $(n_1,n_2,n_3) \equiv (0,4,3) \mod 7$ (respectively, $(4,4,6) \mod 7$). 

Let $B_4 = \overline{B_5} = \emptyset$ and $B_5 = \overline{B_4} = \{b_1\}$.
For $i=4$ (respectively, $i=5$), let $F_1^{(i)} = A(G_1 - \{w,x,y,z_1,z_3,y_1,w_1\} - B_i + w_2z)$, $F_2^{(i)} = A(G_2 - \{z_1\})$, and $F_3^{(i)} = A(G_3 - \{y_1,w_1\} - B_i)$.
Then $|F_1^{(i)}| \geq \lceil (4(|G_1|-7-|B_i|)+3)/7 \rceil$, $|F_2^{(i)}| \geq \lceil (4(|G_2|-1)+3)/7 \rceil$, and $|F_3^{(i)}| \geq \lceil (4(|G_3|-2-|B_i|)+3)/7 \rceil$. 
Now $G[F_1^{(i)} \cup F_2^{(i)} \cup F_3^{(i)} + \{w,x,y\} - (\{z_2,a_1\} \cap (F_1^{(i)} \triangle F_2^{(i)})) - (\{b_1\} \cap (F_1^{(i)} \triangle F_3^{(i)}))]$ is an induced forest in $G$, showing $a(G) \geq |F_1^{(i)}| + |F_2^{(i)}| + |F_3^{(i)}| + 3 - 2 - |\overline{B_i}| \geq \lceil (4n+3)/7 \rceil$, a contradiction.

\medskip
{\it Subcase 6.4}: $(n_1,n_2,n_3) \equiv (4,0,3) \mod 7$ (respectively, $(1,0,6) \mod 7$). 

Let $B_6 = \overline{B_7} = \emptyset$ and $B_7 = \overline{B_6} = \{b_1\}$.
For $i=6$ (resp. $i=7$), let $F_1^{(i)} = A((G_1 - \{w,x,y,z,y_1,w_1\} - B_i)/\{z_1z_2z_3\})$, $F_2^{(i)} = A(G_2)$, and $F_3^{(i)} = A(G_3 - \{y_1,w_1\} - B_i)$. 
Then $|F_1^{(i)}| \geq \lceil (4(|G_1|-8-|B_i|)+3)/7 \rceil$, 
$|F_2^{(i)}| \geq \lceil (4|G_2|+3)/7 \rceil$
and $|F_3^{(i)}| \geq \lceil (4(|G_3|-2-|B_i|)+3)/7 \rceil$. 
Now $G[F_1^{(i)} \cup F_2^{(i)} \cup F_3^{(i)} + \{w,x,y,z\} - (\{z_1,z_2,a_1\}\cap (F_1^{(i)} \triangle F_2^{(i)})) - (\{b_1\}\cap (F_1^{(i)} \triangle F_3^{(i)}))]$ (if $z' \not\in F_1^{(i)}$) or $G[(F_1^{(i)} - z') \cup F_2^{(i)} \cup F_3^{(i)} + \{x,y,z_1,z_2,z_3\}  - (\{z_1,z_2,a_1\}\cap ( (F_1^{(i)} \cup \{z_1,z_2\}) \triangle F_2^{(i)})) - (\{b_1\}\cap (F_1^{(i)} \triangle F_3^{(i)}))]$ (if $z' \in F_1^{(i)}$) is an induced forest in $G$, showing $a(G) \geq |F_1^{(i)}| + |F_2^{(i)}| + |F_3^{(i)}| + 4 - 3 - |\overline{B_i}| \geq \lceil (4n+3)/7 \rceil$, a contradiction.

\medskip
{\it Subcase 6.5}: $(n_1,n_2,n_3) \equiv (3,4,0) \mod 7$. 

Let $F_1^{(8)} = A(G_1 - \{w,x,y,z_1\} + zy_1)$, $F_2^{(8)} = A(G_2 - z_1)$ and $F_3^{(i)} = A(G_3)$. 
Then $|F_1^{(i)}| \geq \lceil (4(|G_1|-4)+3)/7 \rceil$,
$|F_2^{(i)}| \geq \lceil (4(|G_2|-1)+3)/7 \rceil$
and $|F_3^{(i)}| \geq \lceil (4|G_3|+3)/7 \rceil$. 
Now $G[F_1^{(8)} \cup F_2^{(8)} \cup F_3^{(8)} + \{x,y\} - (\{z_2,a_1\} \cap (F_1^{(8)} \triangle F_2^{(8)})) - (\{w_1,y_1,b_1\} \cap (F_1^{(8)} \triangle F_3^{(8)})) \}]$ is an induced forest in $G$, showing $a(G) \geq |F_1^{(8)}| + |F_2^{(8)}| + |F_3^{(8)}| + 2 - 5 \geq \lceil (4n+3)/7 \rceil$, a contradiction.

\medskip
Case 7: $|N(z_1) \cap N(z_2)| \leq 2$, $|N(z_2) \cap N(z_3)| > 2$ and $|N(y_1) \cap N(w_1)| > 2$. 

There exist $c_1 \in N(z_2) \cap N(z_3), b_1 \in N(y_1) \cap N(w_1)$ and subgraphs $G_1,G_2,G_3$ such that 
$G_2$ is the maximal subgraph of $G$ contained in the closed region of the plane bounded by the cycle $zz_2c_1z_3z$ and containing $N(z_2)\cap N(z_3) - \{z\}$,  
$G_3$ is the maximal subgraph of $G$ contained in the closed region of the plane bounded by the cycle $ww_1b_1y_1w$ and containing $N(y_1)\cap N(w_1) - \{w\}$,
and $G_1$ is obtained from $G$ by removing $G_2-\{z_2,c_1,z_3\}$ and $G_3-\{w_1,b_1,y_1\}$. 
Let $B_i = \{b_1\}$ if $i=1,2$ and $B_i =\emptyset$ if $i=3,4$ and $\overline{B_i} = \{b_1\} - B_i$. Let $C_i = \{c_1\}$ if $i=1,3$ and $C_i = \emptyset$ if $i=2,4$ and $\overline{C_i} = \{c_1\} - C_i$. 
For $i = 1,2,3,4$, let $F_1^{(i)} = A(G_1 - \{w,x,y,z,z_1,z_2,z_3,y_1,w_1\} - B_i - C_i)$, $F_2^{(i)} = A(G_2 - \{z_2,z_3\} - C_i)$, and $F_3^{(i)} = A(G_3 - \{y_1,w_1\} - B_i)$. 
Then $|F_1^{(i)}| \geq \lceil (4(|G_1|-9 - |B_i| - |C_i|)+3)/7 \rceil$, $|F_2^{(i)}| \geq \lceil (4(|G_2|-2-|C_i|)+3)/7 \rceil$ and $|F_3^{(i)}| \geq \lceil (4(|G_3|-2-|B_i|)+3)/7 \rceil$. 
Now $G[F_1^{(i)} \cup F_2^{(i)} \cup F_3^{(i)} + \{w,x,y,z\} - (\{c_1\} \cap (F_1^{(i)} \triangle F_2^{(i)})) - (\{b_1\} \cap (F_1^{(i)} \triangle F_3^{(i)}))]$ is an induced forest in $G$, showing $a(G) \geq |F_1^{(i)}| + |F_2^{(i)}| + |F_3^{(i)}|  + 4 - |\overline{C_i}| - |\overline{B_i}|$. 
By Lemma~\ref{ineq2}(5)
 (with $a = |G_1| - 9, a_1 = |G_2| - 2, a_2 = |G_3| - 2, c=4$),
$(n_1,n_2,n_3):=(4(|G_1|-9)+3,4(|G_2|-2)+3,4(|G_3|-2)+3) \equiv (0,0,0)$, $(1,0,0)$,  $(4,0,3)$, $(4,3,0)$, $(3,0,4)$, $(4,0,4)$, $(3,4,0)$, $(4,4,0)$, $(1,6,0)$, $(1,0,6)$, 
 $(0,3,4)$, $(0,4,3)$, $(0,4,4)$, $(6,4,4)$, $(4,4,6)$, $(4,6,4) \mod 7$.

\medskip
{\it Subcase 7.1}: $(n_1,n_2,n_3) \equiv (0,0,0), (1,0,0) \mod 7$. 

Let $F_1^{(1)} = A((G_1 - \{w,x,y,z\})/\{y_1w_1,z_1z_2z_3\})$ with $w'$ (respectively, $z'$) as the identification of $\{y_1,w_1\}$ (respectively, $\{z_1,z_2,z_3\}$), $F_2^{(1)} = A(G_2)$, and $F_3^{(1)} = A(G_3)$. 
Then $|F_1^{(1)}| \geq \lceil (4(|G_1|-7)+3)/7 \rceil$, $|F_2^{(1)}| \geq \lceil (4|G_2|+3)/7 \rceil$, and $|F_3^{(1)}| \geq \lceil (4|G_3|+3)/7 \rceil$. 
Let $F^{(1)} := \overline{F_1}^{(1)} \cup F_2^{(1)} \cup F_3^{(1)} - \{z_3,z_2,c_1\} \cap (\overline{F_1}^{(1)} \triangle F_2^{(1)}) - \{y_1,w_1,b_1\} \cap (\overline{F_1}^{(1)} \triangle F_3^{(1)})$ 
where $\overline{F_1}^{(1)} = F_1^{(1)} + \{w,x,y,z\}$ if $w',z' \not\in F_1^{(1)}$; $\overline{F_1}^{(1)} = F_1^{(1)} + \{w_1,y_1,x,y,z\}$ if $z' \not\in F_1^{(1)}, w' \in F_1^{(1)}$; $\overline{F_1}^{(1)} = F_1^{(1)} + \{x,y,z_1,z_2,z_3\} - \{z'\}$ if $w' \not\in F_1^{(1)}, z' \in F_1^{(1)}$; and $\overline{F_1}^{(1)} = F_1^{(1)} + \{w_1,y_1,x,z_1,z_2,z_3\} - \{w',z'\}$ if $w', z' \in F_1^{(1)}$.
Therefore, $G[F^{(1)}]$ is an induced forest in $G$, showing $a(G) \geq |F_1^{(1)}| + |F_2^{(1)}| + 4 - 6 \geq \lceil (4n+3)/7 \rceil$, a contradiction.

\medskip
{\it Subcase 7.2}: $(n_1,n_2,n_3) \equiv (4,4,0), (3,4,0) \mod 7$.

Let $F_1^{(2)} = A((G_1 - \{w,x,y,z_1,z_3\})/y_1w_1 + \{w'z, w_2z\})$ with $w'$ as the identification of $\{y_1,w_1\}$, $F_2^{(2)} = A(G_2 - \{z_3\})$, and $F_3^{(2)} = A(G_3)$. 
Then $|F_1^{(2)}| \geq \lceil (4(|G_1|-6)+3)/7 \rceil$, $|F_2^{(2)}| \geq \lceil (4(|G_2|-1)+3)/7 \rceil$ and  $|F_3^{(2)}| \geq \lceil (4|G_3|+3)/7 \rceil$. 
Now $G[F_1^{(2)} \cup F_2^{(2)} \cup F_3^{(2)} + \{w,x,y\} - (\{z_2,c_1\} \cap (F_1^{(2)} \triangle F_2^{(2)})) -  (\{w_1,y_1,b_1\} \cap (F_1^{(2)} \triangle F_3^{(2)}))  ]$ (if $w' \not\in F_1^{(2)}$) or $G[(F_1^{(2)} - w') \cup F_2^{(2)} \cup F_3^{(2)} + \{w_1,y_1,x,y\} - (\{z_2,c_1\} \cap (F_1^{(2)} \triangle F_2^{(2)})) -  (\{w_1,y_1,b_1\} \cap ((F_1^{(2)} \cup \{w_1,y_1\}) \triangle F_3^{(2)}))]$ is an induced forest in $G$, showing $a(G) \geq |F_1^{(2)}| + |F_2^{(2)}| + |F_3^{(2)}| + 3 - 5 \geq \lceil (4n+3)/7 \rceil$, a contradiction.

\medskip
{\it Subcase 7.3}: $(n_1,n_2,n_3) \equiv (0,4,3), (0,4,4) \mod 7$ (respectively, $(4,4,6) \mod 7$).

Let $B_3 = \overline{B_4} = \emptyset$ and $B_4 = \overline{B_3} = \{b_1\}$.
For $i=3$ (respectively, $i=4$), let $F_1^{(i)} = A(G_1 - \{w,x,y,z_1,z_3,y_1,w_1\} - B_i + w_2z)$, $F_2^{(i)} = A(G_2 - z_3)$, and $F_3^{(i)} = A(G_3 -\{y_1,w_1\} - B_i)$. 
Then $|F_1^{(i)}| \geq \lceil (4(|G_1|-7-|B_i|)+3)/7 \rceil$, 
$|F_2^{(i)}| \geq \lceil (4(|G_2|-1)+3)/7 \rceil$
and $|F_3^{(i)}| \geq \lceil (4(|G_3|-2-|B_i|)+3)/7 \rceil$. 
Now $G[F_1^{(i)} \cup F_2^{(i)} \cup F_3^{(i)} + \{w,x,y\} - (\{b_1\} \cap (F_1^{(i)} \triangle F_3^{(i)})) - (\{z_2,c_1\} \cap (F_1^{(i)} \triangle F_2^{(i)}))]$ is an induced forest in $G$, giving $a(G) \geq |F_1^{(i)}| + |F_2^{(i)}| + |F_3^{(i)}| + 3 - 2 - |\overline{B_i}| \geq \lceil (4n+3)/7 \rceil$, a contradiction.

\medskip
{\it Subcase 7.4}: $(n_1,n_2,n_3) \equiv (3,0,4),(4,0,4),(6,4,4),(4,6,4) \mod 7$.

Let $F_1^{(5)} = A(G_1 - \{x,y,z,z_1,z_3,y_1\} + z_2w)$, $F_2^{(5)} = A(G_2 - z_3)$, and $F_3^{(5)} = A(G_3 -\{y_1\})$. 
Then $|F_1^{(5)}| \geq \lceil (4(|G_1|-6)+3)/7 \rceil$, 
$|F_2^{(5)}| \geq \lceil (4(|G_2|-1)+3)/7 \rceil$
and $|F_3^{(5)}| \geq \lceil (4(|G_3|-1)+3)/7 \rceil$. 
Now $G[F_1^{(5)} \cup F_2^{(5)} \cup F_3^{(5)} + \{z,x,y\} - (\{b_1,w_1\} \cap (F_1^{(i)} \triangle F_3^{(i)})) - (\{z_2,c_1\} \cap (F_1^{(i)} \triangle F_2^{(i)}))]$ is an induced forest in $G$, showing $a(G) \geq |F_1^{(5)}| + |F_2^{(5)}| + |F_3^{(5)}| + 3 - 4 \geq \lceil (4n+3)/7 \rceil$, a contradiction.

\medskip
{\it Subcase 7.5}: $(n_1,n_2,n_3) \equiv (4,3,0) \mod 7$ (respectively, $(1,6,0) \mod 7$).

Let $C_6 = \overline{C_7} = \emptyset$ and $C_7 = \overline{C_6} = \{c_1\}$.
For $i=6$ (respectively, $i=7$), let $F_1^{(i)} = A((G_1 - \{w,x,y,z,z_1,z_2,z_3\} - C_i)/y_1w_1)$ with $w'$ as the identification of $\{y_1,w_1\}$, $F_2^{(i)} = A(G_2 - \{z_2,z_3\} - C_i)$, and $F_3^{(i)} = A(G_3)$. 
Then $|F_1^{(i)}| \geq \lceil (4(|G_1|-8-|C_i|)+3)/7 \rceil$, $|F_2^{(i)}| \geq \lceil (4(|G_2|-2-|C_i|)+3)/7 \rceil$ 
and $|F_3^{(i)}| \geq \lceil (4|G_3|+3)/7 \rceil$. 
Now $G[F_1^{(i)} \cup F_2^{(i)} \cup F_3^{(i)} + \{w,x,y,z\} - (\{y_1,b_1,w_1\} \cap (F_1^{(i)} \triangle F_3^{(i)})) - (\{c_1\} \cap (F_1^{(i)} \triangle F_2^{(i)}))]$ (if $w' \not\in F_1^{(i)}$) or $G[(F_1^{(i)} - w') \cup F_2^{(i)} \cup F_3^{(i)} + \{w_1,y_1,x,y,z\}  - (\{y_1,b_1,w_1\} \cap ((F_1^{(i)} \cup \{w_1,y_1\}) \triangle F_3^{(i)})) - (\{c_1\} \cap (F_1^{(i)} \triangle F_2^{(i)}))]$  is an induced forest in $G$, showing $a(G) \geq |F_1^{(i)}| + |F_2^{(i)}| + |F_3^{(i)}| + 4 - 3  - |\overline{C_i}| \geq \lceil (4n+3)/7 \rceil$, a contradiction.  

\medskip
{\it Subcase 7.6}: $(n_1,n_2,n_3) \equiv (4,0,3) \mod 7$ (respectively $(1,0,6) \mod 7$).

Let $B_8 = \overline{B_9} = \emptyset$ and $B_9 = \overline{B_8} = \{b_1\}$.
For $i=8$ (respectively, $i=9$), let $F_1^{(i)} = A((G_1 - \{w,x,y,z,y_1,w_1\} - B_i)/\{z_1z_2z_3\})$ with $z'$ as the identification of $\{z_1,z_2,z_3\}$, $F_2^{(i)} = A(G_2)$, and $F_3^{(i)} = A(G_3 - \{y_1,w_1\} - B_i)$. 
Then $|F_1^{(i)}| \geq \lceil (4(|G_1|-8-|B_i|)+3)/7 \rceil$,
$|F_2^{(i)}| \geq \lceil (4|G_2|+3)/7 \rceil$
and $|F_3^{(i)}| \geq \lceil (4(|G_3|-2-|B_i|)+3)/7 \rceil$.
Now $G[F_1^{(i)} \cup F_2^{(i)} \cup F_3^{(i)} + \{w,x,y,z\} - (\{b_1\} \cap (F_1^{(i)} \triangle F_3^{(i)})) - (\{z_2,z_3,c_1\} \cap (F_1^{(i)} \triangle F_2^{(i)}))]$ (if $z' \not\in F_1^{(i)}$) or $G[(F_1^{(i)} - z') \cup F_2^{(i)} \cup F_3^{(i)} + \{x,y,z_1,z_2,z_3\}  - (\{b_1\} \cap (F_1^{(i)} \triangle F_3^{(i)})) - (\{z_2,z_3,c_1\} \cap ((F_1^{(i)} \cup \{z_2,z_3\}) \triangle F_2^{(i)}))]$ is an induced forest in $G$, showing $a(G) \geq |F_1^{(i)}| + |F_2^{(i)}| + |F_3^{(i)}| + 4 - 3 - |\overline{B_i}| \geq \lceil (4n+3)/7 \rceil$, a contradiction.   

\medskip
{\it Subcase 7.7}: $(n_1,n_2,n_3) \equiv (0,3,4) \mod 7$.

Let $F_1^{(10)} = A(G_1 - \{x,y,z,y_1,z_1,z_2,z_3\})$, $F_2^{(10)} = A(G_2 - \{z_2,z_3\})$, and $F_3^{(10)} = A(G_3 -\{y_1\})$. 
Then $|F_1^{(10)}| \geq \lceil (4(|G_1|-7)+3)/7 \rceil$, 
$|F_2^{(10)}| \geq \lceil (4(|G_2|-2)+3)/7 \rceil$
and $|F_3^{(10)}| \geq \lceil (4(|G_3|-1)+3)/7 \rceil$. 
Now $F^{(10)} := G[F_1^{(10)} \cup F_2^{(10)} \cup F_3^{(10)} + \{z,x,y\} - (\{b_1,w_1\} \cap (F_1^{(i)} \triangle F_3^{(i)})) - (\{c_1\} \cap (F_1^{(i)} \triangle F_2^{(i)}))]$ is an induced forest in $G$, showing $a(G) \geq |F_1^{(10)}| + |F_2^{(10)}| + |F_3^{(10)}| + 3 - 3 \geq \lceil (4n+3)/7 \rceil$, a contradiction.  
\qed

\section{Configurations around 5-vertices and 6-vertices}

First, we define certain configurations around a $5$-vertex or $6$-vertex.

\begin{defi}
\label{5structure} Let $x$ be a $5$-vertex in $G$ and $x_1,x_2,x_3,x_4,x_5$ be neighbors of $x$ in cyclic order around $x$. 
\begin{itemize}
\item[(i)] $x$ is of type \textbf{5-2-A} if $\{x_1,x_3\} \subseteq V_3, \{x_2,x_4,x_5\} \subseteq V_{\geq 4}$ such that if $N(x_1) = \{x_1',x_1'',x\}$ and $N(x_3) = \{x_3',x_3'',x\}$, then for $v \in  \{x_1',x_1''\}$, either $v \in V_{\leq 4}$ or $R_{v,\{x_1\}} \neq \emptyset$; and for $u \in  \{x_3',x_3''\}$, either $u \in V_{\leq 4}$ or $R_{u,\{x_3\}} \neq \emptyset$;
\item[(ii)] $x$ is of type \textbf{5-2-B} if $\{x_1,x_3\} \subseteq V_3, \{x_2,x_4,x_5\} \subseteq V_{\geq 4}$ such that if $N(x_1) = \{x_1',x_1'',x\}$ and $N(x_3) = \{x_3',x_3'',x\}$, then for $v \in  \{x_3',x_3''\}$, either $v \in V_{\leq 4}$ or $R_{v,\{x_3\}} \neq \emptyset$; and $x_1' \in V_{\geq 5}$ and $R_{x_1',\{x_1\}} = \emptyset$;
\item[(iii)] $x$ is of type \textbf{5-2-C} if $\{x_1,x_3\} \subseteq V_3, \{x_2,x_4,x_5\} \subseteq V_{\geq 4}$ such that if $N(x_1) = \{x_1',x_1'',x\}$ and $N(x_3) = \{x_3',x_3'',x\}$, then $x_1' \in V_{\geq 5}, R_{x_1',\{x_1\}} = \emptyset$, $x_3' \in V_{\geq 5}$ and $R_{x_3',\{x_3\}} = \emptyset$;
\item[(iv)] $x$ is of type \textbf{5-1-A} if $x_1 \in V_3$, $\{x_2,x_3,x_4,x_5\} \subseteq V_{\geq 4}$ such that if $N(x_1) = \{x_1',x_1'',x\}$, then for $v \in  \{x_1',x_1''\}$, either $v \in V_{\leq 4}$ or $R_{v,\{x_1\}} \neq \emptyset$;
\item[(v)] $x$ is of type \textbf{5-1-B} if $x_1 \in V_3$, $\{x_2,x_3,x_4,x_5\} \subseteq V_{\geq 4}$ such that if $N(x_1) = \{x_1',x_1'',x\}$, then $x_1' \in V_{\geq 5}$ and $ R_{x_1',\{x_1\}} = \emptyset$;
\item[(vi)] $x$ is of type \textbf{5-0} if $\{x_1,x_2,x_3,x_4,x_5\} \subseteq V_{\geq 4}$.
\end{itemize}
\end{defi}

\begin{defi}
\label{6structure} Let $v$ be a $6$-vertex in $G$ and $v_1,v_2,v_3,v_4,v_5,v_6$ be neighbors of $v$ in cyclic order around $v$. 
\begin{itemize}
\item[(i)] $v$ is of type \textbf{6-3} if $\{v_1,v_3,v_5\} \subseteq V_3$ and $\{v_2,v_4,v_6\} \subseteq V_{\geq 4}$;
\item[(ii)] $v$ is of type \textbf{6-2-A} if $\{v_1,v_3\} \subseteq V_3$ and $\{v_2,v_4,v_5,v_6\} \subseteq V_{\geq 4}$;
\item[(iii)] $v$ is of type \textbf{6-2-B} if $\{v_1,v_4\} \subseteq V_3$ and $\{v_2,v_3,v_5,v_6\} \subseteq V_{\geq 4}$;
\item[(iv)] $v$ is of type \textbf{6-1} if $\{v_1\} \subseteq V_3$ and $\{v_2,v_3,v_4,v_5,v_6\} \subseteq V_{\geq 4}$;
\item[(v)] $v$ is of type \textbf{6-0} if $\{v_1,v_2,v_3,v_4,v_5,v_6\} \subseteq V_{\geq 4}$.
\end{itemize}
\end{defi}

\begin{lem}
\label{5-2-B-lemma3}
The following configuration is impossible in $G$: 
$x$ is a $5$-vertex of type 5-2-B with neighbors $x_1,y,x_3,z,x_2$ in cyclic order around $x$, $\{y,z\} \subseteq V_3, x_1 \in V_4$, $N(z) = \{z_1,z_2\}$ with $\{z_1x_2, z_2x_3\} \subseteq E(G)$, $\{z_1,z_2\} \subseteq V_4$, and $xx_2wx_1x$ forms a facial cycle where $w \in V_3$. 
\end{lem}


\pf 
Let $N(w) = \{x_1,w_1,x_2\}$, $N(z_1) = \{z,x_2,s_1,s_2\}$ and $N(z_2) = \{z,x_3,t,s_2\}$. 

First, we claim that $|N(x_1) \cap N(y)| \leq 2$. For otherwise, suppose $N(x_1) \cap N(y) = \{x,p_1,p_2\}$. 
There exists a separation $(G_1,G_2)$ such that $V(G_1 \cap G_2) = \{p_1,p_2\}$, $\{x,y,x_1\} \subseteq V(G_1)$, and $N(p_1) \cap N(p_2) - \{y,x_1\} \subseteq V(G_2)$. 
Let $F_1^{(1)} = A(G_1 - \{x,y,x_1,p_1,p_2\})$, and $F_2^{(1)} = A(G_2 - p_2$).
Then $|F_1^{(1)}| \geq \lceil (4(|G_1|-5)+3)/7 \rceil$, and $|F_2^{(1)}| \geq \lceil (4(|G_2|-1)+3)/7 \rceil$. 
Now $G[F_1^{(1)} \cup F_2^{(1)} + \{x_1,y\}]$  is an induced forest in $G$, showing $a(G) \geq |F_1^{(1)}| + |F_2^{(1)}| + 2 \geq \lceil (4n+3)/7 \rceil$, a contradiction. Thus, let $N(x_1) \cap N(y) = \{x,y_1\}$.


By Lemma \ref{No434Edge}, $z_2x_2 \not\in E(G)$ and  $z_1x_3 \not\in E(G)$.

We also claim that $w_1z_2 \not\in E(G)$. Otherwise, 
there exists a separation $(G_1,G_2)$ such that $V(G_1 \cap G_2) = \{w_1,x_2,z,z_2\}$, $\{x,y,w,x_1\} \subseteq V(G_1)$, and $\{z_1,s_1,s_2\} \subseteq V(G_2)$. 
Let $F_1^{(3)} = A(G_1 - \{w_1,x_2,z,z_2,w,x,x_1,y,y_1\})$, and $F_2^{(3)} = A(G_2 - \{w_1,x_2,z,z_2,z_1,s_2\})$. 
Then $|F_1^{(3)}| \geq \lceil (4(|G_1|-9)+3)/7 \rceil$, and $|F_2^{(3)}| \geq \lceil (4(|G_2|-6)+3)/7 \rceil$. 
Now $G[F_1^{(3)} \cup F_2^{(3)} + \{z,z_1,z_2,w,x_1,y\}]$  is an induced forest in $G$, showing $a(G) \geq |F_1^{(3)}| + |F_2^{(3)}| + 6 \geq \lceil (4n+3)/7 \rceil$, a contradiction.


We further claim that $s_1z_2 \not\in E(G)$. Otherwise, there exists a separation $(G_1,G_2)$ such that $V(G_1 \cap G_2) = \{s_1,z_1,z_2\}$, $\{x,y,w,x_1,z\} \subseteq V(G_1)$, and $s_2 \in V(G_2)$. 
Let $F_1^{(4)} = A(G_1 - \{s_1,z_1,z_2,z\})$, and $F_2^{(4)} = A(G_2 - \{s_1,z_1,z_2,s_2\})$. 
Then $|F_1^{(4)}| \geq \lceil (4(|G_1|-4)+3)/7 \rceil$, and $|F_2^{(4)}| \geq \lceil (4(|G_2|-4)+3)/7 \rceil$. 
Now $G[F_1^{(4)} \cup F_2^{(4)} + \{z_1,z_2\}]$  is an induced forest in $G$, showing $a(G) \geq |F_1^{(4)}| + |F_2^{(4)}| + 2$. 
By Lemma \ref{ineq2}(8)
(with $a = |G_1| - 4, a_1 =  |G_2| - 4, c = 2$)
$(4(|G_1| - 4)+3,4(|G_2| - 4)+3) \equiv (0,0),(0,6),(0,5),(5,0),(6,6),(6,0) \mod 7$.
Let $F_1^{(5)} = A(G_1 - \{s_1,z_1,z_2,x,z\})$, and $F_2^{(5)} = A(G_2 - \{s_1,z_1,z_2\})$. 
Then $|F_1^{(5)}| \geq \lceil (4(|G_1|-5)+3)/7 \rceil$, and $|F_2^{(5)}| \geq \lceil (4(|G_2|-3)+3)/7 \rceil$. 
Define $G[F_1^{(5)} \cup F_2^{(5)} + \{z,z_2\}]$  is an induced forest in $G$, showing $a(G) \geq |F_1^{(5)}| + |F_2^{(5)}| + 2 \geq \lceil (4n+3)/7 \rceil$, a contradiction.

Note that $s_1x \not\in E(G)$. Otherwise, since $G$ is simple, $s_1 \not\in \{x_2,z\}$. $s_1 \not\in \{x_1,y\}$ by Lemma \ref{basic_lemma} ($G$ is a quadrangulation). $s_1 \neq x_3$ by second claim. Similarly, $tx \not\in E(G)$.

\medskip

We now distinguish several cases.

\medskip
Case 1: $|N(w_1) \cap N(x_2)| \leq 2$ and $|N(s_1) \cap N(s_2)| \leq 2$.

Let $F' = A((G - \{w,x,z,z_1\} )/\{x_1y,w_1x_2,s_1s_2\} + z_2u_2)$ with $u_1$ (respectively, $u_2,u_3$) as the identification of  $\{x_1,y\}$ (respectively, $\{w_1,x_2\},\{s_1,s_2\}$). 
Then $|F'| \geq \lceil (4(n-7)+3)/7 \rceil$. 
Note $u_1 \in F'$ by Lemma \ref{basic_lemma} since $|N(u_1)| = 3$. 
Let $F = F' + \{x_1,y,z,z_1,w\} - \{u_1\}$ if $u_2,u_3 \not\in F'$, and otherwise, $F$ obtained from $F' + \{x_1,y,z,z_1,w\} - \{u_1\}$ by deleting $\{u_2,w\}$ (respectively, $\{u_3,z_1\}$) and adding $\{w_1,x_2\}$ (respectively, $\{s_1,s_2\}$) when $u_2 \in F'$ (respectively, $u_3 \in F'$).
Therefore, $G[F]$ is an induced forest in $G$, showing $a(G) \geq |F'| + 4 \geq \lceil (4n+3)/7 \rceil$, a contradiction.

\medskip
Case 2: $|N(w_1) \cap N(x_2)| \leq 2$ and $|N(s_1) \cap N(s_2)| > 2$.

There exist $a_1 \in N(s_1) \cap N(s_2)$ and a separation $(G_1,G_2)$ such that $V(G_1 \cap G_2) = \{s_1,s_2,a_1\}$, $\{x,y,w,z,x_1,x_2,x_3\} \subseteq V(G_1)$, and $N(s_1) \cap N(s_2) - \{z_1\} \subseteq V(G_2)$. 
Let $A_1 = \{a_1\}$ and $A_2 = \emptyset$. 
For $i=1,2$, let $F_1^{(i)} = A((G_1 - \{w,x,z,z_1,s_1,s_2\} - A_i)/\{x_1y,w_1x_2\} + u_2z_2)$ with $u_1$ (respectively, $u_2$) as the identification of $\{x_1,y\}$ (respectively, $\{w_1,x_2\}$), and $F_2^{(i)} = A(G_2 - \{s_1,s_2\} - A_i)$. 
Then $|F_1^{(i)}| \geq \lceil (4(|G_1|-8-|A_i|)+3)/7 \rceil$, and $|F_2^{(i)}| \geq \lceil (4(|G_2|-2-|A_i|)+3)/7 \rceil$. 
Note $u_1 \in F_1^{(i)}$ by Lemma \ref{basic_lemma} since $|N(u_1)| = 3$. 
Let $F = (F_1^{(i)} - u_1) \cup F_2^{(i)} +  \{x_1,y,z,z_1,w\}  - (\{a_1\} \cap (F_1^{(i)} \triangle F_2^{(i)}))$.
Now $G[F]$ (if $u_2 \not\in F_1^{(i)}$) or $G[F - \{u_2,w\} + \{w_1,x_2\}]$ (if $u_2 \in F_1^{(i)}$)
is an induced forest in $G$, showing $a(G) \geq |F_1^{(i)}| + |F_2^{(i)}| + 4 - (1-|A_i|)$. 
By Lemma \ref{ineq2}(2) 
(with $a = |G_1|-8, a_1 = |G_2|-2, c=4$), 
$(4(|G_1|-8)+3,4(|G_2|-2)+3) \equiv (4,0),(0,4) \mod 7$.

\medskip
{\it Subcase 2.1:} $(4(|G_1|-8)+3,4(|G_2|-2)+3) \equiv (4,0) \mod 7$.

Let $F_1^{(3)} = A((G_1 - \{w,x,z,z_1\})/\{x_1y,w_1x_2,s_1s_2\} + z_2u_2)$  with $u_1$ (respectively, $u_2,u_3$) as the identification of  $\{x_1,y\}$ (respectively, $\{w_1,x_2\},\{s_1,s_2\}$), and $F_2^{(3)} = A(G_2)$. 
Then $|F_1^{(3)}| \geq \lceil (4(|G_1|-7)+3)/7 \rceil$ and $|F_2^{(3)}| \geq \lceil (4|G_2|+3)/7 \rceil$. 
Note $u_1 \in F_1^{(4)}$ by Lemma \ref{basic_lemma} since $|N(u_1)| = 3$. 
Let $F^{(3)} := \overline{F_1}^{(3)} \cup F_2^{(3)}  - (\{s_1,s_2,a_1\} \cap (\overline{F_1}^{(3)} \triangle F_2^{(3)}))$,
where $\overline{F_1}^{(3)} = {F_1}^{(3)} + \{x_1,y,z,z_1,w\} - u_1$ if $u_2,u_3 \not\in F_1^{(3)}$, and otherwise, $\overline{F_1}^{(3)}$ obtained from ${F_1}^{(3)} + \{x_1,y,z,z_1,w\} - u_1$ by deleting $\{u_2,w\}$ (respectively, $\{u_3,z_1\}$) and adding $\{w_1,x_2\}$ (respectively, $\{s_1,s_2\}$) when $u_2 \in {F_1}^{(3)}$ (respectively, $u_3 \in {F_1}^{(3)}$).
Therefore, $G[F^{(3)}]$ is an induced forest in $G$, showing $a(G) \geq |F_1^{(3)}| + |F_2^{(3)}| + 4 - 3 \geq \lceil (4n+3)/7 \rceil$, a contradiction.

\medskip
{\it Subcase 2.2:} $(4(|G_1|-8)+3,4(|G_2|-2)+3) \equiv (0,4) \mod 7$.

If $|N(x_3) \cap N(t)| \leq 2$, 
 let $F_1^{(4)} = A((G_1 - \{s_2,x_2,z,z_1,z_2\})/x_3t + s_1x)$ with $u$ as the identification of $\{x_3,t\}$, and $F_2^{(4)}= A(G_2 - \{s_2\})$. 
Then $|F_1^{(4)}| \geq \lceil (4(|G_1|-6)+3)/7 \rceil$, and $|F_2^{(4)}| \geq \lceil (4(|G_2|-1)+3)/7 \rceil$. 
Let $F = F_1^{(4)} \cup F_2^{(4)} + \{z,z_1,z_2\} - (\{s_1,a_1\} \cap (F_1^{(4)} \triangle F_2^{(4)}))$.
Now $G[F]$ (if $u \not\in F_1^{(4)}$) or $G[F - \{u,z_2\}+\{x_3,t\}]$ (if $u \in F_1^{(4)}$)
is an induced forest in $G$, showing $a(G) \geq |F_1^{(4)}| + |F_2^{(4)}| + 3 - 2 \geq \lceil (4n+3)/7 \rceil$, a contradiction.

So, $|N(x_3) \cap N(t)| > 2$. 
There exist $b_1 \in N(x_3) \cap N(t)$ and subgraphs $G_1',G_2',G_3'$ such that
$G_2' = G_2$,
$G_3'$ is the maximal subgraph of $G$ contained in the closed region of the plane bounded by the cycle $z_2x_3b_1tz_2$ and containing $N(x_3) \cap N(t) - \{z_2\}$,
and $G_1'$ is obtained from $G$ by removing $G_2'-\{s_1,s_2,a_1\}$ and $G_3'-\{x_3,b_1,t\}$.
 Let $B_5 = \{b_1\}$ or $B_6 = \emptyset$.
 For $i=5,6$, 
 let $F_1^{(i)} = A(G_1' - \{z,z_1,z_2,s_2,x_3,t,x_2\} - B_i + s_1x)$, $F_2^{(i)} = A(G_2' - \{s_2\})$, and $F_3^{(i)} = A(G_3' - \{x_3,t\} - B_i)$. 
 Then $|F_1^{(i)}| \geq \lceil (4(|G_1'|-7-|B_i|)+3)/7 \rceil$, $|F_2^{(i)}| \geq \lceil (4(|G_2'|-1)+3)/7 \rceil = \lceil (4(|G_2|-1)+3)/7 \rceil = (4(|G_2|-1)+3)/7 + 6/7 $,
  and $|F_3^{(i)}| \geq \lceil (4(|G_3'|-2-|B_i|)+3)/7 \rceil$. 
 Now $G[F_1^{(i)} \cup F_2^{(i)} \cup F_3^{(i)} + \{z,z_1,z_2\} - (\{s_1,a_1\} \cap (F_1^{(i)} \triangle F_2^{(i)})) - (\{b_1\} \cap (F_1^{(i)} \triangle F_3^{(i)}))]$ 
 is an induced forest in $G$, showing $a(G) \geq |F_1^{(i)}| + |F_2^{(i)}| + |F_3^{(i)}| + 3 - 2 - (1-|B_i|) \geq \lceil (4n+3)/7 \rceil$. 
By Lemma \ref{ineq2}(2)
(with $a = |G_1'|-7, a_1 = |G_3'|-2, c = (4(|G_2|-1)+3)/7 + 6/7 + 1$)
, $(4(|G_1'|-7)+3,4(|G_2'|-2)+3,4(|G_3'|-2)+3) \equiv (0,4,4), (4,4,0) \mod 7$.

  If $(4(|G_1'|-7)+3,4(|G_2'|-2)+3,4(|G_3'|-2)+3) \equiv (4,4,0) \mod 7$, 
  let $F_1^{(7)} = A((G_1' - \{z,z_1,z_2,s_2,x_2\})/x_3t + s_1x)$ with $u$ as the identification of $\{x_3,t\}$, $F_2^{(7)} = A(G_2' - \{s_2\})$, and $F_3^{(7)} = A(G_3')$. 
  Then $|F_1^{(7)}| \geq \lceil (4(|G_1'|-6)+3)/7 \rceil$, $|F_2^{(7)}| \geq \lceil (4(|G_2'|-1)+3)/7 \rceil$, and $|F_3{(7)}| \geq \lceil (4|G_3'|+3)/7 \rceil$. 
  Now $G[F_1^{(7)} \cup F_2^{(7)} \cup F_3^{(7)} + \{z,z_1,z_2\}  - (\{s_1,a_1\} \cap (F_1^{(7)} \triangle F_2^{(7)})) - (\{t,b_1,x_3\} \cap (F_1^{(7)} \triangle F_3^{(7)}))]$ (if $u \not\in F_1^{(7)}$)
  or $G[(F_1^{(7)} - u) \cup F_2^{(7)} \cup F_3^{(7)} + \{z,z_1,x_3,t\}  - (\{s_1,a_1\} \cap (F_1^{(7)} \triangle F_2^{(7)})) - (\{t,b_1,x_3\} \cap ((F_1^{(7)} + \{x_3,t\}) \triangle F_3^{(7)}))]$ (if $u \in F_1^{(7)}$)
   is an induced forest in $G$, showing $a(G) \geq |F_1^{(7)}| + |F_2^{(7)}| + |F_3^{(7)}| + 3 - 5 \geq \lceil (4n+3)/7 \rceil$, a contradiction.

 If $(4(|G_1'|-7)+3,4(|G_2'|-2)+3,4(|G_3'|-2)+3) \equiv (0,4,4) \mod 7$,
 let $F_1^{(8)} = A((G_1' - \{z,z_2,s_2,x_3\})/xz_1)$ with $u$ as the identification of $\{x,z_1\}$, $F_2^{(8)} = A(G_2' - \{s_2\})$, and $F_3^{(8)} = A(G_3' - \{x_3\})$. 
 Then $|F_1^{(8)}| \geq \lceil (4(|G_1'|-5)+3)/7 \rceil$, $|F_2^{(8)}| \geq \lceil (4(|G_2'|-1)+3)/7 \rceil$, and $|F_3^{(8)}| \geq \lceil (4(|G_3'|-1)+3)/7 \rceil$. 
Let $F = F_1^{(8)} \cup F_2^{(8)} \cup F_3^{(8)}  + \{z,z_2\}  - (\{s_1,a_1\} \cap (F_1^{(8)} \triangle F_2^{(8)})) - (\{t,b_1\} \cap (F_1^{(8)} \triangle F_3^{(8)}))$. 
 Now $G[F]$ (if $u \not\in F_1^{(8)}$) or $G[F - \{u,z\}+\{x,z_1\}]$ (if $u \in F_1^{(8)}$)
 is an induced forest in $G$, showing $a(G) \geq |F_1^{(8)}| + |F_2^{(8)}| + |F_3^{(8)}| + 2 - 4 \geq \lceil (4n+3)/7 \rceil$, a contradiction.


 \medskip

Case 3: $|N(w_1) \cap N(x_2)| > 2$. 
There exist $c_1 \in N(w_1) \cap N(x_2)$ and a separation $(G_1,G_2)$ such that $V(G_1 \cap G_2) = \{w_1,x_2,c_1\}$, $\{x,y,w,z,x_1,x_3,z_1,z_2,s_1,s_2\} \subseteq V(G_1)$, and $N(w_1) \cap N(x_2) - \{w\} \subseteq V(G_2)$. 
By the fourth claim, $s_1z_2 \not\in E(G)$. 
Let $C_1 = \{c_1\}$ and $C_2 = \emptyset$. 
For $i=1,2$, let $F_1^{(i)} = A((G_1 - \{w,x,w_1,x_2,z,z_1,s_2\} - C_i)/x_1y + s_1z_2)$ with $u$ as the identification of $\{x_1,y\}$, and $F_2^{(i)} = A(G_2 - \{w_1,x_2\} - C_i)$. 
Then $|F_1^{(i)}| \geq \lceil (4(|G_1|-8-|C_i|)+3)/7 \rceil$, and $|F_2^{(i)}| \geq \lceil (4(|G_2|-2-|C_i|)+3)/7 \rceil$. 
Note $u \in F_1^{(i)}$ by Lemma \ref{basic_lemma} since $|N(u)| = 3$. 
Now $G[(F_1^{(i)} - u) \cup F_2^{(i)} + \{w,x_1,y,z,z_1\}  - (\{c_1\} \cap (F_1^{(i)} \triangle F_2^{(i)} ))]$
is an induced forest in $G$, showing $a(G) \geq |F_1^{(i)}| + |F_2^{(i)}| + 4 - (1-|C_i|)$. 
By Lemma \ref{ineq2}(2)
(with $a = |G_1|-8, a_1 = |G_2|-2, c=4$),
 $(4(|G_1|-8)+3,4(|G_2|-2)+3) \equiv (4,0),(0,4) \mod 7$.

\medskip

{\it Subcase 3.1:} $(4(|G_1|-8)+3,4(|G_2|-2)+3) \equiv (4,0) \mod 7$.

 If $|N(s_1) \cap N(s_2)| \leq 2$, then let $F_1^{(3)} = A((G_1 - \{w,x,z,z_1\} )/\{x_1y,w_1x_2,s_1s_2\}+ z_2u_2)$ with $u_1$ (respectively, $u_2,u_3$) as the identification of  $\{x_1,y\}$ (respectively, $\{w_1,x_2\},\{s_1,s_2\}$), and $F_2^{(3)} = A(G_2)$. 
 Then $|F_1^{(3)}| \geq \lceil (4(|G_1|-7)+3)/7 \rceil$, and $|F_2^{(3)}| \geq \lceil (4|G_2|+3)/7 \rceil$. 
 Note $u_1 \in F_1^{(10)}$ by Lemma \ref{basic_lemma} since $|N(u_1)| = 3$. 
 Let $F^{(3)} := \overline{F_1}^{(3)} \cup F_2^{(3)}  - (\{w_1,x_2,c_1\} \cap (\overline{F_1}^{(3)} \triangle F_2^{(3)}))$, 
where $\overline{F_1}^{(3)} = {F_1}^{(3)} + \{x_1,y,z,z_1,w\} - u_1$ if $u_2,u_3 \not\in F_1^{(3)}$, and otherwise, $\overline{F_1}^{(3)}$ obtained from ${F_1}^{(3)} + \{x_1,y,z,z_1,w\} - u_1$ by deleting $\{u_2,w\}$ (respectively, $\{u_3,z_1\}$) and adding $\{w_1,x_2\}$ (respectively, $\{s_1,s_2\}$) when $u_2 \in {F_1}^{(3)}$ (respectively, $u_3 \in {F_1}^{(3)}$).
Therefore, $G[F^{(3)}]$ is an induced forest in $G$, showing $a(G) \geq |F_1^{(3)}| + |F_2^{(3)}| + 4 - 3 \geq \lceil (4n+3)/7 \rceil$, a contradiction.

So $|N(s_1) \cap N(s_2)| >  2$. 
There exist $a_1 \in N(s_1) \cap N(s_2)$ and subgraphs $G_1',G_2',G_3'$ such that
$G_2' = G_2$,
$G_3'$ is the maximal subgraph of $G$ contained in the closed region of the plane bounded by the cycle $z_1s_1a_1s_2z_1$ and containing $N(s_1) \cap N(s_2) - \{z_1\}$,
and $G_1'$ is obtained from $G$ by removing $G_2'-\{w_1,x_2,w\}$ and $G_3'-\{s_1,a_1,s_2\}$.
Let $A_4 = \{a_1\}$ and $A_5 = \emptyset$. 
For $i=4,5$, let $F_1^{(i)} = A((G_1' - \{w,x,z,z_1,s_1,s_2\} - A_i)/\{x_1y,w_1x_2\} + u_2z_2)$ with $u_1$ (respectively, $u_2$) as the identification of $\{x_1,y\}$,  $F_2^{(i)} = A(G_2')$, and $F_3^{(i)} = A(G_3' - \{s_1,s_2\} - A_i)$. 
Then $|F_1^{(i)}| \geq \lceil (4(|G_1'|-8-|A_i|)+3)/7 \rceil$, 
$|F_2^{(i)}| \geq \lceil (4|G_2'|+3)/7 \rceil = \lceil (4|G_2|+3)/7 \rceil =  (4|G_2|+3)/7 + 6/7 $, 
and $|F_3^{(i)}| \geq \lceil (4(|G_3'|-2-|A_i|)+3)/7 \rceil$. 
Note $u_1 \in F_1^{(i)}$ by Lemma \ref{basic_lemma} since $|N(u_1)| = 3$. 
Now $G[(F_1^{(i)} - u_1) \cup F_2^{(i)}  \cup F_3^{(i)} + \{x_1,y,z,z_1,w\} - (\{w_1,x_2,c_1\} \cap (F_1^{(i)} \triangle F_2^{(i)})) - (\{a_1\} \cap (F_1^{(i)} \triangle F_3^{(i)}))]$ (if $u_2 \not\in F_1^{(i)}$ )
or $G[(F_1^{(i)} - \{u_1,u_2\}) \cup F_2^{(i)}  \cup F_3^{(i)} + \{x_1,y,z,z_1,w_1,x_2\} - (\{w_1,x_2,c_1\} \cap ( (F_1^{(i)} \cup \{w_1,x_2\}) \triangle F_2^{(i)})) - (\{a_1\} \cap (F_1^{(i)} \triangle F_3^{(i)}))]$ (if $u_2 \in F_1^{(i)}$ )
is an induced forest in $G$, showing $a(G) \geq |F_1^{(i)}| + |F_2^{(i)}| + 4 - 3 - (1-|A_i|)$. 
By Lemma \ref{ineq2}(2)
(with $a = |G_1'|-8, a_1 = |G_3'|-2, c=(4|G_2|+3)/7 + 6/7 + 1 $),
 $(4(|G_1'|-8)+3,4(|G_2'|-2)+3,4(|G_3'|-2)+3) \equiv (4,0,0),(0,0,4) \mod 7$.

If $(4(|G_1'|-8)+3,4(|G_2'|-2)+3,4(|G_3'|-2)+3) \equiv (4,0,0) \mod 7$,
let $F_1^{(6)} = A((G_1' - \{w,x,z,z_1\} )/\{x_1y,w_1x_2,s_1s_2\} + z_2u_2)$ with $u_1$ (respectively, $u_2,u_3$) as the identification of  $\{x_1,y\}$ (respectively, $\{w_1,x_2\},\{s_1,s_2\}$), $F_2^{(6)} = A(G_2')$, and $F_3^{(6)} = A(G_3')$. 
Then $|F_1^{(6)}| \geq \lceil (4(|G_1'|-7)+3)/7 \rceil$, $|F_2^{(6)}| \geq \lceil (4|G_2'|+3)/7 \rceil$, and $|F_3^{(6)}| \geq \lceil (4|G_3'|+3)/7 \rceil$. 
Note $u_1 \in F_1^{(6)}$ by Lemma \ref{basic_lemma} since $|N(u_1)| = 3$. 
Let $F^{(6)} := \overline{F_1}^{(6)} \cup F_2^{(6)}  \cup F_3^{(6)} - (\{w_1,x_2,c_1\} \cap (\overline{F_1}^{(6)} \triangle F_2^{(6)}))  - (\{s_1,s_2,a_1\} \cap (\overline{F_1}^{(6)} \triangle F_3^{(6)}))$,
where $\overline{F_1}^{(6)} = {F_1}^{(6)} + \{x_1,y,z,z_1,w\} - u_1$ if $u_2,u_3 \not\in F_1^{(6)}$, and otherwise, $\overline{F_1}^{(6)}$ obtained from ${F_1}^{(6)} + \{x_1,y,z,z_1,w\} - u_1$ by deleting $\{u_2,w\}$ (respectively, $\{u_3,z_1\}$) and adding $\{w_1,x_2\}$ (respectively, $\{s_1,s_2\}$) when $u_2 \in {F_1}^{(6)}$ (respectively, $u_3 \in {F_1}^{(6)}$).
Therefore, $F^{(6)}$ is an induced forest in $G$, showing $|F_1^{(6)} | + |F_2^{(6)} | + |F_3^{(6)} |  + 4 - 6 \geq \lceil (4n+3)/7 \rceil$, a contradiction. 

If $(4(|G_1'|-8)+3,4(|G_2'|-2)+3,4(|G_3'|-2)+3) \equiv (0,0,4) \mod 7$,
 let $F_1^{(7)} = A((G_1' - \{w_1,x_2,w,x,z,z_1,c_1,s_2\} + s_1z_2)/x_1y)$ with $u_1$ as the identification of  $\{x_1,y\}$, $F_2^{(7)} = A(G_2' - \{w_1,x_2,c_1\})$, and $F_3^{(7)} = A(G_3' - s_2)$. 
 Then $|F_1^{(7)}| \geq \lceil (4(|G_1'|-9)+3)/7 \rceil$, $|F_2^{(7)}| \geq \lceil (4(|G_2'|-3)+3)/7 \rceil$, and $|F_3^{(7)}| \geq \lceil (4(|G_3'|-1)+3)/7 \rceil$.
 Note $u_1 \in F_1^{(7)}$ by Lemma \ref{basic_lemma} since $|N(u_1)| = 3$. 
 Now $G[(F_1^{(7)} - u_1) \cup F_2^{(7)} \cup F_3^{(7)} + \{x_1,y,z,z_1,w\} - (\{s_1,a_1\} \cap (F_1^{(7)} \triangle F_2^{(7)}) )]$
 is an induced forest in $G$, showing $a(G) \geq |F_1^{(7)} | + |F_2^{(7)} | + |F_3^{(7)} |  + 4 - 2 \geq \lceil (4n+3)/7 \rceil$, a contradiction.

 \medskip

{\it Subcase 3.2:} $(4(|G_1|-8)+3,4(|G_2|-2)+3) \equiv (0,4) \mod 7$. 

If $|N(x_3) \cap N(t)| \leq 2$,  
let $F_1^{(8)} = A((G_1 - \{z,z_1,z_2,x_2,s_2\} + xs_1)/x_3t)$ with $u$ as the identification of $\{x_3,t\}$, and $F_2^{(8)} = A((G_2 -x_2)$. 
Then $|F_1^{(8)}| \geq \lceil (4(|G_1|-6)+3)/7 \rceil$, and $|F_2^{(8)}| \geq \lceil (4(|G_2|-1)+3)/7 \rceil$.
Let $F = F_1^{(8)} \cup F_2^{(8)} + \{z,z_1,z_2\} - (\{w_1,c_1\} \cap (F_1^{(8)} \triangle F_2^{(8)}))$.
Now $G[F]$ (if $u \not\in F_1^{(8)}$) or $G[F-\{u,z_2\} +\{x_3,t\}]$ (if $u \in F_1^{(8)}$)
is an induced forest in $G$, showing $a(G) \geq |F_1^{(8)} | + |F_2^{(8)} | + 3 - 2 \geq \lceil (4n+3)/7 \rceil$, a contradiction.

So $|N(x_3) \cap N(t)| > 2$. 
There exist $b_1 \in N(x_3) \cap N(t)$ and subgraphs $G_1',G_2',G_3'$ such that
$G_2' = G_2$,
$G_3'$ is the maximal subgraph of $G$ contained in the closed region of the plane bounded by the cycle $z_2x_3b_1tz_2$ and containing $N(x_3) \cap N(t) - \{z_2\}$,
and $G_1'$ is obtained from $G$ by removing $G_2'-\{w_1,x_2,c_1\}$ and $G_3'-\{x_3,b_1,t\}$.
Let $B_9 = \emptyset$ and $B_{10} = \{b_1\}$. 
For $i=9,10$, 
let $F_1^{(i)} = A(G_1'- \{z,z_1,z_2,x_2,s_2,x_3,t\} - B_i + xs_1)$, $F_2^{(i)} = A(G_2' - \{x_2\})$, and $F_3^{(i)} = A(G_3' - \{x_3,t\} - B_i)$.
Then $|F_1^{(i)}| \geq \lceil (4(|G_1'|-7-|B_i|)+3)/7 \rceil$, $|F_2^{(i)}| \geq \lceil (4(|G_2'|-1)+3)/7 \rceil = \lceil (4(|G_2|-1)+3)/7 \rceil = (4(|G_2|-1)+3)/7 + 6/7$, and $|F_3^{(i)}| \geq \lceil (4(|G_3'|-2-|B_i|)+3)/7 \rceil$.
Now $G[F_1^{(i)} \cup F_2^{(i)} \cup F_3^{(i)} + \{z,z_1,z_2\} -  (\{c_1,w_1\} \cap (F_1^{(i)} \triangle F_2^{(i)})) - (\{b_1\} \cap (F_1^{(i)} \triangle F_3^{(i)}))]$ 
is an induced forest in $G$, showing $a(G) \geq |F_1^{(i)} | + |F_2^{(i)} | + 3 - 2 - (1-|B_i|)$.
By Lemma \ref{ineq2}(2)
(with $a = |G_1'|-7, a_1 = |G_3'|-2 , c = (4(|G_2|-1)+3)/7 + 6/7 + 1$ ),
 $(4(|G_1'|-7)+3,4(|G_2'|-2)+3,4(|G_3'|-2)+3) \equiv (4,4,0),(0,4,4) \mod 7$.

If $(4(|G_1'|-7)+3,4(|G_2'|-2)+3,4(|G_3'|-2)+3) \equiv (4,4,0) \mod 7$, 
let $F_1^{(11)} = A((G_1' - \{z,z_1,z_2,x_2,s_2\} + xs_1)/x_3t)$ with $u$ as the identification of $\{x_3,t\}$, $F_2^{(11)} = A(G_2' -x_2)$, and $F_3^{(11)} = A(G_3')$. 
Then $|F_1^{(11)}| \geq \lceil (4(|G_1'|-6)+3)/7 \rceil$, $|F_2^{(11)}| \geq \lceil (4(|G_2'|-1)+3)/7 \rceil$, and $|F_3^{(11)}| \geq \lceil (4|G_3'|+3)/7 \rceil$.
Now $G[F_1^{(11)} \cup F_2^{(11)} \cup F_3^{(11)} + \{z,z_1,z_2\} - (\{w_1,c_1\} \cap (F_1^{(11)} \triangle F_2^{(11)})) - (\{x_3,t,b_1\} \cap (F_1^{(11)} \triangle F_3^{(11)}))]$ (if $u \not\in F_1^{(11)}$) 
or $G[(F_1^{(11)} - u) \cup F_2^{(11)} \cup F_3^{(11)} + \{z,z_1,x_3,t\} - (\{w_1,c_1\} \cap (F_1^{(11)} \triangle F_2^{(11)})) - (\{x_3,t,b_1\} \cap ((F_1^{(11)} + \{x_3,t\}) \triangle F_3^{(11)}))]$ (if $u \in F_1^{(11)}$)
is an induced forest in $G$, showing $a(G) \geq |F_1^{(11)} | + |F_2^{(11)} | + |F_3^{(11)} | + 3 - 5 \geq \lceil (4n+3)/7 \rceil$, a contradiction.

So $(4(|G_1'|-7)+3,4(|G_2'|-2)+3,4(|G_3'|-2)+3) \equiv (0,4,4) \mod 7$. 
If $|N(x_2) \cap N(s_1)| \leq 2$,
let $F_1^{(12)} = A((G_1' \cup G_2' - \{z,z_1,z_2,s_2,x_3\})/x_2s_1 + xt)$ with $u$ as the identification of $\{x_2,s_1\}$, and $F_2^{(12)} = A(G_3' - \{x_3\})$.
Then $|F_1^{(12)}| \geq \lceil (4((n+3-|G_3'|)-6)+3)/7 \rceil$, and $|F_2^{(12)}| \geq \lceil (4(|G_3'|-1)+3)/7 \rceil$.
Let $F = F_1^{(12)} \cup F_2^{(12)} + \{z,z_1,z_2\} - (\{b_1,t\} \cap (F_1^{(12)} \triangle F_2^{(12)}))$.
Now $G[F]$ (if $u \not\in F_1^{(12)}$) or $G[F -\{u,z_1\} +\{x_2,s_1\}]$ (if $u \in F_1^{(12)}$)
is an induced forest in $G$, showing $a(G) \geq |F_1^{(12)} | + |F_2^{(12)} | + 3 - 2 \geq \lceil (4n+3)/7 \rceil$, a contradiction.
So $|N(x_2) \cap N(s_1)| > 2$. 
There exist $e_1 \in N(x_2) \cap N(s_1)$ and subgraphs $G_1'',G_2'',G_3'',G_4''$ such that 
$G_2'' = G_2'$,
$G_3'' = G_3'$,
$G_4''$ is the maximal subgraph of $G$ contained in the closed region of the plane bounded by the cycle $z_1x_2e_1s_1z_1$ and containing $N(s_1) \cap N(x_2) - \{z_1\}$,
and $G_1''$ is obtained from $G$ by removing $G_2''-\{w_1,x_2,c_1\}$, $G_3''-\{x_3,b_1,t\}$ and $G_4''-\{s_1,x_2,e_1\}$.
Let $E_{13} = \emptyset$ and $E_{14} = \{e_1\}$.
For $i=13,14$, let $F_1^{(i)} = A(G_1''- \{z,z_1,z_2,x_2,s_1,s_2,x_3\} - E_i + xt)$, $F_2^{(i)} = A(G_2'' - \{x_2\})$, $F_3^{(i)} = A(G_3'' - \{x_3\})$, and $F_4^{(i)} = A(G_4'' - \{s_1,x_2\} - E_i)$.
Then $|F_1^{(i)}| \geq \lceil (4(|G_1''|-7-|E_i|)+3)/7 \rceil$, 
$|F_2^{(i)}| \geq \lceil (4(|G_2''|-1)+3)/7 \rceil = \lceil (4(|G_2'|-1)+3)/7 \rceil = (4(|G_2'|-1)+3)/7 + 6/7$, 
$|F_3^{(i)}| \geq \lceil (4(|G_3''|-1)+3)/7 \rceil = \lceil (4(|G_3'|-1)+3)/7 \rceil = (4(|G_3'|-1)+3)/7 + 6/7$,  
and $|F_4^{(i)}| \geq \lceil (4(|G_4''|-2-|E_i|)+3)/7 \rceil$.
Now $G[F_1^{(i)} \cup F_2^{(i)} \cup F_3^{(i)} + \{z,z_1,z_2\} - (\{w_1,c_1\} \cap (F_1^{(i)} \triangle F_2^{(i)})) - (\{b_1,t\} \cap (F_1^{(i)} \triangle F_3^{(i)})) - (\{e_1\} \cap (F_1^{(i)} \triangle F_4^{(i)}))]$ 
is an induced forest in $G$, showing $a(G) \geq |F_1^{(i)} | + |F_2^{(i)} | + |F_3^{(i)} | + |F_4^{(i)} | + 3 - 2 - (1-|E_i|)$.
By Lemma \ref{ineq2}(1) 
(with $k=1, a=|G_1''|-7, a_1 = |G_4''|-2, L = \{1,2\}, b_1 = |G_2''|-1, b_2 = |G_3''|-1, c=1$), 
$a(G) \geq \lceil (4n+3)/7 \rceil $, a contradiction.
\qed

\begin{lem}
\label{5-1-B}
The following configuration is impossible in $G$: 
$x$ is a $5$-vertex of type 5-1-B in $G$ with neighbors $y_1,y_3,y_4,y_2,y_5$ in cyclic order around $x$. $xy_1y_3'y_3x$, $xy_4y_4'y_2x$, $xy_1z_1y_5x$, $xy_2z_2y_5x$ are facial cycles. $\{y_1,y_2\} \subseteq V_4, \{z_1,z_2,y_3',y_4\} \subseteq V_3$, $N(z_1) = \{y_1,y_5,w_1\}$ and $N(z_2) = \{y_2,y_5,w_2\}$. 
\end{lem}



\pf 
Let $N(y_1) = \{x,z_1,y_3',y_1'\}$ and $N(y_2) = \{x,z_2,y_4',y_2'\}$. 

First, we claim that $w_1x \not\in E(G)$. Otherwise $w_1x \in E(G)$. 
Since $G$ is simple, $w_1 \not\in \{y_1,y_5\}$. 
If $w_1 = y_2$ and $y_1y_4' \not\in E(G)$, then let $F' = A(G - \{z_1,z_2,x,y_5,y_2\} + y_1y_4')$. 
Then $|F'| \geq \lceil (4(n-5)+3)/7 \rceil$. 
Now $G[F' + \{z_1,z_2,y_2\}]$  is an induced forest in $G$, showing $a(G) \geq |F'| + 3 \geq \lceil (4n+3)/7 \rceil$, a contradiction. 
If $w_1 = y_2$ and $y_1y_4' \in E(G)$, 
 let $F' = A(G - \{z_1,z_2,x,y_5,y_2,y_1,y_4'\})$. 
Then $|F'| \geq \lceil (4(n-7)+3)/7 \rceil$. 
Now $G[F' + \{z_1,z_2,y_2,y_1\}]$ is an induced forest in $G$, showing $a(G) \geq |F'| + 4 \geq \lceil (4n+3)/7 \rceil$, a contradiction.  
If $w_1 = y_3$, then since $G$ is plane, there exists a separation $(G_1,G_2)$ such that $V(G_1 \cap G_2) = \{z_1,y_1,y_3\}$, $y_1' \in V(G_1)$, and $\{x,y_5\} \subseteq V(G_2)$. 
Let $F_1^{(1)} = A(G_1 - \{z_1,y_3\})$, and $F_2^{(1)} = A(G_2 - \{z_1,y_1,y_3,y_5\})$. 
Then $|F_1^{(1)}| \geq \lceil (4(|G_1|-2)+3)/7 \rceil$, and $|F_2^{(1)}| \geq \lceil (4(|G_2|-4)+3)/7 \rceil$. 
Now $G[F_1^{(1)} \cup F_2^{(1)} + \{z_1\}]$ is an induced forest in $G$, showing $a(G) \geq |F_1^{(1)}| + |F_2^{(1)}| + 1$.
By Lemma \ref{ineq2}(7) 
(with $k=2,a_1 = |G_1|-2, a_2 = |G_2| - 4, c=1$), 
$(4(|G_1|-2)+3,4(|G_2| - 4)+3) \equiv (0,6), (6,0), (0,0) \mod 7 $. 
Let $F_i^{(2)} = A(G_i - \{z_1,y_1,y_3\})$ for $i=1,2$. 
Then $|F_i^{(2)}| \geq \lceil (4(|G_i|-3)+3)/7 \rceil$. 
Now $G[F_1^{(2)} \cup F_2^{(2)} + \{z_1\}]$  is an induced forest in $G$, showing $a(G) \geq |F_1^{(2)}| + |F_2^{(2)}| + 1 \geq \lceil (4n+3)/7 \rceil$, a contradiction. 
If $w_1 = y_4$, then since $G$ is plane, there exists a separation $(G_1,G_2)$ such that $V(G_1 \cap G_2) = \{z_1,x,y_4\}$, $y_1 \in V(G_1)$, and $\{y_2,y_5\} \subseteq V(G_2)$. 
Let $F_1^{(3)} = A(G_1 - \{z_1,x,y_4,y_1\})$, and $F_2^{(3)} = A(G_2 - \{z_1,y_4\})$. 
Then $|F_1^{(3)}| \geq \lceil (4(|G_1|-4)+3)/7 \rceil$, and $|F_2^{(3)}| \geq \lceil (4(|G_2|-2)+3)/7 \rceil$. 
Now $G[F_1^{(3)} \cup F_2^{(3)} + \{z_1\}]$ is an induced forest in $G$, showing $a(G) \geq |F_1^{(3)}| + |F_2^{(3)}| + 1$. 
By Lemma \ref{ineq2}(7) 
(with $k=2,a_1 = |G_1|-4, a_2 = |G_2| - 2, c=1$), 
$(4(|G_1|-4)+3,4(|G_2| - 2)+3) \equiv (0,6), (6,0), (0,0) \mod 7 $. 
Let $F_1^{(4)}= A(G_1 - \{z_1,x,y_4\})$, and $F_2^{(4)} = A(G_2 - \{z_1,x,y_4\})$. 
Then $|F_1^{(4)}| \geq \lceil (4(|G_1|-3)+3)/7 \rceil$, and $|F_2^{(4)}| \geq \lceil (4(|G_2|-3)+3)/7 \rceil$. 
Now $G[F_1^{(4)} \cup F_2^{(4)} + \{z_1\}]$  is an induced forest in $G$, showing $a(G) \geq |F_1^{(2)}| + |F_2^{(2)}| + 1 \geq \lceil (4n+3)/7 \rceil$, a contradiction.
Similarly, $w_2x \not\in E(G)$.


Secondly, we claim that $|N(y_1') \cap N(y_3')| \leq 2$. Otherwise, 
there exist $a_1 \in N(y_1') \cap N(y_3')$ 
and a separation $(G_1,G_2)$ such that $V(G_1 \cap G_2) = \{y_1',y_3',a_1\}$, $\{x,y_1,y_2,y_3,y_4,y_5\} \subseteq V(G_1)$, and $N(y_1') \cap N(y_3') - \{y_1\} \subseteq V(G_2)$. 
Let $F_1^{(1)} = A(G_1 - \{y_1',y_3',a_1,y_1,z_1,y_5\} + w_1x)$, and $F_2^{(1)} = A(G_2 - \{y_1',y_3',a_1\})$. 
Then $|F_1^{(1)}| \geq \lceil (4(|G_1|-6)+3)/7 \rceil$, and $|F_2^{(1)}| \geq \lceil (4(|G_2|-3)+3)/7 \rceil$. 
Now $G[F_1^{(1)} \cup F_2^{(1)} + \{z_1,y_1,y_3\}]$  is an induced forest in $G$, showing $a(G) \geq |F_1^{(1)}| + |F_2^{(1)}| + 3 \geq \lceil (4n+3)/7 \rceil$, a contradiction.

Now we prove the lemma. 
If $|N(y_2') \cap N(y_4')| \leq 2$, 
let $F' = A((G - \{z_1,z_2,y_1,y_2,y_5\})/$ $\{y_1'y_3',y_2'y_4'\} + \{w_1x,w_2x\})$ with $u_1$ (respectively, $u_2$) as the identification of  $\{y_1',y_3'\}$ (respectively, $\{y_2',y_4'\}$). 
Then $|F'| \geq \lceil (4(n-7)+3)/7 \rceil$. 
Let $F = F' + \{z_1,z_2,y_1,y_2\}$ if $u_1,u_2 \not\in F'$, and otherwise $F$  obtained from $F' + \{z_1,z_2,y_1,y_2\}$ by deleting $\{u_1,y_1\}$ (respectively, $\{u_2,y_2\}$) and adding $\{y_1',y_3'\}$ (respectively, $\{y_2',y_4'\}$) when $u_1 \in F'$ (respectively, $u_2 \in F'$).
Therefore, $G[F]$ is an induced forest in $G$, showing $a(G) \geq |F'| + 4 \geq \lceil (4n+3)/7 \rceil$, a contradiction. 

So, $|N(y_2') \cap N(y_4')| > 2$. 
Let $B_1 = \{b_1\}$ and $B_2 = \emptyset$. 
There exist $b_1 \in N(y_2') \cap N(y_4')$ and a separation $(G_1,G_2)$ such that $V(G_1 \cap G_2) = \{y_2',y_4',b_1\}$, $\{x,y_1,y_2,y_3,y_4,y_5\} \subseteq V(G_1)$, and $N(y_2') \cap N(y_4') - \{y_2\} \subseteq V(G_2)$. 
For $i=1,2$, let $F_1^{(i)} = A((G_1 - \{z_1,z_2,y_1,y_2,y_5,y_2',y_4'\} - B_i)/y_1'y_3' + \{w_1x,w_2x\})$ with $u$ as the identification of $\{y_1',y_3'\}$, and $F_2^{(i)} = A(G_2 - \{y_2',y_4'\} - B_i)$. 
Then $|F_1^{(i)}| \geq \lceil (4(|G_1|-8-|B_i|)+3)/7 \rceil$, and $|F_2^{(i)}| \geq \lceil (4(|G_2|-2-|B_i|)+3)/7 \rceil$. 
Let $F = F_1^{(i)} \cup F_2^{(i)} + \{z_1,z_2,y_1,y_2\} - (\{b_1\} \cap (F_1^{(i)} \triangle F_2^{(i)} ) )$.
Now $G[F]$ ($u \not\in F_1^{(i)}$) or $G[F -\{u,y_1\} +\{y_1',y_3'\}]$ ($u \in F_1^{(i)}$)
is an induced forest in $G$, showing $a(G) \geq |F_1^{(i)}| + |F_2^{(i)}| + 4 - (1-|B_i|)$.
By Lemma \ref{ineq2}(2)
(with $a = |G_1|-8, a_1 = |G_2|-2, c=4$),
 $(4(|G_1|-8)+3,4(|G_2|-2)+3) \equiv (4,0), (0,4) \mod 7$.


If $(4(|G_1|-8)+3,4(|G_2|-2)+3) \equiv (4,0) \mod 7$, 
let $F_1^{(3)} = A((G - \{z_1,z_2,y_1,y_2,y_5\})/$ $\{y_1'y_3',y_2'y_4'\} + \{w_1x,w_2x\})$ with $u_1$ (respectively, $u_2$) as the identification of  $\{y_1',y_3'\}$ (respectively, $\{y_2',y_4'\}$), and $F_2^{(3)} = A(G_2)$. 
Then $|F_1^{(3)}| \geq \lceil (4(|G_1|-7)+3)/7 \rceil$, and $|F_2^{(3)}| \geq \lceil (4|G_2|+3)/7 \rceil$. 
Let $F^{(3)} = \overline{F_1}^{(3)} \cup F_2^{(3)} - (\{b_1,y_4',y_2'\} \cap (\overline{F_1}^{(3)} \triangle F_2^{(3)}))$,
where
$\overline{F_1}^{(3)} = {F_1}^{(3)} + \{z_1,z_2,y_1,y_2\}$ if $u_1,u_2 \not\in F'$, and let $\overline{F_1}^{(3)}$ obtained from ${F_1}^{(3)} + \{z_1,z_2,y_1,y_2\}$ by deleting $\{u_1,y_1\}$ (respectively, $\{u_2,y_2\}$) and adding $\{y_1',y_3'\}$ (respectively, $\{y_2',y_4'\}$) when $u_1 \in F'$ (respectively, $u_2 \in F'$).
Therefore, $G[F^{(3)}]$ is an induced forest in $G$, showing $a(G) \geq |{F_1}^{(3)}| + |{F_2}^{(3)}| + 4 - 3 \geq \lceil (4n+3)/7 \rceil$, a contradiction.

So $(4(|G_1|-8)+3,4(|G_2|-2)+3) \equiv (0,4) \mod 7$. 
If $y_5y_2' \not\in E(G)$, 
let $F_1^{(4)} = A(G_1 - \{x,y_4,y_4',y_2,z_2,w_2\}+ y_5y_2')$, and $F_2^{(4)} = A(G_2-\{y_4'\})$. 
Then $|F_1^{(4)}| \geq \lceil (4(|G_1|-6)+3)/7 \rceil$, and $|F_2^{(4)}| \geq \lceil (4(|G_2|-1)+3)/7 \rceil$. 
Now $G[F_1^{(4)} \cup F_2^{(4)} + \{y_4,y_2,z_2\} - ( \{b_1,y_2'\}\cap (F_1^{(4)} \triangle F_2^{(4)}))]$  is an induced forest in $G$, showing $a(G) \geq |F_1^{(4)}| + |F_2^{(4)}| + 3 - 2 \geq \lceil (4n+3)/7 \rceil$, a contradiction.

So, $y_5y_2' \in E(G)$, then 
there exists a separation $(G_1',G_2')$ such that $V(G_1' \cap G_2') = \{y_2',y_5,z_2\}$, $\{x,z_1,y_2,y_3,y_4\} \subseteq V(G_1')$, and $w_2 \in V(G_2')$. 
Let $F_1^{(5)} = A((G_1' - \{y_2',y_5,z_2,$ $z_1,y_1\})/y_1'y_3'+ w_1x)$ with $u$ as the identification of $\{y_1',y_3'\}$ and $F_2^{(5)} = A(G_2'-\{y_2',y_5,z_2\})$. 
Then $|F_1^{(5)}| \geq \lceil (4(|G_1|-6)+3)/7 \rceil$, and $|F_2^{(5)}| \geq \lceil (4(|G_2|-3)+3)/7 \rceil$. 
Let $F = F_1^{(5)} \cup F_2^{(5)} + \{z_1,z_2,y_1\}$.
Now $G[F]$ (if $u \not\in F_1^{(5)}$) or $G[F - \{u,y_1\} + \{y_1',y_3'\}]$ (if $u \in F_1^{(5)}$)
 is an induced forest of size $a(G) \geq |F_1^{(5)}| + |F_2^{(5)}| + 3 \geq \lceil (4n+3)/7 \rceil$, a contradiction. 
\qed

\begin{lem}
\label{No6-2-A}
The following configuration is impossible in $G$:
$x$ is a $5^+$-vertex in $G$ with neighbors $x_1,x_2,x_3,...,x_m$ in cyclic order around $x$. $\{x_1,x_k\} \subseteq V_3$, $N(x_1) = \{x,z_1,y_1\}$, $N(x_k) = \{x,z_2,y_2\}$, and $\{y_1x_2, y_2x_{k-1}\} \subseteq E(G)$. Moreover, for $v \in \{y_1,z_1\}$, either $v \subseteq V_{\leq 4}$ or $R_{v,\{x_1\}} \neq \emptyset$; and for $v \in \{y_2,z_2\}$, either $v \subseteq V_{\leq 4}$ or $R_{v,\{x_k\}} \neq \emptyset$. 
\end{lem}

\pf By Lemmas \ref{No3RR}, \ref{No4-3R}, we may assume that $\{y_1,z_1,y_2,z_2\} \subseteq V_{4}$. Let $N(z_1)=\{z_1',x_1,x_m,w_1\}$, $N(y_1) = \{x_1,w_1,x_2,y_1'\}$, $N(z_2)=\{z_2',x_k,x_{k+1},w_2\}$, and $N(y_2) = \{x_k,w_2,$ $x_{k-1},y_2'\}$. By Lemma \ref{No434Edge}, $z_1x_2 \not\in E(G)$, $y_1x_m \not\in E(G)$, $z_2x_{k-1} \not\in E(G)$, and $y_2x_{k+1} \not\in E(G)$.

\medskip
\textit{Claim 1: $z_1'x \not\in E(G)$,  $y_1'x \not\in E(G)$,  $z_2'x \not\in E(G)$ and  $y_2'x \not\in E(G)$}.

For, suppose $z_1'x \in E(G)$. 
Then there exists a separation $(G_1,G_2)$ of $G$ such that $V(G_1 \cap G_2) = \{x,z_1,z_1'\}$, $N(x) \cap N(z_1) - \{x_1\} \subseteq V(G_1)$ and $\{x_1,y_1\} \subseteq V(G_2)$. 
If $|N(w_1) \cap N(y_1')| \leq 2$, let $F_1 = A((G_1-\{x,z_1,z_1',x_1,y_1\})/w_1y_1')$ with $u$ as the identification of $w_1$ and $y_1'$, and $F_2 = A(G_2-\{x,z_1,z_1'\})$. 
Then $|F_1| \geq \lceil (4(|G_1|-6)+3)/7 \rceil$, and $|F_2| \geq \lceil (4(|G_2|-6)+3)/7 \rceil$.
Let $F = F_1 \cup F_2 + \{z_1,x_1,y\}$. 
Now, $G[F]$ (if $u \not\in F_1$) or $G[F - \{u,y_1\} + \{w_1,y_1'\}]$ (if $u \in F_1$) is an induced forest in $G$, showing $a(G) \geq |F_1| + |F_2| + 3 \geq \lceil (4n+3)/7 \rceil$, a contradiction.

So $|N(w_1) \cap N(y_1')| > 2$. 
Then there exist $a_1 \in N(w_1) \cap N(y_1')$ and subgraphs $G_1',G_2',G_3'$ of $G$ such that 
$G_2' = G_2$,
$G_3'$ is the maximal subgraph of $G$ contained in the closed region of the plane bounded by the cycle $yw_1a_1y_1'y$ and containing $N(y_1') \cap N(w_1) - \{y_1\}$,
and $G_1'$ is obtained from $G$ by removing $G_2'-\{z_1',z_1,x\}$, and $G_3'-\{y_1',a_1,w_1\}$. 
Let $A_1 = \{a_1\}$ and $A_2 = \emptyset$.
For $i=1,2$, let $F_1^{(i)} = A(G_1'-\{x,z_1,z_1',x_1,y_1,w_1,y_1'\} - A_i)$, $F_2^{(i)} = A(G_2'-\{x,z_1,z_1'\})$, and $F_3^{(i)} = A(G_3'-\{w_1,y_1'\} - A_i)$. 
Then, $|F_1^{(i)}| \geq \lceil (4(|G_1'| - 7 - |A_i|)+3)/7 \rceil$, $|F_2^{(i)}| \geq \lceil (4(|G_2'| - 3)+3)/7 \rceil$, and  $|F_3^{(i)}| \geq \lceil (4(|G_3'| - 2 - |A_i|)+3)/7 \rceil$. 
Now, $F_1^{(i)} \cup F_2^{(i)} \cup F_3^{(i)} + \{z_1,x_1,y_1\} - (\{a_1\} \cap (F_1^{(i)} \triangle F_3^{(i)}))$ is an induced forest in $G$, showing $a(G) \geq |F_1^{(i)}| +  |F_2^{(i)}| + |F_3^{(i)}| + 3 -  (1-|A_i|)$. 
Let $(n_1,n_2,n_3) := (4(|G_1'| - 7)+3,4(|G_2'| - 7)+3,4(|G_3'| - 7)+3)$. 
By Lemma \ref{ineq2}(2),  $(n_1,n_2,n_3) \equiv (4,0,0), (0,0,4) \mod 7$. 
If $(n_1,n_2,n_3) \equiv (4,0,0) \mod 7 $, 
let $F_1^{(3)} = A((G_1'-\{x,z_1,z_1',x_1,y_1\})/w_1y_1')$ with $u$ as the identification of $w_1$ and $y_1'$, $F_2^{(3)} = A(G_2'-\{x,z_1,z_1'\})$, and $F_3^{(3)} = A(G_3')$. 
Then, $|F_1^{(3)}| \geq \lceil (4(|G_1'| - 6)+3)/7 \rceil$, $|F_2^{(3)}| \geq \lceil (4(|G_2'| - 3)+3)/7 \rceil$, and  $|F_3^{(3)}| \geq \lceil (4|G_3'|+3)/7 \rceil$.
Now, $G[F_1^{(3)} \cup F_2^{(3)} \cup F_3^{(3)} + \{z_1,x_1,y_1\} - (\{w_1,y_1',a_1\} \cap (F_1^{(3)} \triangle F_3^{(3)}))]$ (if $u \not\in F_1^{(3)}$) or $G[F_1^{(3)} \cup F_2^{(3)} \cup F_3^{(3)} + \{z_1,x_1,w_1,y_1'\} - (\{w_1,y_1',a_1\} \cap ( (F_1^{(3)} \cup \{w_1,y_1'\}) \triangle F_3^{(3)}))]$ (if $u \in F_1^{(3)}$) is an induced forest in $G$, showing $a(G) \geq |F_1^{(3)}| +  |F_2^{(3)}| + |F_3^{(3)}| + 3 - 3 \geq \lceil (4n+3)/7 \rceil$, a contradiction.
If $(n_1,n_2,n_3) \equiv (0,0,4) \mod 7$, 
let $F_1^{(4)} = A(G_1'-\{x,z_1,z_1',x_1,w_1\})$, $F_2^{(4)} = A(G_2'-\{x,z_1,z_1'\})$, and $F_3^{(4)} = A(G_3' - w_1)$. 
Then, $|F_1^{(4)}| \geq \lceil (4(|G_1'| - 5)+3)/7 \rceil$, $|F_2^{(4)}| \geq \lceil (4(|G_2'| - 3)+3)/7 \rceil$, and  $|F_3^{(4)}| \geq \lceil (4(|G_3'|-1)+3)/7 \rceil$.
Now, $G[F_1^{(4)} \cup F_2^{(4)} \cup F_3^{(4)} + \{z_1,x_1\} - (\{y_1',a_1\} \cap (F_1^{(4)} \triangle F_3^{(4)}))]$ is an induced forest in $G$, showing $a(G) \geq |F_1^{(4)}| +  |F_2^{(4)}| + |F_3^{(4)}| + 2 - 2 \geq \lceil (4n+3)/7 \rceil$, a contradiction.
By symmetry, we have $y_1'x \not\in E(G)$,  $z_2'x \not\in E(G)$ and  $y_2'x \not\in E(G)$.

\medskip
\textit{Claim 2: If $|N(y_1') \cap N(w_1)| > 2$, $|N(z_1') \cap N(w_1)| > 2$, and  there exist $a_1 \in N(y_1') \cap N(w_1)$ and a separation $(G_1,G_2)$ such that $V(G_1 \cap G_2) = \{w_1,y_1',a_1\}$, $\{x_1,y_1\} \subseteq V(G_1)$, and $N(y_1') \cap N(w_1) - \{y_1\} \subseteq V(G_2)$, then $4(|G_2|-2)+3 \not\equiv 4 \mod 7$.}

For, suppose $4(|G_2|-2)+3 \equiv 4 \mod 7$. There exist $a_1 \in N(y_1') \cap N(w_1), b_1 \in N(z_1') \cap N(w_1)$ and subgraphs $G_1',G_2',G_3'$ such that 
$G_2' = G_2,$
$G_3'$ is the maximal subgraph of $G$ contained in the closed region of the plane bounded by the cycle $z_1z_1'b_1w_1z_1$ and containing $N(z_1') \cap N(w_1) - \{z_1\}$,
and $G_1'$ is obtained from $G$ by removing $G_2'-\{y_1',a_1,w_1\}$, and $G_3'-\{z_1',b_1,w_1\}$.
Let $B_1 = \{b_1\}$ and $B_2 = \emptyset$. 
If  $|N(y_1') \cap N(x_2)| \leq 2$, 
for $i=1,2$, let $F_1^{(i)} = A((G_1'-\{x_1,y_1,z_1,w_1,x_m\} - B_i)/y_1'x_2 + z_1'x)$ with $u$ as the identification of $\{y_1',x_2\}$, $F_2^{(i)} = A(G_2' - w_1)$, and $F_3^{(i)} = A(G_3' - w_1 - B_i)$. 
Then $|F_1^{(i)}| \geq \lceil (4(|G_1'|-6-|B_i|)+3)/7 \rceil$, $|F_2^{(i)}| \geq \lceil (4(|G_2'|-1)+3)/7 \rceil$, and $|F_3^{(i)}| \geq \lceil (4(|G_3'|-1-|B_i|)+3)/7 \rceil$. 
Now $G[F_1^{(i)} \cup F_2^{(i)} \cup F_3^{(i)} \cup F_4^{(i)} + \{x_1,z_1,y_1\} - (\{y_1',a_1\} \cap (F_1^{(i)} \triangle F_2^{(i)})) - (\{z_1',b_1\} \cap (F_1^{(i)} \triangle F_3^{(i)}))]$ (if $u \not\in F_1^{(i)}$)
or $G[(F_1^{(i)} - u) \cup F_2^{(i)} \cup F_3^{(i)} \cup F_4^{(i)} + \{x_1,z_1,y_1',x_2\} - (\{y_1',a_1\} \cap ((F_1^{(i)} + y_1') \triangle F_2^{(i)})) - (\{z_1',b_1\} \cap (F_1^{(i)} \triangle F_3^{(i)}))]$ (if $u \not\in F_1^{(i)}$)
is an induced forest in $G$, showing $a(G) \geq |F_1^{(i)}| + |F_2^{(i)}| + 3 - 3 - (1-|B_i|)$. 
By Lemma \ref{ineq2}(1) (with $k=1$), $a(G) \geq \lceil (4n+3)/7 \rceil$, a contradiction.
Thus, $|N(y_1') \cap x_2| > 2$. Similarly, $|N(z_1') \cap N(x_m)| > 2$.

So $|N(y_1') \cap N(x_2)| > 2$ and $|N(z_1') \cap N(x_m)| > 2$.
There exist $c_1 \in N(y_1') \cap N(x_2), d_1 \in N(z_1') \cap N(x_m)$ and subgraphs $G_1''',G_2''',G_3''',G_4''',G_5'''$ of $G$ such that 
$G_2''' = G_2'$,
$G_3''' = G_3'$,
$G_4'''$ is the maximal subgraph of $G$ contained in the closed region of the plane bounded by the cycle $z_1z_1'd_1x_mz_1$ and containing $N(z_1') \cap N(x_m) - \{z_1\}$,
$G_5'''$ is the maximal subgraph of $G$ contained in the closed region of the plane bounded by the cycle $y_1y_1'c_1x_2y_1$ and containing $N(y_1') \cap N(x_2) - \{y_1\}$,
and $G_1''$ is obtained from $G$ by removing $G_2''-\{y_1',a_1,w_1\}$, $G_3''-\{z_1',b_1,w_1\}$, $G_4''-\{z_1',d_1,x_m\}$ and $G_5''-\{y_1',c_1,x_2\}$.
Let $B_i \subseteq \{b_1\}$, $C_i \subseteq \{c_1\}$, and $D_i \subseteq \{d_1\}$.
For each choice of $B_i,C_i,D_i$, let $F_1^{(i)} = A(G_1'''-\{x_1,y_1,z_1,w_1,z_1',y_1',x_2,x_m\} - B_i-C_i-D_i)$,
$F_2^{(i)} = A(G_2''' - \{w_1,y_1'\})$,
 $F_3^{(i)} = A(G_3''' - \{w_1,z_1'\} - B_i)$,
  $F_4^{(i)} = A(G_4''' - \{x_m,z_1'\} - D_i)$, 
  and $F_5^{(i)} = A(G_5''' - \{x_2,y_1'\} - C_i)$. 
  Then $|F_1^{(i)}| \geq \lceil (4(|G_1'''|-8-|B_i|-|C_i|-|D_i|)+3)/7 \rceil$,
   $|F_2^{(i)}| \geq \lceil (4(|G_2'''|-2)+3)/7 \rceil = \lceil (4(|G_2'|-2)+3)/7 \rceil = (4(|G_2'|-2)+3)/7 + 3/7$,
    $|F_3^{(i)}| \geq \lceil (4(|G_3'''|-2-|B_i|)+3)/7 \rceil$,
     $|F_4^{(i)}| \geq \lceil (4(|G_4'''|-2-|D_i|)+3)/7 \rceil$,
     and $|F_5^{(i)}| \geq \lceil (4(|G_5'''|-2-|C_i|)+3)/7 \rceil$. 
 Now $G[F_1^{(i)} \cup F_2^{(i)} \cup F_3^{(i)} \cup F_4^{(i)} \cup F_5^{(i)} + \{x_1,z_1,y_1\} -  (\{a_1\} \cap (F_1^{(i)} \triangle F_2^{(i)})) - (\{b_1\} \cap (F_1^{(i)} \triangle F_3^{(i)}))  - (\{c_1\} \cap (F_1^{(i)} \triangle F_4^{(i)}))  - (\{d_1\} \cap (F_1^{(i)} \triangle F_5^{(i)}))]$
 is an induced forest in $G$, showing $a(G) \geq |F_1^{(i)}| + |F_2^{(i)}| + |F_3^{(i)}| + |F_4^{(i)}| + |F_5^{(i)}| + 3 - 1 - (1-|B_i|) - (1-|C_i|) - (1-|D_i|)$.
By Lemma \ref{ineq2}(1) 
(with $k=3, a = |G_1'''| - 8, a_2 = |G_3'''| - 2, a_3 = |G_4'''| - 2, a_4 = |G_5'''| - 2, L = \{1\}, b_1 = |G_2'''| - 2, c=2$), 
$a(G) \geq \lceil (4n+3)/7 \rceil$, a contradiction.


\medskip

Now we distinguish several cases.

\medskip

Case 1: either $|N(z_1') \cap N(w_1)| \leq 2$ or $|N(y_1') \cap N(w_1)| \leq 2$; and either $|N(z_2') \cap N(w_2)| \leq 2$ or $|N(y_2') \cap N(w_2)| \leq 2$. 

We may assume that $|N(y_1') \cap N(w_1)| \leq 2$ and $|N(y_2') \cap N(w_2)| \leq 2$. Let $F' = A((G-\{x_1,x_k,y_1,y_2,x\})/\{y_1'w_1,y_2'w_2\} + \{z_1x_2,z_2x_{k-1}\})$ with $u_1$ (respectively, $u_2$) as the identification of $\{y_1',w_1\}$ (respectively, $\{y_2',w_2\}$). 
Then $|F'| \geq \lceil (4(n-7)+3)/7 \rceil$. 
Let $F = F' + \{x_1,x_k,y_1,y_2\}$ if $u_1,u_2 \not\in F'$, and otherwise $F$ obtained from $F' + \{x_1,x_k,y_1,y_2\}$ by deleting $\{y_1,u_1\}$ (respectively,  $\{y_2,u_2\}$) and adding $\{y_1',w_1\}$ (respectively,  $\{y_2',w_2\}$). 
Therefore, $F$ is an induced forest in $G$, showing $a(G) \geq |F'|+4 \geq \lceil (4n+3)/7 \rceil$, a contradiction.

Then, we have (both $|N(z_1') \cap N(w_1)| \geq 3$ and $|N(y_1') \cap N(w_1)| \geq 3$) or (both $|N(z_2') \cap N(w_2)| \geq 3$ and $|N(y_2') \cap N(w_2)| \geq 3$). Suppose $|N(z_1') \cap N(w_1)| \geq 3$ and $|N(y_1') \cap N(w_1)| \geq 3$. 

\medskip

Case 2: $|N(z_2') \cap N(w_2)| \leq 2$ or $|N(y_2') \cap N(w_2)| \leq 2$.

We may assume $|N(y_2') \cap N(w_2)| \leq 2$. 
There exist $a_1 \in N(y_1') \cap N(w_1)$ and a separation $(G_1,G_2)$ of $G$ such that $V(G_1 \cap G_2) = \{w_1,y_1',a_1\}$, $\{x,x_1,x_2,x_3,z_1,z_2\} \subseteq V(G_1)$, and $N(y_1') \cap N(w_1) - \{y_1\} \subseteq V(G_2)$. 
Let $A_1 = \emptyset$ and $A_2 = \{a_1\}$.
For $i=1,2$, let $F_1^{(i)} = A((G_1-\{x_1,x_k,y_1,y_2,x,y_1',w_1\} - A_i)/y_2'w_2 + z_2x_{k-1})$ with $u$ as the identification of $\{y_2',w_2\}$, and $F_2^{(i)} = A(G_2 - \{y_1',w_1\} - A_i)$. 
Then $|F_1^{(i)}| \geq \lceil (4(|G_1|-8 - |A_i|)+3)/7 \rceil$, and $|F_2^{(i)}| \geq \lceil (4(|G_2|-2  - |A_i|)+3)/7 \rceil$. 
Let $F = F_1^{(i)} \cup F_2^{(i)} + \{x_1,x_k,y_1,y_2\} - (\{a_1\} \cap (F_1^{(i)} \triangle F_2^{(i)}))$.
Now $G[F]$ (if $u \not\in F_1^{(i)}$) or $G[F - \{u,y_2\} +\{y_2',w_2\}]$ (if $u \in F_1^{(i)}$)
is an induced forest in $G$, showing $a(G) \geq |F_1^{(i)}| + |F_2^{(i)}| + 4 - (1-|A_i|)$. 
By Lemma \ref{ineq2}(2)
(with $a = |G_1|-8, a_2 = |G_2|-2, c=4$), 
$(4(|G_1|-8)+3,4(|G_2|-2)+3) \equiv (4,0),(0,4) \mod 7$. By \textit{Claim 2}, we have $(4(|G_1|-8)+3,4(|G_2|-2)+3) \equiv (4,0) \mod 7$. So assume it's the case. 
Let $F_1^{(3)} = A((G_1-\{x_1,x_k,y_1,y_2,x\})/\{y_1'w_1,y_2'w_2\} + \{z_1x_2,z_2x_{k-1}\})$ with $u_1$ (respectively, $u_2$) as the identification of $\{y_1',w_1\}$ (respectively, $\{y_2',w_2\}$), and $F_2^{(3)} = A(G_2)$. 
Then $|F_1^{(3)}| \geq \lceil (4(|G_1|-7)+3)/7 \rceil$, and $|F_2^{(3)}| \geq \lceil (4|G_2|+3)/7 \rceil$. 
Let $F^{(3)} = G[\overline{F_1}^{(3)} \cup F_2^{(3)} - (\{w_1,y_1',a_1\} \cap (\overline{F_1}^{(3)} \triangle F_2^{(3)}))]$,
 where 
$\overline{F_1}^{(3)} = {F_1}^{(3)}  + \{x_1,x_k,y_1,y_2\}$ if $u_1,u_2 \not\in F_1^{(3)}$, and otherwise, $\overline{F_1}^{(3)}$ obtained from ${F_1}^{(3)}  + \{x_1,x_k,y_1,y_2\}$ by deleting $\{y_1,u_1\}$ (respectively,  $\{y_2,u_2\}$) and adding $\{y_1',w_1\}$ (respectively,  $\{y_2',w_2\}$). 
Therefore, $F^{(3)}$ is an induced forest in $G$, showing $a(G) \geq |F_1^{(3)}|+|F_2^{(3)}|+4-3 \geq \lceil (4n+3)/7 \rceil$, a contradiction.

\medskip

Case 3: $|N(z_2') \cap N(w_2)| \geq 3$ and $|N(y_2') \cap N(w_2)| \geq 3$.

There exist $a_1 \in N(y_1') \cap N(w_1)$, $c_1 \in N(y_2') \cap N(w_2)$, $b_1 \in N(z_1') \cap N(w_1)$, $d_1 \in N(z_2') \cap N(w_2)$ and subgraphs 
$G_1,G_2,G_3,G_4,G_5$ such that 
$G_2$ is the maximal subgraph of $G$ contained in the closed region of the plane bounded by the cycle $y_1y_1'a_1w_1y_1$ and containing $N(y_1') \cap N(w_1) - \{y_1\}$,
$G_3$ is the maximal subgraph of $G$ contained in the closed region of the plane bounded by the cycle $y_2y_2'c_1w_2y_2$ and containing $N(y_2') \cap N(w_2) - \{y_2\}$,
$G_4$ is the maximal subgraph of $G$ contained in the closed region of the plane bounded by the cycle $z_1z_1'b_1w_1z_1$ and containing $N(z_1') \cap N(w_1) - \{z_1\}$,
$G_5$ is the maximal subgraph of $G$ contained in the closed region of the plane bounded by the cycle $z_2z_2'd_1w_2z_2$ and containing $N(z_2') \cap N(w_2) - \{z_2\}$,
and $G_1$ is obtained from $G$ by removing $G_2-\{y_1',a_1,w_1\}$, $G_3-\{y_2',c_1,w_2\}$, $G_4-\{z_1',b_1,w_1\}$ and $G_5-\{z_2',d_1,w_2\}$.
Let $A_i \subseteq \{a_1\}$, $B_i \subseteq \{b_1\}$, $C_i \subseteq \{c_1\}$ and $D_i \subseteq \{d_1\}$.
Let $G_1' = G_1 \cup G_4 \cup G_5$. 
For each choice of $A_i,C_i$, let $F_1^{(i)} = A(G_1'-\{x_1,y_1,y_1',w_1,x,x_k,y_2',w_2,y_2\} - A_i - C_i + \{z_1x_2,z_2x_{k-1}\})$,
 $F_2^{(i)} = A(G_2 - \{w_1,y_1'\} - A_i)$, and $F_3^{(i)} = A(G_3 - \{w_2,y_2'\} - C_i)$. 
 Then $|F_1^{(i)}| \geq \lceil (4(n+6-|G_2|-|G_3|-9-|A_i|-|C_i|)+3)/7 \rceil$,
  $|F_2^{(i)}| \geq \lceil (4(|G_2|-2 -|A_i|)+3)/7 \rceil$, and $|F_3^{(i)}| \geq \lceil (4(|G_3|-2-|C_i|)+3)/7 \rceil$. 
  Now $G[F_1^{(i)} \cup F_2^{(i)} \cup F_3^{(i)} + \{x_1,y_1,x_k,y_2\} - (\{a_1\} \cap (F_1^{(i)} \triangle F_2^{(i)})) - (\{c_1\} \cap (F_1^{(i)} \triangle F_3^{(i)}))]$
  is an induced forest in $G$, showing $a(G) \geq |F_1^{(i)}| + |F_2^{(i)}| + |F_3^{(i)}| + 4 - (1-|A_i|) - (1-|C_i|)$. 
By Lemma \ref{ineq2}(5)
(with $a = n+6-|G_2|-|G_3|-9, a_1 = |G_2|-2, a_2 = |G_3|-2, c=4$), 
$4(|G_2|-2)+3 \equiv 0,3,4,6 \mod 7$ and $4(|G_3|-2)+3 \equiv 0,3,4,6 \mod 7$.
By \textit{Claim 2} and by symmetry, we have $4(|G_i|-2)+3 \equiv 0,3,6 \mod 7$ for $i=2,3,4,5$ and if $4(|G_2|-2)+3 \equiv 3,6 \mod 7$ or $4(|G_4|-2)+3 \equiv 3,6 \mod 7$, then $4(|G_3|-2)+3 \equiv 4(|G_5|-2)+3 \equiv 0 \mod 7$ and vice versa. 

If $4(|G_2|-2)+3 \equiv 3,6 \mod 7$, then $4(|G_3|-2)+3 \equiv 0 \mod 7$. Let $A_1 = \{a_1\}$ if $4(|G_2|-2)+3 \equiv 6 \mod 7$ and  $A_1=\emptyset$ otherwise. 
Let $F_1^{(1)} = A((G_1'-\{x_1,y_1,y_1',w_1,x,x_k,y_2\} - A_1)/y_2'w_2 + \{z_1x_2, z_2x_{k-1}\})$ with $u$ as the identification of $\{y_2',w_2\}$,
 $F_2^{(1)} = A(G_2 - \{w_1,y_1'\} - A_1)$, 
 and $F_3^{(1)} = A(G_3)$. 
 Then $|F_1^{(1)}| \geq \lceil (4(n+6-|G_2|-|G_3|-8-|A_1|)+3)/7 \rceil$, $|F_2^{(1)}| \geq \lceil (4(|G_2|-2 -|A_1|)+3)/7 \rceil$, and $|F_3^{(1)}| \geq \lceil (4|G_3|+3)/7 \rceil$. 
 Now $G[F_1^{(1)} \cup F_2^{(1)} \cup F_3^{(1)} + \{x_1,y_1,x_k,y_2\} - (\{a_1\} \cap (F_1^{(1)} \triangle F_2^{(1)})) - (\{w_2,y_2',c_1\} \cap (F_1^{(1)} \triangle F_3^{(1)})) ]$ (if $u \not\in F_1^{(1)}$)
 or $G[(F_1^{(1)} - u) \cup F_2^{(1)} \cup F_3^{(1)} + \{x_1,y_1,x_k,y_2',w_2\} - (\{a_1\} \cap (F_1^{(1)} \triangle F_2^{(1)})) - (\{w_2,y_2',c_1\} \cap ((F_1^{(1)} \cup \{y_2',w_2\}) \triangle F_3^{(1)})) ]$ (if $u \in F_1^{(1)}$)
 is an induced forest in $G$, showing  $a(G) \geq |F_1^{(1)}| + |F_2^{(1)}| + |F_3^{(1)}| + 4 - (1-|A_1|) - 3 \geq \lceil (4n+3)/7 \rceil$, a contradiction.

So $4(|G_2|-2)+3 \equiv 4(|G_3|-2)+3 \equiv 0 \mod 7$ by symmetry. Let $F_1^{(2)} = A((G_1'-\{x_1,y_1,x,x_k,y_2\} )/\{y_1'w_1,y_2'w_2\} + \{z_1x_2, z_2x_{k-1}\})$ with $u_1$,$u_2$ as the identification of $\{y_1',w_1\}$, $\{y_2',w_2\}$ respectively, and $F_2^{(2)} = A(G_2)$ and $F_3^{(2)} = A(G_3)$. 
Then $|F_1^{(2)}| \geq \lceil (4(|G_1|-7)+3)/7 \rceil$, $|F_2^{(2)}| \geq \lceil (4|G_2|+3)/7 \rceil$, and $|F_3^{(2)}| \geq \lceil (4|G_3|+3)/7 \rceil$. 
Let $F^{(2)} = G[\overline{F_1}^{(2)} \cup F_2^{(2)} \cup F_3^{(2)} - (\{w_1,y_1',a_1\} \cap (\overline{F_1}^{(2)} \triangle F_2^{(2)})) - (\{w_2,y_2',c_1\} \cap (\overline{F_1}^{(2)} \triangle F_3^{(2)}))]$,
 where 
$\overline{F_1}^{(2)} = {F_1}^{(2)}  + \{x_1,x_k,y_1,y_2\}$ if $u_1,u_2 \not\in F_1^{(2)}$, and otherwise, $\overline{F_1}^{(2)}$ obtained from ${F_1}^{(3)}  + \{x_1,x_k,y_1,y_2\}$ by deleting $\{y_1,u_1\}$ (respectively,  $\{y_2,u_2\}$) and adding $\{y_1',w_1\}$ (respectively,  $\{y_2',w_2\}$) when $u_1 \in F_1^{(2)}$ (respectively,  $u_2 \in F_1^{(2)}$).
Therefore, $F^{(2)}$ is an induced forest in $G$, showing $a(G) \geq |F_1^{(2)}|+|F_2^{(2)}|+|F_3^{(2)}|+4-6 \geq \lceil (4n+3)/7 \rceil$, a contradiction.
\qed

\begin{lem}
\label{No5-2-C-5-2-B}
The following configuration is impossible in $G$:
$v$ is a $5$-vertex of type 5-2-C with neighbors $v_1,v_2,v_3,v_4,v_5$ in cyclic order, $\{v_1,v_3\} \subseteq V_3$, $vv_4xv_5$ is a facial cycle, $v_4$ is a $5$-vertex of type 5-2-B with neighbors $x,v,v_3',v_4',v_4''$ in cyclic order, $v_5 \in V_4, \{x,v_4'\} \subseteq V_3$, $v_4v_3'yv_4'v$, $v_4v_4'zv_4''v$, $v_4v_4''x'xv$ are facial cycles and $x' \in V_4$. 
\end{lem}

\pf Let $N(v_5) = \{x,v,v_5',v_1'\}$ where $x'v_5' \in E(G)$ and $v_1v_1' \in E(G)$.

\medskip

Case 1: $N(v_4'') \cap N(v_5') = \{x'\}$ and $|N(z) \cap N(y)| \leq 2$. 

Let $F' = A((G-\{x',x,v_5,v_1',v_1,v_4,v,v_4',v_3',v_3\})/\{v_4''v_5',yz\})$ with $v',y'$ as the identification of $\{v_4'',v_5'\},\{y,z\}$ respectively. 
Then $|F'| \geq \lceil (4(n-12)+3)/7 \rceil$. 
Let $F = F' + \{v_1,v_3,v_4,v_5,x,x',v_4'\}$ if $v',y' \not\in F'$, and otherwise, $F$ obtained from $F' + \{v_1,v_3,v_4,v_5,x,x',$ $v_4'\}$ by deleting $\{y',v_4'\}$ (respectively, $\{v',x'\}$) and adding $\{y,z\}$ (respectively, $\{v_4'',v_5'\}$) when $y' \in F'$ (respectively, $v' \in F'$ ).
Therefore, $F$ is an induced forest in $G$, showing $a(G) \geq |F'|+7 \geq \lceil (4n+3)/7 \rceil$, a contradiction.

\medskip

Case 2: $N(v_4'') \cap N(v_5') = \{x'\}$ and $|N(z) \cap N(y)| > 2$.

There exist $a_1 \in N(z) \cap N(y)$ and a separation $(G_1,G_2)$ such that $V(G_1 \cap G_2) = \{y,z,a_1\}$, $\{v,v_1,v_2,v_3,v_4,v_5,x,x',v_4,v_4'',v_3'\} \subseteq V(G_1)$, and $N(z) \cap N(y) - \{v_4'\} \subseteq V(G_2)$. 
Let $A_1 = \emptyset$ and $A_2 = \{a_1\}$.
For $i=1,2$, let $F_1^{(i)} = A((G_1 - \{x',x,v_5,v_1',v_1,v_4,v,v_3',v_3,y,z,$ $v_4'\} - A_i)/v_4''v_5')$ with $v'$ as the identification of $\{v_4'',v_5'\}$, and $F_2^{(i)} = A(G_2 - \{y,z\} - A_i)$. 
Then $|F_1^{(i)}| \geq \lceil (4(|G_1|-13-|A_i|)+3)/7 \rceil$, and $|F_2^{(i)}| \geq \lceil (4(|G_2|-2-|A_i|)+3)/7 \rceil$. 
Let $F = F_1^{(i)} \cup F_2^{(i)} + \{v_1,v_3,v_4,v_5,x,x',v_4'\} - (\{a_1\} \cap (F_1^{(i)} \triangle F_2^{(i)}))$.
Now $G[F]$ (if $v' \not\in F_1^{(i)}$) or $G[F-\{v',x'\}+\{v_4'',v_5'\}]$ (if $v' \in F_1^{(i)}$)
is an induced forest in $G$, showing $a(G) \geq |F_1^{(i)}| + |F_2^{(i)}| + 7 - (1-|A_i|)$. 
By Lemma \ref{ineq2}(1) 
(with $k=1, a= |G_1|-13, a_1 = |G_2|- 2, L = \emptyset, c = 7 $), 
$a(G) \geq \lceil (4n+3)/7 \rceil$, a contradiction.


\medskip

Case 3: $|N(v_4'') \cap N(v_5')| > 1$.

There exist $b_1 \in N(v_4'') \cap N(v_5')$ and a separation $(G_1,G_2)$ such that $V(G_1 \cap G_2) = \{v_4'',v_5,x',b_1\}$, $\{v,v_1,v_2,v_3,v_4,v_5,x,x',v_4,v_4'',v_3'\} \subseteq  V(G_1)$, and $N(v_4'') \cap N(v_5') - \{b_1\} \subseteq V(G_2)$. 
Let $F_1^{(1)} = A(G_1 - \{x',x,v_5,v_1',v_1,v_4,v,v_3',v_3,v_4'',v_5',b_1\})$, and $F_2^{(1)} = A(G_2 - \{x',v_4'',v_5',b_1\})$. 
Then $|F_1^{(1)}| \geq \lceil (4(|G_1|-12)+3)/7 \rceil$, and $|F_2^{(1)}| \geq \lceil (4(|G_2|-4)+3)/7 \rceil$. 
Then $G[F_1^{(1)} \cup F_2^{(1)} + \{v_1,v_3,v_4,v_5,x,x'\}]$ is an induced forest in $G$, showing $a(G) \geq |F_1^{(1)}| + |F_2^{(1)}| + 6$. 
If $v_4'b_1 \not\in E(G)$, 
let $F_1^{(2)} = A(G_1 - \{x',x,v_5,v_1',v_1,v_4,v,v_3',v_3,v_4'',v_5'\} + v_4'b_1)$,
 and $F_2^{(2)} = A(G_2 - \{x',v_4'',v_5'\})$. 
 Then $|F_1^{(2)}| \geq \lceil (4(|G_1|-11)+3)/7 \rceil$, and $|F_2^{(2)}| \geq \lceil (4(|G_2|-3)+3)/7 \rceil$. 
 Now $G[F_1^{(2)} \cup F_2^{(2)} + \{v_1,v_3,v_4,v_5,x,x'\} - (\{b_1\} \cap (F_1^{(2)} \triangle F_2^{(2)}))]$
 is an induced forest in $G$, showing $a(G) \geq |F_1^{(2)}| + |F_2^{(2)}| + 6 - 1$.
 By Lemma \ref{ineq2}(1) 
 (with $k=1, a=|G_1|-11, a_1 = |G_2|-3, L = \emptyset, c=6 $), 
 $a(G) \geq \lceil (4n+3)/7 \rceil$, a contradiction.
 So $v_4'a_1 \in E(G)$. 
 Let $F_1^{(3)} = A(G_1 - \{x',x,v_5,v_1',v_1,v_4,v,v_3',v_3,v_4'',v_5',b_1,v_4'\})$, and $F_2^{(3)} = A(G_2 - \{x',v_4'',v_5',b_1\})$. 
 Then $|F_1^{(3)}| \geq \lceil (4(|G_1|-13)+3)/7 \rceil$, and $|F_2^{(3)}| \geq \lceil (4(|G_2|-4)+3)/7 \rceil$. 
 Now $G[F_1^{(3)} \cup F_2^{(3)} + \{v_1,v_3,v_4,v_5,x,x',v_4'\}]$
  is an induced forest in $G$, showing $a(G) \geq |F_1^{(3)}| + |F_2^{(3)}| + 7 \geq \lceil (4n+3)/7 \rceil$, a contradiction.
\qed

\begin{lem}
\label{No5-2-C-5-1-A}
The following configuration is impossible in $G$: 
$v$ is a $5$-vertex of type 5-2-C with neighbors $v_1,v_2,v_3,v_4,v_5$ in cyclic  order, $\{v_1,v_3\} \subseteq V_3$; $vv_4xv_5v$ is a facial cycle. $v_4$ is a $5$-vertex of type 5-1-A with neighbors $x,v,v_3',v_4',v_4''$ in cyclic order, $v_5 \in V_4, x \in V_3$; $v_4v_3'yv_4'v_4$, $v_4v_4'zv_4''v_4$, $v_4v_4''x'xv_4$ are facial cycles, $\{y,z\} \subseteq V_3$, and $x' \in V_{\leq 4}$. 
\end{lem}

\pf Let $xx'v_5'v_5$, $vv_5v_1'v_1$ bound $4$-faces. Let $F' = A(G-\{v_5',v_1',x',v_5,x,v_1,v_4'',v_4,v,z,v_4',$ $y,v_3',v_3\})$. 
Then $|F'| \geq \lceil (4(n-14)+3)/7 \rceil$. 
Now $G[F' + \{x',x,v_5,v_1,v_3,v_4,y,z\}]$ is an induced forest in $G$, showing $a(G) \geq |F'|+8 \geq \lceil (4n+3)/7 \rceil$, a contradiction.
\qed

\begin{lem}
\label{No5-2-C-5-2-B_second}
The following configuration is impossible in $G$: 
$v$ is a $5$-vertex of type 5-2-B with neighbors $v_1,v_2,v_3,v_4,v_5$ in cyclic  order, $\{v_1,v_3\} \subseteq V_3$, $N(v_1) = \{v,v_1',v_1''\}$, $N(v_3) = \{v,v_3',v_3''\}$, $v_1'' \in V_{\geq 5}$ and $\{v_3',v_3''\} \subseteq V_4$;
$vv_1v_1'v_2v$, $vv_2v_3'v_3v$, $vv_3v_3''v_4v$, $vv_5v_1''v_1v$ are facial cycles.
$v_2$ is a $5$-vertex of type 5-2-C with neighbors $v,v_1',v_2',v_2'',v_3'$ in cyclic order, $\{v_1',v_2''\} \subseteq V_3$. 
\end{lem}

\pf Let $t \in N(v_2'') \cap N(v_3') $ and $v_2v_2''tv_3'v_2$ bound a $4$-face. 
Let $N(v_3'') = \{v_3,v_4,s_1,s_2\}$ and $v_3v_3's_2v_3''v_3$ bound a $4$-face. 
Let $w \in N(v_1'') \cap N(v_1') $ and $v_1v_1''wv_1'v_1$ bound a $4$-face. 

By Lemma \ref{No434Edge}, $v_3'v_4 \not\in E(G)$. 
We claim that $v_2v_1'' \not\in E(G)$. 
Since $G$ is simple, $v_1'' \not\in \{v,v_1'\} $. 
Since $v_1'' \in V_{\geq 5}$, $v_1'' \not\in \{v_3',v_2''\} $. 
If $v_2' = v_1''$, then since $G$ is a quadrangulation, $v_2v_1''wv_1'v_2$ bound a $4$-face and thus $N(w) = \{v_1'',v_1'\}$. But this contradicts Lemma \ref{No333}.

If $|N(s_1) \cap N(s_2)| \leq 2$, then 
let $F' = A(G-\{v,v_1,v_1',v_3,v_3'',w\}/s_1s_2 + \{v_3'v_4, v_2v_1''\})$ with $s'$ the identification of $\{s_1,s_2\}$. 
Then $|F'| \geq \lceil (4(n-7)+3)/7 \rceil$. 
Now $G[F' + \{v_3,v_3'',v_1,v_1'\}]$ (if $s' \not \in F'$ ) 
or $G[F' - s' + \{v_3,s_1,s_2,v_1,v_1'\}]$ (if $s' \in F'$ ) 
 is an induced forest in $G$, showing $a(G) \geq |F'|+4 \geq \lceil (4n+3)/7 \rceil$, a contradiction.
 
So $|N(s_1) \cap N(s_2)| > 2$.
There exist $a_1 \in N(s_1) \cap N(s_2)$ and a separation $(G_1,G_2)$ such that $V(G_1 \cap G_2) = \{s_1,s_2,a_1\}$, $\{v,v_1,v_2,v_3,v_4,v_5,v_1',v_1'',v_2',v_2'',t,v_3',v_3''\} \subseteq V(G_1)$, and $N(s_1) \cap N(s_2) - \{v_3''\} \subseteq V(G_2)$. 
Let $A_1 = \emptyset$ and $A_2 = \{a_1\}$.
For $i=1,2$, let $F_1^{(i)} = A(G_1-\{v,v_1,v_1',v_3,v_3'',w,s_1,s_2\} - A_i + \{v_3'v_4, v_2v_1''\})$ with $s'$ the identification of $\{s_1,s_2\}$, and $F_2^{(i)} = A(G_2 - \{s_1,s_2\} - A_i)$. 
Then $|F_1^{(i)}| \geq \lceil (4(|G_1|-8-|A_i|)+3)/7 \rceil$, and $|F_2^{(i)}| \geq \lceil (4(|G_2|-2-|A_i|)+3)/7 \rceil$. 
Now $F_1^{(i)} \cup F_2^{(i)} + \{v_3,v_3'',v_1,v_1'\} - (\{a_1\} \cap (F_1^{(i)} \triangle F_2^{(i)}))$
is an induced forest in $G$, showing $a(G) \geq |F_1^{(i)}| + |F_2^{(i)}| + 4 - (1-|A_i|)$. 
By Lemma \ref{ineq2}(2) 
(with $a= |G_1|-8, a_1 = |G_2|- 2, L = \emptyset, c = 4 $), 
$(4(|G_1|-8)+3, 4(|G_2|- 2)+3) \equiv (0,4), (4,0) \mod 7 $. 

If $(4(|G_1|-8)+3, 4(|G_2|- 2)+3) \equiv (4,0) \mod 7 $, then 
let $F_1^{(3)} = A(G_1-\{v,v_1,v_1',v_3,v_3'',$ $w\}/s_1s_2 + \{v_3'v_4, v_2v_1''\})$ with $s'$ the identification of $\{s_1,s_2\}$, and $F_2^{(3)} = A(G_2)$. 
Then $|F_1^{(3)}| \geq \lceil (4(|G_1|-7)+3)/7 \rceil$, and $|F_2^{(3)}| \geq \lceil (4|G_2|+3)/7 \rceil$. 
Now $G[F_1^{(3)} \cup F_2^{(3)} + \{v_3,v_3'',v_1,v_1'\} - (\{a_1,s_1,s_2\} \cap (F_1^{(3)} \triangle F_2^{(3)}))]$ (if $s' \not \in F_1^{(3)}$ ) 
or $G[F_1^{(3)} \cup F_2^{(3)} - s' + \{v_3,s_1,s_2,v_1,$ $v_1'\} - (\{a_1,s_1,s_2\} \cap ((F_1^{(3)} + \{s_1,s_2\}) \triangle F_2^{(3)}))]$ (if $s' \in F_1^{(3)}$ ) 
 is an induced forest in $G$, showing $a(G) \geq |F_1^{(3)}| + |F_2^{(3)}| + 4 - 3 \geq \lceil (4n+3)/7 \rceil$, a contradiction.

So $(4(|G_1|-8)+3, 4(|G_2|- 2)+3) \equiv (0,4) \mod 7 $.
First, we claim that $vt \not \in E(G)$. $t \not \in \{v_2,v_3\}$ since $G$ is simple. 
$t \neq v_4$ by Lemma \ref{No434Edge}.
$t \neq v_1$ since $v_1'' \in V_{\geq 5}$ and $v_1' \in V_3$. 
Suppose $t = v_5$. Since $G$ is a quadrangulation, $v_2''v_5 \in E(G)$. 
let $F_4 = A(G-\{v_3,v_3',v_3'',s_2,v_4,v,v_2,v_2'',v_1,v_1',v_5,w\})$. 
Then $|F_4| \geq \lceil (4(n-12)+3)/7 \rceil$. 
Now $G[F_4 + \{v_3'',v_3,v_3',v,v_2'',v_1,v_1'\}]$ is an induced forest in $G$, showing $a(G) \geq |F_4|+7 \geq \lceil (4n+3)/7 \rceil$, a contradiction.
Secondly, suppose $|N(s_1) \cap N(v_4)| \leq 2$. 
Then let $F_1^{(5)} = A(G_1-\{v_3'',s_2,v_3,v_3',v_2\}/s_1v_4 + \{vt\})$ with $s'$ the identification of $\{s_1,v_4\}$, and $F_2^{(5)} = A(G_2-s_2)$. 
Then $|F_1^{(5)}| \geq \lceil (4(|G_1|-6)+3)/7 \rceil$, and $|F_2^{(5)}| \geq \lceil (4(|G_2|-1)+3)/7 \rceil$. 
Now $G[F_1^{(5)} \cup F_2^{(5)} + \{v_3,v_3'',v_3'\} - (\{a_1,s_1\} \cap (F_1^{(5)} \triangle F_2^{(5)}))]$ (if $s' \not \in F_1^{(5)}$ ) 
or $G[F_1^{(5)} \cup F_2^{(5)} - s' + \{s_1,v_4,v_3,v_3'\} - (\{a_1,s_1\} \cap ((F_1^{(5)}+s_1) \triangle F_2^{(5)}))]$ (if $s'  \in F_1^{(5)}$ )
 is an induced forest in $G$, showing $a(G) \geq |F_1^{(5)}| + |F_2^{(5)}| + 3 - 2 \geq \lceil (4n+3)/7 \rceil$, a contradiction.
 Now, $|N(s_1) \cap N(v_4)| > 2$. 
 There exist $a_1 \in N(s_1) \cap N(s_2)$, $b_1 \in N(s_1) \cap N(v_4)$ and subgraphs 
$G_1', G_2', G_3'$ such that $G_2' = G_2$,
$G_3'$ is the maximal subgraph of $G$ contained in the closed region of the plane bounded by the cycle $v_3''s_1b_1v_4v_3''$ and containing $N(s_1) \cap N(v_4) - \{v_3''\}$,
and $G_1'$ is obtained from $G$ by removing $G_2'-\{s_1,a_1,s_2\}$ and $G_3'-\{s_1,b_1,v_4\}$.
Let $B_7  = \{b_1\}$ and $B_8  = \emptyset$. 
For $i=7,8$, let $F_1^{(i)} = A((G_1'-\{v_3,v_3'',s_1,v_4\} - B_i)/vv_3')$ with $v'$ as the identification of $\{v,v_3'\}$,
 $F_2^{(i)} = A(G_2' - \{s_1\})$, and $F_3^{(i)} = A(G_3' - \{s_1,v_4\} - B_i)$. 
 Then $|F_1^{(i)}| \geq \lceil (4(|G_1'|-5-|B_i|)+3)/7 \rceil$,
  $|F_2^{(i)}| \geq \lceil (4(|G_2'|-1)+3)/7 \rceil = \lceil (4(|G_2|-1)+3)/7 \rceil = (4(|G_2|-1)+3)/7 + 6/7$,
   and $|F_3^{(i)}| \geq \lceil (4(|G_3'|-2-|B_i|)+3)/7 \rceil$. 
Now $G[F_1^{(i)} \cup F_2^{(i)} \cup F_3^{(i)} + \{v_3,v_3''\} - (\{a_1,s_2\} \cap (F_1^{(i)} \triangle F_2^{(i)})) - (\{b_1\} \cap (F_1^{(i)} \triangle F_3^{(i)}))]$
(if $v' \not\in F_1^{(i)}$)
or $G[F_1^{(i)} \cup F_2^{(i)} \cup F_3^{(i)} - v' + \{v,v_3',v_3''\} - (\{a_1,s_2\} \cap (F_1^{(i)} \triangle F_2^{(i)})) - (\{b_1\} \cap (F_1^{(i)} \triangle F_3^{(i)}))]$
(if $v' \in F_1^{(i)}$)
  is an induced forest in $G$, showing $a(G) \geq |F_1^{(i)}| + |F_2^{(i)}| + |F_3^{(i)}| + 2 - 2 - (1-|B_i|)$. 
  By Lemma \ref{ineq2}(1) 
  (with $k=1, a = |G_1'|-5, a_1 = |G_3'|-2, L = \{1\}, b_1= |G_2'|-1, c=0$), 
  $a(G) \geq \lceil (4n+3)/7 \rceil$, a contradiction.
\qed


\section{Proof of Theorem \ref{MainTheorem}}

We define the discharging rules as follow: For each $v \in V(G)$, let $ch(v):= |N(v)| - 4$. Let $\mathcal{F}$ be the set of all the faces of graph $G$. For each $f\in \mathcal{F}$, let $ch(f) := |f| - 4$. Then, by Euler's Formula, the total charge of graph $G$ is $$\sum_{v \in V(G)} ch(v) + \sum_{f \in \mathcal{F}} ch(f) = \sum_{v \in V(G)} (N(v) - 4) + \sum_{f \in \mathcal{F}} (|f| - 4) = 4|E(G)| - 4|V(G)| - 4|\mathcal{F}| = -8$$

\begin{defi}
\label{new_charge}
For $v \in V(G)$, suppose $|N(v)| \geq 5$. We redistribute the charges as follow:
\begin{itemize}
\item[(i)] Suppose $R_{v,U} \neq \emptyset$ for some $U \subseteq V(G)$. If $R_{v,U} = \{\{r\}\}$, then $v$ sends charge $|N(v)| - 4$ to $r$; If $R_{v,U} = \{\{r_1,r_2\}\}$, then $v$ sends charge $(|N(v)| - 4)/2$ to both $r_1$ and $r_2$; If $R_{v,U} = \{R_1,R_2\}$, then $v$ sends charge $|N(v)| - 4$ to $R_1 \cap R_2$; 
\item[(ii)] Suppose $R_{v,\{u\}} = \emptyset $ and $vu \in E(G)$ for some $u \in V_3$. Let $N(u) = \{v,u_1,u_2\}$. If for both $w \in \{u_1,u_2\}$, either $w \in V_{\leq 4}$ or $R_{w,\{u\}} \neq \emptyset$, then $v$ sends charge $1$ to $u$; If $u_2 \in V_{\geq 5}$ and $R_{u_2,\{u\}} = \emptyset$, then $v$ sends charge $1/2$ to $u$;
\item[(iii)] Suppose $R_{v,\{u\}} = \emptyset $ and $xwyvx$ is a facial cycle such that $\{x,y\} \subseteq N(v), x \in V_{\geq 5}$, $w \in V_3$ and $y \in V_4$. If neither $v$ nor $x$ is of type 5-2-C, then $v$ sends charge $1/4$ to $x$. 
%
\end{itemize}
\end{defi}
We denote the new charge of $v$ as $ch'(v)$. We remark that if $v$ sends charge $1/4$ to $x$ in both faces bounded by $xw_1y_1vx$ and $xw_2y_2vx$ by Definition \ref{new_charge} (iii), then $v$ sends charge $1/2$ to $x$.


We show that for $v \in V(G)$, $ch'(v) \geq 0$. If $|N(v)| = 2$, then by Lemma \ref{2summary} and Definition \ref{new_charge} (i), $v$ either receives at least $1$ from $\{v_5,v_5'\} \subseteq V_{\geq 5} \cap N(v)$ where $R_{v_5,\{v\}} = R_{v_5',\{v\}} =\emptyset$ or at least $2$ from $v_6 \in V_{\geq 6} \cap N(v)$ where $R_{v_6,\{v\}} = \emptyset$. Hence, $ch'(v) \geq ch(v)+2 = 0$. 
Suppose $|N(v)| = 3$ with $N(v) = \{v_1,v_2,v_3\}$. If $R_{v_3,\{v\}} \neq \emptyset$, then by Lemmas \ref{No3RR}, \ref{No4-3R}, $\{v_1,v_2\} \subseteq V_{\geq 5}$, $R_{v_1,\{v\}} = R_{v_2,\{v\}} = \emptyset$; thus, by Definition \ref{new_charge}(ii), $v$ receives $1/2$ from each of $v_1$ and $v_2$,  and $ch'(v) = ch(v) + 1/2 + 1/2 = 0$. 
Now, assume $R_{v_i,\{v\}} = \emptyset$ for $i=1,2,3$. By Corollary \ref{3summary}, there exists $v_1 \in N(v) \cap V_{\geq 5}$. 
By Definition \ref{new_charge}(ii), $v$ receives at least $1$ from $N(v)$ and thus $ch'(v) \geq ch(v) + 1 = 0$. 
If $|N(v)| = 4$, then $v$ does not receive or send charge to other vertices. Therefore, $ch'(v) = ch(v) = 0$. 
If $|N(v)| \geq 5$ and $R_{v,\emptyset} \neq \emptyset$, then by Definition \ref{new_charge}(i) and Lemma \ref{NoTwoRNeighbors}, $v$ sends $|N(v)| - 4$ to $R_{v,\emptyset}$ only. Therefore, $ch'(v) = ch(v) - (|N(v)| - 4) = 0$. 

Next, assume $|N(v)| \geq 5$ and $R_{v,\emptyset} = \emptyset$. We distinguish the cases by Definition \ref{5structure}: 

\begin{itemize}
\item $v$ is of type 5-2-A. By Lemma \ref{No6-2-A}, $v$ does not exist in $G$, a contradiction;

\item $v$ is of type 5-2-B. Let $N(v) = \{v_1,v_2,v_3,v_4,v_5\}$ in order. Let $\{v_1,v_3\} \subseteq V_3$, $\{v_2,v_4,v_5\} \subseteq V_{\geq 4}$, $N(v_1) = \{v_1',v_1'',v\}$, $N(v_3) = \{v_3',v_3'',v\}$, $u \in V_{\leq 4}$ or $R_{u,\{v_3\}} \neq \emptyset$ for $u \in  \{v_3',v_3''\}$ and $v_1' \in V_{\geq 5}$, $R_{v_1',\{v_1\}} = \emptyset$. By Definition \ref{new_charge}(ii), $v$ sends $1/2$ to $v_1$, $1$ to $v_3$.
By Lemma \ref{No3RR}, \ref{No4-3R}, \ref{No3345},  $R_{u,\{v_3\}} = \emptyset$ and $u \in V_4$ for $u \in  \{v_3',v_3''\}$.
By Lemma \ref{5-2-B-combined}, $\{v_2,v_4\} \subseteq V_{\geq 5}$ and $R_{v_2,\{v\}} = R_{v_4,\{v\}} = \emptyset$. By Lemma \ref{No5-2-C-5-2-B}, $v_4$ is not of type 5-2-C. By Lemma \ref{No5-2-C-5-2-B_second}, $v_2$ is not of type 5-2-C.
 Hence, $v$ receives $1/4$ from each of $v_2$ and $v_4$  by Definition \ref{new_charge}(iii). In addition, by Lemma \ref{5-2-B-lemma3}, $v$ does not send charge to $v_4$. So $ch'(v) = ch(v)-1-1/2+1/4+1/4 = 0$; 

\item $v$ is of type 5-2-C. Let $N(v) = \{v_1,v_2,v_3,v_4,v_5\}$ in order. Let $\{v_1,v_3\} \subseteq V_3$, $\{v_2,v_4,v_5\} \subseteq V_{\geq 4}$, $N(v_1) = \{v_1',v_1'',v\}$, $N(v_3) = \{v_3',v_3'',v\}$, $v_1' \in V_{\geq 5}$, $R_{v_1',\{v_1\}} = \emptyset$, $v_3' \in V_{\geq 5}$, and $R_{v_3',\{v_3\}} = \emptyset$. By Definition \ref{new_charge}(ii)(iii), $v$ sends $1/2$ to both $v_1$ and $v_3$. So $ch'(v) = ch(v)-1/2-1/2 = 0$; 

\item $v$ is of type 5-1-A. Let $N(v) = \{v_1,v_2,v_3,v_4,v_5\}$ in order. Let $v_1 \in V_3, \{v_2,v_3,v_4,v_5\} \subseteq V_{\geq 4}$, $N(v_1) = \{v_1',v_1'',v\}$, and  $u \in V_{\leq 4}$ or $R_{u,\{v_1\}} \neq \emptyset$ for $u \in  \{v_1',v_1''\}$.
Let $vv_2v_2'v_3v$, $vv_3v_3'v_4v$, $vv_4v_4'v_5v$ be facial cycles. 
By Lemmas \ref{No3RR}, \ref{No4-3R}, \ref{No3345}, $R_{u,\{v_1\}} = \emptyset$ and $u \in V_4$  for $u \in  \{v_1',v_1''\}$.
By Lemma \ref{5-2-B-combined}, $\{v_2,v_5\} \subseteq V_{\geq 5}$ and $R_{v_2,\{v\}} = R_{v_5,\{v\}} = \emptyset$. 
If $v_2,v_5$ are not of type 5-2-C, then by Definition \ref{new_charge}(iii) $v$ receives $1/4$ from each of $v_2$ and $v_5$. By Definition \ref{new_charge}(ii)(iii), $v$ sends $1$ to $v_1$ and $1/4$ to at most two of $\{v_2,v_3,v_4,v_5\}$. So $ch'(v) \geq ch(v)-1-1/4 \times 2+1/4+1/4 = 0$; 
If both $v_2$ and $v_5$ are of type 5-2-C, then $\{v_2',v_4'\} \subseteq V_3$ and by Lemma \ref{No5-2-C-5-1-A} $v_3' \not\in V_3$. By Definition \ref{new_charge}(ii)(iii), $v$ sends $1$ to $v_1$ and no charge to $\{v_2,v_3,v_4,v_5\}$. So $ch'(v) \geq ch(v)-1 = 0$; 
If exactly one of $v_2$ and $v_5$ is of type 5-2-C, say $v_2$, then by Lemma \ref{No5-2-C-5-1-A}, $|\{v_3',v_4'\} \cap V_3| \leq 1$. By Definition \ref{new_charge}(ii)(iii), $v$ sends $1$ to $v_1$ and $1/4$ to at most one of $\{v_3,v_4,v_5\}$ and $v$ receives $1/4$ from $v_5$. So $ch'(v) \geq ch(v) + 1/4 -1 - 1/4 = 0$;

\item $v$ is of type 5-1-B. Let $N(v) = \{v_1,v_2,v_3,v_4,v_5\}$ in order. Let $v_1 \in V_3, \{v_2,v_3,v_4,v_5\} \subseteq V_{\geq 4}$, $N(v_1) = \{v_1',v_1'',v\}$, $v_1' \in V_{\geq 5}$, and $R_{v_1',\{v_1\}} =  \emptyset$. By Lemma \ref{5-1-B} and Definition \ref{new_charge}(ii)(iii), $v$ sends $1/2$ to $v_1$ and $1/4$ to at most two of $\{v_2,v_3,v_4,v_5\}$. So $ch'(v) \geq ch(v)-1/2-1/4 \times 2 = 0$;

\item $v$ is of type 5-0. By Definition \ref{new_charge}(iii), $ch'(v) \geq ch(v)-1/4 \times 4 = 0$. 
\end{itemize}

Suppose $|N(v)| = 6$ and $R_{v,\emptyset} = \emptyset$. We distinguish the cases by Definition \ref{6structure}:

\begin{itemize}
\item $v$ is of type 6-3. By Lemma \ref{No6-2-A} and Definition \ref{new_charge}(ii), $v$ sends $1$ to at most one of $\{v_1,v_3,v_5\}$. By Definition \ref{new_charge}(iii), $v$ sends no charge to $\{v_2,v_4,v_6\}$. So $ch'(v) \geq ch(v)-1-1/2\times 2 = 0$; 

\item $v$ is of type 6-2-A. By Lemma \ref{No6-2-A} and Definition \ref{new_charge}(ii), $v$ sends $1$ to at most one of $\{v_1,v_3\}$. By Definition \ref{new_charge}(iii), $v$ sends no charge to $v_2$ and $1/4$ to at most two of $\{v_4,v_5,v_6\}$. So $ch'(v) \geq ch(v)-1-1/2-2 \times 1/4 = 0$; 

\item $v$ is of type 6-2-B. Let $N(v_1) = \{v_1',v_1'',v\}$ and $N(v_3) = \{v_3',v_3'',v\}$. 
By Lemma \ref{No6-2-A} and Definition \ref{new_charge}(ii), $v$ sends $1$ to at most one of $\{v_1,v_4\}$. By Definition \ref{new_charge}(iii), $v$ sends $1/4$ to at most one of $\{v_2,v_3\}$ and to at most one of $\{v_5,v_6\}$. So $ch'(v) \geq ch(v)-1-1/2-2 \times 1/4 = 0$;

\item $v$ is of type 6-1. By Definition \ref{new_charge}(ii), $v$ sends at most $1$ to $v_1$. By Definition \ref{new_charge}(iii), $v$ sends $1/4$ to at most four of $\{v_2,v_3,v_4,v_5,v_6\}$. So $ch'(v) \geq ch(v)-1-4 \times 1/4 = 0$;

\item $v$ is of type 6-0. By Definition \ref{new_charge}(iii), $ch'(v) \geq ch(v)-6 \times 1/4 = 1/2$.
\end{itemize}

Suppose $|N(v)| = 7$ and $R_{v,\emptyset} = \emptyset$. Let $N(v):=\{v_1,v_2,v_3,v_4,v_5,v_6,v_7\}$ in order.
 If $|N(v) \cap V_3| = 3$, then we may assume that they are $v_1,v_3,v_5$. 
 By Lemma \ref{No6-2-A} and Definition \ref{new_charge}(ii), $v$ sends 1 to at most one of $N(v) \cap V_3$. 
So $ch'(v) \geq ch(v)-1-2 \times 1/2-1/4 > 0$. If $|N(v) \cap V_3| \leq 2$, then by Definition \ref{new_charge}(iii), $v$ sends $1/4$ to at most three of $N(v)$. So, $ch'(v) \geq ch(v)-1-1-3\times 1/4 > 0$. Suppose $|N(v)| \geq 8$. We observe that if we amortize the redistribution of charge to all the faces which $v$ is incident with, then $v$ sends at most $1/2$ in each face. So $ch'(v) \geq ch(v)-1/2 \times |N(v)|  = |N(v)|/2 - 4 \geq 0$.

Therefore, $ch'(v) \geq 0$ for $v \in V(G)$. Since $G$ is a quadrangulation by Lemma \ref{basic_lemma}, $ch'(f) := ch(f) = 0$. Then, the total charge after redistribution is $\sum_{v \in V(G)} ch'(v) + \sum_{f \in \mathcal{F}} ch'(f) \geq 0$, which contradicts Euler's Formula. 

To conclude, the minimum counterexample $G$ does not exist. This completes the proof of Theorem \ref{MainTheorem}. 



\end{document}